\documentclass[a4paper,oneside,11pt]{article}
\usepackage[a4paper]{geometry}
\usepackage{aeguill}                 
\usepackage{graphicx}
\usepackage{caption}
\usepackage{amsmath,amsfonts,amssymb,amsthm}
\usepackage{mathabx}
\usepackage{pstricks,comment} 
\usepackage{pst-plot}
\usepackage{pstricks-add}
\usepackage{multido}
\usepackage{url}
\usepackage{epstopdf,enumerate}
\usepackage{srcltx}
\usepackage[utf8]{inputenc} 
\usepackage[T1]{fontenc} 
\usepackage[english]{babel} 
\usepackage{hyperref}
\usepackage[all]{xy} 
\usepackage{titlesec}
\usepackage{mathabx}
\usepackage{tikz}
\usepackage{fancyhdr}
\usepackage{mathrsfs}
\usepackage[us,12hr]{datetime}
\fancypagestyle{plain}{
\fancyhf{}
\rfoot{ {\ddmmyyyydate\today}}
\lfoot{\thepage}
}

    \newtheoremstyle{TheoremNum}
        {\topsep}{\topsep}              
        {\itshape}                      
        {}                              
        {\bfseries}                     
        {.}                             
        { }                             
        {\thmname{#1}\thmnote{ \bfseries #3}}

{\theoremstyle {definition} \newtheorem {defi} {Definition}[section]}
{\theoremstyle {plain}  \newtheorem {theo} [defi] {Theorem}}
{\theoremstyle {plain}  \newtheorem {coro} [defi] {Corollary}}
{\theoremstyle {plain} \newtheorem {prop} [defi] {Proposition}}
{\theoremstyle {plain} \newtheorem {lem}[defi] {Lemma}}
{\theoremstyle {plain} \newtheorem {rmq}[defi] {Remark}}
{\theoremstyle {plain} }
{\theoremstyle{TheoremNum} }
{\theoremstyle{TheoremNum} }
{\theoremstyle{TheoremNum} }

\newcommand{\Aut}{\mathrm{Aut}}
\newcommand{\Out}{\mathrm{Out}}

\newcommand{\Mod}{\mathrm{Mod}}

\newcommand{\Stab}{\mathrm{Stab}}

\newcommand{\Curr}{\mathrm{Curr}}
\newcommand{\PCurr}{\mathbb{P}\mathrm{Curr}}
\newcommand{\Supp}{\mathrm{Supp}}

\newcommand{\ZZ}{\mathbb{Z}}

\newcommand{\NN}{\mathbb{N}}
\newcommand{\RR}{\mathbb{R}}

\newcommand{\dem}{\noindent{\bf Proof. }}

\title{North-South type dynamics of
relative atoroidal automorphisms of free groups on a relative space of currents}
\author{Yassine Guerch}
\date{\today}
\begin{document}
\maketitle
\renewcommand*\labelenumi{(\theenumi)}

\begin{abstract}
This paper, which is the second of a series of three papers, studies dynamical properties of elements of $\Out(F_{\tt n})$, the outer automorphism group of a nonabelian free group $F_{\tt n}$. We prove that, for every exponentially growing outer automorphism of $F_{\tt n}$, there exists a preferred compact topological space, the space of currents relative to a malnomal subgroup system, on which $\phi$ acts by homeomorphism with a North-South dynamics behavior.
\end{abstract}

\section{Introduction}\label{Section Introduction}
Let ${\tt n} \geq 2$. This paper is the second of a sequence of three papers where we study the growth of the conjugacy classes of elements of $F_{\tt n}$ under iterations of elements of $\Out(F_{\tt n})$, the outer automorphism group of a nonabelian free group of rank ${\tt n}$. An outer automorphism $\phi \in \Out(F_{\tt n})$ is \emph{exponentially growing} if there exist $g \in F_{\tt n}$, a representative $\Phi$ of $\phi$, a free basis $\mathfrak{B}$ of $F_{\tt n}$ and a constant $K>0$ such that, for every $m \in \NN^*$, we have $$\ell_{\mathfrak{B}}(\Phi^m(g)) \geq e^{Km},$$ where $\ell_{\mathfrak{B}}(\Phi^m(g))$ denotes the length of $\Phi^m(g)$ in the basis $\mathfrak{B}$. Such an element $g$ is said to be \emph{exponentially growing under iteration of $\phi$} and the set of elements of $F_{\tt n}$ which have exponential growth under iteration of $\phi$ is the \emph{pure exponential part of $\phi$}. It is known, using for instance the train track technology of Bestvina and Handel (see~\cite{BesHan92}), that every element $g$ of $F_{\tt n}$ which is not exponentially growing under iteration of $\phi$ is \emph{polynomially growing under iteration of $\phi$}, that is, there exist $\Phi \in \phi$ and an integer $K \in \NN$ such that, for every $m \in \NN^*$, we have $\ell_{\mathfrak{B}}(\Phi^m(g)) \leq (m+1)^K$.

Initiated by Švarc, Milnor and Wolf, and particularly developped by Guivarc'h, Gromov and Grigorchuk, growth problems in groups is a major field of study in geometric and dynamical group theory, see for instance~\cite{LubotzkySegal2003,Mann12,Helfgott2015}. Many works study the subfield of the element growths
under iteration of group automorphisms (see for instance~\cite{BesFeiHan00,Levitt09,clay2019atoroidal}), for instance in the context of hyperbolic groups. See
in particular~\cite{Coulon2019} for examples of intermediate growth rates. As another example, Dahmani and Krishna~\cite{DahmaniKrish2020} found a sufficient condition for the suspension of an automorphism of a hyperbolic group to be relatively hyperbolic, and this condition is linked with the structure of the set of all elements of the hyperbolic group which have polynomial growth under iterations of the considered automorphism. Such exponentially growing outer automorphisms of $F_{\tt n}$ were already studied in distinct contexts. For instance, Bestvina, Feighn and Handel~\cite{BesFeiHan00} used them to prove the Tits alternative for $\Out(F_{\tt n})$. 

If $\phi \in \Out(F_{\tt n})$, we denote by $\mathrm{Poly}(\phi)$ the set of elements $g$ of $F_{\tt n}$ such that $g$ is polynomially growing under iteration of $\phi$. Let $\mathrm{Poly}(H)=\bigcap_{\phi \in H} \mathrm{Poly}(\phi)$. The aim of this series of papers is to prove the following theorem.

\begin{theo}\label{Theo intro 2}
Let ${\tt n} \geq 3$ and let $H$ be a subgroup of $\Out(F_{\tt n})$. There exists $\phi \in H$ such that $\mathrm{Poly}(\phi)=\mathrm{Poly}(H)$.
\end{theo}

Informally, Theorem~\ref{Theo intro 2} shows that the exponential growth of a subgroup $H$ of $\Out(F_{\tt n})$ is encaptured by the exponential growth of a single element of $H$. Indeed, if $g \in F_{\tt n}$ has exponential growth for some element $\psi \in H$, then $g$ has exponential growth for an element $\phi \in H$ given by Theorem~\ref{Theo intro 2}. The proof relies on dynamical properties of the action of outer automorphisms on some preferred topological space. In this article, we study the dynamical properties of the elements of the subgroup $H$ of $F_{\tt n}$ that will be used in~\cite{Guerch2021Polygrowth} in order to construct an element $\phi \in H$ given by Theorem~\ref{Theo intro 2}. 

Let $\phi \in \Out(F_{\tt n})$ be an exponentially growing outer automorphism. In this article, we construct natural (compact, metrizable) topological spaces $X$ on which a subgroup of $\Out(F_{\tt n})$ containing $\phi$ acts by homeomorphisms with the additional property that $\phi$ acts with \emph{North-South dynamics}: there exist two proper disjoint closed subsets of $X$ such that every point of $X$ which is not contained in these subsets converges to one of the two subsets under positive or negative iteration of $\phi$. North-South dynamics are preferred tools to apply ping-pong arguments similar to the ones of Tits~\cite{Tits72} and are used to obtain structural properties of some groups.

The topological space $X$ that we use in the proof of Theorem~\ref{Theo intro 2} is constructed in such a way that it allows us to create a dictionnary between dynamical properties of the action of $\phi$ on $X$ and growth properties of elements of $F_{\tt n}$ under iterations of $\phi$. In order to construct $X$, we first need to detect all the elements $g$ of $F_{\tt n}$ such that the length of $[g]$ with respect to any basis of $F_{\tt n}$ grows at most polynomially fast fast under iteration of $\phi$. Levitt~\cite{Levitt09} proved that there exist finitely many finitely generated subgroups $H_1,\ldots,H_k$ of $F_{\tt n}$ such that the conjugacy class of an element $g$ of $F_{\tt n}$ is not exponentially growing under iteration of $\phi$ if and only if $g$ is contained in a conjugate of some $H_i$ for $i \in \{1,\ldots,k\}$. Moreover, the set $\mathcal{A}(\phi)=\{[H_1],\ldots, [H_k]\}$ is a \emph{malnormal subgroup system}: for every $i \in \{1,\ldots,k\}$, the group $H_i$ is a malnormal subgroup of $F_{\tt n}$ and for every distinct subgroups $A$ and $B$ such that $[A], [B] \in \mathcal{A}(\phi)$, we have $A \cap B=\{e\}$. Every element of $F_{\tt n}$ which is contained in a conjugate of some $H_i$ with $i \in \{1,\ldots,k\}$ has polynomial growth under iteration of $\phi$. Moreover, we have $\mathrm{Poly}(\phi)=\bigcup_{i=1}^r\bigcup_{g \in F_{\tt n}} gH_ig^{-1}$.

In~\cite{Guerch2021currents}, we construct a compact, metrizable space, called the space of projectivised currents relative to $\mathcal{A}(\phi)$, denoted by $\PCurr(F_{\tt n},\mathcal{A}(\phi))$, which is the space of projectivised Radon measures on the double boundary of $F_{\tt n}$ relative to $\mathcal{A}(\phi)$, equipped with the weak-star topology (see~Section~\ref{Section relative currents} for precise definitions). In~\cite{Guerch2021currents}, we proved that the set of currents associated with \emph{$\mathcal{A}(\phi)$-nonperipheral} conjugacy classes of elements of $g$ of $F_{\tt n}$, that is, such that $g$ is not contained in the conjugacy class of some $H_i$ with $i \in \{1,\ldots,k\}$, is dense in $\PCurr(F_{\tt n},\mathcal{A}(\phi))$. Thus, the set of conjugacy classes of elements of $F_{\tt n}$ whose length grows exponentially fast under iteration of $\phi$ is dense in $\PCurr(F_{\tt n},\mathcal{A}(\phi))$. If we denote by $\Out(F_{\tt n},\mathcal{A}(\phi))$ the subgroup of $\Out(F_{\tt n})$ consisting in every element $\psi \in \Out(F_{\tt n})$ such that $\psi(\mathcal{A}(\phi))=\mathcal{A}(\phi)$, the group $\Out(F_{\tt n},\mathcal{A}(\phi))$ acts by homeomorphisms on $\PCurr(F_{\tt n},\mathcal{A}(\phi))$ by pushing forward the measures. In this article, we prove the following theorem.

\begin{theo}[see~Theorem~\ref{Theo North south dynamics relative atoroidal}]\label{Theo intro 1}
Let ${\tt n} \geq 3$ and let $\phi$ be an exponentially growing outer automorphism. The outer automorphism $\phi$ acts with North-South dynamics on the space $\PCurr(F_{\tt n},\mathcal{A}(\phi))$.
\end{theo}

In fact, we prove a slightly stronger result since we prove a \emph{uniform North-South dynamics result}, that is, the convergence in the North-South dynamics statement can be made uniform on compact subsets of $\PCurr(F_{\tt n},\mathcal{A}(\phi))$. As explained above, North-South dynamics results given by Theorem~\ref{Theo intro 1} will be a key point in the proof of Theorem~\ref{Theo intro 2}. 

Such dynamical results already appear in similar contexts. For instance, Tits proved in~\cite{Tits72} its alternative for linear groups using North-South dynamics and ping-pong arguments. In the context of the mapping class group $\Mod(S)$ of a compact connected orientable surface $S$ of genus at least $2$, pseudo-Anosov elements acts with North-South dynamics on the space of projectivised measured foliations (\cite{Thurston88}, see also the work of Ivanov~\cite{Ivanov92}) or the curve complex~\cite{MasMin99}. Using this North-South dynamics, Ivanov~\cite{Ivanov92} (see also the work of McCarthy~\cite{McCarthy85}) later proved a Tits alternative for subgroups of $\Mod(S)$. Similarly, North-South dynamics results were obtained for certain classes of outer automorphisms of $F_{\tt n}$. For instance, \emph{fully irreducible outer automorphisms} act on the compactified Outer space~\cite{LevLus03} or the space of projectivised currents~(\cite{Martin95}, see also the work of Uyanik \cite{Uyanik2015}) with a North-South dynamics and \emph{atoroidal outer automorphisms} act on the space of projectivsed currents with a North-South dynamics~\cite{LustigUyanik2019,Uyanik2019}. Clay and Uyanik~\cite{clay2019atoroidal} applied this result in the proof of the fact that, for every subgroup $H$ of $\Out(F_{\tt n})$, either $H$ contains an atoroidal outer automorphism or there exists a nontrivial element $g$ of $F_{\tt n}$ such that, for every element $\phi \in H$, there exists $k \in \NN^*$ such that we have $\phi^k([g])=[g]$. Such dynamical results were later extended to relative contexts by Gupta~\cite{gupta2017relative,gupta18}.  

In order to prove Theorem~\ref{Theo intro 2}, we will need a slightly stronger result than Theorem~\ref{Theo intro 1}. Indeed, let $\phi \in \Out(F_{\tt n})$ and let $\mathcal{A}(\phi)=\{[H_1],\ldots,[H_k]\}$. Suppose that $\phi$ preserves the conjugacy class of a corank one free factor $A$ of $F_{\tt n}$. Let $\mathcal{A}(\phi) \wedge A$ be the malnormal subgroup system consisting in the conjugacy classes of the intersection of the conjugates of the subgroups $H_i$ with $i \in \{1,\ldots,k\}$ with $A$. Note that, by Theorem~\ref{Theo intro 1}, there exist closed disjoint subsets $\Delta_{\pm}(\phi|_A)$ such that the outer automorphism $\phi|_A \in \Out(A,\mathcal{A}(\phi) \wedge A)$ acts with North-South dynamics on $\PCurr(A,\mathcal{A}(\phi) \wedge A)$ with respect to $\Delta_{\pm}(\phi|_A)$. There is a canonical embedding $\PCurr(A,\mathcal{A}(\phi) \wedge A) \hookrightarrow \PCurr(F_{\tt n},\mathcal{A}(\phi) \wedge A)$, and we denote by $\Delta_{\pm}(\phi)$ the image of $\Delta_{\pm}(\phi|_A)$ in $\PCurr(F_{\tt n},\mathcal{A}(\phi) \wedge A)$. We will need to understand the dynamics of $\phi$ on the space $\PCurr(F_{\tt n},\mathcal{A}(\phi) \wedge A)$. As there might exist elements in $F_{\tt n}$ which have polynomial growth under iterations of $\phi$ and which are not contained in a conjugate of $A$, one cannot apply Theorem~\ref{Theo intro 1} to obtain a North-South dynamics result. However, we obtain the following result.

\begin{theo}[see~Theorem~\ref{Theo North-South dynamics almost atoroidal}]\label{Theo intro 3}
Let ${\tt n} \geq 3$ and let $\phi \in \Out(F_{\tt n})$ be an exponentially growing outer automorphism which preserves a corank one free factor $A$. There exist two convex compact subsets $\widehat{\Delta}_{\pm}(\phi)$ of $\PCurr(F_{\tt n},\mathcal{A}(\phi) \wedge A)$ such that the following holds. Let $U_{\pm}$ be open neighborhoods of $\Delta_{\pm}(\phi)$ in $\PCurr(F_{\tt n},\mathcal{A}(\phi) \wedge A)$ and $\widehat{V}_{\pm}$ be open neighborhoods of $\widehat{\Delta}_{\pm}(\phi)$ in $\PCurr(F_{\tt n},\mathcal{A}(\phi) \wedge A)$. There exists $M \in \NN^*$ such that for every $n \geq M$, we have 
$$
\phi^{\pm n}(\PCurr(F_{\tt n},\mathcal{A}(\phi) \wedge A)-\widehat{V}_{\mp}) \subseteq U_{\pm}.
$$
\end{theo} 

In~\cite[Theorem~4.15]{clay2019atoroidal}, Clay and Uyanik proved an analogue of Theorem~\ref{Theo intro 3} in the context of atoroidal outer automorphisms of $F_{\tt n}$. In Theorem~\ref{Theo intro 3}, the two convex subsets $\widehat{\Delta}_{\pm}(\phi)$ have nonempty intersection, so that Theorem~\ref{Theo intro 3} is not a North-South dynamics result as defined above. However, Theorem~\ref{Theo intro 3} gives a sufficiently precise description of the dynamics of $\phi$ for our considerations. The intersection $\widehat{\Delta}_{+}(\phi) \cap \widehat{\Delta}_{-}(\phi)$ corresponds informally to the polynomial growth part of $\phi$. This intersection, denoted by $K_{PG}$ in the rest of the article, is the closure in $\PCurr(F_{\tt n},\mathcal{A}(\phi) \wedge A)$ of the $(\mathcal{A}(\phi) \wedge A)$-nonperipheral elements of $F_{\tt n}$ which have polynomial growth under iteration of $\phi$. In Section~\ref{Section poly growing currents}, we give a complete study of the subspace $K_{PG}$ in a more general context. 

In fact, Section~\ref{Section PG subgraph of a CT map} is devoted to the study of the polynomial growth of an exponentially growing outer automorphism. Following the works of Bestvina, Feighn and Handel ~\cite{BesFeiHan00,BesFeiHan04}, of Feighn and Handel~\cite{FeiHan06} and of Handel and Mosher~\cite{HandelMosher20}, we use appropriate relative train track representatives of a power of an exponentially growing outer automorphism $\phi$ in order to describe $\mathcal{A}(\phi)$ geometrically. It gives rise to a (not necessarily connected) topological graph $G^{\ast}$ such that the fundamental group of every connected component $G_c^{\ast}$ of $G^{\ast}$ injects into $F_{\tt n}$ and such that the set $\{[\pi_1(G_c^{\ast})]\}_{G_c^{\ast} \in \pi_0(G^{\ast})}$ where $\pi_1(G_c^{\ast})$ is viewed as a subgroup of $F_{\tt n}$ is equal to $\mathcal{A}(\phi)$ (see~Proposition~\ref{Prop circuits in Gpg are elements in poly subgroup}). We then use this characterization of $\mathcal{A}(\phi)$ in Section~\ref{Section poly growing currents} in order to describe the subset $K_{PG}$.

We now sketch a proof of Theorem~\ref{Theo intro 1}. The proofs of Theorem~\ref{Theo intro 1} and Theorem~\ref{Theo intro 3} given in this paper are long and quite technical, this is why we postpone the proof of Theorem~\ref{Theo intro 2} in~\cite{Guerch2021Polygrowth}. Let $\phi \in \Out(F_{\tt n})$ be exponentially growing. The first step is to construct the closed subsets $\Delta_{\pm}(\phi)$ associated with $\phi$ as defined in Therorem~\ref{Theo intro 1}. This is done in Section~\ref{Section Stable unstable currents}. In order to construct them, we use as inspiration the construction given by Lustig and Uyanik in \cite{LustigUyanik2019} (see also \cite{Uyanik2019,gupta2017relative}). We choose an appropriate relative train track representative $f \colon G \to G$ of a power of $\phi$, where $G$ is a graph whose fundamental group is isomorphic to $F_{\tt n}$. A current of $\Delta_+(\phi)$ is then constructed by considering occurrences of paths in $\lim_{m \to \infty} f^m(e)$, where $e$ is an edge in $G$ whose length grows exponentially fast under iteration of $f$ (see Proposition~\ref{Prop existence relative currents atoroidal automorphisms}). Currents of $\Delta_-(\phi)$ are then defined similarly using a representative of a power of $\phi^{-1}$. We then prove Theorem~\ref{Theo intro 1} in Section~\ref{Section North South dyn 1}. Let $[\mu] \in \PCurr(F_{\tt n},\mathcal{A}(\phi))-\Delta_{\pm}(\phi)$ be the current associated with a $\mathcal{A}(\phi)$-nonperipheral conjugacy class $[w] \in F_{\tt n}$. Then $[w]$ is represented by a circuit $\gamma_w$ in the graph $G$. In order to show that we have $\lim_{m \to \infty}\phi^m([\mu]) \in \Delta_+(\phi)$, we prove that the proportion of the path $f^m(\gamma_g)$ which grows exponentially fast under iteration of $f$ tends to $1$ as $m$ goes to infinity. This fact is sufficient to prove that $$\lim_{m \to \infty}\phi^m([\mu]) \in \Delta_+(\phi)$$ (see~Lemma~\ref{Lem conversion goodness closeness PCurr}). We then conclude the proof using the density of currents associated with nonperipheral elements in $F_{\tt n}$ proved in \cite{Guerch2021currents}. Theorem~\ref{Theo intro 3} is then proved in Section~\ref{Section North SOuth dyn almost} using a combination of Theorem~\ref{Theo intro 1} and the description of the space $K_{PG}$.

\medskip

{\small{\bf Acknowledgments. } I warmly thank my advisors, Camille Horbez and Frédéric Paulin, for their precious advices and for carefully reading the different versions of this article. }

\section{Preliminaries}

\subsection{Malnormal subgroup systems of $F_{\tt n}$}\label{Subsection malnormal}

Let ${\tt n}$ be an integer greater than $1$ and let $F_{\tt n}$ be a free group of rank ${\tt n}$. A \emph{subgroup system of $F_{\tt n}$} is a finite (possibly empty) set $\mathcal{A}$ whose elements are conjugacy classes of nontrivial (that is distinct from $\{1\}$) finite rank subgroups of $F_{\tt n}$. There exists a partial order on the set of subgroup systems of $F_{\tt n}$, where $\mathcal{A}_1 \leq \mathcal{A}_2$ if for every subgroup $A_1$ of $F_{\tt n}$ such that $[A_1] \in \mathcal{A}_1$, there exists a subgroup $A_2$ of $F_{\tt n}$ such that $[A_2] \in \mathcal{A}_2$ and $A_1$ is a subgroup of $A_2$. The \emph{stabilizer in $\Out(F_{\tt n})$ of a subgroup system} $\mathcal{A}$, denoted by $\Out(F_{\tt n},\mathcal{A})$, is the set of all elements $\phi \in \Out(F_{\tt n})$ such that $\phi(\mathcal{A})=\mathcal{A}$. 

Recall that a subgroup $A$ of $F_{\tt n}$ is \emph{malnormal} if for every element $x \in F_{\tt n}-A$, we have $xAx^{-1} \cap A=\{e\}$. A subgroup system $\mathcal{A}$ is said to be \emph{malnormal} if every subgroup $A$ of $F_{\tt n}$ such that $[A] \in \mathcal{A}$ is malnormal and, for all subgroups $A_1,A_2$ of $F_{\tt n}$ such that $[A_1],[A_2] \in \mathcal{A}$, if $A_1 \cap A_2$ is nontrivial then $A_1=A_2$. An element $g \in F_{\tt n}$ is \emph{$\mathcal{A}$-peripheral} (or simply \emph{peripheral} if there is no ambiguity) if it is trivial or conjugate into one of the subgroups of $\mathcal{A}$, and \emph{$\mathcal{A}$-nonperipheral} otherwise.

An important class of examples of malnormal subgroup systems is given by the \emph{free factor systems}. A \emph{free factor system of $F_{\tt n}$} is a (possibly empty) set $\mathcal{F}$ of conjugacy classes $\{[A_1],\ldots,[A_r]\}$ of nontrivial subgroups $A_1,\ldots,A_r$ of $F_{\tt n}$ such that there exists an integer $k \in \NN$ with $F_{\tt n}=A_1 \ast \ldots \ast A_r \ast F_k$. The free factor system $\mathcal{F}$ is \emph{sporadic} if $(k+r,k) \leq (2,1)$ for the lexicographic order, and is \emph{nonsporadic} otherwise. Therefore, the sporadic free factor systems are those of the form $\{[C]\}$ where $C$ has rank at least equal to $n-1$ and those of the form $\{[A],[B]\}$ with $F_{\tt n}=A \ast B$. An ascending sequence of free factor systems $\mathcal{F}_1 \leq \ldots \leq \mathcal{F}_i=\{[F_{\tt n}]\}$ of $F_{\tt n}$ is called a \emph{filtration of $F_{\tt n}$}.

Given a free factor system $\mathcal{F}$ of $F_{\tt n}$, a \emph{free factor of $(F_{\tt n},\mathcal{F})$} is a subgroup $A$ of $F_{\tt n}$ such that there exists a free factor system $\mathcal{F}'$ of $F_{\tt n}$ with $[A] \in \mathcal{F}'$ and $\mathcal{F} \leq \mathcal{F}'$. When $\mathcal{F}=\varnothing$, we say that $A$ is a \emph{free factor of $F_{\tt n}$}. A free factor of $(F_{\tt n},\mathcal{F})$ is \emph{proper} if it is nontrivial, not equal to $F_{\tt n}$ and if its conjugacy class does not belong to $\mathcal{F}$. 

Another class of examples of malnormal subgroup systems is the following one. An outer automorphism $\phi \in \Out(F_{\tt n})$ is \emph{exponentially growing} if there exists $g \in F_{\tt n}$ such that the length of the conjugacy class $[g]$ of $g$ in $F_{\tt n}$ with respect to some basis of $F_{\tt n}$ grows exponentially fast under iteration of $\phi$. If $\phi \in \Out(F_{\tt n})$ is not exponentially growing, then $\phi$ is \emph{polynomially growing}. For an automorphism $\alpha \in \Aut(F_{\tt n})$, we say that $\alpha$ is exponentially growing if there exists $g \in F_{\tt n}$ such that the length of $g$ grows exponentially fast under iteration of $\phi$. Otherwise, $\alpha$ is polynomially growing. Let $\phi \in \Out(F_{\tt n})$ be exponentially growing. A subgroup $P$ of $F_{\tt n}$ is a \emph{polynomial subgroup} of $\phi$ if there exist $k \in \NN^*$ and a representative $\alpha$ of $\phi^k$ such that $\alpha(P)=P$ and $\alpha|_P$ is polynomially growing. By~\cite[Proposition~1.4]{Levitt09}, there exist finitely many conjugacy classes $[H_1],\ldots,[H_k]$ of maximal polynomial subgroups of $\phi$. Moreover, the proof of~\cite[Proposition~1.4]{Levitt09} implies that the set $\mathcal{H}=\{[H_1],\ldots,[H_k]\}$ is a malnormal subgroup system. Indeed, Levitt shows that there exists a nontrivial $\RR$-tree $T$ in the boundary of Culler and Vogtmann Outer space~\cite{Vogtmann1986} on which $F_{\tt n}$ acts with trivial arc stabilizers, such that $\phi$ preserves the homothety class of $T$ and such that the groups $H_1 \ldots, H_k$ are elliptic in $T$. If two distinct subgroups $A,B$ of $F_{\tt n}$ such that $[A],[B] \in \mathcal{H}$ fix distinct points in $T$, then their intersection is trivial. If $A$ and $B$ fix the same point $x$ in $T$, then (up to taking a power of $\phi$) $\phi$ preserves $[\Stab(x)]$ and an inductive argument on the rank using $\phi|_{\Stab(x)}$ (the rank of $\Stab(x)$ is less than ${\tt n}$ by a result of Gaboriau-Levitt~\cite{GabLev95}) shows that the intersection of $A$ and $B$ is trivial. We denote this malnormal subgroup system by $\mathcal{A}(\phi)$. Note that, if $H$ is a subgroup of $F_{\tt n}$ such that $[H] \in \mathcal{A}(\phi)$, there exists $\Phi^{-1} \in \phi^{-1}$ such that $\Phi^{-1}(H)=H$ and $\Phi^{-1}|_H$ is polynomially growing. Hence we have $\mathcal{A}(\phi) \leq \mathcal{A}(\phi^{-1})$. By symmetry, we have
\begin{equation}\label{Equation p6}
\mathcal{A}(\phi)=\mathcal{A}(\phi^{-1}).
\end{equation}

Let $\mathcal{A}$ be a malnormal subgroup system and let $\phi \in \Out(F_{\tt n},\mathcal{A})$ be a relative outer automorphism. We say that $\phi$ is \emph{atoroidal relative to $\mathcal{A}$} if, for every $k \in \NN^*$, the element $\phi^k$ does not preserve the conjugacy class of any $\mathcal{A}$-nonperipheral element. We say that $\phi$ is \emph{expanding relative to $\mathcal{A}$} if $\mathcal{A}(\phi) \leq \mathcal{A}$. Note that an expanding outer automorphism relative to $\mathcal{A}$ is in particular atoroidal relative to $\mathcal{A}$. When $\mathcal{A}=\varnothing$, then the outer automorphism $\phi$ is expanding relative to $\mathcal{A}$ if and only if for every nontrivial element $g \in F_{\tt n}$, the length of the conjugacy class $[g]$ of $g$ in $F_{\tt n}$ with respect to some basis of $F_{\tt n}$ grows exponentially fast under iteration of $\phi$. Therefore, by a result of Levitt \cite[Corollary~1.6]{Levitt09}, the outer automorphism $\phi$ is expanding relative to $\mathcal{A}=\varnothing$ if and only if $\phi$ is atoroidal relative to $\mathcal{A}=\varnothing$.

Let $\mathcal{A}=\{[A_1],\ldots,[A_r]\}$ be a malnormal subgroup system and let $\mathcal{F}$ be a free factor system. Let $i \in \{1,\ldots,r\}$. By~\cite[Theorem~3.14]{ScoWal79} for the action of $A_i$ on one of its Cayley graphs, there exist finitely many subgroups $A_i^{(1)},\ldots,A_i^{(k_i)}$ of $A_i$ such that:

\medskip

\noindent{$(1)$ } for every $j \in \{1,\ldots,k_i\}$, there exists a subgroup $B$ of $F_{\tt n}$ such that $[B] \in \mathcal{F}$ and $A_i^{(j)}=B \cap A_i$;

\medskip

\noindent{$(2)$ } for every subgroup $B$ of $F_{\tt n}$ such that $[B] \in \mathcal{F}$ and $B \cap A_i \neq \{e\}$, there exists $j \in \{1,\ldots,k_i\}$ such that $A_i^{(j)}=B \cap A_i$;

\medskip

\noindent{$(3)$ } the subgroup $A_i^{(1)} \ast \ldots \ast A_i^{(k_i)}$ is a free factor of $A_i$.

\bigskip

Thus, one can define a new subgroup system as $$\mathcal{F} \wedge \mathcal{A}=\bigcup_{i=1}^r\{[A_i^{(1)}],\ldots,[A_i^{(k_i)}]\}.$$ Since $\mathcal{A}$ is malnormal, and since, for every $i \in \{1,\ldots,r\}$, the group $A_i^{(1)} \ast \ldots \ast A_i^{(k_i)}$ is a free factor of $A_i$, it follows that the subgroup system $\mathcal{F} \wedge \mathcal{A}$ is a malnormal subgroup system of $F_{\tt n}$. We call it the \emph{meet of $\mathcal{F}$ and $\mathcal{A}$}.

\subsection{Graphs, markings and filtrations}

Let ${\tt n} \geq 2$. A \emph{marked graph} is a pointed (at a vertex $*$), connected, finite graph $G$ (in the sense of~\cite{Serre83}) whose fundamental group is isomorphic to $F_{\tt n}$ which is equipped with a \emph{marking}, that is an isomorphism $\rho \colon F_{\tt n} \to \pi_1(G,\ast)$.

We denote by $VG$ (resp. $\vec{E}G$) the set of vertices (resp. edges) of $G$. Given an edge $e$ of $G$, we denote by $o(e)$ the \emph{origin} of $e$, by $t(e)$ the \emph{terminal point} of $e$ and by $e^{-1}$ the edge of $G$ such that $o(e^{-1})=t(e)$ and $t(e^{-1})=o(e)$. An \emph{edge path $\gamma$ of length $m$} is a concatenation of $m$ edges $\gamma=e_1e_2\ldots e_m$ such that for every $i \in \{1,\ldots,m-1\}$, we have $t(e_i)=o(e_{i+1})$. The length of $\gamma$ is denoted by $\ell(\gamma)$. The edge path $\gamma$ is \emph{reduced} if for every $i \in \{1,\ldots,m-1\}$, we have $e_i \neq e_{i+1}^{-1}$. A reduced edge path is \emph{cyclically reduced} if $t(e_m)=o(e_1)$ and $e_m \neq e_1^{-1}$. A cyclically reduced edge path is also called a \emph{circuit}. For any edge path $\gamma$, there exists a unique reduced edge path homotopic to $\gamma$ relatively to endpoints, we denote it by $[\gamma]$.

Let $G$ and $G'$ be two marked graphs. A \emph{graph map} is a pointed homotopy equivalence $f \colon G \to G'$ such that $f(VG) \subseteq VG'$ and such that the restriction of $f$ to the interior of an edge is an immersion. Thus, for every edge $e \in \vec{E}G$, the image $f(e)$ determines a reduced edge path $[f(e)]$. Given $\phi \in \Out(F_{\tt n})$ and $(G,\rho)$ a marked graph, a \emph{topological representative} of $\phi$ is a graph map $f \colon G \to G$ such that the outer automorphism class of $\rho^{-1} \circ f_{\ast} \circ \rho \in \Aut(F_{\tt n})$ is $\phi$.

Let $f \colon G \to G$ be a topological representative. Let $w \in F_{\tt n}$. We denote by $\gamma_w$ the unique circuit in $G$ which represents the conjugacy class of $w$.

Let $f \colon G \to G$ be a topological representative. A \emph{filtration} for $G$ is an increasing sequence of $f$-invariant (not necessarily connected) subgraphs $\varnothing=G_0 \subsetneq G_1 \subsetneq \ldots \subsetneq G_k=G$. Let $r \in \{1,\ldots,k\}$. The \emph{$r$-th stratum} in this filtration, denoted by $H_r$ is the (not necessarily connected) closure of $G_r-G_{r-1}$. For every $r \in \{1,\ldots,k\}$, there exists a square matrix $M_r$ associated with the stratum $H_r$ called the \emph{transition matrix} of $H_r$. The rows and columns of $M_r$ are indexed by the nonoriented edges in $H_r$ and the entry associated with the pair of nonoriented edges defined by $(e,e') \in \left(EH_r\right)^2$ is the number of occurrences of $e'$ and $e'^{-1}$ in $[f(e)]$.

Recall that a nonnegative square matrix $M=(M_{i,j})_{i,j}$ is \emph{irreducible} if for every $(i,j)$, there exists $p=p(i,j)$ such that $M_{i,j}^p>0$ and that $M$ is \emph{primitive} if there exists $p \in \NN^*$ such that every entry of $M^p$ is positive. For $r \in \{1,\ldots,k\}$, we say that the stratum $H_r$ is \emph{irreducible} if its associated matrix is irreducible and we say that $H_r$ is \emph{primitive} if its associated matrix is primitive. Let $r \in \{1,\ldots,k\}$ and suppose that $M_r$ is irreducible. Then it has a unique real eigenvalue $\lambda_r \geq 1$ called the \emph{Perron-Frobenius} eigenvalue. Let $H_r$ be an irreducible stratum. Then $H_r$ is \emph{exponentially growing (EG)} if $\lambda_r >1$ and is \emph{nonexponentially growing (NEG)} otherwise. Finally, if the matrix associated with the stratum $H_r$ is the zero matrix, then $H_r$ is called a \emph{zero stratum}.  

Let $G$ be a marked graph of $F_{\tt n}$ and let $K$ be a (possibly disconnected) subgraph of $G$. The subgraph $K$ determines a free factor system $\mathcal{F}(K)$ of $F_{\tt n}$ as follows. Let $C_1,\ldots,C_k$ be the noncontractible connected components of $K$. Then, for every $i \in \{1,\ldots,k\}$, the connected component $C_i$ determines the conjugacy class $[A_i]$ of a subgroup $A_i$ of $\pi_1(G)$. Then the set $\{[A_1],\ldots,[A_k]\}$ is a free factor system $\mathcal{F}(K)$ of $F_{\tt n}$. 

Let $\mathcal{F}_1 \leq \ldots \leq \mathcal{F}_i=\{[F_{\tt n}]\}$ be a filtration of $F_{\tt n}$. A \emph{geometric realization of the filtration} is a marked graph $G$ equipped with an increasing sequence $$\varnothing=G_0 \subsetneq G_1 \subsetneq \ldots \subsetneq G_{j}=G$$ of subgraphs of $G$ such that for every $k \in \{1,\ldots,i\}$ there exists \mbox{$\ell \in \{1,\ldots,j\}$} such that $\mathcal{F}_k=\mathcal{F}(G_{\ell})$.

\subsection{Train tracks and CTs}

In this section we introduce the technology of \emph{train tracks}. Train tracks are a type of graph maps introduced by Bestvina and Handel (\cite{BesHan92}). Even though there exist outer automorphisms of $F_{\tt n}$ which do not have a topological representative which is a train track, every outer automorphism has a power which has a topological representative called a
\emph{completely split train track map} (CT). CT maps were introduced by Feighn and Handel (\cite{FeiHan06}). The definition of a CT map being quite technical, we will only state the relevant properties needed for the rest of the article. First we need some preliminary definitions.

Let $G$ be a marked graph of $F_{\tt n}$ and let $f \colon G \to G$ be a graph map. The map $f$ induces a \emph{derivative map} $Df \colon \vec{E}G \to \vec{E}G$ on the set of edges as follows. For every $e \in \vec{E}G$, the map $Df(e)$ is equal to the first edge of the edge path $f(e)$. A \emph{turn} in $G$ is an unordered pair $\{e_1,e_2\}$ of edges in $G$ with $o(e_1)=o(e_2)$. A turn $\{e_1,e_2\}$ is \emph{degenerate} if $e_1=e_2$, and is \emph{nondegenerate} otherwise. A turn $\{e_1,e_2\}$ is \emph{illegal} if there exists $k \in \NN^*$ such that $\{(Df)^k(e_1),(Df)^k(e_2)\}$ is degenerate, and is \emph{legal} otherwise. An edge path $\gamma=e_1e_2\ldots e_i$ is \emph{legal} if for every $j \in \{1,\ldots,i\}$, the turn $\{e_j^{-1},e_{j+1}\}$ is legal.

In order to deal with relative outer automorphisms, we also need a notion of relative legal paths. Let $\varnothing=G_0 \subsetneq G_1 \subsetneq \ldots \subsetneq G_{j}=G$ be the geometric realization of some filtration of $F_{\tt n}$ which is $f$-invariant and let $r \in \{1,\ldots,j\}$. We say that a turn $\{e_1,e_2\}$ is contained in the stratum $H_r$ if $\{e_1,e_2\} \subseteq \vec{E}H_r$. An edge path $\gamma$ of $G$ is \emph{$r$-legal} if every turn in $\gamma$ that is contained in $H_r$ is legal. A \emph{connecting path} for $H_r$ is a nontrivial reduced path $\gamma$ in $G_{r-1}$ whose endpoints are in $G_{r-1} \cap H_r$. A path $\gamma$ in $G$ is $r$-\emph{taken} (or \emph{taken} if $\gamma$ is $r$-taken for some $r$) if it is contained in the reduced image of an iterate of an edge $e \in \vec{E}H_r$, where $H_r$ is an irreducible stratum. The \emph{height of a path $\gamma$} is the maximal $r$ such that $\gamma$ contains an edge of $H_r$. We can now define the notion of a \emph{relative train track map} due to Bestvina and Handel (\cite{BesHan92}).

\begin{defi}\label{Defi relative train track}
Let ${\tt n} \geq 3$. Let $G$ be a marked graph and let $f \colon G \to G$ be a graph map equipped with a $f$-invariant filtration $\varnothing=G_0 \subsetneq G_1 \subsetneq \ldots \subsetneq G_{j}=G$. The map $f$ is a \emph{relative train track map} if, for each exponentially growing stratum $H_r$, the following holds:

\medskip

\noindent{$(1)$ } for every edge $e \in \vec{E}H_r$ and every $k \in \NN^*$, we have $(Df)^k(e) \in \vec{E}H_r$;

\medskip

\noindent{$(2)$ } for every connecting path $\gamma$ for $H_r$, the reduced path $[f(\gamma)]$ is also a connecting path for $H_r$;

\medskip

\noindent{$(3)$ } if $\gamma$ is a height $r$ reduced edge path which is $r$-legal, then so is $[f(\gamma)]$.
\end{defi}

In order to explain the properties of CT maps that we will use in this paper, we will need some further definitions regarding edge paths in a graph.

\begin{defi}
Let ${\tt n} \geq 3$ and let $G$ be a marked graph of $F_{\tt n}$ equipped with an $f$-invariant filtration $\varnothing=G_0 \subsetneq G_1 \subsetneq \ldots \subsetneq G_{j}=G$. Let $\gamma$ be an edge path of $G$.

\medskip

\noindent{$(1)$ } The path $\gamma$ is a \emph{periodic Nielsen path} if there exists $k \in \NN^*$ such that $[f^k(\gamma)]=\gamma$. The minimal such $k$ is the \emph{period}, and if $k=1$, then $\gamma$ is a \emph{Nielsen path}.

\medskip

\noindent{$(2)$ } A \emph{(periodic) indivisible Nielsen path ((p)INP)} is a (periodic) Nielsen path that cannot be written as a nontrivial concatenation of (periodic) Nielsen paths.

\medskip

\noindent{$(3)$ } The path $\gamma$ is an \emph{exceptional path} if there exist a cyclically reduced Nielsen path $w$, edges $e_1,e_2 \in \vec{E}G$ and integers $d_1,d_2,p \in \ZZ^*$ such that for every $i \in \{1,2\}$, we have $f(e_i)=e_iw^{d_i}$ and $\gamma=e_1w^pe_2^{-1}$. The value $|p|$ is called the \emph{width of $\gamma$}.
\end{defi}

\begin{defi}
Let $n \geq 3$, let $G$ be a marked graph of $F_{\tt n}$ and let $f \colon G \to G$ be a relative train track map equipped with a filtration $\varnothing=G_0 \subsetneq G_1 \subsetneq \ldots \subsetneq G_{j}=G$. Let $\gamma$ be a reduced edge path or a circuit of $G$. 

\medskip

\noindent{$(1)$ } A \emph{splitting} of $\gamma$ is a decomposition of $\gamma$ into edge subpaths $\gamma=\gamma_1\gamma_2\ldots\gamma_i$ such that for every $k \in \NN^*$, we have $$[f^k(\gamma)]=[f^k(\gamma_1)]\ldots[f^k(\gamma_i)],$$ that is one can tighten the image of $f^k(\gamma)$ by tightening the image of every $f^k(\gamma_j)$ (where $o(\gamma)$ is the base point in the case where $\gamma$ is a circuit).

\medskip

\noindent{$(2)$ } Let $\gamma$ be a circuit. A \emph{circuital splitting} is a splitting $\gamma=\gamma_1\ldots\gamma_i$ of $\gamma$ such that for every $k \in \NN^*$, the concatenation $[f^k(\gamma_1)]\ldots[f^k(\gamma_i)]$ defines a path whose initial and terminal directions are distinct.

\medskip

\noindent{$(3)$ } Let $\gamma=\gamma_1\gamma_2\ldots\gamma_i$ be a splitting of $\gamma$. The splitting is \emph{complete} if for every \mbox{$j \in \{1,\ldots,i\}$}, the subpath $\gamma_j$ is one of the following:

\medskip

\noindent{$\bullet$} an edge in an irreducible stratum;

\medskip

\noindent{$\bullet$} an INP;

\medskip

\noindent{$\bullet$} an exceptional path;

\medskip

\noindent{$\bullet$} a connecting path in a zero stratum that is both maximal (for the inclusion in $\gamma$) and taken.

\medskip

\end{defi}

Let ${\tt n} \geq 2$, let $G$ be a marked graph of $F_{\tt n}$ and let $f \colon G \to G$ be a relative train track map with respect to a filtration $\varnothing=G_0 \subsetneq G_1 \subsetneq \ldots \subsetneq G_{j}=G$. Let $\gamma$ be an edge path of $G$. Such paths in the above list are called \emph{splitting units}. When $\gamma$ has a complete splitting, we say that $\gamma$ is \emph{completely split}. 

\begin{defi}\cite[Fact~2.16]{HandelMosher20}
Let $p \in \{0,\ldots,j\}$. Let $\gamma=\gamma_1\gamma_2\ldots\gamma_i$ be a splitting of $\gamma$. This splitting is \emph{complete relatively to $G_p$}, or \emph{relatively complete} if there is no ambiguity, if for every \mbox{$j \in \{1,\ldots,i\}$}, the subpath $\gamma_j$ is one of the following:

\medskip

\noindent{$\bullet$} a splitting unit of height at least equal to $p+1$;

\medskip

\noindent{$\bullet$ } a subpath in $G_p$.
\end{defi}

We now describe some properties of CT maps whose complete definition can be found in~\cite[Definition~4.7]{FeiHan06}.

\begin{prop}\label{Prop definition CT}
Let ${\tt n} \geq 3$ and let $G$ be a marked graph of $F_{\tt n}$. Let $f \colon G \to G$ be a completely split train track (CT) map. Then $f$ satisfies the following properties.

\medskip

\noindent{$(1)$ } The map $f$ is a relative train track map and every stratum in $G$ is either irreducible or a zero stratum (\cite[Definition~4.7]{FeiHan06}).

\medskip

\noindent{$(2)$ } If $H_r$ is an NEG stratum, then $H_r$ consists of a single edge $e_r$. Moreover, either $e_r$ is fixed by $f$ or $f(e_r)=e_ru_r$ where $u_r$ is a nontrivial completely split circuit in $G_{r-1}$. The terminal endpoint of each NEG stratum is fixed (\cite[Lemma~4.21]{FeiHan06}).

\medskip

\noindent{$(3)$ } For every filtration element $G_r$, the stratum $H_r$ is a zero stratum if and only if $H_r$ is a contractible component of $G_r$ (\cite[Lemma~4.15]{FeiHan06}).

\medskip

\noindent{$(4)$ } For every zero stratum $H_r$, there exists a unique $\ell>r$ such that $H_{\ell}$ is an EG stratum and, for every vertex $v \in VH_r$, we have $v \in VH_r \cap VH_{\ell}$ and the link of $v$ is contained in $VH_r \cup VH_{\ell}$ (\cite[Definition~4.7]{FeiHan06}).

\medskip

\noindent{$(5)$ } Every periodic Nielsen path has period one (\cite[Lemma~4.13]{FeiHan06}).

\medskip

\noindent{$(6)$ } For every edge $e$ in an irreducible stratum, the reduced path $f(e)$ is completely split. For every taken connecting path $\gamma$ in a zero stratum, $[f(\gamma)]$ is completely split.

\medskip

\noindent{$(7)$ } Every completely split path or circuit has a unique complete splitting.

\medskip

\noindent{$(8)$ } If $\gamma$ is an edge path, there exists $k_0 \in \NN^*$ such that for every $k \geq k_0$, the reduced path $[f^k(\gamma)]$ is completely split (\cite[Lemma~4.25]{FeiHan06}).

\medskip

\noindent{$(9)$ } If $H_r$ is an EG stratum, there is at most one INP $\rho_r$ of height $r$. The initial edges of $\rho_r$ and $\rho_r^{-1}$ are distinct oriented edges in $H_r$ (\cite[Corollary~4.19]{FeiHan06}).

\medskip

\noindent{$(10)$ } If $H_r$ is a zero stratum, no Nielsen path intersects $H_r$ in at least one edge (\cite[Fact~I.1.43]{HandelMosher20}).

\medskip

\noindent{$(11)$ } Let $H_r$ be an NEG stratum such that $H_r=\{e_r\}$, such that $f(e_r)=e_ru_r$ and such that $u_r$ is not trivial. There exists an INP $\sigma$ which intersects $H_r$ nontrivially if and only if $u_r$ is a Nielsen path and there exists $s \in \ZZ$ such that $\sigma=e_ru_r^se_r^{-1}$ (\cite[Definition~4.7]{FeiHan06}).
\end{prop}

\begin{defi}
Let ${\tt n} \geq 2$ and let $G$ be a marked graph of $F_{\tt n}$. Let $f \colon G \to G$ be a completely split train track (CT) map. Let $H_r$ be an NEG stratum and let $e_r$ be the edge of $H_r$. Let $u_r$ be such that $f(e_r)=e_ru_r$. The edge $e_r$ is called a \emph{fixed edge} if $u_r$ is trivial, a \emph{linear edge} if $u_r$ is a Nielsen path and a \emph{superlinear edge} otherwise.
\end{defi}

\begin{lem}\cite[Fact~1.39]{HandelMosher20}\label{Lem Nielsen path Ct concatenation INP and fixed edges}
Let ${\tt n} \geq 2$ and let $G$ be a marked graph of $F_{\tt n}$. Let $f \colon G \to G$ be a CT map. Let $\gamma$ be a Nielsen path. Then $\gamma$ is completely split, and all terms in the complete splitting of $\gamma$ are fixed edges and INPs.
\end{lem}

\begin{lem}\cite[Fact~1.41]{HandelMosher20}\label{Lem INP in EG stratum impossible with zero stratum}
Let ${\tt n} \geq 2$ and let $G$ be a marked graph of $F_{\tt n}$. Let $f \colon G \to G$ be a CT map.

\medskip

\noindent{$(1)$ } Let $H_r$ be a zero stratum and let $H_{\ell}$ be the EG stratum given by Proposition~\ref{Prop definition CT}~$(4)$. There does not exist an INP of height $\ell$.

\medskip

\noindent{$(2)$ } Let $H_r$ be an EG stratum and let $\rho_r$ be an INP of height $r$. Then $\rho_r$ has a decomposition $\rho_r=a_0b_1a_1\ldots b_ka_k$ where, for every $i \in \{0,\ldots,k\}$, the subpath $a_i$ is a nontrivial path contained in $H_r$ and for every $i \in \{1,\ldots,k\}$, the subpath $b_i$ is a Nielsen path contained in $G_{r-1}$. 
\end{lem}

An INP is an \emph{EG INP} if the maximal stratum it intersects is an EG stratum and is an \emph{NEG INP} otherwise. Note that, by Proposition~\ref{Prop definition CT}~$(9)$, there exists only finitely many EG INPs.

\begin{lem}\label{Lem No zero path and Nielsen path adjacent}
Let ${\tt n} \geq 2$. Let $\phi \in \Out(F_{\tt n})$. Suppose that there exists a CT map $f \colon G \to G$ representing a power of $\phi$. Let $\gamma'$ be a nontrivial path in a zero stratum. There does not exist a reduced edge path $\gamma=\alpha\gamma'$ where $\alpha$ is either an INP or a fixed edge. 
\end{lem}

\dem Suppose towards a contradiction that such a path $\gamma=\alpha\gamma'$ exists. Let $H_r$ be the zero stratum containing $\gamma'$. Note that, by Proposition~\ref{Prop definition CT}~$(10)$, the path $\alpha$ does not contain edges in $H_r$. By Proposition~\ref{Prop definition CT}~$(4)$, there exists $\ell >r$ such that $H_{\ell}$ is an EG stratum and such that any edge adjacent to a vertex in $H_r$ and not contained in $H_r$ is in $H_{\ell}$. Hence $\alpha$ has height at least $\ell$. Since $H_{\ell}$ is an EG stratum, the path $\alpha$ is not a fixed edge. Hence $\alpha$ is an INP. By Lemma~\ref{Lem INP in EG stratum impossible with zero stratum}~$(1)$, the height of $\alpha$ is not equal to $\ell$. Let $j>\ell$ be the height of $\alpha$. We distinguish between three cases according to the nature of the stratum $H_j$. By Proposition~\ref{Prop definition CT}~$(10)$, the stratum $H_j$ is not a zero stratum. Hence, by Proposition~\ref{Prop definition CT}~$(1)$, the stratum $H_j$ is irreducible. By Proposition~\ref{Prop definition CT}~$(11)$, if $H_j$ is an NEG stratum, then $\alpha$ is of the form $\alpha=e_jw^ke_j^{-1}$, where $e_j \in H_j$, $k$ is an integer and $w$ is a closed Nielsen path in $G_{j-1}$. But then $e_j^{-1}$ is adjacent to a vertex in $H_r$. This contradicts Proposition~\ref{Prop definition CT}~$(4)$ since $j>\ell$. If $H_j$ is an EG stratum, then by Lemma~\ref{Lem INP in EG stratum impossible with zero stratum}~$(2)$, the path $\alpha$ is the concatenation of subpaths in $H_j$ and Nielsen paths of height at most $j-1$, and $\alpha$ ends with an edge in $H_j$. By Proposition~\ref{Prop definition CT}~$(4)$, we see that $j=\ell$. This contradicts Lemma~\ref{Lem INP in EG stratum impossible with zero stratum}~$(1)$.
\hfill\qedsymbol

\bigskip

The next theorem due to Feighn and Handel is the main existence theorem of the CT maps.

\begin{theo}\cite[Theorem~4.28, Lemma~4.42]{FeiHan06}\label{Theo existence CT}
Let ${\tt n} \geq 3$. There exists a uniform constant $M=M(n) \geq 1$ such that for every $\phi \in \Out(F_{\tt n})$ and every $\phi^M$-invariant filtration $\mathcal{C}$ of $F_{\tt n}$, there exists a CT map $f \colon G \to G$ that represents $\phi^M$ and realizes $\mathcal{C}$.
\end{theo}

\subsection{Relative currents}\label{Section relative currents}

In this section, we define the notion of \emph{currents of $F_{\tt n}$ relative to a malnormal subgroup system}. The section follows~~\cite{Guerch2021currents} (see the work of Gupta~\cite{gupta2017relative} for the particular case of free factor systems and Guirardel and Horbez \cite{Guirardelhorbez19laminations} in the context of free products of groups). It is closely related to the notion of conjugacy classes of $\mathcal{A}$-nonperipheral elements of $F_{\tt n}$. 

Let $\partial_{\infty}F_{\tt n}$ be the Gromov boundary of $F_{\tt n}$. The \emph{double boundary of $F_{\tt n}$} is the quotient topological space $$\partial^2F_{\tt n}=\left(\partial_{\infty} F_{\tt n} \times \partial_{\infty} F_{\tt n} \setminus \Delta \right)/\sim,$$ where $\sim$ is the equivalence relation generated by the flip relation $(x,y)\sim(y,x)$ and $\Delta$ is the diagonal, endowed with the diagonal action of $F_{\tt n}$. We denote by $\{x,y\}$ the equivalence class of $(x,y)$.

Let $T$ be the Cayley graph of $F_{\tt n}$ with respect to a free basis $\mathfrak{B}$. The boundary of $T$ is naturally homeomorphic to $\partial_{\infty}F_{\tt n}$ and the set $\partial^2F_{\tt n}$ is then identified with the set of unoriented bi-infinite geodesics in $T$. Let $\gamma$ be a finite geodesic path in $T$. The path $\gamma$ determines a subset in $\partial^2F_{\tt n}$ called the \emph{cylinder set of $\gamma$}, denoted by $C(\gamma)$, which consists in all unoriented bi-infinite geodesics in $T$ that contain $\gamma$. Such cylinder sets form a basis for a topology on $\partial^2 F_{\tt n}$, and in this topology, the cylinder sets are both open and closed, hence compact. The action of $F_{\tt n}$ on $\partial^2F_{\tt n}$ has a dense orbit.

For every nontrivial subgroup $A$ of $F_{\tt n}$, let $T_A$ be the minimal $A$-invariant subtree of $T$. Let $\mathcal{A}=\{[A_1],\ldots,[A_r]\}$ be a malnormal subgroup system of $F_{\tt n}$. By malnormality of $\mathcal{A}$, there exists $L \in \NN^*$ such that for all distinct subgroups $A,B$ of $F_{\tt n}$ such that $[A],[B] \in \mathcal{A}$, the diameter of the intersection $T_{A} \cap T_{B}$ is at most $L$ (see for instance~\cite[Section~I.1.1.2]{HandelMosher20}). Let $i \in \{1,\ldots,r\}$. Let $\Gamma_i$ be the set of subgroups $B$ of $F_{\tt n}$ such that there exists $g_B \in F_{\tt n}$ such that $B=g_BA_ig_B^{-1}$ and the tree $T_B$ contains the base point $e$ of $T$. Note that, by malnormality of $\mathcal{A}$, for every $i \in \{1,\ldots,r\}$, the set $\Gamma_i$ is finite. For an element $w \in F_{\tt n}$, let $\widehat{\gamma_w}$ be the geodesic path in $T$ starting at $e$ and labeled by $w$. Let $C_i$ be the set of elements $w$ of $F_{\tt n}$ such that the length of $\widehat{\gamma_w}$ is equal to $L+2$ and, for every $B \in \Gamma_i$, the path $\widehat{\gamma_w}$ is not contained in $T_{B}$. Let $\mathscr{C}=\bigcap_{i=1}^r C_i$. Since we are looking at geodesic paths of length equal to $L+2$, the set $\mathscr{C}$ is finite. Moreover, it only depends on the choice of $\mathcal{A}$, $\mathfrak{B}$ and $L$.

\begin{lem}\cite[Lemma~2.3]{Guerch2021currents}\label{Lem finite set of words determines vertex group system}
Let $\mathfrak{B}$, $T$, $\mathcal{A}=\{[A_1],\ldots,[A_r]\}$, $L \in \NN^*$, $\Gamma_1,\ldots,\Gamma_r$, $\mathscr{C}$ be as above. The finite set $\mathscr{C}=\mathscr{C}(A_1,\ldots, A_k)$ is nonempty. Moreover, it satisfies the following properties:

\medskip

\noindent{$(1)$ } every $\mathcal{A}$-nonperipheral cyclically reduced element $g \in F_{\tt n}$ has a power which contains an element of $\mathscr{C}$ as a subword; 

\medskip

\noindent{$(2)$ } for every $\mathcal{A}$-nonperipheral cyclically reduced element $g \in F_{\tt n}$, if $c_g$ is the geodesic ray in $T$ starting from $e$ obtained by concatenating infinitely many edge paths labeled by $g$, there exists an edge path in $c_g$ labeled by a word in $\mathscr{C}$ at distance at most $L+2$ from $\bigcup_{i=1}^r \bigcup_{B \in \Gamma_i} T_B$;

\medskip

\noindent{$(3)$ } if $\gamma$ is a path in $T$ which contains a subpath labeled by an element of $\mathscr{C}$, then for every $i \in \{1,\ldots,r\}$ and every $g \in F_{\tt n}$, the path $\gamma$ is not contained in $T_{gA_ig^{-1}}$.

\end{lem}

Let $A$ be a nontrivial subgroup of $F_{\tt n}$ of finite rank. The induced $A$-equivariant inclusion $\partial_{\infty} A \hookrightarrow \partial_{\infty} F_{\tt n}$ induces an inclusion $\partial^2 A \hookrightarrow \partial^2 F_{\tt n}$. Let $$\partial^2\mathcal{A}= \bigcup_{i=1}^r \bigcup_{g \in F_{\tt n}} \partial^2 \left(gA_ig^{-1}\right).$$ Let $\partial^2(F_{\tt n},\mathcal{A})=\partial^2F_{\tt n} -\partial^2\mathcal{A}$ be the \emph{double boundary of $F_{\tt n}$ relative to $\mathcal{A}$}. This subset is invariant under the action of $F_{\tt n}$ on $\partial^2F_{\tt n}$ and inherits the subspace topology of $\partial^2F_{\tt n}$.

\begin{lem}\cite[Lemma~2.5]{Guerch2021currents}\label{Lem double boundary open}
Let $\mathrm{Cyl}(\mathscr{C})$ be the set of cylinder sets of the form $C(\gamma)$, where the element of $F_{\tt n}$ determined by the geodesic edge path $\gamma$ contains an element of $\mathscr{C}$ as a subword. We have $$\partial^2(F_{\tt n},\mathcal{A})=\bigcup_{C(\gamma) \in \mathrm{Cyl}(\mathscr{C})}C(\gamma).$$ In particular, the space $\partial^2(F_{\tt n},\mathcal{A})$ is an open subset of $\partial^2 F_{\tt n}$.
\end{lem}

\begin{lem}\cite[Lemma~2.6, Lemma~2.7]{Guerch2021currents}\label{Lem Properties of relative boundary}
Let ${\tt n} \geq 3$ and let $\mathcal{A}$ be a malnormal subgroup system of $F_{\tt n}$. The space $\partial^2(F_{\tt n},\mathcal{A})$ is locally compact and the action of $F_{\tt n}$ on $\partial^2(F_{\tt n},\mathcal{A})$ has a dense orbit.
\end{lem}

We can now define a \emph{relative current}. Let ${\tt n} \geq 3$ and let $\mathcal{A}$ be a malnormal subgroup system of $F_{\tt n}$. A \emph{relative current on $(F_{\tt n},\mathcal{A})$} is a (possibly zero) $F_{\tt n}$-invariant Radon measure $\mu$ on $\partial^2(F_{\tt n},\mathcal{A})$. The set $\Curr(F_{\tt n},\mathcal{A})$ of all 
relative currents on $(F_{\tt n},\mathcal{A})$ is equipped with the weak-$\ast$ topology: a sequence $(\mu_n)_{n \in \NN}$ in $\Curr(F_{\tt n},\mathcal{A})^{\NN}$ converges to a current \mbox{$\mu \in \Curr(F_{\tt n},\mathcal{A})$} if and only if for all disjoint clopen subsets $S,S' \subseteq \partial^2(F_{\tt n},\mathcal{A})$, the sequence $(\mu_n(S \times S'))_{n \in \NN}$ converges to $\mu(S \times S')$.

The group $\Out(F_{\tt n},\mathcal{A})$ acts on $\Curr(F_{\tt n},\mathcal{A})$ as follows. Let $\phi \in \Out(F_{\tt n},\mathcal{A})$, let $\Phi$ be a representative of $\phi$, let $\mu \in \Curr(F_{\tt n},\mathcal{A})$ and let $C$ be a Borel subset of $\partial^2(F_{\tt n},\mathcal{A})$. Then, since $\phi$ preserves $\mathcal{A}$, we see that $\Phi^{-1}(C) \in \partial^2(F_{\tt n},\mathcal{A})$. Then we set $$\phi(\mu)(C)=\mu(\Phi^{-1}(C)),$$ which is well-defined since $\mu$ is $F_{\tt n}$-invariant.

\bigskip

Every conjugacy class of nonperipheral element $g \in F_{\tt n}$ determines a relative current $\eta_{[g]}$ as follows. Suppose first that $g$ is \emph{root-free}, that is $g$ is not a proper power of any element in $F_{\tt n}$. Let $\gamma$ be a finite geodesic path in the Cayley graph $T$. Then $\eta_{[g]}(C(\gamma))$ is the number of axes in $T$ of conjugates of $g$ that contain the path $\gamma$. If $g=h^k$ with $k \geq 2$ and $h$ root-free, we set $\eta_{[g]}=k \;\eta_{[h]}$. Such currents are called \emph{rational currents}.

\bigskip

Let $G$ be a pointed connected graph whose fundamental group is isomorphic to $F_{\tt n}$. Let $\widetilde{G}$ be the universal cover of $G$. There exists a (nonunique, but fixed) $F_{\tt n}$-equivariant quasi-isometry $\widetilde{m} \colon \widetilde G \to T$ which extends uniquely to a homeomorphism $\widehat{m} \colon \partial_{\infty} G \to \partial_{\infty} F_{\tt n}$. Therefore, if $\widetilde{\gamma}$ is a reduced edge path in $\widetilde{G}$, we can define the cylinder set in $\partial^2 F_{\tt n}$ defined by $\widetilde{\gamma}$ as 

$$C_{\widetilde{m}}(\widetilde{\gamma})=C([\widetilde{m}(\widetilde{\gamma})]).$$

Let $\gamma$ be a reduced edge path in $G$ and let $\widetilde{\gamma}$ be a lift of $\gamma$ in $\widetilde{G}$. Let $\mu \in \Curr(F_{\tt n},\mathcal{A})$. We define the \emph{number of occurrences of $\gamma$ in $\mu$} as 
\begin{equation}\label{Equation defi evaluation current}
\left\langle \gamma,\mu \right\rangle_{\widetilde{m}}=\mu(C_{\widetilde{m}}(\widetilde{\gamma})).
\end{equation}
For every such graph $G$, we fix once and for all the quasi-isometry $\widetilde{m}\colon \widetilde{G} \to T$. Therefore, when the graph $G$ is fixed, we will generally omit the mention of $\widetilde{m}$. We also define the \emph{simplicial length of $\mu$} as: $$ \lVert \mu \rVert=\sum_{e \in \vec{E}G} \left\langle e,\mu \right\rangle.$$ For any given reduced edge path $\gamma$, the functions $\left\langle \gamma,. \right\rangle$ and $\lVert . \rVert$ are continuous, linear functions of $\Curr(F_{\tt n},\mathcal{A})$. 

Let $\mu \in \Curr(F_{\tt n},\mathcal{A})$. The \emph{support of $\mu$}, denoted by $\Supp(\mu)$, is the support of the Borel measure $\mu$ on $\partial^2(F_{\tt n},\mathcal{A})$. We recall that $\Supp(\mu)$ is a closed subset of $\partial^2(F_{\tt n},\mathcal{A})$.

In the rest of the article, rather than considering the space of relative currents itself, we will consider the set of \emph{projectivized relative currents}, denoted by $$\mathbb{P}\Curr(F_{\tt n},\mathcal{A})=(\Curr(F_{\tt n},\mathcal{A})-\{0\})/\sim,$$ where $\mu \sim \nu$ if there exists $\lambda \in \RR_+^*$ such that $\mu=\lambda \nu$. The projective class of a current $\mu \in \Curr(F_{\tt n},\mathcal{A})$ will be denoted by $[\mu]$. We have the following properties.

\begin{lem}\cite[Lemma~3.3]{Guerch2021currents}\label{Lem PCurr compact}
Let ${\tt n} \geq 3$ and let $\mathcal{A}$ be a malnormal subgroup system of $F_{\tt n}$. The space $\PCurr(F_{\tt n},\mathcal{A})$ is compact.
\end{lem}

\begin{prop}\cite[Theorem~1.1]{Guerch2021currents}\label{Prop density rational currents}
Let ${\tt n} \geq 3$ and let $\mathcal{A}$ be a malnormal subgroup system of $F_{\tt n}$. The set of projectivised rational currents about nonperipheral elements of $F_{\tt n}$ is dense in $\PCurr(F_{\tt n},\mathcal{A})$.
\end{prop}

\section{The polynomially growing subgraph of a CT map}\label{Section PG subgraph of a CT map}

In this section, let ${\tt n} \geq 3$ and let $\mathcal{F}$ be a free factor system of $F_{\tt n}$. Let $\phi \in \Out(F_{\tt n},\mathcal{F})$. Let $f \colon G \to G$ be a CT map with filtration $\varnothing=G_0 \subsetneq G_1 \subsetneq \ldots \subsetneq G_k=G$ representing a power of $\phi$ and such that there exists $p \in \{1,\ldots,k-1\}$ such that $\mathcal{F}(G_p)=\mathcal{F}$. 

We construct a subgraph of $G$, called the \emph{polynomially growing subgraph} of $G$ and denoted by $G_{PG}$, which encaptures the information regarding polynomial growth in the graph $G$. We then define a notion of length relative to $G_{PG}$, called the \emph{exponential length}, which measures the time spent by an edge path outside of $G_{PG}$. Finally, we construct a subspace of $\PCurr(F_{\tt n},\mathcal{F})$ which consists in the currents whose support maps to $G_{PG}$.

\subsection{Definitions and first properties}

We define in this section the \emph{polynomially growing subgraph} $G_{PG}$ of $G$ and proves some of its properties.

\begin{defi}\label{Defi Gpg}

\noindent{$(1)$ }  Let $G_{PG}$ be the (not necessarily connected) subgraph of $G$ whose edges are the edges $e$ of $G$ in an NEG stratum such that for every $k \in \NN^*$, the path $[f^k(e)]$ does not contain a splitting unit which is an edge in an EG stratum. 

\medskip

\noindent{$(2)$ } Let $\mathcal{N}_{PG}'$ be the set of all Nielsen paths in $G$.  

\medskip

\noindent{$(3)$ } Let $\mathcal{N}_{PG}$ be the subset of $\mathcal{N}_{PG}'$ consisting in all Nielsen paths which are either EG INPs or concatenations of (at least $2$) nonclosed EG INPs.

\medskip

\noindent{$(4)$ } Let $\mathcal{Z}$ be the subgraph of $G$ whose edges are the edges contained in a zero stratum.
\end{defi}

\begin{figure}[ht]
\centering
\begin{tikzpicture}
\draw (0,0) ellipse (2cm and 1cm);
\draw (7,0) ellipse (2cm and 1cm);
\draw (2,0) -- (5,0);
\draw (2,0) node {$\bullet$};
\draw (5,0) node {$\bullet$};
\draw (0,1) node[above] {$G_{PG}$};
\draw (7,1) node[above] {$G_{PG}$};
\draw (3.5,0) node [above] {$\gamma$};
\end{tikzpicture}
\caption{A path $\gamma$ in $\mathcal{N}_{PG}$ between two connected components of $G_{PG}$.}\label{Figure path in Npg}
\end{figure}
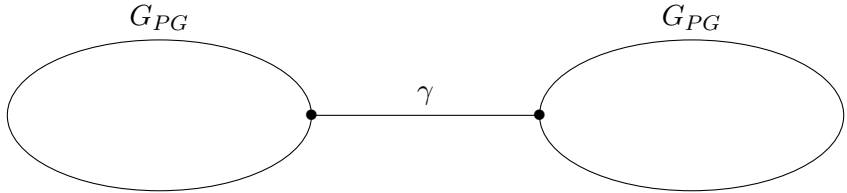

Note that, by Lemma~\ref{Lem Nielsen path Ct concatenation INP and fixed edges}, every path in $\mathcal{N}_{PG}'$ (and hence every path in $\mathcal{N}_{PG}$) has a complete splitting consisting in fixed edges and INPs. Since a complete splitting is unique by Proposition~\ref{Prop definition CT}~$(7)$, if $\gamma$ is a reduced path in $\mathcal{N}_{PG}$, then the splitting of $\gamma$ given in Definition~\ref{Defi Gpg}~$(3)$ is the complete splitting of $\gamma$. Moreover, $\gamma$ is either an EG INP or the complete splitting of $\gamma$ has at least two splitting units and all of them are nonclosed EG INPs. In particular, the set $\mathcal{N}_{PG}$ does not contain Nielsen paths such that one of their splitting units is either a fixed edge or an NEG INP. Moreover, a Nielsen path which is a concatenation of at least $2$ splitting units and such that one of them is a closed EG INP is not in $\mathcal{N}_{PG}$. Excluding such paths from $\mathcal{N}_{PG}$ ensures a finiteness result for $\mathcal{N}_{PG}$ (see~Lemma~\ref{Lem Nielsen paths in NPG properties}~$(1)$). Informally, paths in $\mathcal{N}_{PG}$ play the role of low-dynamics bridges between connected components of $G_{PG}$ (see~Figure~\ref{Figure path in Npg}). We will see in Proposition~\ref{Prop circuits in Gpg are elements in poly subgroup} that a cycle in $G$ has polynomial growth under iteration of $f$ if and only if is a concatenation of paths in $G_{PG}$ and paths in $\mathcal{N}_{PG}$.

Note that, with $p$ defined at the beginning of Section~\ref{Section PG subgraph of a CT map}, one can similarly define the \emph{polynomially growing subgraph of $G_p$}, denoted by $G_{PG,\mathcal{F}}$, which is the subgraph $G_{PG} \cap G_p$. We can also define similarly $\mathcal{N}_{PG,\mathcal{F}}'$, $\mathcal{N}_{PG,\mathcal{F}}$ and $\mathcal{Z}_{\mathcal{F}}$ by considering the paths of $\mathcal{N}_{PG}'$, $\mathcal{N}_{PG}$ and $\mathcal{Z}$ contained in $G_p$.

We now recall a lemma due to Bestvina and Handel regarding $r$-legal paths.

\begin{lem}\cite[Lemma~5.8]{BesHan92}\label{Lem r legal paths and concatenation}
Let $f \colon G \to G$ be a relative train track map. Let $H_r$ be an EG stratum. Suppose that $\sigma=a_1b_1a_2\ldots a_{\ell}b_{\ell}$ is the decomposition of an $r$-legal path into subpaths $a_j \subseteq H_r$ and $b_j \subseteq G_{r-1}$ (where $a_1$ and $b_{\ell}$ might be trivial). Then for every $i \in \{1,\ldots,\ell\}$, the path $f(a_{\ell})$ is a reduced edge path and $$[f(\sigma)]=f(a_1)[f(b_1)]f(a_2)\ldots f(a_{\ell})[f(b_{\ell})].$$
\end{lem}

Note that, if $H_r$ is an EG stratum and if $\sigma=a_1b_1a_2\ldots a_{\ell}b_{\ell}$ is an $r$-legal path as in Lemma~\ref{Lem r legal paths and concatenation}, then for every $i \in \{1,\ldots,\ell\}$, as $a_{i} \subseteq H_r$, the path $a_i$ grows exponentially fast under iteration of $f$. Hence, by Lemma~\ref{Lem r legal paths and concatenation} the path $\sigma$ grows exponentially fast under iteration of $f$. We now prove some results regarding paths in $\mathcal{N}_{PG}$.

\begin{lem}\label{Lem EG INP is not cdc}
Let $\sigma$ be an EG INP. 

\medskip

\noindent{$(1)$ } There do not exist nontrivial subpaths $c,d$ of $\sigma$ such that $\sigma=cdc$.

\medskip

\noindent{$(2)$ } Let $\gamma \in \{\sigma^{\pm 1}\}$. There do not exist paths $\gamma_1,\gamma_2,\gamma_3$ such that $\gamma_2$ is nontrivial, $\gamma_1$ or $\gamma_3$ is nontrivial and $\sigma=\gamma_1\gamma_2$ and $\gamma=\gamma_2\gamma_3$.
\end{lem}

\dem $(1)$ Let $r$ be the height of $\sigma$. Suppose towards a contradiction that such a decomposition $\sigma=cdc$ exists. By~\cite[Lemma~5.11]{BesHan92}, there exist two distinct $r$-legal paths $\alpha$ and $\beta$ such that $\sigma=\alpha\beta$ and such that the turn $\{Df(\alpha^{-1}),Df(\beta)\}$ is the only height $r$ illegal turn. Moreover, there exists a path $\tau$ such that $[f(\alpha)]=\alpha\tau$ and $[f(\beta)]=\tau^{-1}\beta$. Hence $c$ is contained in $\alpha$ and in $\beta$ and is $r$-legal. Thus, there exist two paths $d_1$ and $d_2$ such that $\alpha=cd_1$ and $\beta=d_2c$.

First we claim that for every $k \in \NN^*$, there exists a path $\tau_k$ such that $[f^k(\alpha)]=\alpha \tau_k$ and $[f^k(\beta)]=\tau_k^{-1}\beta$. The proof is by induction on $k$. The base case follows from the existence of $\tau$. Suppose now that $\tau_{k-1}$ exists. We have: $$[f^k(\alpha)]=[f(\alpha\tau_{k-1})]=[f(\alpha)][f(\tau_{k-1})]=\alpha\tau[f(\tau_{k-1})]=\alpha \tau_k,$$ where the second equality comes from the fact that $\alpha$ is $r$-legal, that $\alpha$ ends with an edge in $H_r$ and from Lemma~\ref{Lem r legal paths and concatenation}. Similarly, we have $[f^k(\beta)]=\tau_k^{-1}\beta$. This proves the claim.

We now claim that, up to taking a power of $f$, there exists a cycle $e$ such that $[f(c)]=\alpha e \beta$. Indeed, by Proposition~\ref{Prop definition CT}~$(9)$, the path $\sigma$ starts and ends with an edge in $H_r$. Hence the path $c$ starts and ends with an edge in $H_r$. Since $c$ is $r$-legal, we see that the length of $[f^k(c)]$ goes to infinity as $k$ goes to infinity by Lemma~\ref{Lem r legal paths and concatenation}. But, for every $k \in \NN^*$, there exists a path $\tau_k$ such that $[f^k(\alpha)]=\alpha \tau_k$ and $[f^k(\beta)]=\tau_k^{-1}\beta$. By Lemma~\ref{Lem r legal paths and concatenation}, since $c$ is the initial segment of $\alpha$ and since $\alpha$ is $r$-legal, there is no identification between $[f(c)]$ and $[f(d_1)]$. Thus, there exists $k_1 \in \NN^*$ such that $[f^{k_1}(c)]$ starts with $\alpha$. Similarly, there exists $k_2 \in \NN^*$ such that $[f^{k_2}(c)]$ ends with $\beta$. Thus, up to taking a power of $f$, and since the paths $\alpha$ and $\beta$ are $r$-legal, we may suppose that there exists a (reduced) cycle $e$ such that $[f(c)]=\alpha e \beta$.

Finally, we claim that the cycle $e$ is trivial. Indeed, since the paths $\alpha$ and $\beta$ are $r$-legal, and since $c$ starts and ends with an edge in $H_r$, we see that $$[f(\alpha)]=[f(c)][f(d_1)]=\alpha e \beta [f(d_1)]$$ and $$[f(\beta)]=[f(d_2)][f(c)]=[f(d_2)]\alpha e \beta.$$ Recall that there exists $k \in \NN^*$ such that $[f(\alpha)]=\alpha\tau_k$ and $[f(\beta)]=\tau_k^{-1}\beta$. This implies that $\tau_k=e \beta [f(d_1)]$ and that $\tau_k^{-1}=[f(d_2)]\alpha e$, that is $\tau_k=e^{-1}\alpha^{-1} [f(d_2)]^{-1}$. This shows that $e=e^{-1}$, that is, $e$ is trivial. This proves the claim.

Therefore, we see that $[f(c)]=\alpha\beta=\sigma$. But $\sigma$ contains a height $r$ illegal turn, whereas $c$ is an $r$-legal path. This contradicts Proposition~\ref{Prop definition CT}~$(1)$ and Definition~\ref{Defi relative train track}~$(3)$. This concludes the proof of $(1)$.

\bigskip

\noindent{$(2)$ } Let $\sigma,\gamma$ be as in the assertion of the lemma. Suppose towards a contradiction that there exist three paths $\gamma_1,\gamma_2,\gamma_3$ such that $\gamma_2$ is nontrivial and $\sigma=\gamma_1\gamma_2$ and $\gamma=\gamma_2\gamma_3$. Suppose first that $\gamma=\sigma$. Then either a nontrivial initial segment of $\gamma_2$ is its terminal segment or there exists a path $\gamma_4$ such that $\sigma=\gamma_2\gamma_4\gamma_2$. The first case is not possible as otherwise $\sigma$ would contain two illegal turns. This contradicts the fact that $\sigma$ contains a unique illegal turn (see~\cite[Lemma~5.11]{BesHan92}). The second case is not possible by Lemma~\ref{Lem EG INP is not cdc}~$(1)$. Suppose now that $\gamma=\sigma^{-1}$. But $\sigma^{-1}=\gamma_2^{-1}\gamma_1^{-1}$. Therefore we see that $\gamma_2^{-1}=\gamma_2$, that is, $\gamma_2$ is trivial. This leads to a contradiction. This concludes the proof.
\hfill\qedsymbol

\begin{lem}\label{Lem Nielsen paths in NPG properties}
\noindent{$(1)$ } There are only finitely many paths in $\mathcal{N}_{PG}$.

\medskip

\noindent{$(2)$ } Let $\gamma,\gamma'$ be paths in $\mathcal{N}_{PG}$. Suppose that $\gamma$ has a decomposition $\gamma=\gamma_1\gamma_2$ such that $\gamma_2$ is an initial segment of $\gamma'$. Then $\gamma_1,\gamma_2 \in \mathcal{N}_{PG}$ and $\gamma_1\gamma' \in \mathcal{N}_{PG}$.

\medskip

\noindent{$(3)$ } Let $\gamma,\gamma'$ be paths in $\mathcal{N}_{PG}$. Suppose that $\gamma' \subseteq \gamma$. Then one of the following holds:

{$(a)$ } there exist (possibly trivial) paths $\gamma_1, \gamma_2 \in \mathcal{N}_{PG}$ such that $\gamma=\gamma_1\gamma'\gamma_2$;

{$(b)$ } there exists an INP $\sigma$ in the complete splitting of $\gamma$ such that $\gamma' \subsetneq \sigma$ and $\gamma'$ is not an initial or a terminal segment of $\sigma$.

\medskip

\noindent{$(4)$ } Let $\gamma,\gamma'$ be two paths in $\mathcal{N}_{PG}$. Suppose that there exist three paths $\gamma_1$, $\gamma_2$ and $\gamma_3$ such that $\gamma=\gamma_1\gamma_2$, $\gamma'=\gamma_2^{-1}\gamma_3$ and the path $\gamma_1\gamma_3$ is reduced. Then $\gamma_2 \in \mathcal{N}_{PG}$ and $\gamma_1\gamma_3 \in \mathcal{N}_{PG}$.
\end{lem}

\dem $(1)$ First note that, since there are only finitely many EG strata in $G$, there are only finitely many EG INPs by Proposition~\ref{Prop definition CT}~$(9)$. Let $\gamma$ be a path in $\mathcal{N}_{PG}$ which is a concatenation of at least $2$ nonclosed EG INPs. Let $\gamma=\sigma_1\ldots\sigma_k$ be the complete splitting of $\gamma$ given by Lemma~\ref{Lem Nielsen path Ct concatenation INP and fixed edges}. As $\gamma$ is a concatenation of nonclosed EG INPs, every splitting unit of $\gamma$ is a nonclosed EG INP. By Proposition~\ref{Prop definition CT}~$(9)$, an INP contained in the complete splitting of $\gamma$ is entirely determined by the highest stratum $H_r$ such that $\gamma$ contains an edge of $H_r$. For every $i \in \{1,\ldots,k\}$, let $r_i$ be the height of $\sigma_i$. Let $i \in \{2,\ldots,k\}$. Since $\sigma_i$ is not closed, by \cite[Fact~1.42(1)(a)]{HandelMosher20}, one of the endpoints of $\sigma_i$ is not contained in $G_{r_i-1}$. Since there exists a unique INP of height $r_i$ by Proposition~\ref{Prop definition CT}~$(9)$, either $r_{i-1}<r_i$ or $r_i<r_{i-1}$. We treat the case $r_1<r_2$, the case $r_2<r_1$ being similar. We claim that, for every $i \in \{1,\ldots,k-1\}$, we have $r_{i+1} >r_i$. The proof is by induction on $i$. The base case is true by hypothesis. Let $i \in \{2,\ldots,k-1\}$. Since $r_{i-1}<r_i$, the origin of $\sigma_i$ is contained in $G_{r_i-1}$ and the terminal point of $\sigma_i$ is not contained in $G_{r_i-1}$. Thus, the first edge of $\sigma_{i+1}$ is contained in $\overline{G-G_{r_i-1}}$. Since there exists a unique INP of height $r_i$ we necessarily have $r_i <r_{i+1}$. Thus, the sequence of maximal heights of INPs in $\gamma$ is (strictly) monotonic. Since there are only finitely many EG strata, there are only finitely many paths in $\mathcal{N}_{PG}$. This concludes the proof of $(1)$.

\bigskip

\noindent{$(2)$ } Let $\gamma' \in \mathcal{N}_{PG}$ and let $\gamma=\gamma_1\gamma_2$ be as in the assertion of the lemma. We claim that $\gamma_2 \in \mathcal{N}_{PG}$ and that the splitting units of $\gamma_2$ are splitting units of both $\gamma$ and $\gamma'$. This will conclude the proof of Assertion~$(2)$ because $\gamma_1$ will be a concatenation of splitting units of $\gamma$, that is, it will be either an EG INP or a concatenation of nonclosed EG INPs (cf Definition~\ref{Defi Gpg}~$(3)$). Hence we will have $\gamma_1 \in \mathcal{N}_{PG}$ and $\gamma_1\gamma' \in \mathcal{N}_{PG}$. We show that $\gamma_2$ is a concatenation of INPs which are splitting units of $\gamma'$. A similar proof will show that the splitting units of $\gamma_2$ will also be splitting units of $\gamma$. Indeed, the path $\gamma'$ has a splitting $\gamma'=\sigma_1'\sigma_2'\ldots\sigma_k'$ which consists in EG INPs. Let $r'$ be the height of $\sigma_1'$. By Proposition~\ref{Prop definition CT}~$(9)$, there exists a unique unoriented INP of height $r'$ and this INP starts and ends with an edge in $H_{r'}$. Let $\sigma$ be the INP of $\gamma$ which has a decomposition $\sigma=\sigma_1\sigma_2$, where $\sigma_2$ is a nontrivial initial segment of $\gamma'$. As every splitting unit of $\gamma$ is an EG INP, so is $\sigma$. Let $r$ be the height of $\sigma$. Since the first edge of $\sigma_1'$ is of height $r'$, we cannot have $r'>r$. If $r=r'$, then by the uniqueness statement in Proposition~\ref{Prop definition CT}~$(9)$, we have $\sigma_1' \in \{\sigma,\sigma^{-1}\}$. Note that, if $\sigma_1$ is nontrivial, there exist reduced paths $\tau_1,\tau_2$ such that $\sigma=\sigma_1\tau_1$ and $\sigma_1'=\tau_1\tau_2$. This contradicts Lemma~\ref{Lem EG INP is not cdc}~$(2)$ applied to $\sigma$ and $\sigma_1'$. Thus, we see that $\sigma=\sigma_1'$ and $\sigma_1' \subseteq \gamma_2$. If $r'<r$, then by Lemma~\ref{Lem INP in EG stratum impossible with zero stratum}~$(2)$, the path $\sigma$ has a decomposition $\sigma=a_1b_1\ldots b_{k-1}a_k$ such that, for every $i \in \{1,\ldots,k\}$, the path $a_i$ is a path in $H_r$ and for every $i \in \{1,\ldots,k-1\}$, the path $b_i$ is a Nielsen path in $G_{r-1}$. Hence there exists $i \in \{1,\ldots,k-1\}$ such that $\sigma_1'$ is contained in $b_i$. Therefore, we see that $\sigma_1' \subseteq \sigma \subseteq \gamma$. As $\sigma_1' \subseteq \gamma'$, we see that $\sigma_1' \subseteq \gamma \cap  \gamma'=\gamma_2$. If $\gamma_2=\sigma_1'$, then we are done. Otherwise, the path $\gamma_2$ contains an edge of $\sigma_2'$. As $\sigma_2'$ is an EG INP, the same argument as for $\sigma_1'$ shows that $\sigma_2' \subseteq \gamma_2$, and an inductive argument shows that $\gamma_2$ is a concatenation of INPs in the splitting of $\gamma'$. Hence $\gamma_2$ is a Nielsen path. Therefore, we see that $\gamma_2 \in \mathcal{N}_{PG}$ and that $\gamma_2$ is composed of splitting units of $\gamma'$. Similarly, we see that $\gamma_2$ is composed of splitting units which are splitting units of both $\gamma$ and $\gamma'$. Hence $\gamma_1$ is composed of splitting units of $\gamma$. This concludes the proof of $(2)$.

\medskip

\noindent{$(3)$ } Let $\gamma$, $\gamma'$ be as in the assertion of the lemma. Let $\gamma=\sigma_1\ldots\sigma_k$ be the complete splitting of $\gamma$ and let $\gamma'=\sigma_1'\ldots\sigma_m'$ be the complete splitting of $\gamma'$, which exist by Lemma~\ref{Lem Nielsen path Ct concatenation INP and fixed edges}. Recall that every splitting unit of both $\gamma$ and $\gamma'$ is an EG INP. There exists $i \in \{1,\ldots,k\}$ such that $\sigma_i$ contains an initial segment of $\sigma_1'$. We claim that $\sigma_1'$ is either equal to $\sigma_i$ or $\gamma'$ is strictly contained in $\sigma_i$. Indeed, let $r$ be the height of $\sigma_i$ and let $r'$ be the height of $\sigma_1'$.  Since the first edge of $\sigma_1'$ is of height $r'$, we cannot have $r'>r$. 

Suppose first that $r'<r$. By Lemma~\ref{Lem INP in EG stratum impossible with zero stratum}~$(2)$, the path $\sigma_i$ has a decomposition $\sigma_i=a_1b_1\ldots b_{p-1}a_p$ such that, for every $i \in \{1,\ldots,p\}$, the path $a_i$ is a path in $H_r$ and for every $j \in \{1,\ldots,p-1\}$, the path $b_j$ is a Nielsen path in $G_{r-1}$. Hence there exists $j \in \{1,\ldots,p-1\}$ such that $\sigma_1'$ is contained in $b_j$. We claim that, for every $\ell \in \{1,\ldots,m\}$, the splitting unit $\sigma_{\ell}'$ is contained in $b_j$. The proof is by induction on $\ell$. For the base case, we already know that $\sigma_1' \subseteq b_j$. Suppose that for some $\ell \in \{2,\ldots,m\}$, the path $\sigma_{\ell-1}'$ is contained in $b_j$. By Proposition~\ref{Prop definition CT}~$(9)$, the path $\sigma_i$ ends with an edge in $H_r$. Hence the path $a_p$ is nontrivial. Since $\sigma_{\ell-1}'$ is contained in $b_j$, the path $\sigma_{\ell}'$ intersects $\sigma_i$ nontrivially. Let $r_{\ell}$ be the height of $\sigma_{\ell}'$. Recall that $\sigma_{\ell}'$ is an EG INP. By Proposition~\ref{Prop definition CT}~$(9)$, the path $\sigma_{\ell}'$ starts with an edge in $H_{r_{\ell}}$. Hence $r_{\ell} \leq r$. Suppose towards a contradiction that $r_{\ell}=r$. Then, by the uniqueness statement of Proposition~\ref{Prop definition CT}~$(9)$, we see that $\sigma_{\ell}' \in \{\sigma_i^{\pm 1}\}$. As $\sigma_i$ contains an initial segment of $\sigma_{\ell}'$, there exist three paths $\gamma_1$, $\gamma_2$ and $\gamma_3$ of $G$ such that $\gamma_2$ is nontrivial and $\sigma_i=\gamma_1\gamma_2$ and $\sigma_{\ell}'=\gamma_2\gamma_3$. Since $\sigma_{\ell-1}'$ is contained in $\sigma_i$, the path $\gamma_1$ is nontrivial. This contradicts Lemma~\ref{Lem EG INP is not cdc}~$(2)$. Therefore we have $r_{\ell}<r$. But then $\sigma_{\ell}'$ cannot intersect $a_{j+1}$. This implies that $\sigma_{\ell}'$ is contained in $b_j$. This proves the claim and the fact that $\gamma' \subsetneq \sigma_i$ and $\gamma'$ is not an initial or a terminal segment of $\sigma_i$.

Suppose now that $r=r'$. By the uniqueness statement of Proposition~\ref{Prop definition CT}~$(9)$, we see that $\sigma_{1}' \in \{\sigma_i^{\pm 1}\}$. As $\sigma_i$ contains an initial segment of $\sigma_1'$, there exist three paths $\gamma_1$, $\gamma_2$ and $\gamma_3$ of $G$ such that $\gamma_2$ is nontrivial and $\sigma_i=\gamma_1\gamma_2$ and $\sigma_{1}'=\gamma_2\gamma_3$. By Lemma~\ref{Lem EG INP is not cdc}~$(2)$, we necessarily have that $\gamma_1$ and $\gamma_3$ are trivial. Thus, we see that $\sigma_i=\sigma_1'$. Therefore, $\gamma'$ is an initial segment of $\sigma_i\ldots\sigma_k$ and is a Nielsen path. By~\cite[Corollary~4.12]{FeiHan06}, for every $j \in \{1,\ldots,m\}$, we have $\sigma_{i+j-1}=\sigma_j'$. Thus,
there exist (possibly trivial) paths $\gamma_1, \gamma_2 \in \mathcal{N}_{PG}$ such that $\gamma=\gamma_1\gamma'\gamma_2$. This concludes the proof of $(3)$.

\bigskip

\noindent{$(4)$ } Let $\gamma$, $\gamma'$, $\gamma_1$, $\gamma_2$ and $\gamma_3$ be as in the assertion of the lemma. Let $\gamma=\alpha_1\ldots\alpha_{k}$ and $\gamma'=\beta_1\ldots\beta_{\ell}$ be the complete splittings of $\gamma$ and $\gamma'$ given by Lemma~\ref{Lem Nielsen path Ct concatenation INP and fixed edges}. By definition of $\mathcal{N}_{PG}$, every splitting unit of $\gamma$ and $\gamma'$ is an EG INP. Let $i \in \{1,\ldots,k\}$ be such that $\alpha_i$ contains the first edge of $\gamma_2$. Let $j \in \{1,\ldots,\ell\}$ be such that $\beta_j$ contains the last edge of $\gamma_2^{-1}$. We claim that $\alpha_i \subseteq \gamma_2$ and that $\beta_j \subseteq \gamma_2^{-1}$. By~\cite[Corollary~4.12]{FeiHan06} applied to $\gamma_2^{-1}$ and $\gamma^{-1}$, there exists a path $\delta_i$ contained in $\alpha_i$ such that the decomposition $\gamma_2=\delta_i\alpha_{i+1}\ldots\alpha_k$ is a splitting of $\gamma_2$. Similarly, there exists a path $\delta_j'$ in $\beta_j$ such that $\gamma_2^{-1}=\beta_1\ldots\beta_{j-1}\delta_j'$ is a splitting of $\gamma_2^{-1}$. By Proposition~\ref{Prop definition CT}~$(9)$, an EG INP starts with an edge of highest height and an EG INP is entirely determined by its height. Hence $\alpha_k=\beta_1^{-1}$. Note that the paths $\delta_i\alpha_{i+1}\ldots\alpha_{k-1}$ and $\beta_2\ldots\beta_{j-1}\delta_j'$ satisfy the same hypotheses as $\delta_i\alpha_{i+1}\ldots\alpha_k$ and $\beta_1\ldots\beta_{j-1}\delta_j'$. Applying the same arguments, we see that $i=j$ and for every $s \in \{1,\ldots,j-1\}$, we have $\beta_s=\alpha_{k-s+1}^{-1}$. Hence we see that $\delta_i=\delta_j'^{-1}$. Let $r$ be the height of $\alpha_i$ and let $r'$ be the height of $\beta_j$. Note that by Proposition~\ref{Prop definition CT}~$(9)$ applied to $\alpha_i$ and $\beta_j$, the path $\delta_i$ ends with an edge in $H_{r}$ and $\delta_j'^{-1}$ ends with an edge in $H_{r'}$. Therefore, we see that $r=r'$. By uniqueness of EG INPs of height $r_i$ given by Proposition~\ref{Prop definition CT}~$(9)$, and since $\gamma_1\gamma_3$ is reduced, we see that $\alpha_i=\beta_j^{-1}$, that $\alpha_i \subseteq \gamma_2$ and that $\beta_j \subseteq \gamma_2^{-1}$. This shows that $\gamma_2$ is a path in $\mathcal{N}_{PG}$. By Assertion~$(2)$ applied to $\gamma$ and $\gamma_2$, the path $\gamma_1$ is contained in $\mathcal{N}_{PG}$. Similarly, we see that the path $\gamma_3$ is contained in $\mathcal{N}_{PG}$. Since the path $\gamma_1\gamma_3$ is reduced, we see that $\gamma_1\gamma_3 \in \mathcal{N}_{PG}$. This concludes the proof.

\hfill\qedsymbol

\begin{lem}\label{Lem concatenation paths Gpg Npg}
Let $\gamma$ and $\gamma'$ be two reduced edge paths in $G$ which are concatenations of paths in $G_{PG}$ and $\mathcal{N}_{PG}$. Suppose that there exist three paths $\gamma_1$, $\gamma_2$ and $\gamma_3$ such that $\gamma=\gamma_1\gamma_2$, $\gamma'=\gamma_2^{-1}\gamma_3$ and $\gamma_1\gamma_3$ is reduced. Then $\gamma_2$ and $\gamma_1\gamma_3$ are concatenations of paths in $G_{PG}$ and $\mathcal{N}_{PG}$.
\end{lem}

\dem Let $\gamma=b_0a_1b_1\ldots a_kb_k$ be the decomposition of the path $\gamma$ such that for every $i \in \{0,\ldots,k\}$, the path $b_i$ is in $G_{PG}$ and for every $i \in \{1,\ldots,k\}$, the path $a_i$ is a maximal subpath of $\gamma$ contained in $\mathcal{N}_{PG}$. The existence of the paths $a_i$ follows from Lemma~\ref{Lem Nielsen paths in NPG properties}~$(2)$. Let $\gamma'=d_0c_1d_1\ldots c_{\ell}d_{\ell}$ be the similar decomposition of $\gamma'$. Let $e$ be the initial edge of $\gamma_2$.

\medskip

\noindent{\bf Claim. } There exists $i \in \{0,\ldots,k\}$ such that $b_i$ contains $e$ if and only if there exists $j \in \{0,\ldots,\ell\}$ such that the edge $e^{-1}$ is contained in $d_j$. 

\medskip

\dem The proof of the two directions being similar, we only prove one direction. Suppose that there exists $i \in \{0,\ldots,k\}$ such that $b_i$ contains $e$. Suppose towards a contradiction that there exists $j \in \{1,\ldots,\ell\}$ such that $e^{-1}$ is contained in $c_j$. It follows that there exists an EG INP $\sigma$ of $c_j$ such that $e^{-1}$ is contained in $\sigma$. Let $r$ be the height of $\sigma$. Let $\delta^{-1}$ be the subpath of $\sigma$ contained in $\gamma_2^{-1}$. Note that, as $\gamma_2^{-1}$ is an initial segment of $\gamma'$, the path $\delta^{-1}$ is an initial segment of $\sigma$. By Proposition~\ref{Prop definition CT}~$(9)$, the path $\delta^{-1}$ starts with an edge in $H_r$. As $\delta$ is contained in $\gamma$, the terminal edge of $\delta$ is an edge in an EG stratum. Since every edge in $G_{PG}$ is contained in an NEG stratum, there exists $s \in \{1,\ldots,k\}$ such that $a_s$ contains a terminal segment of $\delta$. Since the initial edge $e$ of $\gamma_2$ is not contained in $a_s$ by hypothesis, the path $\delta$ contains the initial segment $\delta'$ of $a_s$. Hence the terminal segment $\delta'^{-1}$ of $a_s^{-1}$ is the initial segment $\delta'^{-1}$ of $\sigma$. By Lemma~\ref{Lem Nielsen paths in NPG properties}~$(2)$ applied to $a_s^{-1}$ and $\sigma$ and~\cite[Corollary~4.12]{FeiHan06}, the path $\delta'^{-1}$ is contained in $\mathcal{N}_{PG}$ and is a concatenation of splitting units of $\sigma$. As $\sigma$ contains a unique splitting unit, this implies that $\delta'=\sigma$. As $\delta' \subseteq \delta^{-1} \subseteq \sigma$, we see that $\delta^{-1}=\sigma$. Note that the edge $\delta^{-1}$ ends with $e^{-1}$. But $\sigma$ ends with an edge in an EG stratum by Proposition~\ref{Prop definition CT}~$(9)$, that is, $e^{-1}$ is an edge in an EG stratum. But every edge in $b_i$ is contained in an NEG stratum by definition of $G_{PG}$. This contradicts the fact that $e \subseteq b_i$. This concludes the proof of the claim.
\hfill\qedsymbol

\medskip

Suppose first that there exists $i \in \{1,\ldots,k\}$, such that $e$ is contained in $b_i$. By the above claim, there exists $j \in \{0,\ldots,\ell\}$ such that $e^{-1}$ is contained in $d_j$. Let $\tau$ and $\tau'$ be such that $\gamma=b_0a_1b_1\ldots a_i\tau\gamma_2$ and $\gamma'=\gamma_2^{-1}\tau'c_{j+1}\ldots d_{\ell}$. Note that $\tau \subseteq b_i$ and $\tau' \subseteq d_j$. Then we have $\gamma_1=b_0a_1b_1\ldots a_i\tau$ and $\gamma_3=\tau'c_{j+1}\ldots d_{\ell}$. Since the path $\gamma_1\gamma_3$ is reduced, so is $\tau\tau'$. Moreover the reduced edge path $\tau\tau'$ is contained in $G_{PG}$ and $\gamma_1\gamma_3=b_0a_1b_1\ldots a_i\tau\tau'c_{j+1}\ldots d_{\ell}$ is a concatenation of paths in $G_{PG}$ and in $\mathcal{N}_{PG}$. Moreover, let $\delta''$ be the maximal subpath of $b_i$ contained in $\gamma_2$. Then $\gamma_2=\delta''a_{i+1}\ldots b_k$ is a concatenation of paths in $G_{PG}$ and in $\mathcal{N}_{PG}$.

Suppose now that there exists $i \in \{1,\ldots,k\}$ such that the initial edge $e$ of $\gamma_2$ is contained in $a_i$. By the above claim, there exists $j \in \{1,\ldots,\ell\}$ such that $e^{-1}$ is contained in $c_j$. Let $\delta'$ be the terminal segment of $a_i$ contained in $\gamma_2$. By Proposition~\ref{Prop definition CT}~$(9)$, the terminal edge $e'$ of $\delta'$ is an edge in an EG stratum. Since $G_{PG}$ does not contain any edge in an EG stratum, there exists $s \leq j$ such that $c_s$ contains $e'^{-1}$. We claim that $s=j$. Indeed, suppose towards a contradiction that $s<j$. Let $\delta^{-1}$ be the terminal segment of $c_s$ whose first edge is $e'^{-1}$. Then $\delta$ is a terminal segment of $a_i$ and $\delta$ is an initial segment of $c_s^{-1}$. By Lemma~\ref{Lem Nielsen paths in NPG properties}~$(2)$ applied to $a_i$ and $c_s^{-1}$, the path $\delta$ is a concatenation of splitting units of $a_i$ and $c_s^{-1}$. If $\delta$ is properly contained in $\delta'$, there exists an EG INP $\sigma$ which is a splitting unit of $a_i$ and such that the last edge of $\sigma$ is the last edge of $\delta'$ not contained in $\delta$. But, by Proposition~\ref{Prop definition CT}~$(9)$, the terminal edge $e_{\sigma}$ of $\sigma$ is in an EG stratum. However, the first edge of $d_s$ (which is the edge $e_{\sigma}^{-1}$) is in $G_{PG}$. This leads to a contradiction. Hence $\delta=\delta'$. But $\delta$ intersects $c_j$ nontrivially. Hence we have $s=j$.

Therefore, $\delta'^{-1}$ is contained in $c_j$. We claim that $\delta'^{-1}$ is an initial segment of $c_j$. Indeed, otherwise let $\epsilon'$ be the initial segment of $c_j$ whose endpoint is the origin of $\delta'^{-1}$. By Proposition~\ref{Prop definition CT}~$(9)$, the first edge of $\epsilon'$ is an edge in an EG stratum. Hence there exists $p >i$ such that $a_p$ contains the terminal edge of $\epsilon'^{-1}$. Let $\epsilon^{-1}$ be the subpath of $\epsilon'^{-1}$ contained in $a_p$. Then $\epsilon^{-1}$ is an initial segment of $a_p$ and $\epsilon$ is an initial segment of $c_j$. By Lemma~\ref{Lem Nielsen paths in NPG properties}~$(2)$ applied to $a_p^{-1}$ and $c_j$, the path $\epsilon$ is a concatenation of splitting units of $a_p^{-1}$ and $c_j$. But since $\epsilon$ is properly contained in $c_j$ as it does not intersect $\delta'^{-1}$, the path $\epsilon$ is adjacent to a splitting unit of $c_j$. Since an EG INP starts with an edge in an EG stratum by Proposition~\ref{Prop definition CT}~$(9)$, the path $b_{p-1}$ ends with an edge in an EG stratum. This contradicts the fact that $b_{p-1}$ is contained in $G_{PG}$.

Hence $\delta'^{-1}$ is an initial segment of $c_j$ and $\delta'$ is a terminal segment of $a_i$. Let $\tau$ and $\tau'$ be two paths such that $a_i=\tau \delta'$ and $c_j=\delta'^{-1}\tau'$. By Lemma~\ref{Lem Nielsen paths in NPG properties}~$(4)$ applied to $a_i$ and $c_j$, the path $\delta'$ is in $\mathcal{N}_{PG}$ and the path $\tau\tau'$ is in $\mathcal{N}_{PG}$. Hence $\gamma_2=\tau b_ia_{i+1}\ldots b_k$ and $\gamma_1\gamma_3=b_0a_1b_1\ldots a_i\tau\tau'c_{j+1}\ldots d_{\ell}$ are concatenations of paths in $G_{PG}$ and in $\mathcal{N}_{PG}$. This concludes the proof.
\hfill\qedsymbol

\begin{lem}\label{Lem closed nielsen paths in Npg'}
Let $\gamma$ be a closed Nielsen path of $G$. Then $\gamma$ is a concatenation of paths in $G_{PG}$ and in $\mathcal{N}_{PG}$.
\end{lem}

\dem Let $\gamma$ be a closed Nielsen path of $G$. We prove the result by induction on the height $r$ of $\gamma$. If $r=0$, there is nothing to prove. Assume that $r \geq 1$. By Lemma~\ref{Lem Nielsen path Ct concatenation INP and fixed edges}, the path $\gamma$ is completely split, and every splitting unit in its complete splitting is either an INP or a fixed edge. Let $\gamma=\sigma_1\ldots\sigma_k$ be the complete splitting of $\gamma$. For every $i \in \{1,\ldots,k\}$, let $r_i$ be the height of $\sigma_i$. We prove that for every $i \in \{1,\ldots,k\}$, the path $\sigma_i$ is a concatenation of paths in $G_{PG}$ and in $\mathcal{N}_{PG}$. Let $i \in \{1,\ldots,k\}$. If $\sigma_i$ is a fixed edge, it is contained in $G_{PG}$. Suppose that $\sigma_i$ is an NEG INP. By Proposition~\ref{Prop definition CT}~$(11)$, there exists an edge $e_{r_i} \in EH_{r_i}$, a Nielsen path $w$ in $G_{r_i-1}$ and an integer $s \in \ZZ^*$ such that $\sigma_i=e_{r_i}w^se_{r_i}^{-1}$. Moreover, we have $f(e_{r_i})=e_{r_i}w$. Hence for every $j \in \NN^*$, we have $[f^j(e_{r_i})]=e_{r_i}w^j$. Since $w$ is a Nielsen path, by Lemma~\ref{Lem Nielsen path Ct concatenation INP and fixed edges}, the path $w$ is completely split and its complete splitting is made of fixed edges and INPs. Thus, for every $j \in \NN^*$, the complete splitting of $[f^j(e_{r_i})]$ does not contain splitting units which are edges in $EG$ strata. By definition of $G_{PG}$, we have $e_{r_i} \in \vec{E}G_{PG}$. Moreover, by the induction hypothesis, the path $w^s$ is a concatenation of paths in $G_{PG}$ and in $\mathcal{N}_{PG}$. Hence $\sigma_i$ is a concatenation of paths in $G_{PG}$ and in $\mathcal{N}_{PG}$. Finally, if $\sigma_i$ is an EG INP, then it is contained in $\mathcal{N}_{PG}$. Hence $\gamma$ is a concatenation of paths in $G_{PG}$ and in $\mathcal{N}_{PG}$.
\hfill\qedsymbol

\begin{lem}\label{Lem NEG INP in Npg}
Let $\gamma$ be either an NEG INP or an exceptional path. Then $\gamma$ is a concatenation of paths in $G_{PG}$ and in $\mathcal{N}_{PG}$.
\end{lem}

\dem We claim that there exist edges $e_1,e_2$ and a closed Nielsen path $w$ such that $\gamma=e_1we_2^{-1}$ and, for every $i \in \{1,2\}$, we have $f(e_i)=e_iw^{d_i}$ for some $d_i \in \ZZ^*$. If $\gamma$ is an exceptional path, it follows from the definition. If $\gamma$ is an NEG INP, let $r$ be the height of $\gamma$. Then $H_r$ is an NEG stratum. As $\gamma$ is a Nielsen path, we can apply Proposition~\ref{Prop definition CT}~$(11)$ to conclude the proof of the claim. Since $e_1$ and $e_2$ are linear edges, for every $k \in \NN^*$, the paths $[f^k(e_1)]$ and $[f^k(e_1)]$ do not contain splitting units which are edges in EG strata. Thus $e_1$ and $e_2$ are contained in $G_{PG}$. By Lemma~\ref{Lem closed nielsen paths in Npg'}, the path $w$ is a concatenation of paths in $G_{PG}$ and in $\mathcal{N}_{PG}$. Hence $\gamma$ is a concatenation of paths in $G_{PG}$ and in $\mathcal{N}_{PG}$. This concludes the proof.
\hfill\qedsymbol

\begin{lem}\label{Lem Nielsen paths in Npg}
Let $\gamma$ be a Nielsen path in $G$. Then $\gamma$ is a concatenation of paths in $G_{PG}$ and in $\mathcal{N}_{PG}$.
\end{lem}

\dem By Lemma~\ref{Lem Nielsen path Ct concatenation INP and fixed edges}, the path $\gamma$ is completely split, and every splitting unit in its complete splitting is either an INP or a fixed edge. Let $\gamma=\sigma_1\ldots\sigma_k$ be the complete splitting of $\gamma$. Let $i \in \{1,\ldots,k\}$. If $\sigma_i$ is a fixed edge, then $\sigma_i$ is contained in $G_{PG}$. If $\sigma_i$ is an NEG INP then, by Lemma~\ref{Lem NEG INP in Npg}, the path $\sigma_i$ is a concatenation of paths in $G_{PG}$ and in $\mathcal{N}_{PG}$. If $\sigma_i$ is an EG INP then, by definition, we have $\sigma_i \in \mathcal{N}_{PG}$. Hence $\gamma$ is a concatenation of paths in $G_{PG}$ and in $\mathcal{N}_{PG}$.
\hfill\qedsymbol

\begin{lem}\label{Lem iterate of a path in GPG}
\noindent{$(1)$ } Let $\gamma$ be an edge in $G_{PG}$ (resp. an edge in $G_{PG,\mathcal{F}}$). The path $[f(\gamma)]$ is a concatenation of paths in $G_{PG}$ and in $\mathcal{N}_{PG}$ (resp. a concatenation of paths in $G_{PG,\mathcal{F}}$ and in $\mathcal{N}_{PG,\mathcal{F}}$).

\medskip

\noindent{$(2)$ } Let $\gamma$ be an edge path contained in $G_{PG}$ (resp. an edge path in $G_{PG,\mathcal{F}}$). The path $[f(\gamma)]$ is a concatenation of paths in $G_{PG}$ and in $\mathcal{N}_{PG}$ (resp. a concatenation of paths in $G_{PG,\mathcal{F}}$ and in $\mathcal{N}_{PG,\mathcal{F}}$).

\medskip

\noindent{$(3)$ } Let $\gamma$ be an edge path which is a concatenation of paths in $G_{PG}$ and in $\mathcal{N}_{PG}$ (resp. a concatenation of paths in $G_{PG,\mathcal{F}}$ and in $\mathcal{N}_{PG,\mathcal{F}}$). The path $[f(\gamma)]$ is a concatenation of paths in $G_{PG}$ and in $\mathcal{N}_{PG}$ (resp. a concatenation of paths in $G_{PG,\mathcal{F}}$ and in $\mathcal{N}_{PG,\mathcal{F}}$). 
\end{lem}

\dem We prove Assertions~$(1)$, $(2)$, $(3)$ for paths in $G_{PG}$ and in $\mathcal{N}_{PG}$, the proofs for paths in $G_{PG,\mathcal{F}}$ and $\mathcal{N}_{PG,\mathcal{F}}$ being similar, using the fact that $f(G_p)=G_p$. 

\medskip

\noindent{$(1)$ } Let $\gamma$ be an edge of $G_{PG}$. By definition of $G_{PG}$, the edge $\gamma$ is an edge in an NEG stratum. By Proposition~\ref{Prop definition CT}~$(6)$, the path $[f(\gamma)]$ is completely split. Let $[f(\gamma)]=\gamma_1\ldots\gamma_m$ be the complete splitting of $[f(\gamma)]$. Since $\gamma$ is an edge in an NEG stratum, by Proposition~\ref{Prop definition CT}~$(2)$, we have $\gamma_1=\gamma$. Suppose towards a contradiction that $[f(\gamma)]$ is not a concatenation of paths in $G_{PG}$ and in $\mathcal{N}_{PG}$. It follows that there exists $i \in \{1,\ldots,m\}$ and an edge $e$ of $\gamma_i$ which is not contained in $G_{PG}$ and is not contained in a subpath of $[f(\gamma)]$ contained in $\mathcal{N}_{PG}$. Hence $\gamma_i$ is not an EG INP nor a fixed edge. By Lemma~\ref{Lem NEG INP in Npg}, the path $\gamma_i$ cannot be an NEG INP or an exceptional path. Hence $\gamma_i$ is either an edge in an irreducible stratum or a maximal taken connecting path in a zero stratum. Suppose first that $\gamma_i$ is a maximal taken connecting path in a zero stratum. By Proposition~\ref{Prop definition CT}~$(4)$, the path $\gamma_i$ cannot be adjacent to an edge in an NEG stratum nor an edge in a zero stratum. As $\gamma_1=\gamma$, we see that $i \geq 3$ and that $\gamma_{i-1}$ ends with an edge in an EG stratum. By Lemma~\ref{Lem No zero path and Nielsen path adjacent} (applied to $\gamma=\gamma_{i-1}\gamma_i$), the path $\gamma_{i-1}$ is not an EG INP. Therefore we see that $\gamma_{i-1}$ is an edge in an EG stratum. This contradicts the definition of the edges in $G_{PG}$. Hence we are reduced to the case where $\gamma_i$ is an edge in an irreducible stratum. Therefore, we have $\gamma_i=e$. By definition of $G_{PG}$ and as $e \notin \vec{E}G_{PG}$, there exists $k \in \NN^*$ such that $[f^k(\gamma_i)]$ contains a splitting unit which is an edge in an EG stratum. This contradicts the fact that $\gamma$ is contained in $G_{PG}$. This concludes the proof of $(1)$. 

\bigskip

\noindent{$(2)$ } Let $\gamma$ be a path in $G_{PG}$. We prove by induction on the length of $\gamma$ that $[f(\gamma)]$ is a concatenation of paths in $G_{PG}$ and in $\mathcal{N}_{PG}$. The case where $\gamma$ is an edge follows from $(1)$. Suppose now that the length of $\gamma$ is at least equal to $2$. Let $e$ be the last edge of $\gamma$ and let $\gamma'$ be an edge path such that $\gamma=\gamma' e$. Hence $\gamma'$ and $e$ are paths in $G_{PG}$. By the induction hypothesis, the paths $[f(\gamma')]$ and $[f(e)]$ are concatenations of paths in $G_{PG}$ and in $\mathcal{N}_{PG}$. It remains to show that identifications between $[f(\gamma')]$ and $[f(e)]$ do not create paths which are not concatenations of paths in $G_{PG}$ and in $\mathcal{N}_{PG}$. Let $\alpha$, $\beta$ and $\sigma$ be paths such that $[f(\gamma')]=\alpha\sigma$, $[f(e')]=\sigma^{-1}\beta$ and $\alpha\beta$ is reduced. By Lemma~\ref{Lem concatenation paths Gpg Npg} applied to $[f(\gamma')]$ and $[f(e')]$, the path $[f(\gamma)]$ is a concatenation of paths in $G_{PG}$ and in $\mathcal{N}_{PG}$. This concludes the proof of $(2)$.

\bigskip

\noindent{$(3)$ } Let $\gamma$ be a concatenation of paths in $G_{PG}$ and in $\mathcal{N}_{PG}$. Let $\gamma=\gamma_0'\gamma_1\gamma_1'\ldots\gamma_k\gamma_k'$ be a decomposition of $\gamma$ such that for every $i \in \{1,\ldots,k\}$, the path $\gamma_i$ is a maximal subpath of $\gamma$ in $\mathcal{N}_{PG}$ and for every $i \in \{0,\ldots,k\}$, the path $\gamma_i'$ is a path in $G_{PG}$. Such a decomposition is possible by Lemma~\ref{Lem Nielsen paths in NPG properties}~$(2)$. We prove the result by induction on $k$. If $k=0$, the proof follows from Assertion~$(2)$. Suppose that the result is true for $k'<k$. Then the paths $\gamma'=\gamma_0'\gamma_1\gamma_1'\ldots\gamma_{k-1}\gamma_{k-1}'$ and $\gamma''=\gamma_k \gamma_k'$ satisfy the induction hypothesis. Hence the paths $[f(\gamma')]$ and $[f(\gamma'')]$ are concatenations of paths in $G_{PG}$ and in $\mathcal{N}_{PG}$. Let $\alpha$, $\beta$ and $\sigma$ be three paths such that $[f(\gamma')]=\alpha\beta$, $[f(\gamma'')]=\beta^{-1}\sigma$ and $\alpha\beta$ is reduced. By Lemma~\ref{Lem concatenation paths Gpg Npg}, the path $[f(\gamma)]=\alpha\sigma$ is a concatenation of paths in $G_{PG}$ and in $\mathcal{N}_{PG}$. This concludes the proof.
\hfill\qedsymbol

\bigskip

For the next lemma, we recall a definition due to Bestvina, Feighn and Handel (\cite[Section~6]{BesFeiHan00}, see also~\cite[Definition~III.1.2]{HandelMosher20}). Let $H_{r_+}$ be the EG stratum of $G$ of maximal height $r_+$. By Proposition~\ref{Prop definition CT}~$(9)$, there exists at most one unoriented INP $\rho_{r_+}$ of height $r_+$ (we suppose that $\rho_{r_+}$ is a point if such a nontrivial INP does not exist). Following~\cite[Definition~III.1.2]{HandelMosher20}, let $Z_{r_+}$ be the subgraph of $G$ consisting in all edges $e'$ such that for every $m \in \NN^*$ and every splitting unit $\sigma$ of $[f^m(e')]$, the path $\sigma$ is not an edge in $H_{r_+}$. Let $\left\langle Z_{r_+}, \rho_{r_+} \right\rangle$ be the set consisting in the following paths:

\noindent{$(i)$ } paths in $Z_{r_+}$;

\noindent{$(ii)$ } paths in $\{\rho_{r_+},\rho_{r_+}^{-1}\}$;

\noindent{$(iii)$ } concatenations of paths in $Z_{r_+}$ and in $\{\rho_{r_+},\rho_{r_+}^{-1}\}$.

\noindent Note that $\left\langle Z_{r_+}, \rho_{r_+} \right\rangle$ contains every path in $G_{r_+-1}$.

\begin{lem}\label{Lem Gpg in Zr}
The set $\left\langle Z_{r_+}, \rho_{r_+} \right\rangle$ contains every path which is a concatenation of paths in $G_{PG}$ and in $\mathcal{N}_{PG}$.
\end{lem}

\dem It suffices to prove that $\left\langle Z_{r_+}, \rho_{r_+} \right\rangle$ contains every edge of $G_{PG}$ and every EG INP. Let $e$ be an edge in $G_{PG}$. By definition of $G_{PG}$, for every $k \in \NN^*$, the complete splitting of $[f^k(e)]$ does not contain a splitting unit which is an edge in an EG stratum. In particular, for every $k \in \NN^*$, the complete splitting of $[f^k(e)]$ does not contain a splitting unit which is an edge in $H_{r_+}$. Hence $e \subseteq Z_{r_+}$ and $G_{PG}$ is a subgraph of $Z_{r_+}$. Let $\rho$ be an EG INP and let $r$ be the height of $\rho$. By definition of $r_+$, we have $r \leq r_+$. If $r=r_+$, by Proposition~\ref{Prop definition CT}~$(9)$, we have $\rho \in \{\rho_{r_+},\rho_{r_+}^{-1}\}$, hence we have $\rho \in \left\langle Z_{r_+}, \rho_{r_+} \right\rangle$. If $r<r_+$, then $\rho$ is contained in $G_{r_+-1}$. Hence $\rho$ is contained in $\left\langle Z_{r_+}, \rho_{r_+} \right\rangle$ by the above remark.
\hfill\qedsymbol

\bigskip

We now define a graph which will be used in the proof of Lemma~\ref{Lem bijection on circuits Gpg}. Let $G^{\ast}$ be the finite, not necessarily connected, graph defined as follows:

\noindent{$(a)$ } vertices of $G^{\ast}$ are the vertices in $G_{PG}$ and the endpoints of EG INPs in $G$ which are not in $G_{PG}$;

\noindent{$(b)$ } we add one edge between two vertices corresponding to vertices in $G_{PG}$ if there exists an edge in $G_{PG}$ between the corresponding vertices of $G_{PG}$;

\noindent{$(c)$ } we add one edge between two vertices corresponding to the endpoints of an EG INP.

Note that we have a natural continuous application $p_{G^{\ast}} \colon G^{\ast} \to G$ which sends an edge as defined in $(b)$ to the corresponding edge in $G_{PG}$ and which sends an edge as defined in $(c)$ to the corresponding EG INP in $G$. Let $x \in VG^{\ast}$.

\begin{lem}\label{Lem Injectivity pi1}
\noindent{$(1)$ } If $\gamma$ is a nontrivial reduced path in $G^{\ast}$, so is $p_{G^{\ast}}(\gamma)$.

\medskip

\noindent{$(2)$ } The homomorphism $$p_{G^{\ast}}' \colon \pi_1(G^{\ast},x) \to \pi_1(G,p_{G^{\ast}}(x))$$ induced by $p_{G^{\ast}}$ is injective.
\end{lem}

\dem $(1)$ Let $\gamma$ be a reduced path in $G^{\ast}$. Suppose towards a contradiction that $p_{G^{\ast}}(\gamma)$ is not a reduced path in $G$. Thus, there exist an edge $e \in \vec{E}G$ and two paths $a$ and $b$ such that $p_{G^{\ast}}(\gamma)=aee^{-1}b$. Let $e^{\ast}$ be an arc in $\gamma$ such that $p_{G^{\ast}}(e^{\ast})=ee^{-1}$. Note that, by definition of $p_{G^{\ast}}$, the application $p_{G^{\ast}}$ sends edges of $G^{\ast}$ to reduced edge paths in $G$. In particular, the path $e^{\ast}$ is not contained in a single edge of $G^{\ast}$. As the image of an edge in $G^{\ast}$ by $p_{G^{\ast}}$ is either an edge in $G$ or an edge path, we see that the path $e^{\ast}$ is contained in at most two edges of $G^{\ast}$. Let $e_1,e_2 \in G^{\ast}$ be such that $e^{\ast} \subseteq e_1e_2$. Suppose first that $p_{G^{\ast}}(e_1)$ and $p_{G^{\ast}}(e_2)$ are edges in $G_{PG}$. Then $p_{G^{\ast}}(e_1)=e$ and $p_{G^{\ast}}(e_2)=e^{-1}$. But, as $\gamma$ is reduced, we have $e_1 \neq e_2^{-1}$. Thus we have $p_{G^{\ast}}(e_1) \neq p_{G^{\ast}}(e_2)^{-1}$. Suppose now that $p_{G^{\ast}}(e_1)$ is an edge in $G_{PG}$ and $p_{G^{\ast}}(e_2)$ is an EG INP. By Proposition~\ref{Prop definition CT}~$(9)$, the first edge of $p_{G^{\ast}}(e_2)$ is an edge in an EG stratum. By definition, every edge in $G_{PG}$ is an edge in an NEG stratum. Hence the turn $\{p_{G^{\ast}}(e_1)^{-1},p_{G^{\ast}}(e_2)\}$ is nondegenerate. Therefore, we see that $p_{G^{\ast}}(e^{\ast}) \neq ee^{-1}$. Finally, suppose that $p_{G^{\ast}}(e_1)$ and $p_{G^{\ast}}(e_2)$ are EG INPs. for every $i \in \{1,2\}$, let $r_i$ be the height of $p_{G^{\ast}}(e_i)$. By Proposition~\ref{Prop definition CT}~$(9)$, the last edge of $p_{G^{\ast}}(e_1)$ is in $H_{r_1}$ whereas the first edge of $p_{G^{\ast}}(e_2)$ is in $H_{r_2}$. Hence if $r_1 \neq r_2$, there is no identification between $p_{G^{\ast}}(e_1)$ and $p_{G^{\ast}}(e_2)$. Hence $p_{G^{\ast}}(e^{\ast}) \neq ee^{-1}$. If $r_1=r_2$, then by the uniqueness statement in Proposition~\ref{Prop definition CT}~$(9)$, we have $p_{G^{\ast}}(e_2) \in \{p_{G^{\ast}}(e_1),p_{G^{\ast}}(e_1)^{-1}\}$. Hence $e_2 \in \{e_1,e_1^{-1}\}$. As $\gamma$ is a reduced path, we see that $e_2=e_1$. Hence $e_1$ is a loop and $p_{G^{\ast}}(e_1)$ is a closed EG INP. By Proposition~\ref{Prop definition CT}~$(9)$, the initial and terminal edges of $p_{G^{\ast}}(e_1)$ are distinct unoriented edges. Hence the path $p_{G^{\ast}}(e_1)p_{G^{\ast}}(e_2)$ is a reduced path and $p_{G^{\ast}}(e^{\ast}) \neq ee^{-1}$. As we have ruled out every case, we see that such a path $e^{\ast}$ does not exist. This concludes the proof of Assertion~$(1)$.

\bigskip

\noindent{$(2)$ } Let $\gamma$ be a nontrivial reduced closed path in $G^{\ast}$ based at $x$. By Assertion~$(1)$, the path $p_{G^{\ast}}(\gamma)$ is a nontrivial reduced closed path in $G$. Hence the kernel of $p_{G^{\ast}}'$ is trivial.

\hfill\qedsymbol

\begin{lem}\label{Lem bijection on circuits Gpg}
The application $[f]$ which sends a circuit $\alpha$ in $G$ to $[f(\alpha)]$ preserves the set of circuits which are concatenations of paths in $G_{PG}$ and in $\mathcal{N}_{PG}$. Moreover, $[f]$ restricts to a bijection on the set of circuits which are concatenations of paths in $G_{PG}$ and in $\mathcal{N}_{PG}$.
\end{lem}

\dem The first part follows from Lemma~\ref{Lem iterate of a path in GPG}~$(3)$. By~\cite[Lemma~III.1.6~$(2),(5)$]{HandelMosher20}, the application $[f]$ preserves $\left\langle Z_{r_+}, \rho_{r_+} \right\rangle$ and restricts to a bijection on the set of circuits of $\left\langle Z_{r_+}, \rho_{r_+} \right\rangle$. By Lemma~\ref{Lem Gpg in Zr} concatenations of paths in $G_{PG}$ and in $\mathcal{N}_{PG}$ are contained in $\left\langle Z_{r_+}, \rho_{r_+} \right\rangle$. By Lemma~\ref{Lem iterate of a path in GPG}, the application $[f]$ preserves concatenations of paths in $G_{PG}$ and in $\mathcal{N}_{PG}$. In particular, this shows that $[f]$ is injective when restricted to the set of paths which are concatenations of paths in $G_{PG}$ and in $\mathcal{N}_{PG}$. 

For surjectivity, let $\alpha$ be a circuit in $G$ which is a concatenation of paths in $G_{PG}$ and in $\mathcal{N}_{PG}$ and let $x$ be a vertex in $\alpha$ which is either an endpoint of an edge in $G_{PG}$ or an endpoint of an EG INP contained in $\alpha$. Note that by Proposition~\ref{Prop definition CT}~$(2)$, the endpoint of every edge in $G_{PG}$ is fixed by $f$. Moreover, the endpoint of every EG INP is fixed by $f$. Therefore, $f$ fixes $x$. 
The circuit $\alpha$ naturally corresponds to a circuit $\alpha'$ in $G^{\ast}$. Let $x'$ be the vertex of $\alpha'$ corresponding to $x$ (which exists by the choices made on $x$). Since $[f]$ preserves concatenations of paths in $G_{PG}$ and in $\mathcal{N}_{PG}$ by Lemma~\ref{Lem iterate of a path in GPG}, the application $[f]$ induces an application $$[f]_{G^{\ast}} \colon \pi_1(G^{\ast},x') \to \pi_1(G^{\ast},x').$$ Note that, by Lemma~\ref{Lem Injectivity pi1}, the group $\pi_1(G^{\ast},x')$ is naturally identified with a subgroup of $\pi_1(G,x)$. By \cite[Lemma~6.0.6]{BesFeiHan00}, the application $[f]_{G^{\ast}}$ is a bijection. Hence there exists a closed path $\beta'$ in $G^{\ast}$ such that $[f]_{G^{\ast}}([\beta'])=\alpha'$. Let $\beta$ be the circuit corresponding to $\beta'$ in $G$. Then $\beta$ is a concatenation of paths in $G_{PG}$ and in $\mathcal{N}_{PG}$ and $[f(\beta)]=\alpha$. This concludes the proof.
\hfill\qedsymbol

\begin{prop}\label{Prop circuits in Gpg are elements in poly subgroup}
Let ${\tt n} \geq 3$. Let $\phi \in \Out(F_{\tt n},\mathcal{F})$ be an exponentially growing outer automorphism, let $f \colon G \to G$ be a CT map representing a power of $\phi$. Let $w \in F_{\tt n}$. There exists a subgroup $A$ of $F_{\tt n}$ such that $[A] \in \mathcal{A}(\phi)$ and $w \in A$ if and only if the circuit $\gamma_w$ of $G$ associated with $w$ is a concatenation of paths in $G_{PG}$ and in $\mathcal{N}_{PG}$.
\end{prop}

\dem Suppose first that $\gamma_w$ is a concatenation of paths in $G_{PG}$ and in $\mathcal{N}_{PG}$. We claim that $\gamma_w$ has polynomial growth under iteration of $f$. By Proposition~\ref{Prop definition CT}~$(8)$, there exists $m \in \NN^*$ such that $[f^m(\gamma_w)]$ is completely split. By Lemma~\ref{Lem iterate of a path in GPG}~$(3)$, the path $[f^m(\gamma_w)]$ is a concatenation of paths in $G_{PG}$ and in $\mathcal{N}_{PG}$. Hence every splitting unit of $[f^m(\gamma_w)]$ is either an edge of $G_{PG}$ or an INP. Let $[f^m(\gamma_w)]=\gamma_1\ldots \gamma_k$ be the complete splitting of $[f^m(\gamma_w)]$. For every $i \geq m$, we have $$\ell[f^i(\gamma_w)])=\sum_{j=1}^k \ell([f^i(\gamma_j)]).$$ Therefore, it suffices to prove that, for every $j \in \{1,\ldots,k\}$, there exists a polynomial $P_j \in \ZZ[X]$ such that for every $i \in \NN^*$, we have $$\ell([f^i(\gamma_j)])=\operatorname{O}(P(i)).$$

\noindent{\bf Claim. } There exists a polynomial $P \in \ZZ[X]$ such that for every edge $e \in \vec{E}G_{PG}$ and every $i \in \NN^*$, we have $$\ell([f^i(e)])=\operatorname{O}(P(i)).$$

\dem Since there are finitely many edges in $G_{PG}$, it suffices to prove the claim for a single edge $e \in \vec{E}G_{PG}$. Let $e \in \vec{E}G_{PG}$. By Proposition~\ref{Prop definition CT}~$(2)$, there exists a cyclically reduced, completely split circuit $w$ of height less than the one of $e$ and such that $f(e)=ew$. By Lemma~\ref{Lem iterate of a path in GPG}~$(1)$, the path $w$ is a concatenation of paths in $G_{PG}$ and in $\mathcal{N}_{PG}$. We prove the claim by induction on the height of $e$. Suppose first that $e$ has minimal height in $G_{PG}$. By minimality of $e$, the path $w$ does not contain a splitting unit which is an edge in $G_{PG}$. Hence $w$ is either trivial or a path in $\mathcal{N}_{PG}$, that it, a closed Nielsen path. If $w$ is trivial then $e$ is a fixed edge and $P=1$ satisfies the claim. Suppose that $w$ is a closed Nielsen path. For every $i \in \NN^*$, we have $[f^i(e)]=ew^i$. Hence $\ell([f^i(e)]) \leq i\ell(w)+1$. Then the polynomial $P(i)=i\ell(w)+1$ satisfies the assertion of the claim. This proves the base case. Suppose now that $e$ has height $r$. Let $w=w_1\ldots w_k$ be the complete splitting of $w$. Recall that, for every reduced path $x$ in $G$, we have $[f([f(x)])]=[f^2(x)]$. Thus, for every $i \in \NN^*$. we have $$[f^i(e)]=ew_1\ldots w_k[f(w_1)]\ldots [f(w_k)] \ldots [f^{i-1}(w_1)]\ldots[f^{i-1}(w_k)].$$ Hence, for every $i \in \NN^*$, we have $$\ell([f^i(e)])=1+\sum_{\ell=1}^k\sum_{j=0}^{i-1} \ell([f^j(w_{\ell})]).$$ Hence it suffices, for every $\ell \in \{1,\ldots,k\}$, to find a polynomial $P_{\ell} \in \ZZ[X]$ such that, for every $i \in \NN^*$, we have $$\ell([f^i(w_{\ell})])=\operatorname{O}(P_{\ell}(i)).$$ Let $\ell \in \{1,\ldots,k\}$. As $w$ is a concatenation of paths in $G_{PG}$ and in $\mathcal{N}_{PG}$, every splitting unit of $w$ is either an edge in $G_{PG}$ or an INP. If $w_{\ell}$ is an edge in $G_{PG}$, the polynomial $P_{\ell}$ exists using the induction hypothesis. If $w_{\ell}$ is an INP, then the polynomial $P_{\ell}(i)=\ell(w_{\ell})$ satisfies the conclusion of the claim. This proves the existence of the polynomial $P$.
\hfill\qedsymbol

\bigskip

Let $j \in \{1,\ldots,k\}$. If $\gamma_k$ is an edge in $G_{PG}$ which is a splitting unit of $[f^m(\gamma_w)]$, by the above claim, the polynomial $P_j$ exists. If $\gamma_{j}$ is an INP, then the polynomial $P_{\ell}(x)=\ell(\gamma_{j})$ satisfies the conclusion. Thus, the path $\gamma_w$ has polynomial growth under iteration of $[f]$. Therefore, $[w]$ has polynomial growth under iterates of $\phi$. By the definition of $\mathcal{A}(\phi)$, there exists a subgroup $A$ of $F_{\tt n}$ such that $[A] \in \mathcal{A}(\phi)$ and $w \in A$.

Conversely, suppose that there exists a subgroup $A$ of $F_{\tt n}$ such that $[A] \in \mathcal{A}(\phi)$ and $w \in A$. Let $m \in \NN^*$ be such that $[f^m(\gamma_w)]$ is completely split, which exists by Proposition~\ref{Prop definition CT}~$(7)$. Since $[w]$ has polynomial growth under iteration of $\phi$, there does not exist a splitting unit of $[f^m(\gamma_w)]$ which is an edge in an EG stratum or a superlinear edge with exponential growth. Suppose towards a contradiction that a splitting unit $\sigma$ of $[f^m(\gamma_w)]$ is contained in a zero stratum. By Proposition~\ref{Prop definition CT}~$(3)$, every zero stratum of $G$ is contractible. As $[f^m(\gamma_w)]$ is a cycle, it is not contained in a zero stratum. By Proposition~\ref{Prop definition CT}~$(4)$, every edge adjacent to $\sigma$ and not contained in the same stratum as $\sigma$ is in an EG stratum. Hence there exists a splitting unit $\sigma'$ of $[f^m(\gamma_w)]$ such that $\sigma\sigma' \subseteq [f^m(\gamma_w)]$ and $\sigma'$ the first edge of $\sigma$ is in an EG stratum. Hence $\sigma'$ is either an edge in an EG stratum or an INP. But, by Lemma~\ref{Lem No zero path and Nielsen path adjacent}, the path $\sigma'$ is not an INP. Hence $\sigma'$ is an edge in an EG stratum. This contradicts the fact that $[w]$ has polynomial growth under iteration of $\phi$. Hence every splitting unit of $[f^m(\gamma_w)]$ is either an INP, an exceptional path or an edge in an NEG stratum whose iterates by $f$ do not contain splitting units which are edges in EG strata. Edges in the last category are precisely the edges in $G_{PG}$. By Lemma~\ref{Lem NEG INP in Npg} and Lemma~\ref{Lem Nielsen paths in Npg} every INP and every exceptional path is a concatenation of paths in $G_{PG}$ and in $\mathcal{N}_{PG}$. Thus, the path $[f^m(\gamma_w)]$ is a concatenation of paths in $G_{PG}$ and in $\mathcal{N}_{PG}$. By Lemma~\ref{Lem bijection on circuits Gpg}, the circuit $\gamma_w$ is a concatenation of paths in $G_{PG}$ and in $\mathcal{N}_{PG}$.
\hfill\qedsymbol

\bigskip

Let $\mathcal{F}$ be a nonsporadic free factor system of $F_{\tt n}$ and let $\phi \in \Out(F_{\tt n},\mathcal{F})$. We say that $\phi$ is \emph{fully irreducible relative to $\mathcal{F}$} if no power of $\phi$ preserves a proper free factor system $\mathcal{F}'$ of $F_{\tt n}$ such that $\mathcal{F} < \mathcal{F}'$. The following corollary will be used in~\cite{Guerch2021Polygrowth}. It is a well-known result but we did not find a precise statement in the literature.

\begin{coro}\label{Coro fully irreducible}
Let ${\tt n } \geq 3$ and let $\mathcal{F}$ be a nonsporadic free factor system of $F_{\tt n}$. Let $\phi \in \Out(F_{\tt n},\mathcal{F})$ be a fully irreducible outer automorphism relative to $\mathcal{F}$. There exists at most one (up to taking inverse) conjugacy class $[g]$ of root-free $\mathcal{F}$-nonperipheral element of $F_{\tt n}$ which has polynomial growth under iteration of $\phi$. Moreover, the conjugacy class $[g]$ is $\phi$-periodic.
\end{coro}

\dem Let $f \colon G \to G$ be a CT map representing a power of $\phi$ and let $G'$ be a subgraph of $G$ such that $\mathcal{F}(G')=\mathcal{F}$. Since $\phi$ is irreducible relative to $\mathcal{F}$ and since $\mathcal{F}$ is nonsporadic, we see that $\overline{G-G'}$ is an EG stratum $H_r$. Let $[g]$ be the conjugacy class of a root-free $\mathcal{F}$-nonperipheral element $g$ of $F_{\tt n}$. Then $\gamma_g$ has height $r$. Suppose that $[g]$ has polynomial growth with respect to $\phi$. By Proposition~\ref{Prop circuits in Gpg are elements in poly subgroup}, the circuit $\gamma_g$ is a concatenation of paths in $G_{PG}$ and in $\mathcal{N}_{PG}$. Since $\gamma_g$ has height $r$ and since $H_r$ is an $EG$ stratum, every subpath $\alpha$ of $\gamma_g$ contained in $H_r$ is contained in a concatenation of INPs of height $r$. By Proposition~\ref{Prop definition CT}~$(9)$, there exists at most one INP $\sigma$ of height $r$. Moreover, one of its endpoints is not contained in $G'=G_{r-1}$ (see~\cite[I.Fact~1.42]{HandelMosher20}). Hence $\sigma$ is necessarily a closed EG INP. Since the endpoint of $\sigma$ is not in $G_{r-1}$ and since $\gamma_g$ is a concatenation of paths in $G_{PG}$ and $\mathcal{N}_{PG}$, we see that $\gamma_g$ is an iteration of the closed path $\sigma$. Since $g$ is root-free, we have $\gamma_g=\sigma^{\pm 1}$. This concludes the proof. 
\hfill\qedsymbol

\subsection{The exponential length of a CT map}\label{Section exponential length}

In this section, we define the \emph{exponential length function} $\ell_{exp}$, and its relative version $\ell_{\mathcal{F}}$, of paths in CT maps. We compute its value for some paths in $G$. Let $G_{PG}'=G_{PG} \cup \mathcal{Z}$ (see Definition~\ref{Defi Gpg}) and let $G_{PG,\mathcal{F}}'=G_{PG,\mathcal{F}} \cup \mathcal{Z}_{\mathcal{F}}$.

Let $\gamma$ be a reduced edge path in $G$. By Lemma~\ref{Lem Nielsen paths in NPG properties}~$(2)$, every path of $\mathcal{N}_{PG}$ which is contained in $\gamma$ is contained in a unique maximal subpath of $\gamma$ contained in $\mathcal{N}_{PG}$. Thus, the path $\gamma$ has a unique decomposition into edge paths $\gamma=\gamma_0\gamma_1'\gamma_1\ldots \gamma_k\gamma_k'$ where: 

\medskip

\noindent{$(1)$ } for every $i \in \{0,\ldots,k\}$, the path $\gamma_i$ is a maximal path in $\mathcal{N}_{PG}$  contained in $\gamma$ (where $\gamma_0$ and $\gamma_k'$ might be trivial);

\medskip

\noindent{$(2)$ } for every $\gamma' \in \mathcal{N}_{PG}$ contained in $\gamma$, there exists $i \in \{1,\ldots,k\}$ such that $\gamma' \subseteq \gamma_i$. 

\medskip 

Such a decomposition of $\gamma$ is called the \emph{exponential decomposition of $\gamma$}. Note that the exponential decomposition of $\gamma$ is not necessarily a splitting of $\gamma$. We denote by $\mathcal{N}_{PG}^{\max}(\gamma)$ the set consisting in  all paths $\gamma_i$, with $i \in \{0,\ldots,k\}$. Similarly, $\gamma$ has a decomposition $\gamma=\alpha_0\alpha_1'\alpha_1\ldots \alpha_m\alpha_m'$, where for every $i \in \{0,\ldots,m\}$, the path $\alpha_i$ is a maximal path in $\mathcal{N}_{PG,\mathcal{F}}$ and for every $\gamma' \in \mathcal{N}_{PG,\mathcal{F}}$ contained in $\gamma$, there exists $i \in \{1,\ldots,k\}$ such that $\gamma' \subseteq \alpha_i$.  Such a decomposition is called the \emph{$\mathcal{F}$-exponential decomposition of $\gamma$}. We denote by $\mathcal{N}_{PG,\mathcal{F}}^{\max}(\gamma)$ the set consisting in all paths $\alpha_i$, with $i \in \{0,\ldots,m\}$.

\begin{defi}
\noindent{$(1)$ } Let $\gamma$ be a reduced edge path in $G$. The \emph{exponential length of $\gamma$}, denoted by $\ell_{exp}(\gamma)$ is: $$\ell_{exp}(\gamma)=\ell\left(\gamma \cap \overline{G-G_{PG}'}\right)-\sum_{\alpha \in \mathcal{N}_{PG}^{\max}(\gamma)} \ell\left(\alpha \cap \overline{G-G_{PG}'}\right).$$

\medskip

\noindent{$(2)$ }  Let $\gamma$ be a reduced edge path in $G$.
The \emph{$\mathcal{F}$-exponential length of $\gamma$}, denoted by $\ell_{\mathcal{F}}(\gamma)$ is: $$\ell_{\mathcal{F}}(\gamma)=\ell\left(\gamma \cap \overline{G-G_{PG,\mathcal{F}}'}\right)-\sum_{\alpha \in \mathcal{N}_{PG,\mathcal{F}}^{\max}(\gamma)} \ell\left(\alpha \cap \overline{G-G_{PG,\mathcal{F}}'}\right).$$

\medskip

\noindent{$(3)$ } Let $\gamma$ be a reduced edge path in $G$ and let $\gamma=\gamma_0\gamma_1'\gamma_1\ldots\gamma_k'\gamma_k$ be the exponential decomposition of $\gamma$. A \emph{$PG$-relative complete splitting} of the path $\gamma$ is a splitting $\gamma=\delta_1\ldots\delta_m$ such that for every $i \in \{1,\ldots,m\}$, the path $\delta_i$ is one of the following paths:

\medskip

\noindent{$\bullet$ } a splitting unit of positive exponential length not contained in some $\gamma_i$ for $i \in \{0,\ldots,k\}$;

\medskip

\noindent{$\bullet$ } a maximal taken connecting path in a zero stratum;

\medskip

\noindent{$\bullet$ } a subpath of $\gamma$ which is a concatenation of subpaths contained in $G_{PG}$ and Nielsen paths in $\mathcal{N}_{PG}$.

\medskip

We call the above paths \emph{$PG$-relative splitting units}. If $\gamma$ is a circuit, a \emph{$PG$-relative circuital complete splitting} of $\gamma$ is a circuital splitting of $\gamma$ which is a $PG$-relative complete splitting of $\gamma$.

\medskip

\noindent{$(4)$ } A \emph{factor} of a $PG$-relative completely split edge path $\gamma$ is a concatenation of $PG$-relative splitting units of some given $PG$-relative complete splitting of $\gamma$.
\end{defi}

Note that if $\gamma$ is an edge path of $G$, then $\ell_{exp}(\gamma) \geq 0$. Indeed, two paths $\gamma_1$ and $\gamma_2$ contained in $\mathcal{N}_{PG}^{\max}(\gamma)$ are either equal or disjoint. Let $\gamma=\gamma_0\gamma_1'\gamma_1\ldots\gamma_k'\gamma_k$ be the exponential decomposition of $\gamma$. For every $i \in \{1,\ldots,k\}$, we have \mbox{$\ell_{exp}(\gamma_i')=\ell(\gamma_i' \cap \overline{G-G_{PG}'})$} and $$\ell_{exp}(\gamma)=\sum_{i=1}^k \ell_{exp}(\gamma_i').$$ We prove the existence of $PG$-relative complete splittings in Lemma~\ref{Lem completely split and paths in Npg}. Note that a $PG$-relative complete splitting of a reduced edge path $\gamma$ is not necessarily unique. Indeed, it might be possible that one can split a $PG$-relative splitting unit of $\gamma$ which is a concatenation of paths in $G_{PG}$ and in $\mathcal{N}_{PG}$ into two $PG$-relative splitting units which are concatenations of paths in $G_{PG}$ and in $\mathcal{N}_{PG}$. 

In the rest of the section, we describe some properties of the exponential length.

\begin{lem}\label{Lem Exponential length less exp length subpaths}
Let $\gamma$ be a reduced edge path in $G$ and let $\gamma=\gamma_1\gamma_2$ be a decomposition of $\gamma$ into two edge paths. We have:
$$\ell_{exp}(\gamma) \leq \ell_{exp}(\gamma_1)+\ell_{exp}(\gamma_2).$$ 
\end{lem}

\dem It is immediate that $$\ell(\gamma \cap \overline{G-G_{PG}'})=\ell(\gamma_1 \cap \overline{G-G_{PG}'})+ \ell(\gamma_2 \cap \overline{G-G_{PG}'}).$$ Let $i \in \{1,2\}$. Let $\gamma' \in \mathcal{N}_{PG}^{\max}(\gamma_i)$. Then there exists $\gamma'' \in \mathcal{N}_{PG}^{\max}(\gamma)$ such that $\gamma' \subseteq \gamma''$. In particular, we have $$\sum_{\gamma'' \in \mathcal{N}_{PG}^{\max}(\gamma)}\ell(\gamma'' \cap \overline{G-G_{PG}'}) \geq \sum_{\gamma' \in \mathcal{N}_{PG}^{\max}(\gamma_1)}\ell(\gamma' \cap \overline{G-G_{PG}'})+\sum_{\gamma' \in \mathcal{N}_{PG}^{\max}(\gamma_2)}\ell(\gamma' \cap \overline{G-G_{PG}'}).$$ By definition of the exponential length, this concludes the proof.
\hfill\qedsymbol

\bigskip

Note that we do not necessarily have equality in Lemma~\ref{Lem Exponential length less exp length subpaths}. Indeed, let $\gamma=\gamma_1\gamma_2$ be as in Lemma~\ref{Lem Exponential length less exp length subpaths}. Suppose that the endpoint of $\gamma_1$ is contained in a path $\gamma'$ of $\mathcal{N}_{PG}^{\max}(\gamma)$. Then $\gamma'$ is not necessarily a concatenation of paths in $\mathcal{N}_{PG}^{\max}(\gamma_1)$ and $\mathcal{N}_{PG}^{\max}(\gamma_2)$. Therefore, we might have: $$\sum_{\gamma' \in \mathcal{N}_{PG}^{\max}(\gamma)}\ell(\gamma' \cap \overline{G-G_{PG}'}) > \sum_{\gamma' \in \mathcal{N}_{PG}^{\max}(\gamma_1)}\ell(\gamma' \cap \overline{G-G_{PG}'})+\sum_{\gamma' \in \mathcal{N}_{PG}^{\max}(\gamma_2)}\ell(\gamma' \cap \overline{G-G_{PG}'}),$$ and a strict inequality in Lemma~\ref{Lem Exponential length less exp length subpaths}. In particular, a proper subpath of $\gamma$ might have greater exponential length than $\gamma$ itself. For instance, if $\gamma$ is a reduced path in $G$ such that $\ell_{exp}(\gamma)=0$, it is possible that there exists a proper subpath $\gamma'$ of $\gamma$ such that $\ell_{exp}(\gamma')>0$. However, there exists a bound, depending only on $G$, on the difference of the exponential length of a subpath of $\gamma$ and the exponential length of $\gamma$ (see~Lemma~\ref{Lem bound exponential length subpath}).

\bigskip

If $\gamma$ is a path in $G$ such that $\ell_{exp}(\gamma)=0$, we do not necessarily have $\ell_{exp}([f(\gamma)])=0$. Indeed, if $\gamma$ is an edge in a zero stratum such that $[f(\gamma)]$ contains a splitting unit which is an edge in an EG stratum, we have $\ell_{exp}([f(\gamma)])>0$. However, the following lemma describes an important situation where the map $f$ preserves the property of having zero exponential length.

\begin{lem}\label{Lem exponential length paths in Gpg}
Let $\gamma$ be a reduced edge path which is a concatenation of paths in $G_{PG}$ and in $\mathcal{N}_{PG}$. For every $n \in \NN$, we have $\ell_{exp}([f^n(\gamma)])=0$.
\end{lem}

\dem Since the $[f]$-image of a concatenation of paths in $G_{PG}$ and in $\mathcal{N}_{PG}$ is a concatenation of paths in $G_{PG}$ and in $\mathcal{N}_{PG}$ by Lemma~\ref{Lem iterate of a path in GPG}, it suffices to prove the result for $n=0$. Let $\gamma$ be a concatenation of paths in $G_{PG}$ and in $\mathcal{N}_{PG}$. Let $\gamma=\gamma_0\gamma_1'\gamma_1\ldots\gamma_k\gamma_k'$ be the exponential decomposition of $\gamma$: for every $i \in \{1,\ldots,k\}$, the path $\gamma_i$ is a maximal subpath of $\gamma$ in $\mathcal{N}_{PG}$ and for every $i \in \{0,\ldots,k\}$, the path $\gamma_i'$ is a path in $G_{PG}$. Note that for every $i \in \{1,\ldots,k\}$, we have $\gamma_i \in \mathcal{N}_{PG}^{\max}(\gamma)$. By definition of the exponential length, we have $\ell_{exp}(\gamma)=\sum_{i=0}^k \ell_{exp}(\gamma_i')=0$.
\hfill\qedsymbol

\begin{coro}\label{Lem paths in Npg' have zero exponential length}
Let $\gamma$ be a path of $\mathcal{N}_{PG}'$. Then $\ell_{exp}(\gamma)=0$. In particular, if $\gamma$ is either a closed Nielsen path, an NEG INP or an exceptional path, we have $\ell_{exp}(\gamma)=0$.
\end{coro}

\dem By Lemma~\ref{Lem Nielsen paths in Npg}, the path $\gamma$ is a concatenation of paths in $G_{PG}$ and in $\mathcal{N}_{PG}$. By Lemma~\ref{Lem exponential length paths in Gpg}, we have $\ell_{exp}(\gamma)=0$. The second assertion follows from Lemmas~\ref{Lem closed nielsen paths in Npg'} and~\ref{Lem NEG INP in Npg}.
\hfill\qedsymbol

\begin{lem}\label{Lem completely split and paths in Npg}
Let $\gamma$ be a completely split edge path and let $\gamma=\gamma_1\ldots\gamma_m$ be its complete splitting. Let $\gamma' \in \mathcal{N}_{PG}^{\max}(\gamma)$. Then either $\gamma'$ is a concatenation of splitting units of $\gamma$ or there exists $i \in \{1,\ldots,m\}$ such that $\gamma' \subsetneq \gamma_i$. Moreover, the complete splitting of $\gamma$ is a $PG$-relative complete splitting of $\gamma$.
\end{lem}

\dem Let $e$ be the first edge of $\gamma'$ and let $i \in \{1,\ldots,m\}$ be such that $e$ is contained in $\gamma_i$. Let $\sigma$ be the splitting unit of $\gamma'$ containing $e$. By Proposition~\ref{Prop definition CT}~$(9)$, the edge $e$ is in an EG stratum. Hence $\gamma_i$ is either an edge in an EG stratum, an exceptional path or an INP. Since $\gamma'$ is a Nielsen path, and since $\gamma_i$ is a splitting unit of $\gamma$, we see that $\gamma_i$ is not an edge in an EG stratum. If $\gamma_i$ is either an NEG INP or an exceptional path, then Proposition~\ref{Prop definition CT}~$(11)$ implies that $\gamma_i$ starts and ends with edges in NEG strata whose height are strictly higher than the one of $e$. Since the height of $e$ is equal to the height of $\sigma$, we see that $\gamma_i$ contains $\sigma$. An inductive argument shows that $\gamma'$ is contained in $\gamma_i$. 

Suppose now that $\gamma_i$ is an EG INP. By Lemma~\ref{Lem Nielsen paths in NPG properties}~$(2)$ applied to $\gamma_i$ and $\gamma'$, either $\gamma'$ is contained in $\gamma_i$ or $\gamma_i$ is the initial segment of $\gamma'$. If $\gamma'$ is contained in $\gamma_i$, by maximality of $\gamma'$, we see that $\gamma'=\gamma_i$. Suppose that $\gamma'$ is the initial segment of the completely split edge path $\gamma_i\ldots\gamma_k$. Then \cite[Corollary~4.12]{FeiHan06} implies that $\gamma'$ is a factor of $\gamma$. 

The last assertion of the lemma follows from the following observations. Every splitting unit of $\gamma$ which is either an INP or an exceptional path is a concatenation of paths in $G_{PG}$ and in $\mathcal{N}_{PG}$ by Lemma~\ref{Lem NEG INP in Npg}. Moreover, by the first assertion of the lemma, every splitting unit of $\gamma$ which is an edge in an irreducible stratum not contained in $G_{PG}$ does not intersect a path in $\mathcal{N}_{PG}^{\max}(\gamma)$. Hence the complete splitting of $\gamma$ is a $PG$-relative complete splitting.
\hfill\qedsymbol

\bigskip

$PG$-relative completely split edge paths are well-adapted to the computation of the exponential length as explained by the following lemma.

\begin{lem}\label{Lem compute exp length completely split}
Let $\gamma$ be a $PG$-relative completely split edge path and let $\gamma=\alpha_1\ldots \alpha_{\ell}$ be a $PG$-relative complete splitting.

\medskip 

\noindent{$(1)$ } For every path $\gamma' \in \mathcal{N}_{PG}^{\max}(\gamma)$, there exists a minimal concatenation of $PG$-relative splitting units  $\delta$ of $\gamma$ such that $\gamma' \subseteq \delta$; every $PG$-relative splitting unit of $\delta$ is a concatenation of paths in $G_{PG}$ and in $\mathcal{N}_{PG}$; for every $PG$-relative splitting unit $\delta'$ of $\delta$, the intersection $\delta' \cap \gamma'$ is an element of $\mathcal{N}_{PG}^{\max}(\delta')$.

\medskip

\noindent{$(2)$ } We have $\ell_{exp}(\gamma)=\sum_{i=1}^{\ell} \ell_{exp}(\alpha_i)$ and $\ell_{\mathcal{F}}(\gamma)=\sum_{i=1}^{\ell} \ell_{\mathcal{F}}(\alpha_i)$.
\end{lem}

\dem $(1)$ Let $\gamma=\gamma_0\gamma_1'\gamma_1\ldots\gamma_k'\gamma_k$ be the exponential decomposition of $\gamma$ where, for every $i \in \{0,\ldots,k\}$, we have $\gamma_i \in \mathcal{N}_{PG}^{\max}(\gamma)$. Let $i \in \{0,\ldots,k\}$. Let $j \in \{1,\ldots,\ell\}$ be such that $\alpha_j$ contains an initial segment of $\gamma_i$. By Proposition~\ref{Prop definition CT}~$(10)$, the splitting unit $\alpha_j$ is not contained in a zero stratum. Moreover, by definition of the $PG$-relative splitting units, if $\alpha_j$ is an edge in an irreducible stratum of positive exponential length, it is not contained in $\gamma_i$. Hence, by the description of $PG$-relative splitting units, the path $\alpha_j$ is a concatenation of paths in $G_{PG}$ and in $\mathcal{N}_{PG}$. By Proposition~\ref{Prop definition CT}~$(9)$, the path $\gamma_i$ starts with an edge in an EG stratum. Hence there exists a path $\beta_j$ in $\mathcal{N}_{PG}^{\max}(\alpha_j)$ which contains an initial segment of $\gamma_i$. By maximality of $\gamma_i$, we see that $\beta_j \subseteq \gamma_i$. Suppose first that $\beta_j=\gamma_i$. Then setting $\delta=\alpha_j$ proves the first assertion. Suppose now that $\beta_j \subsetneq \gamma_i$. By Lemma~\ref{Lem Nielsen paths in NPG properties}~$(2)$ applied to $\gamma=\gamma_i^{-1}$ and $\gamma'=\beta_j^{-1}$, the path $[\beta_j^{-1}\gamma_i]$ is a path in $\mathcal{N}_{PG}$. Therefore, by Proposition~\ref{Prop definition CT}~$(9)$, the path $[\beta_j^{-1}\gamma_i]$ starts with an edge in an EG stratum. Note that, as $\alpha_j$ is a concatenation of paths in $G_{PG}$ and in $\mathcal{N}_{PG}$, if $\alpha_j$ contains the first edge $e$ of $[\beta_j^{-1}\gamma_i]$, then $e$ would be contained in an EG INP contained in $\alpha_j$. Since $\beta_j$ is a maximal subpath of $\alpha_j$ in $\mathcal{N}_{PG}$, we see that $[\beta_j^{-1}\gamma_i]$ is contained in $\gamma''=\alpha_{j+1}\ldots\alpha_{\ell}$ and is in $\mathcal{N}_{PG}^{\max}(\gamma'')$. We can thus apply the same arguments to the paths $[\beta_j^{-1}\gamma_i]$ and $\gamma''$. This concludes the proof of $(1)$. 

The proof of $(2)$ follows as the exponential length and the $\mathcal{F}$-length are computed by removing paths in $G_{PG}$ and in $\mathcal{N}_{PG}$. As all subpaths in $G_{PG}$ are contained in a splitting unit of $\gamma$ and as subpaths in $\mathcal{N}_{PG}$ are obtained by concatenating paths in $\amalg_{j=1}^{\ell} \mathcal{N}_{PG}^{\max}(\alpha_j)$, we see that $\ell_{exp}(\gamma)=\sum_{i=1}^{\ell} \ell_{exp}(\alpha_i)$ and $\ell_{\mathcal{F}}(\gamma)=\sum_{i=1}^{\ell} \ell_{\mathcal{F}}(\alpha_i)$. 
\hfill\qedsymbol

\bigskip

The following property of the exponential length allows us to pass, if needed, to a further iterate of the CT map $f$.

\begin{lem}\label{Lem exponential length goes to infinity}
For every edge $e$ of $\overline{G-G_{PG}'}$, we have $$\lim_{n \to \infty} \ell_{exp}([f^n(e)])=\infty \text{ and } \lim_{n \to \infty} \ell_{\mathcal{F}}([f^n(e)])=\infty.$$ Moreover, the sequences $(\ell_{exp}([f^n(e)]))_{n \in \NN}$ and $(\ell_{\mathcal{F}}([f^n(e)]))_{n \in \NN}$ grows exponentially fast. 
\end{lem}

\dem We prove the result concerning $\ell_{exp}$, the proof of the result concerning $\ell_{\mathcal{F}}$ follows from the fact that for every reduced edge path $\gamma$ in $G$, we have $\ell_{exp}(\gamma) \leq \ell_{\mathcal{F}}(\gamma)$. Let $e$ be an edge of $\overline{G-G_{PG}'}$. Since every iterate of $e$ is completely split by Proposition~\ref{Prop definition CT}~$(6)$ and since there exists an iterate of $e$ which contains a splitting unit which is an edge in an EG stratum, we may suppose that $e$ is an edge in an EG stratum $H_r$. Since $H_r$ is an EG stratum, the number of edges in $[f^n(e)] \cap H_r$ grows exponentially fast as $n$ goes to infinity. Therefore the number of splitting units of $[f^n(e)]$ which are edges of $H_r$ grows exponentially fast and $\lim_{n \to \infty} \ell_{exp}([f^n(e)])=\infty$. 
\hfill\qedsymbol

\begin{lem}\label{Lem exponential length monotonic for completely split edge path}
Let $\gamma$ be a $PG$-relative completely split edge path. There exists $n_0 \in \NN^*$ such that for every $k \geq n_0$, we have $\ell_{exp}([f^k(\gamma)]) \geq \ell_{exp}(\gamma)$.
\end{lem}

\dem Let $\gamma=\gamma_1\ldots\gamma_k$ be a $PG$-relative complete splitting of $\gamma$. By Lemma~\ref{Lem compute exp length completely split}, it suffices to prove the assertion for every subpath $\gamma_i$, with $i \in \{1,\ldots,k\}$. Let $i \in \{1,\ldots,k\}$. If $\gamma_i$ is a concatenation of paths in $G_{PG}$ and in $\mathcal{N}_{PG}$, then $\ell_{exp}([f(\gamma_i)])=\ell_{exp}(\gamma_i)=0$ by Lemma~\ref{Lem exponential length paths in Gpg}. If $\gamma_i$ is a maximal taken connecting path in a zero stratum, we have $\ell_{exp}(\gamma_i)=0$. Hence $\ell_{exp}([f(\gamma_i)]) \geq \ell_{exp}(\gamma_i)$. In the other cases, $\gamma_i$ is an edge in an irreducible stratum which is not contained in $G_{PG}$. By Lemma~\ref{Lem exponential length goes to infinity}, we have $\lim_{n \to \infty}\ell_{exp}([f^n(\gamma_i)])=\infty$. Hence there exists $n_0 \in \NN^*$ such that, for every $k \geq n_0$, we have $\ell_{exp}([f^k(\gamma_i)]) \geq \ell_{exp}(\gamma_i)$, and $n_0$ may be chosen to be independent of $\gamma_i$ with $i \in \{1,\ldots,k\}$.
\hfill\qedsymbol

\bigskip

The last lemma in this section shows that the exponential length of a $PG$-relative completely split edge path encaptures the splitting units which are edges with exponential growth under iterates of $f$.

\begin{lem}\label{Lem splitting units positive exp length}
Let $\gamma$ be a $PG$-relative completely split edge path, let $\gamma=\gamma_1\ldots\gamma_k$ be a $PG$-relative complete splitting and let $i \in \{1,\ldots,k\}$. Then $\ell_{exp}(\gamma_i)>0$ if and only if $\gamma_i$ is an edge in an irreducible stratum not contained in $G_{PG}$. In particular, the value $\ell_{exp}(\gamma)$  is the number of splitting units which are edges in $\overline{G-G_{PG}'}$.
\end{lem}

\dem Suppose first that $\gamma_i$ is either a concatenation of paths in $G_{PG}$ and in $\mathcal{N}_{PG}$ or a maximal taken connecting path in a zero stratum. By Lemma~\ref{Lem exponential length paths in Gpg}, we have $\ell_{exp}(\gamma_i)=0$. Suppose that $\gamma_i$ is an edge in an irreducible stratum which is not contained in $G_{PG}$. Since there does not exist an EG INP of length $1$, by definition of the exponential length, we have $\ell_{exp}(\gamma_i)=1>0$. This concludes the proof of the first part of the lemma. The computation of $\ell_{exp}(\gamma)$ follows from Lemma~\ref{Lem compute exp length completely split}~$(2)$.
\hfill\qedsymbol

\subsection{The space of polynomially growing currents}\label{Section poly growing currents}

In this section, let $\mathcal{F}$ be a free factor system and let $\phi \in \Out(F_{\tt n},\mathcal{F})$ be an exponentially growing outer automorphism. Recall the definition of $\mathcal{A}(\phi)$ and $\mathcal{F} \wedge \mathcal{A}(\phi)$ from Section~\ref{Subsection malnormal}. We define a subspace of $\PCurr(F_{\tt n},\mathcal{F} \wedge \mathcal{A}(\phi))$, called the \emph{space of polynomially growing currents}. It consists in the currents whose support is contained in $\partial^2\mathcal{A}(\phi)$ (see Lemma~\ref{Lem polynomially growing currents empty}). In order to define it, we first need to show that the exponential length extends to a continuous function $\Psi \colon \PCurr(F_{\tt n},\mathcal{F} \wedge \mathcal{A}(\phi)) \to \RR$. The space of polynomially growing currents will then be defined as a level set of $\Psi$.

We first need some preliminary results concerning paths in $\mathcal{N}_{PG}$. For a path $\gamma \in \mathcal{N}_{PG}$, let $\mathcal{N}_{PG}^{++}(\gamma)$ be the subset of $\mathcal{N}_{PG}$ which consists in all paths $\gamma' \in \mathcal{N}_{PG}$ such that $\gamma \subsetneq \gamma'$ and $\gamma'$ is minimal for this property. Let $\gamma' \in \mathcal{N}_{PG}^{++}(\gamma)$. By Lemma~\ref{Lem Nielsen paths in NPG properties}~$(3)$, either $\gamma$ is properly contained in an INP $\sigma$ of the complete splitting of $\gamma'$, or there exist (possibly trivial) paths $\gamma_1,\gamma_2 \in \mathcal{N}_{PG}$ such that $\gamma'=\gamma_1\gamma\gamma_2$. By minimality, either $\gamma_1$ or $\gamma_2$ is trivial. Moreover, a result of Feighn and Handel (\cite[Corollary~4.12]{FeiHan06}) shows that, in this case, splitting units of the complete splittings of $\gamma_1$, $\gamma_2$ and $\gamma$ are splitting units of $\gamma'$. Thus the set $\mathcal{N}_{PG}^{++}(\gamma)$ can be partitioned into three disjoint subsets: $$\mathcal{N}_{PG}^{++}(\gamma)=\mathcal{N}_{PG,INP}^{++}(\gamma) \amalg \mathcal{N}_{PG,left}^{++}(\gamma) \amalg \mathcal{N}_{PG,right}^{++}(\gamma),$$ where $\mathcal{N}_{PG,INP}^{++}(\gamma)$ is the set of paths in $\mathcal{N}_{PG}^{++}(\gamma)$ such that one of their splitting units properly contains $\gamma$, $\mathcal{N}_{PG,left}^{++}(\gamma)$ is the set of paths $\gamma' \in \mathcal{N}_{PG}^{++}(\gamma)$ such that $\gamma'=\gamma_1\gamma$ and $\mathcal{N}_{PG,right}^{++}(\gamma)$ is the set of paths $\gamma' \in \mathcal{N}_{PG}^{++}(\gamma)$ such that $\gamma'=\gamma\gamma_2$. One can also define similarly the three sets $\mathcal{N}_{PG,INP,\mathcal{F}}^{++}(\gamma)$, $\mathcal{N}_{PG,left,\mathcal{F}}^{++}(\gamma)$ and $\mathcal{N}_{PG,right,\mathcal{F}}^{++}(\gamma)$ as the restriction to the paths in $\mathcal{N}_{PG,INP}^{++}(\gamma)$, $\mathcal{N}_{PG,left}^{++}(\gamma)$ and $\mathcal{N}_{PG,right}^{++}(\gamma)$ contained in $G_p$. We emphasize on the fact that a path in $\mathcal{N}_{PG,INP}^{++}(\gamma)$ might contain several occurrences of the path $\gamma$. However, a path in $\mathcal{N}_{PG,left}^{++}(\gamma)$ or in $\mathcal{N}_{PG,right}^{++}(\gamma)$ contains a unique occurrence of $\gamma$. Indeed, let $\gamma' \in \mathcal{N}_{PG,left}^{++}(\gamma)$ (the proof for $\mathcal{N}_{PG,right}^{++}(\gamma)$ being similar). Then $\gamma'=\gamma_1\gamma_2$ with $\gamma_1 \in \mathcal{N}_{PG}$ and $\gamma_2=\gamma$. Let $\gamma_3$ be an occurrence of $\gamma$ which contains an edge of $\gamma_1$. By Lemma~\ref{Lem EG INP is not cdc}~$(2)$, the path $\gamma_3$ cannot intersect $\gamma_2$ nontrivially. Hence $\gamma_3 \subseteq \gamma_1$. Hence $\gamma_1 \in \mathcal{N}_{PG}$ and $\gamma_1$ contains an occurrence of $\gamma$. This contradicts the minimality of $\gamma'$.

\begin{lem}\label{Lem intersection of nielsen paths in Npgmin}
Let $\gamma$ be a path in $\mathcal{N}_{PG}$. Let $\gamma_1,\gamma_2$ be two distinct paths in $ \mathcal{N}_{PG}^{++}(\gamma)$. Suppose that there exist three paths $\mu_1,\mu_2,\mu_3$ such that $\gamma_1=\mu_1\mu_2$, $\gamma_2=\mu_2\mu_3$ and $\gamma$ is contained in $\mu_2$. Then $\gamma_1 \in \mathcal{N}_{PG,left}^{++}(\gamma)$, $\gamma_2 \in \mathcal{N}_{PG,right}^{++}(\gamma)$ and $\mu_2=\gamma$.
\end{lem}

\dem By Lemma~\ref{Lem Nielsen paths in NPG properties}~$(2)$, the path $\mu_2$ belongs to $\mathcal{N}_{PG}$ and contains $\gamma$. Since $\gamma_1$ and $\gamma_2$ are minimal paths of $\mathcal{N}_{PG}$ for the property of properly containing $\gamma$, we have $\mu_2=\gamma$. Therefore, we see that $\gamma_1=\mu_1\gamma$ and $\gamma_2=\gamma\mu_3$. This shows that $\gamma_1 \in \mathcal{N}_{PG,left}^{++}(\gamma)$ and that $\gamma_2 \in \mathcal{N}_{PG,right}^{++}(\gamma)$.
\hfill\qedsymbol

\bigskip

Lemma~\ref{Lem intersection of nielsen paths in Npgmin} implies that an occurrence of $\gamma$ in the intersection of paths in $\mathcal{N}_{PG}^{++}(\gamma)$ is well-controlled. Following Lemma~\ref{Lem intersection of nielsen paths in Npgmin}, we then define $\mathcal{N}_{PG,lr}^{++}(\gamma)$ to be the set of paths of the form $\gamma_1\gamma\gamma_2$, where $\gamma_1\gamma \in \mathcal{N}_{PG,left}^{++}(\gamma)$ and $\gamma\gamma_2 \in \mathcal{N}_{PG,right}^{++}(\gamma)$. We define similarly the set $\mathcal{N}_{PG,lr,\mathcal{F}}^{++}(\gamma)$ to be the set of all paths in $\mathcal{N}_{PG,lr}^{++}(\gamma)$ contained in $G_p$. As for $\mathcal{N}_{PG,left}^{++}(\gamma)$ and $\mathcal{N}_{PG,right}^{++}(\gamma)$, a path in $\mathcal{N}_{PG,lr}^{++}(\gamma)$ contains a unique occurrence of $\gamma$.

\bigskip

Given two paths $\gamma$ and $\gamma'$ of $G$ let $N(\gamma',\gamma)$ be the number of occurrences of $\gamma$ and $\gamma^{-1}$ in $\gamma'$. Using the finiteness of $\mathcal{N}_{PG}$ (see Lemma~\ref{Lem Nielsen paths in NPG properties}~$(1)$), we denote by $$\Psi_0' \colon \Curr(F_{\tt n},\mathcal{F} \wedge \mathcal{A}(\phi)) \to \RR$$ the continuous function $$\Psi_0'(\nu)=\sum\limits_{\gamma \in \mathcal{N}_{PG}}\Big(\left\langle \gamma,\nu \right\rangle-\sum\limits_{\gamma' \in \mathcal{N}_{PG}^{++}(\gamma)} \left\langle \gamma',\nu \right\rangle N(\gamma',\gamma) + \sum\limits_{\gamma' \in \mathcal{N}_{PG,lr}^{++}(\gamma)} \left\langle \gamma',\nu \right\rangle \Big)\ell\left(\gamma \cap \overline{G-G_{PG}'}\right),$$ and by $\Psi_0  \colon \Curr(F_{\tt n},\mathcal{F} \wedge \mathcal{A}(\phi)) \to \RR$ the continuous linear function $$\begin{array}{ccl}
\Psi_0(\nu)&=&\frac{1}{2}\Big(\sum_{e \in \vec{E}(\overline{G-G_{PG}'})}\left\langle e,\nu \right\rangle-\Psi_0'(\nu)\Big) \\
{} &=& \frac{1}{2}\Big(\sum_{e \in \vec{E}(\overline{G-G_{PG}'})}\Big(\left\langle e,\nu \right\rangle- \\
{} & {} & \sum\limits_{\gamma \in \mathcal{N}_{PG},e \subseteq \gamma} \Big(\left\langle \gamma,\nu \right\rangle-\sum\limits_{\gamma' \in \mathcal{N}_{PG}^{++}(\gamma)} \left\langle \gamma',\nu \right\rangle N(\gamma',\gamma) + \sum\limits_{\gamma' \in \mathcal{N}_{PG,lr}^{++}(\gamma)} \left\langle \gamma',\nu \right\rangle \Big)N(\gamma,e)\Big)\Big).
\end{array}
.$$

\begin{defi}\label{Defi Kpg}
The space of polynomially growing currents, denoted by $K_{PG}(f)$, is the compact subset of $\PCurr(F_{\tt n},\mathcal{F} \wedge \mathcal{A}(\phi))$ consisting in all projective classes of currents $[\nu] \in \PCurr(F_{\tt n},\mathcal{F} \wedge \mathcal{A}(\phi))$ such that: $$\Psi_0(\nu)=0.$$
\end{defi}

Finally, we define the {\it $\mathcal{F}$-simplicial length} function $\lVert.\rVert_{\mathcal{F}} \colon \Curr(F_{\tt n},\mathcal{F} \wedge \mathcal{A}(\phi)) \to \RR$ as 
$$
\begin{array}{ccl}
\lVert \nu \rVert_{\mathcal{F}}&=&\frac{1}{2}\Big(\sum_{e \in \vec{E}(\overline{G-G_{PG,\mathcal{F}}'})}\left\langle e,\nu \right\rangle- \\
{} & {} & \sum\limits_{\gamma \in \mathcal{N}_{PG,\mathcal{F}}}\Big(\left\langle \gamma,\nu \right\rangle-\sum\limits_{\gamma' \in \mathcal{N}_{PG,\mathcal{F}}^{++}(\gamma)} \left\langle \gamma',\nu \right\rangle N(\gamma',\gamma) + \sum\limits_{\gamma' \in \mathcal{N}_{PG,lr,\mathcal{F}}^{++}(\gamma)} \left\langle \gamma',\nu \right\rangle \Big)\ell\left(\gamma \cap \overline{G-G_{PG}'}\right)\Big).
\end{array}
$$

\begin{lem}\label{Lem psi0 equal exp length}
Let $w \in F_{\tt n}$ be a nonperipheral element with conjugacy class $[w]$, associated rational current $\eta_{[w]}$ and associated reduced edge path $\gamma_w$ in $G$. Then $$\begin{array}{c}
\Psi_0(\eta_{[w]})=\ell_{exp}(\gamma_w); \\
\lVert \eta_{[w]} \rVert_{\mathcal{F}}=\ell_{\mathcal{F}}(\gamma_w).
\end{array}$$ 
Therefore $\eta_{[w]} \in K_{PG}(f)$ if and only if $$\ell_{exp}(\gamma_w)=0.$$ 

In particular, there exist a basis $\mathfrak{B}$ of $F_{\tt n}$ and a constant $C>0$ such that, for every $\mathcal{F} \wedge \mathcal{A}(\phi)$-nonperipheral element $g \in F_{\tt n}$, we have $\lVert \eta_{[g]} \rVert_{\mathcal{F}} \in \NN^*$ and $$\ell_{\mathfrak{B}}([g]) \geq C\;\lVert \eta_{[g]} \rVert_{\mathcal{F}}.$$
\end{lem}

\dem We prove the result for $\Psi_0$, the proof for $\lVert \eta_{[w]} \rVert_{\mathcal{F}}$ being similar. First note  that $$\sum_{e \in \vec{E}(\overline{G-G_{PG}'})}\left\langle e, \eta_{[w]}\right\rangle=2\ell(\gamma_w \cap \overline{G-G_{PG}'}),$$ where the factor $2$ follows from the fact that the sum on the left hand side is over oriented edges. Therefore, it remains to prove that 

\begin{equation}\label{Equation Kpoly section 3}
\Psi_0'(\eta_{[w]})=\sum_{\gamma \in \mathcal{N}_{PG}^{\max}(\gamma_w)} \ell\left(\gamma \cap \overline{G-G_{PG}'}\right).
\end{equation} 
Let $\gamma \in \mathcal{N}_{PG}$. Then the value $$\left\langle \gamma,\eta_{[w]} \right\rangle-\sum\limits_{\gamma' \in \mathcal{N}_{PG}^{++}(\gamma)}\left\langle \gamma', \eta_{[w]} \right\rangle N(\gamma',\gamma)+ \sum\limits_{\gamma' \in \mathcal{N}_{PG,lr}^{++}(\gamma)} \left\langle \gamma',\eta_{[w]} \right\rangle$$ measures the number of occurrences of $\gamma$ or $\gamma^{-1}$ in $\gamma_w$ which are not induced by an occurrence of a path $\gamma' \in \mathcal{N}_{PG}$ containing properly $\gamma$ or $\gamma^{-1}$ and contained in $\gamma_w$. Indeed, an occurrence of $\gamma$ in a path $\gamma' \in \mathcal{N}_{PG}$ containing properly $\gamma$ will be counted in $\sum_{\gamma' \in \mathcal{N}_{PG}^{++}(\gamma)}\left\langle \gamma', \eta_{[w]} \right\rangle N(\gamma',\gamma)$. Moreover, if an occurrence of $\gamma$ is contained in two distinct paths $\gamma_1, \gamma_2 \in \mathcal{N}_{PG}^{++}(\gamma)$, Lemma~\ref{Lem intersection of nielsen paths in Npgmin} ensures that this occurrence is contained in a path $\gamma_3 \in \mathcal{N}_{PG,lr}^{\max}(\gamma)$. Therefore, the value $$-\sum_{\gamma' \in \mathcal{N}_{PG}^{++}(\gamma)}\left\langle \gamma', \eta_{[w]} \right\rangle N(\gamma',\gamma)+\sum_{\gamma' \in \mathcal{N}_{PG,lr}^{++}(\gamma)} \left\langle \gamma',\eta_{[w]} \right\rangle$$ measures an occurrence of $\gamma$ or $\gamma^{-1}$ in a larger path, and each such occurrence will be counted exactly once. Therefore, the equation below Equation~\eqref{Equation Kpoly section 3} measures the number of occurrences of $\gamma$ and $\gamma^{-1}$ in $\mathcal{N}_{PG}^{\max}(\gamma_w)$. Thus, the equality~\eqref{Equation Kpoly section 3} holds. The last assertions of Lemma~\ref{Lem psi0 equal exp length} then follows by definition of $K_{PG}(f)$ and of $\ell_{\mathcal{F}}$.
\hfill\qedsymbol

\bigskip

Note that in the proof of Lemma~\ref{Lem psi0 equal exp length}, we show that, for every edge $e \in \vec{E}(\overline{G-G_{PG}'})$ and every nonperipheral element $w \in F_{\tt n}$, the value: $$\left\langle e,\eta_{[w]} \right\rangle-\sum\limits_{\gamma \in \mathcal{N}_{PG},e \subseteq \gamma}\Big(\left\langle \gamma,\eta_{[w]} \right\rangle-\sum\limits_{\gamma' \in \mathcal{N}_{PG}^{++}(\gamma)} \left\langle \gamma',\eta_{[w]} \right\rangle N(\gamma',\gamma) + \sum\limits_{\gamma' \in \mathcal{N}_{PG,lr}^{++}(\gamma)} \left\langle \gamma',\eta_{[w]} \right\rangle \Big)N(\gamma,e)$$ measures the number of occurrences of $e$ in $\gamma_w$ which are not contained in a path of $\mathcal{N}_{PG}^{\max}(\gamma_w)$. Thus, for every nonperipheral element and every edge $e \in \vec{E}(\overline{G-G_{PG}'})$, we have:

$$\left\langle e,\eta_{[w]} \right\rangle-\sum\limits_{\gamma \in \mathcal{N}_{PG},e \subseteq \gamma}\Big(\left\langle \gamma,\eta_{[w]} \right\rangle-\sum\limits_{\gamma' \in \mathcal{N}_{PG}^{++}(\gamma)} \left\langle \gamma',\eta_{[w]} \right\rangle N(\gamma',\gamma) + \sum\limits_{\gamma' \in \mathcal{N}_{PG,lr}^{++}(\gamma)} \left\langle \gamma',\eta_{[w]} \right\rangle \Big)N(\gamma,e) \geq 0.$$ The density of rational currents given by Proposition~\ref{Prop density rational currents} and the continuity of $\left\langle .,. \right\rangle$ then shows that for every current $\nu \in \Curr(F_{\tt n},\mathcal{F} \wedge \mathcal{A}(\phi))$ and every edge $e \in \vec{E}(\overline{G-G_{PG}'})$, we have : $$\left\langle e,\nu \right\rangle-\sum\limits_{\gamma \in \mathcal{N}_{PG},e \subseteq \gamma}\Big(\left\langle \gamma,\nu \right\rangle-\sum\limits_{\gamma' \in \mathcal{N}_{PG}^{++}(\gamma)} \left\langle \gamma',\nu \right\rangle N(\gamma',\gamma) + \sum\limits_{\gamma' \in \mathcal{N}_{PG,lr}^{++}(\gamma)} \left\langle \gamma',\nu \right\rangle \Big)N(\gamma,e) \geq 0.$$

\begin{lem}\label{Lem polynomially growing currents empty}
Let ${\tt n} \geq 3$ and let $\mathcal{F}$ be a free factor system. Let $\phi \in \Out(F_{\tt n},\mathcal{F})$ be an exponentially growing outer automorphism. Let $f \colon G \to G$ be a CT map representing a power of $\phi$. 

\medskip

\noindent{$(1)$ } If $[\nu] \in K_{PG}(f)$, then $\Supp(\nu) \subseteq \partial^2(F_{\tt n},\mathcal{F} \wedge \mathcal{A}(\phi)) \cap \partial^2 \mathcal{A}(\phi)$. In particular, if $\phi$ is expanding relative to $\mathcal{F}$, then $K_{PG}(f)=\varnothing$.

\medskip

\noindent{$(2)$ } Conversely, if $\nu \in \Curr(F_{\tt n},\mathcal{F} \wedge \mathcal{A}(\phi))$ is such that the support $\Supp(\nu)$ of $\nu$ is contained in $\partial^2(F_{\tt n},\mathcal{F} \wedge \mathcal{A}(\phi)) \cap \partial^2 \mathcal{A}(\phi)$, then $[\nu] \in K_{PG}(f)$. Thus we have $$K_{PG}(f)=\{[\mu] \in \PCurr(F_{\tt n},\mathcal{F} \wedge \mathcal{A})\;|\; \Supp(\mu) \subseteq \partial^2(F_{\tt n},\mathcal{F} \wedge \mathcal{A}(\phi)) \cap \partial^2 \mathcal{A}(\phi)\}.$$

\medskip

\noindent{$(3)$ } If $\nu \in \Curr(F_{\tt n},\mathcal{F} \wedge \mathcal{A}(\phi))$, we have $\lVert \nu \rVert_{\mathcal{F}}=0$ if and only if $\nu=0$.

\end{lem}

\dem The proof of $(3)$ being identical to the proof of $(1)$ and $(2)$ replacing $G_{PG}'$ and $\mathcal{N}_{PG}$ by $G_{PG,\mathcal{F}}'$ and $\mathcal{N}_{PG,\mathcal{F}}$, we only prove $(1)$ and $(2)$. For the proof of both $(1)$ and $(2)$, let $\mathcal{B}$ be a free basis of $F_{\tt n}$ and let $T$ be the Cayley graph of $F_{\tt n}$ associated with $\mathcal{B}$. Let $\mathscr{C}(\mathcal{A}(\phi))$ be the set of elements of $F_{\tt n}$ associated with $\mathcal{A}(\phi)$ given by Lemma~\ref{Lem finite set of words determines vertex group system}. Recall that $\mathrm{Cyl}(\mathscr{C}(\mathcal{A}(\phi)))$ is the set of cylinder subsets of the form $C(\gamma)$, where $\gamma$ is a geodesic edge path in $T$ starting at the base point whose associated element $w \in F_{\tt n}$ contains a word of $\mathscr{C}(\mathcal{A}(\phi))$ as a subword.

\medskip

\noindent{$(1)$ } Let $\nu \in \Curr(F_{\tt n},\mathcal{F} \wedge \mathcal{A}(\phi))$ nonzero be such that $\Supp(\nu)$ is not contained in $\partial^2(F_{\tt n},\mathcal{F} \wedge \mathcal{A}(\phi)) \cap \partial^2 \mathcal{A}(\phi)$. Then $\Supp(\nu) \cap \partial^2(F_{\tt n},\mathcal{A}(\phi)) \neq \varnothing$. Hence the restriction of $\nu$ to $\partial^2(F_{\tt n},\mathcal{A}(\phi))$ induces a nonzero current $\nu' \in \Curr(F_{\tt n},\mathcal{A}(\phi))$. By Lemma~\ref{Lem double boundary open} applied to $\mathcal{A}=\mathcal{A}(\phi)$ and $\nu'$, there exists $C(\gamma) \in \mathscr{C}(\mathcal{A}(\phi))$ such that $\nu(C(\gamma))>0$. Let $w$ be the element of $F_{\tt n}$ associated with $\gamma$, and let $\gamma_w'$ be the reduced circuit in $G$ associated with the conjugacy class of $w$. Up to taking a larger geodesic edge path $\gamma''\supseteq \gamma$ in $T$ such that $\nu(C(\gamma''))>0$ (which exists by additivity of $\nu$), we may suppose that $w$ is cyclically reduced. By Lemma~\ref{Lem finite set of words determines vertex group system}~$(3)$, the path $\gamma$ is not contained in any tree $T_A$ such that $[A] \in \mathcal{A}(\phi)$. As $w$ is cyclically reduced, the translation axis in $T$ of $w$ contains $\gamma$. Hence $\{w^{+\infty},w^{-\infty}\} \notin \partial^2 \mathcal{A}(\phi)$ and $w$ is not contained in any subgroup $A$ such that $[A] \in \mathcal{A}(\phi)$. By Proposition~\ref{Prop circuits in Gpg are elements in poly subgroup}, the circuit $\gamma_w'$ is not a concatenation of paths in $G_{PG}$ and in $\mathcal{N}_{PG}$. Therefore, there exists an edge $e$ of $G$ such that $$\left\langle e,\nu \right\rangle-\sum_{\gamma \in \mathcal{N}_{PG},e \subseteq \gamma} \Big(\left\langle \gamma,\nu \right\rangle-\sum\limits_{\gamma' \in \mathcal{N}_{PG}^{++}(\gamma)} \left\langle \gamma',\nu \right\rangle N(\gamma',\gamma)+ \sum\limits_{\gamma' \in \mathcal{N}_{PG,lr}^{++}(\gamma)} \left\langle \gamma',\nu \right\rangle \Big)N(\gamma,e) > 0.$$ Thus, we see that $\Psi_0(\nu)>0$ and that $[\nu] \notin K_{PG}(f)$. The second part of $(1)$ follows from the fact that, if $\phi$ is expanding relative to $\mathcal{F}$, then $\partial^2 \mathcal{A}(\phi) \subseteq \partial^2 \mathcal{F}$. This proves $(1)$.

\medskip

\noindent{$(2)$ } Let $\nu \in \Curr(F_{\tt n},\mathcal{F} \wedge \mathcal{A}(\phi))$ be such that $\Supp(\nu) \subseteq \partial^2(F_{\tt n},\mathcal{F} \wedge \mathcal{A}(\phi)) \cap \partial^2 \mathcal{A}(\phi)$. Let $e$ be an edge such that $\left\langle e,\nu \right\rangle >0$. By Lemma~\ref{Lem Nielsen paths in NPG properties}~$(1)$, there exists a constant $C_1>0$ such that, for every path $\gamma' \in \mathcal{N}_{PG}$, we have $\ell(\gamma') \leq C_1$. Recall the definition of the graph $G^{\ast}$ and the application $p_{G^{\ast}} \colon G^{\ast} \to G$. from Lemma~\ref{Lem Injectivity pi1}. Let $C_2$ be the length of a maximal path in a maximal forest of $p_{G^{\ast}}(G^{\ast})$. Let $C=\max\{2C_1,C_2\}$.

\medskip

\noindent{\bf Claim. } Let $\gamma$, $\delta_1$ and $\delta_2$ be reduced paths such that $\gamma=\delta_1 e \delta_2$, $\ell(\delta_1),\ell(\delta_2) \geq 2C$ and $\left\langle \gamma, \nu \right\rangle >0$. Let $\gamma=\gamma_0\gamma_1'\gamma_1\ldots \gamma_k\gamma_k'$ be the exponential decomposition of $\gamma$ (where, for every $i \in \{0,\ldots,k\}$, the path $\gamma_i$ is contained in $\mathcal{N}_{PG}$). Either $e \in \vec{E}G_{PG}'$ or $e$ is contained in an EG stratum and there exists $i \in \{0,\ldots,k\}$ such that $e \subseteq \gamma_i$.

\medskip

\dem Since $\Supp(\nu) \subseteq \partial^2(F_{\tt n},\mathcal{F} \wedge \mathcal{A}(\phi)) \cap \partial^2 \mathcal{A}(\phi)$, there exists a subgroup $A$ of $F_{\tt n}$ such that $[A] \in \mathcal{A}(\phi)$, and two elements $a$ and $b$ of $A$ such that the geodesic path in $\widetilde{G}$ representing $\{a^{+\infty},b^{+\infty}\} \in \partial^2 A$ contains a lift of $\gamma$. If $b=a^{-1}$, then $\gamma$ is contained in an iterate of $a$ and, by Proposition~\ref{Prop circuits in Gpg are elements in poly subgroup}, $\gamma$ is contained in a concatenation of paths in $G_{PG}$ and $\mathcal{N}_{PG}$. The claim follows in this case. So we may assume that $b \neq a^{-1}$. Suppose first that the axes $\mathrm{Ax}(a)$ and $\mathrm{Ax}(b)$ of $a$ and $b$ are disjoint. Then $\gamma$ is contained in the axis of $a^{-1}b$. Thus,  by Proposition~\ref{Prop circuits in Gpg are elements in poly subgroup}, $\gamma$ is contained in a concatenation of paths in $G_{PG}$ and $\mathcal{N}_{PG}$ and the claim follows in this case.

Suppose now that $\mathrm{Ax}(a) \cap \mathrm{Ax}(b)\neq \varnothing$. Let $\gamma_a'$ and $\gamma_b'$ be the reduced circuit in $G$ associated with $a$ and $b$. Then $\gamma$ is contained in the union of $\gamma_a' \cup \gamma_b'$. Recall that, by Proposition~\ref{Prop circuits in Gpg are elements in poly subgroup}, the paths $\gamma_a'$ and $\gamma_b'$ are concatenation of paths in $G_{PG}$ and $\mathcal{N}_{PG}$. Hence there exist reduced circuits $\alpha$ and $\beta$ in $G^{\ast}$ and reduced arcs $\tau,\tau_e$ in $G^{\ast}$ such that $p_{G^{\ast}}(\alpha)=\gamma_a'$ and $p^{\ast}(\beta)=\gamma_b'$ and such that $p_{G^{\ast}}(\tau)=\gamma$ and $p_{G^{\ast}}(\tau_e)=e$. By the choice of $C$, and as $\ell(\delta_1),\ell(\delta_2) \geq 2C$, one can remove an initial and a terminal segment of $\tau$ so that the resulting path $\tau'$ is nontrivial, is contained in a subgraph of $G^{\ast}$ with no leaf and is such that $\ell(p_{G^{\ast}}(\tau')) \geq 2C+1$. Thus, there exist subpaths $\tau_1'$, $\tau_1'',\tau_2'$, $\tau_2''$ of $\tau$ and a reduced circuit $\delta$ of $G^{\ast}$ such that:

\noindent{$(i)$ } $\ell(p_{G^{\ast}}(\tau_1')),\ell(p_{G^{\ast}}(\tau_2')) \geq C$,

\noindent{$(ii)$ } $\tau=\tau_1''\tau_1'e\tau_2'\tau_2''$,

\noindent{$(iii)$ } $\tau'=\tau_1'e\tau_2' \subseteq \delta$. 

By Lemma~\ref{Lem Injectivity pi1}~$(1)$, the path $p_{G^{\ast}}(\delta)$ is a reduced ciruit which contains $e$. Since $\ell(p_{G^{\ast}}(\tau_1')),\ell(p_{G^{\ast}}(\tau_2')) \geq C \geq 2C_1$, if $\gamma' \in \mathcal{N}_{PG}^{\max}(p_{G^{\ast}}(\delta))$ is such that $e \subseteq \gamma'$, then $\gamma' \subseteq \tau_1'e\tau_2'$. Hence it suffices to prove the claim for $\gamma=p_{G^{\ast}}(\delta)$. As $\delta$ is a concatenation of paths in $G_{PG}$ and in $\mathcal{N}_{PG}$, the claim follows.
\hfill\qedsymbol

\medskip

Suppose towards a contradiction that there exists an edge $e \in \overline{G-G_{PG}'}$ such that:

\begin{equation}\label{Equation Kpg}
\left\langle e,\nu \right\rangle-\sum_{\gamma \in \mathcal{N}_{PG},e \subseteq \gamma} \Big(\left\langle \gamma,\nu \right\rangle-\sum\limits_{\gamma' \in \mathcal{N}_{PG}^{++}(\gamma)} \left\langle \gamma',\nu \right\rangle N(\gamma',\gamma)+ \sum\limits_{\gamma' \in \mathcal{N}_{PG,lr}^{++}(\gamma)} \left\langle \gamma',\nu \right\rangle \Big)N(\gamma,e) > 0.
\end{equation} 

By additivity of $\nu$, there exists a reduced path $\gamma$ of length $4C+1$ such that the path $\gamma$ has a decomposition $\gamma=\gamma_1 e \gamma_2$, where for every $i \in \{1,2\}$, the path $\gamma_i$ has length equal to $2C$ and we have $\nu(C(\gamma)) > 0$. By Equation~\ref{Equation Kpg}, we can choose $\gamma$ such that if $\gamma' \in \mathcal{N}_{PG}^{\max}(\gamma)$, then $\gamma'$ does not contain $e$. Hence $e \notin G_{PG}'$ and $e$ is not contained in a subpath of $\mathcal{N}_{PG}^{\max}(\gamma)$. This contradicts the above claim and this concludes the proof.
\hfill\qedsymbol

\bigskip

Let $\mathcal{F}$ be a free factor system and let $\phi \in \Out(F_{\tt n},\mathcal{F})$ be an exponentially growing outer automorphism. Note that, by Lemma~\ref{Lem polynomially growing currents empty} and since for every $k \in \NN^*$, we have $\mathcal{A}(\phi)=\mathcal{A}(\phi^k)$, the space $K_{PG}(f)$ does not depend on the CT map $f$ and does not depend on the chosen power of $\phi$. Therefore, we will simply write $K_{PG}(\phi)$ instead. Moreover, since $\mathcal{A}(\phi)=\mathcal{A}(\phi^{-1})$, we see that $K_{PG}(\phi)=K_{PG}(\phi^{-1})$. 

\bigskip

For the next proposition, let $C_1>0$ be a constant such that for every path $\gamma \in \mathcal{N}_{PG}$, we have $\ell(\gamma) \leq C_1$. It exists since $\mathcal{N}_{PG}$ is finite by Lemma~\ref{Lem Nielsen paths in NPG properties}~$(1)$. Let $L$ be the malnormality constant associated with $\mathcal{A}(\phi)$ as defined above Lemma~\ref{Lem finite set of words determines vertex group system} and let $C_0=\max\{C_1,L\}$. Let $\mathscr{C}$ be the set of elements of $F_{\tt n}$ associated with $\mathcal{F} \wedge \mathcal{A}(\phi)$ given above Lemma~\ref{Lem finite set of words determines vertex group system}. Let $\mathcal{P}(\mathcal{F} \wedge \mathcal{A}(\phi))$ be the set of reduced paths $\gamma$ in $G$ such that $C(\gamma) \in \mathrm{Cyl}(\mathscr{C})$, $\ell(\gamma)>C_0$ and $\gamma$ is not contained in a concatenation of paths in $G_{PG,\mathcal{F}}$ and $\mathcal{N}_{PG,\mathcal{F}}$.

\begin{lem}\label{Lem path positive exp length cover double boundary}
Let ${\tt n} \geq 3$, let $\mathcal{F}$ be a free factor system of $F_{\tt n}$ and let $\phi \in \Out(F_{\tt n},\mathcal{F})$ be an exponentially growing outer automorphism. We have $$\partial^2(F_{\tt n},\mathcal{F} \wedge \mathcal{A}(\phi))=\bigcup_{\gamma \in \mathcal{P}(\mathcal{F} \wedge \mathcal{A}(\phi))}C(\gamma).$$
\end{lem}

\dem Let $A_1,\ldots,A_r$ be subgroups of $F_{\tt n}$ such that $\mathcal{F} \wedge \mathcal{A}(\phi)=\{[A_1],\ldots,[A_r]\}$ and $\mathscr{C}=\mathscr{C}(A_1,\ldots,A_r)$. By Lemma~\ref{Lem double boundary open}, we have $$\partial^2(F_{\tt n},\mathcal{F} \wedge \mathcal{A}(\phi))=\bigcup_{C(\gamma) \in \mathrm{Cyl}(\mathscr{C})}C(\gamma).$$ Note that, for every path $\gamma \subseteq G$, we have $$C(\gamma)=\bigcup_{e \in \vec{E}G, \; \ell(\gamma e)> \ell(\gamma)} C(\gamma e).$$ Hence we have $$\partial^2(F_{\tt n},\mathcal{F} \wedge \mathcal{A}(\phi))=\bigcup_{C(\gamma) \in \mathrm{Cyl}(\mathscr{C}), \; \ell(\gamma)>C_0}C(\gamma).$$ So it suffices to prove that we can restrict our considerations to paths $\gamma$ which are not contained in a concatenation of paths in $G_{PG,\mathcal{F}}$ and $\mathcal{N}_{PG,\mathcal{F}}$. Let $\gamma$ be a path such that $C(\gamma) \in \mathrm{Cyl}(\mathscr{C})$ and $\ell(\gamma)>C_0$. By Lemma~\ref{Lem finite set of words determines vertex group system}~$(3)$, the path $\gamma$ is not contained in any tree $T_{gA_ig^{-1}}$ with $g \in F_{\tt n}$ and $i \in \{1,\ldots,r\}$. Moreover, it is not contained in any path of $\mathcal{N}_{PG}$ since $\ell(\gamma) > C_1$. Suppose that $\gamma$ is contained in a concatenation of paths in $G_{PG,\mathcal{F}}$ and $\mathcal{N}_{PG,\mathcal{F}}$. Suppose first that there does not exist a circuit which contains $\gamma$ and which is a concatenation of paths in $G_{PG,\mathcal{F}}$ and $\mathcal{N}_{PG,\mathcal{F}}$. Recall the definition of $G^{\ast}$ and $p_{G^{\ast}}$ from Lemma~\ref{Lem Injectivity pi1} and let $G_{\mathcal{F}}^{\ast}=p_{G^{\ast}}^{-1}(G_p)$. By assumption, either there does not exist an immersed path (not necessarily an edge path) $\gamma^{\ast}$ in $G_{\mathcal{F}}^{\ast}$ such that $ p_{G^{\ast}}(\gamma^{\ast})=\gamma$ or there exists an immersed path $\gamma^{\ast}$ in $G_{\mathcal{F}}^{\ast}$ such that $p_{G^{\ast}}(\gamma^{\ast})=\gamma$ and $\gamma^{\ast}$ is not contained in a circuit of $G_{\mathcal{F}}^{\ast}$ (recall that $G^{\ast}$ might contain univalent vertices). In the first case, we have $\ell_{\mathcal{F}}(\gamma)>0$. In the second case, since $G^{\ast}$ is finite, by Lemma~\ref{Lem Injectivity pi1}, up to considering $\gamma^{-1}$, there exists $d \in \NN^*$ such that for every path of $\gamma'$ such that $\gamma\gamma'$ is a reduced path in $G$ and $\ell(\gamma\gamma')=\ell(\gamma)+d$, the path $\gamma\gamma'$ is not the image by $p_{G^{\ast}}$ of an immersed path in $G_{\mathcal{F}}^{\ast}$. Thus we have $\ell_{\mathcal{F}}(\gamma\gamma')>0$. Using the fact that $$C(\gamma)=\bigcup_{e \in \vec{E}G, \ell(\gamma e)> \ell(\gamma)} C(\gamma e),$$ we can replace $\gamma$ by paths $\gamma''$ such that $\gamma \subseteq \gamma''$ and $\gamma''$ is not contained in a concatenation of paths in $G_{PG,\mathcal{F}}$ and $\mathcal{N}_{PG,\mathcal{F}}$. This concludes the proof.
\hfill\qedsymbol

\bigskip

Let $\nu$ be a nonzero current in $\Curr(F_{\tt n},\mathcal{F} \wedge \mathcal{A}(\phi))$. By Lemma~\ref{Lem polynomially growing currents empty}~$(3)$, we have $\lVert \nu \rVert_{\mathcal{F}}\neq 0$. The following result characterizes limits in \mbox{$\PCurr(F_{\tt n},\mathcal{F} \wedge \mathcal{A}(\phi))$}. The result is due to Kapovich~\cite[Lemma~3.5]{Kapovich2006} for a nonrelative context.

\begin{lem}\label{Lem Caracterisation limites PCurr}
Let ${\tt n} \geq 3$ and let $\mathcal{F}$ be a free factor system of $F_{\tt n}$. Let $\phi \in \Out(F_{\tt n},\mathcal{F})$ be an exponentially growing outer automorphism. Let $([\mu_n])_{n \in \NN}$ be a sequence in $\PCurr(F_{\tt n},\mathcal{F} \wedge \mathcal{A}(\phi))$ and let $[\mu] \in \PCurr(F_{\tt n},\mathcal{F} \wedge \mathcal{A}(\phi))$. Let $G$ be a graph whose fundamental group is isomorphic to $F_{\tt n}$ and such that there exists a subgraph $G_p$ of $G$ such that $\mathcal{F}(G_p)=\mathcal{F}$. Then $\lim\limits_{n \to \infty} [\mu_n]=[\mu]$ if and only if, for every reduced edge path $\gamma \in \mathcal{P}(\mathcal{F} \wedge \mathcal{A}(\phi))$, we have 
\begin{equation}\label{Equation lem convergence}
\lim\limits_{n \to \infty} \frac{\left\langle \gamma,\mu_n \right\rangle}{\lVert \mu_n \rVert_{\mathcal{F}}}=\frac{\left\langle \gamma,\mu \right\rangle}{\lVert \mu \rVert_{\mathcal{F}}}.
\end{equation}
\end{lem}

\dem Suppose first that $\lim\limits_{n \to \infty} [\mu_n]=[\mu]$. Thus there exists a sequence $(\lambda_n)_{n \in \NN^*}$ of positive real numbers such that $\lim\limits_{n \to \infty} \lambda_n\mu_n=\mu$. By continuity of  $\lVert. \rVert_{\mathcal{F}}$, we have $\lim\limits_{n \to \infty} \lVert \lambda_n\mu_n\rVert_{\mathcal{F}}=\lVert\mu\rVert_{\mathcal{F}}$. By linearity of $\lVert. \rVert_{\mathcal{F}}$ and $\left\langle .,. \right\rangle$ in the second variable, for every reduced edge path $\gamma \in \mathcal{P}(\mathcal{F} \wedge \mathcal{A}(\phi))$, we have 
$$
\lim\limits_{n \to \infty} \frac{\left\langle \gamma,\lambda_n\mu_n \right\rangle}{\lVert \lambda_n\mu_n \rVert_{\mathcal{F}}}=\lim\limits_{n \to \infty} \frac{\left\langle \gamma,\mu_n \right\rangle}{\lVert \mu_n \rVert_{\mathcal{F}}}=\frac{\left\langle \gamma,\mu \right\rangle}{\lVert \mu \rVert_{\mathcal{F}}}.
$$

Suppose now that for every reduced edge path $\gamma \in \mathcal{P}(\mathcal{F} \wedge \mathcal{A}(\phi))$, Equation~\eqref{Equation lem convergence} holds. By Lemma~\ref{Lem path positive exp length cover double boundary}, for every Borel subset $B$ of $\partial^2(F_{\tt n},\mathcal{F} \wedge \mathcal{A}(\phi))$ such that $\mu(\partial B)=0$, we have $$\lim\limits_{n \to \infty} \frac{\mu_n(B)}{\lVert \mu_n \rVert_{\mathcal{F}}}=\frac{\mu(B)}{\lVert \mu \rVert_{\mathcal{F}}}.$$ Hence we have $\lim\limits_{n \to \infty} [\mu_n]=[\mu]$.
\hfill\qedsymbol

\section{Stable and unstable currents for relative atoroidal outer automorphisms}\label{Section Stable unstable currents}

Let ${\tt n} \geq 3$ and let $\mathcal{F}$ be a free factor system of $F_{\tt n}$. Let $\phi \in \Out(F_{\tt n},\mathcal{F})$ be an atoroidal outer automorphism relative to $\mathcal{F}$. In this section, under additional hypotheses on $\phi$, we construct two $\phi$-invariant convex subsets of $\PCurr(F_{\tt n},\mathcal{F})$. We will then show in the following section that, with respect to these convex subsets, the outer automorphism $\phi$ acts with \emph{generalized north-south dynamics}. 

In order to define the extremal points of these simplices, we need some results regarding substitution dynamics.

\subsection{Substitution dynamics}

Let $A$ be a finite set with cardinality at least equal to $2$. Let $\zeta$ be a \emph{substitution on $A$}, that is, a map from $A$ to the set of nonempty finite words on $A$. The substitution $\zeta$ induces a map on the set of all finite words on $A$ by concatenation, which we still denote by $\zeta$. We can therefore iterate the substitution $\zeta$. For a word $w$ on $A$, we will denote by $|w|$ the length of $w$ on the alphabet $A$.

To the substitution $\zeta$ one can associate its \emph{transition matrix} $M$, which is a square matrix whose rows and columns are indexed by letters in $A$ and, for all $a,b \in A$, $M(a,b)$ is the number of occurrences of $a$ in $\zeta(b)$. Likewise, for $n \geq 1$, the matrix $M^n$ is the transition matrix for $\zeta^n$. We say that a substitution $\zeta$ is \emph{irreducible} if its transition matrix is irreducible, and that the substitution is \emph{primitive} if its transition matrix is.

Let $\ell \in \NN^*$ and let $A_{\ell}$ be the set of words on $A$ of length $\ell$. As defined in~\cite[Section~5.4.1]{Queffelec1987}, the substitution $\zeta$ induces a substitution $\zeta_{\ell}$ on $A_{\ell}$ as follows. Let $w=x_1\ldots x_{\ell} \in A_{\ell}$. Then $\zeta_{\ell}(w)=w_1w_2\ldots w_{|\zeta(x_1)|}$, where, for every $i \in \{1,\ldots,|\zeta(x_1)|\}$, the word $w_i$ is the subword of $\zeta(w)$ of length $\ell$ starting at the $i^{th}$ position of $\zeta(x_1)$. Therefore, $\zeta_{\ell}$ is the concatenation of the $|\zeta(x_1)|$ first subwords of $\zeta(w)$ of length $\ell$. Note that the number of $i \in \{1,\ldots,|\zeta(x_1)|\}$ such that $w_i$ that is not contained in $\zeta(x_1)$ is bounded by $\ell-1$. Let $|\cdot|_{\ell}$ be the length of words on $A_{\ell}$. Then $|\zeta_{\ell}(w)|_{\ell}=|\zeta(x_1)|$. Denote by $M_{\ell}$ the transition matrix of $\zeta_{\ell}$. Note that, for every $n,\ell \geq 1$, we have $(\zeta^n)_{\ell}=(\zeta_{\ell})^n$ as applications on the set of words on $A_{\ell}$ and thus $(M^n)_{\ell}=(M_{\ell})^n$.

Consider now a partition of the alphabet $A= \coprod_{i=0}^k B_i$. Suppose that the transition matrix associated with the substitution $\zeta$ is lower block triangular with respect to this partition. Therefore, for every $i \in \{0,\ldots,k\}$, for every $x \in B_i$ and for every $j <i$, the word $\zeta(x)$ does not contain letters in $B_j$. In the remainder of the article, for every $i \in \{0,\ldots,k\}$ the diagonal block in $M$ corresponding to the block $B_i$ will be denoted by  $M_{B_i}$.

The partition of $A$ induces a partition of $A_{\ell}$ as follows. For every $i \in \{0,\ldots,k\}$, let $\widetilde{B}_i \subseteq A_{\ell}$ be the set of all words on $A$ of length $\ell$ which start with a letter in $B_i$ and which, for every $j<i$ do not contain a letter in $B_j$. Let $\overline{B}_i$ be the set of all words $w$ on $A$ of length $\ell$ which start with a letter in $B_i$ and such that there exists $j<i$ such that $w$ contains a letter in $B_j$ (note that $\overline{B}_0$ is empty). Then $\widetilde{B}_i \cup \overline{B}_i$ is the set of all words on $A$ of length $\ell$ which starts with a letter in $B_i$. The hypothesis on the substitution $\zeta$ implies that the transition matrix $M_{\ell}$ is lower block triangular with respect to the partition $$\widetilde{B}_0 \amalg \overline{B}_1 \amalg \widetilde{B}_1 \amalg \ldots \amalg \overline{B}_k \amalg \widetilde{B}_k$$ of $A_{\ell}$. As before, for every $i \in \{0,\ldots,k\}$, we will denote by $M_{\ell,\overline{B}_i}$ the diagonal block in $M_{\ell}$ corresponding to $\overline{B}_i$ and by $M_{\ell,\widetilde{B}_i}$ the diagonal block in $M_{\ell}$ corresponding to $\widetilde{B}_i$.

\begin{lem}\cite[Lemma~8.8]{gupta2017relative}\label{Lem eigenvalues of substitution}
Let $A$ be a finite alphabet equipped with a partition $A= \amalg_{i=0}^k B_i$. Let $\zeta$ be a substitution and let $M$ be its transition matrix. Let $\ell \in \NN^*$.

\medskip

\noindent{$(1)$ } The eigenvalues of $M_{\ell,\widetilde{B}_i}$ are those of $M_{B_i}$ with possibly additional eigenvalues of absolute value at most equal to $1$. 

\medskip

\noindent{$(2)$ } The eigenvalues of $M_{\ell,\overline{B}_i}$ have absolute value at most equal to $1$.
\end{lem}

Fix an integer $p \in \{0,\ldots,k\}$. For every $i \geq p$, let $\overline{B}_i^{(p)}$ be the subset of $\overline{B}_i$ consisting in all words $w$ of length $\ell$ which start with a letter in $B_i$ and such that there exists $j<p$ such that $w$ contains a letter in $B_j$. Then, for every $i \geq p$, the block $M_{\ell,\overline{B}_i}$ decomposes into a lower triangular block matrix where the columns and rows corresponding to $\overline{B}_i^{(p)}$ are on the top left. Let $M_{\ell,\overline{B}_i^{(p)}}$ be the corresponding block matrix. By Lemma~\ref{Lem eigenvalues of substitution}~$(2)$, the eigenvalues of $M_{\ell,\overline{B}_i^{(p)}}$ have absolute value at most $1$. Moreover, for every $i,j \geq p$, for every word $w$ contained in $ \widetilde{B}_j \cup \overline{B}_j-\overline{B}_j^{(p)}$, the word $\zeta_{\ell}(w)$ considered as a word on $A_{\ell}$ does not contain any word of $\overline{B}_i^{(p)}$. Let $M_{\ell}(p)$ be the matrix obtained from $M_{\ell}$ by deleting, for every $i \geq p$, every row and column corresponding to elements in $\widetilde{B}_i$, and every row and columns corresponding to elements of $\overline{B}_i$ which do not belong to $\overline{B}_i^{(p)}$. Note that, by Lemma~\ref{Lem eigenvalues of substitution}~$(1)$, the eigenvalues of $M_{\ell}(p)$ are those of every block $M_{B_j}$ with $j<p$ with possibly additional eigenvalues of absolute value at most $1$. 

We can now prove a result concerning the number of occurrences of words in iterates of a letter. For words $w,v$ on $A$, we denote by $(w,v)$ the number of occurrences of $w$ in $v$, so that $M=((a,\zeta(b))_{a,b \in A}$. For a word $w$ on $A$, we denote by $\lvert\lvert w \rvert\rvert_{(p)}$ the number of letters in $w$ which are contained in some $B_j$ for $j<p$.

\begin{prop}\label{Prop number of occurences substitution}
Let $A$ be an alphabet equipped with a partition $A= \amalg_{i=0}^k B_i$. Let $\zeta$ be a substitution on $A$ and let $M$ be its transition matrix. Suppose that $M$ is lower triangular by block with respect to the partition of $A$. Let $p \in \NN^*$. Let $a \in \bigcup_{j<p} B_j$ be such that $\zeta(a)$ starts with $a$. Suppose that there exists $j<p$ such that $M_{B_j}$ is a primitive block whose Perron-Frobenius eigenvalue is greater than $1$ and such that there exists $n_j \geq 1$ such that $\zeta^{n_j}(a)$ contains a letter of $B_j$. Let $w$ be a word such that $w$ contains a letter in $B_j$. Then

$$\lim_{n \to  \infty} \frac{(w,\zeta^n(a))}{\lvert\lvert\zeta^n(a)\rvert\rvert_{(p)}}$$ exists and is finite. Furthermore there exists a word $w$ containing a letter in some $B_j$ with $j<p$ such that this limit is positive.
\end{prop}

\dem The proof follows \cite[Lemma~8.9]{gupta2017relative} (see also~\cite{LustigUyanik2017} for similar statements). First, up to replacing $A$ by the smallest $\zeta$-invariant subalphabet of $A$ containing $a$ (which still satisfies the hypotheses of Proposition~\ref{Prop number of occurences substitution}), we may suppose that, for every letter $x \in A$, there exists $n_x \geq 1$ such that $\zeta^{n_x}(a)$ contains the letter $x$. Let $\alpha$ be a word on $A$ with length $\ell \geq 1$ that starts with $a$. Note that, since $a \in \cup_{j<p} B_j$, the word $\alpha$ defines a column and a row in $M_{\ell}(p)$. Recall that for every $n$ the number of occurrences of a word $w$ in $\zeta^n(a)$ differs from the number of occurrences of the letter $w \in A_{\ell}$ in $\zeta_{\ell}^n(\alpha)$ by at most $\ell-1$. Moreover, we have $(w,\zeta_{\ell}^n(\alpha))=M_{\ell}^n(p)(w,\alpha)$. 

Let $S$ be the set of all $s<p$ such that $M_{B_s}$ is a primitive block with associated Perron-Frobenius eigenvalue greater than $1$. By assumption, the set $S$ is a nonempty finite set. Let $S'$ be the subset of $S$ consisting in all such $B_s$ such that the associated Perron-Frobenius eigenvalue is maximal. Call this eigenvalue $\lambda$. By Lemma~\ref{Lem eigenvalues of substitution}, the eigenvalue $\lambda$ is also the maximal eigenvalue of the matrix $M_{\ell}(p)$. Let $d_{\lambda}$ be the size of the maximal Jordan block of $M_{\ell}(p)$ associated with $\lambda$. Then the growth under iterates of the maximal Jordan block of $\frac{M_{\ell}(p)}{\lambda}$ is polynomial of degree $d_{\lambda}$. Therefore, we have

$$\lim_{n \to  \infty} \frac{(w,\zeta^n(a))}{\lambda^n n^{d_{\lambda}}}=\lim_{n \to  \infty} \frac{(w,\zeta_{\ell}^n(\alpha))}{\lambda^n n^{d_{\lambda}}}=\lim_{n \to  \infty} \frac{M_{\ell}^n(p)(w,\alpha)}{\lambda^n n^{d_{\lambda}}}=d_{w,a},$$ 
where $d_{w,a}$ is a real number. Moreover, the limit does not depend on the choice of $\alpha$ since, for any $n$, and for any two columns of $M_{\ell}^n(p)$ corresponding to words starting with the same letter, the sum of the values of each column differ by at most $\ell-1$ (see~\cite[Lemma~8.6]{gupta2017relative}). Moreover, there exists a word $w$ such that the limit is positive since we quotiented by the growth of the iterates of the Jordan block with maximal eigenvalue.

Let $||\cdot||$ be the $L_1$-norm on $\RR^{|A_{\ell}|}$. By \cite[Remark~4.1]{LustigUyanik2017}, since $\lim_{n \to  \infty} \frac{M_{\ell}^n(p)(w,\alpha)}{\lambda^n n^{d_{\lambda}}}$ exists, so does $$\lim_{n \to  \infty} \frac{M_{\ell}^n(p)(w,\alpha)}{||M_{\ell}^n(p)(\alpha)||},$$ where $||M_{\ell}^n(p)(\alpha)||$ is the norm of the column of $M_{\ell}^n(p)$ corresponding to $\alpha$.

\medskip

\noindent{\bf Claim. } Suppose that there exists $C \geq 1$ such that for every $n \in \NN$, we have $$\lvert\lvert\zeta^n(a)\rvert\rvert_{(p)} \leq ||M_{\ell}^n(p)(\alpha)|| \leq C\lvert\lvert\zeta^n(a)\rvert\rvert_{(p)}.$$ 
Then $$\lim_{n \to  \infty} \frac{(w,\zeta^n(a))}{\lvert\lvert\zeta^n(a)\rvert\rvert_{(p)}}$$ exists for all words $w$ on $A$ and is positive for some word $w$.

\medskip

\dem Recall that two sequences $(u_n)_{n \in \NN}$ and $(v_n)_{n \in \NN}$ with values in $\RR$ are equivalent if there exists a sequence $(\epsilon_n)_{n \in \NN}$ tending to zero such that $u_n=(1+\epsilon_n)v_n$. Recall that there exists $C' >0$ such that the sequence $(\lVert M_{\ell}^n(p)(\alpha) \rVert)_{n \in \NN}$ is equivalent to $(C'\lambda^nn^{d_{\lambda}})_{n \in \NN}$. Recall also that for every $n$, the value of $\lVert \zeta^n(a) \rVert_{(p)}$ is the norm of $M^n(p)(v_a)$, where $v_a$ is the vector whose coordinates is $1$ on the coordinate associated with $a$ and $0$ otherwise. Hence, since the matrix $M^n(p)$ is nonnegative and not the zero matrix, there exist $C_a, \lambda_{a} \in \RR_+^*$ and $d_{a} \in \NN$ such that the sequence $(\lVert \zeta^n(a) \rVert_{(p)})_{n \in \NN}$ is equivalent to $(C_a\lambda_a^n n^{d_a})_{n \in \NN}$.  Thus, by the assumption of the claim, since the limit $$\lim_{n \to  \infty} \frac{M_{\ell}^n(p)(w,\alpha)}{||M_{\ell}^n(p)(\alpha)||}$$ exists, and is not equal to zero for some $w$, the same is true for $$\lim_{n \to  \infty} \frac{(w,\zeta^n(a))}{\lvert\lvert\zeta^n(a)\rvert\rvert_{(p)}}.$$ This proves the claim.
\hfill\qedsymbol

\medskip

Therefore, in order to conclude the proof of the proposition, it remains to prove that the hypothesis of the claim is true in our context. Let $\zeta^n(a)=x_1\ldots x_{|\zeta^n(a)|}$ and let $$\zeta_{\ell}^n(\alpha)=w_1\ldots w_{|\zeta^n(a)|}.$$ Let $X^n(a)$ be the list $x_1,\ldots,x_{|\zeta^n(a)|}$ and let $X_{<p}^n(a)$ be the sublist of $X^n(a)$ consisting in all letters in $\cup_{i=1}^{p-1} B_i$. Let $X^{(\ell,n)}(\alpha)$ be the list $w_1,\ldots,w_{|\zeta^n(a)|}$ and let $X_{<p}^{(\ell,n)}(\alpha)$ be the sublist of $X^{(\ell,n)}(\alpha)$ which consists in all elements of $X^{(\ell,n)}(\alpha)$ that do not belong to $\cup_{i \leq p} \widetilde{B}_i \cup \overline{B}_i-\overline{B}_i^{(p)}$. Note that $|X_{<p}^{(\ell,n)}(\alpha)|=||M_{\ell}^n(p)(\alpha)||$ and that $|X_{<p}^n(a)|=\lvert\lvert\zeta^n(a)\rvert\rvert_{(p)}$. The fact that $\lvert\lvert\zeta^n(a)\rvert\rvert_{(p)} \leq ||M_{\ell}^n(p)(\alpha)||$ follows from the fact that we have an injection from $X_{<p}^n(a)$ to $X_{<p}^{(\ell,n)}(\alpha)$ by sending the letter $x_i \in X_{<p}^n(a)$ to $w_i \in X_{<p}^{(\ell,n)}(\alpha)$. Since every word of length $\ell$ contained in $X_{<p}^{(\ell,n)}(\alpha)$ contains a letter in $X_{<p}^n(a)$, we have an application from $X_{<p}^{(\ell,n)}(\alpha)$ to $ X_{<p}^n(a)$ defined as follows. Let $w \in X_{<p}^{(\ell,n)}(\alpha)$ and let $j_w \in \{1,\ldots,|\zeta^n(a)|\}$ be the minimal integer such that $x_{j_w} \in X_{<p}^n(a)$ and $x_{j_w}$ is a letter in $w$. Then the application sends $w$ to $x_{j_w}$. By construction, the cardinal of the preimage of any $x \in X_{<p}^n(a)$ is at most equal to $\ell$. Therefore, we have $$\lvert\lvert\zeta^n(a)\rvert\rvert_{(p)} \leq ||M_{\ell}^n(p)(\alpha)|| \leq \ell\lvert\lvert\zeta^n(a)\rvert\rvert_{(p)}.$$ This concludes the proof.
\hfill\qedsymbol

\bigskip

\subsection{Construction of the attractive and repulsive currents for relative almost atoroidal automorphisms}

Let ${\tt n} \geq 3$ and let $\mathcal{F}=\{[A_1],\ldots,[A_k]\}$ be a free factor system of $F_{\tt n}$. We first define a class of outer automorphisms of $F_{\tt n}$ which we will study in the rest of the article. If $\phi \in \Out(F_{\tt n},\mathcal{F})$ and $\phi$ preserves the conjugacy class of every $A_i$ with $i \in \{1,\ldots,k\}$, we denote by $\phi|_{\mathcal{F}}$ the element $([\phi_1|_{A_1}],\ldots,[\phi_k|_{A_k}])$, where, for every $i \in \{1,\ldots,k\}$, the element $\phi_i$ is a representative of $\phi$ such that $\phi_i(A_i)=A_i$ and $[\phi_i|_{\mathcal{A}_i}]$ is an element of $\Out(A_i)$. Note that the outer class of $\phi_i|_{A_i}$ in $\Out(A_i)$ does not depend on the choice of $\phi_i$.

\begin{defi}\label{Defi almost atoroidal outer automorphism}
Let ${\tt n} \geq 3$ and let $\mathcal{F}=\{[A_1],\ldots,[A_k]\}$ be a free factor system of $F_{\tt n}$. Let $\phi \in \Out(F_{\tt n},\mathcal{F})$. The outer automorphism $\phi$ is \emph{almost atoroidal relative to} $\mathcal{F}$ if $\phi$ preserves the conjugacy class of every $A_i$ with $i \in \{1,\ldots,k\}$ and $\phi$ is one of the following:

\medskip

\noindent{$(1)$ } an atoroidal outer automorphism relative to $\mathcal{F}$. 

\medskip

\noindent{$(2)$ } an outer automorphism which preserves a sequence of free factor systems $\mathcal{F} \leq \mathcal{F}_1 \leq \{F_{\tt n}\}$ with $\mathcal{F}_1=\{[B_1],\ldots,[B_{\ell}]\}$ and such that:

\noindent{$(a)$ } $\mathcal{F}_1 \leq \{F_{\tt n}\}$ is sporadic, 

\noindent{$(b)$ } $\phi$ preserves the conjugacy class of every $B_i$ with $i \in \{1,\ldots,\ell\}$, the element $\phi|_{\mathcal{F}_1}$ is an expanding atoroidal outer automorphism relative to $\mathcal{F}$ and $\phi$ is not expanding relative to $\mathcal{F}$ ($\mathcal{F}$ might be equal to $\mathcal{F}_1$).
\end{defi}

The main example of an almost atoroidal automorphism is the following. Suppose that $\mathcal{F}_1=[A]$ and let $\phi \in \Out(F_{\tt n},\mathcal{F})$ be such that $\phi([A])=[A]$. Then $\phi$ is almost atoroidal if $\phi|_{[A]}$ is expanding relative to $\mathcal{F}$. Indeed, either $\phi$ is expanding relative to $\mathcal{F}$ and in this case $\phi$ satisfies~$(1)$ or $\phi$ is not expanding relative to $\mathcal{F}$ and $\phi$ satisfies ~$(2)$. Almost atoroidality allows us to deal with sporadic extensions. 

Let $\phi \in \Out(F_{\tt n},\mathcal{F})$ be an almost atoroidal outer automorphism relative to $\mathcal{F}$.  In this section, we construct a nontrivial convex compact subset in $\PCurr(F_{\tt n},\mathcal{F}\wedge  \mathcal{A}(\phi))$ associated with $\phi$. We follow the construction of \cite{Uyanik2019} in the context of atoroidal automorphisms. By Theorem~\ref{Theo existence CT}, there exists $M \geq 1$ such that $\phi^M$ is represented by a CT map $f \colon G \to G$ with filtration $\varnothing=G_0 \subsetneq G_1 \subsetneq \ldots \subsetneq G_k=G$ and such that there exists $p \in \{1,\ldots,k\}$ such that $\mathcal{F}(G_p)=\mathcal{F}$. For a splitting unit $\sigma$ in $G$, we say that $\sigma$ is \emph{expanding} if $\lim_{m \to \infty} \ell_{exp}([f^m(\sigma)])= +\infty$. Note that, by Lemma~\ref{Lem splitting units positive exp length}, this is equivalent to saying that there exists $N \in \NN^*$ such that $[f^N(\sigma)]$ contains a splitting unit which is an edge in an EG stratum. Moreover, a splitting unit $\sigma$ which is an expanding splitting unit is either an edge in $\overline{G-G_{PG}'}$ or a maximal taken connecting path in zero stratum such that a reduced iterate of $\sigma$ contains an edge in $\overline{G-G_{PG}'}$ as a splitting unit. In particular, there are finitely many expanding splitting units by Proposition~\ref{Prop definition CT}~$(3)$.

Let $\gamma$ and $\gamma'$ be two finite reduced subpaths of $G$. We denote by $\#(\gamma,\gamma')$ the number of occurrences of $\gamma$ in $\gamma'$ and by $\left\langle \gamma,\gamma' \right\rangle$ the sum \begin{equation}\label{Equation definition crochet}
\left\langle \gamma,\gamma' \right\rangle=\#(\gamma,\gamma')+\#(\gamma^{-1},\gamma'). 
\end{equation}
The next proposition shows the existence of relative currents associated with relative atoroidal outer automorphisms. Once we have constructed these currents for relative atoroidal outer automorphisms, we will also be able to construct attractive and repulsive simplices for every almost atoroidal outer automorphisms relative to $\mathcal{F}$. The proposition and its proof are inspired by the same result in the absolute context due to Uyanik (\cite[Proposition~3.3]{Uyanik2019}) and by the proof due to Gupta in the relative fully irreducible context (\cite[Proposition~~8.13]{gupta2017relative}). Recall the definition of $\mathcal{P}(\mathcal{F} \wedge \mathcal{A}(\phi))$ before Lemma~\ref{Lem path positive exp length cover double boundary} and $\mathscr{C}$ before Lemma~\ref{Lem finite set of words determines vertex group system}.

\begin{prop}\label{Prop existence relative currents atoroidal automorphisms}
Let ${\tt n} \geq 3$ and let $\mathcal{F}$ be a free factor system of $F_{\tt n}$. Let $\phi \in \Out(F_{\tt n},\mathcal{F})$ be an atoroidal outer automorphism relative to $\mathcal{F}$. Let $f \colon G \to G$ be a CT map that represents a power of $\phi$ with filtration $\varnothing=G_0 \subsetneq G_1 \subsetneq \ldots \subsetneq G_k=G$ and such that there exists $p \in \{1,\ldots,k\}$ such that $\mathcal{F}(G_p)=\mathcal{F}$. Let $\gamma \in \mathcal{P}(\mathcal{F} \wedge \mathcal{A}(\phi))$ and let $\sigma$ be an expanding splitting unit with fixed initial direction. 

\medskip

\noindent{$(1)$ } The limit $$\sigma_{\gamma}=\lim_{m \to \infty} \frac{\left\langle \gamma,[f^{m}(\sigma)] \right\rangle}{\ell_{\mathcal{F}}([f^{m}(\sigma)])}$$ exists and is finite.

\medskip

\noindent{$(2)$ } There exists a unique current $\eta_{\sigma} \in \Curr(F_{\tt n},\mathcal{F} \wedge \mathcal{A}(\phi))$ such that, for every finite reduced edge path $\gamma \in \mathcal{P}(\mathcal{F} \wedge \mathcal{A}(\phi))$, we have: $$\eta_{\sigma}(C(\gamma))=\sigma_{\gamma}.$$
\end{prop}

\dem $(1)$ We may suppose that $\gamma$ occurs in a reduced iterate of $\sigma$ as otherwise $\sigma_{\gamma}=0$. We first treat the case where $\sigma$ is an expanding splitting unit which is an edge in an irreducible stratum. Let $r$ be the height of $\sigma$. In order to prove the proposition in this case, we want to apply Proposition~\ref{Prop number of occurences substitution} to the CT map $f$ seen as a substitution on the set of splitting units contained in iterates of $\sigma$. However, the set of splitting units might be infinite since exceptional paths may have arbitrarily large widths and INPs arbitrarily large lengths. Instead, we construct a finite alphabet $A_{\gamma}$ \emph{depending on $\gamma$}. The alphabet is constructed as follows by associating a letter to every splitting unit occurring in a reduced iterate of $\sigma$. However some letters will correspond to infinitely many splitting units.

\medskip

\noindent{$(a)$ } We add one letter for each of the finitely many edges in irreducible strata that are contained in a reduced iterate of $\sigma$.

\medskip

\noindent{$(b)$ } We add one letter for each reduced maximal taken connecting path in a zero stratum contained in a reduced iterate of $\sigma$.

\medskip

\noindent{$(c)$ } We add one letter for each INP contained in a reduced iterate of $\sigma$ and such that the stratum of maximal height it intersects is an EG stratum.

\medskip

\noindent{$(d)$ } Let $\delta$ be an INP such that the stratum of maximal height it intersects is an NEG stratum and such that it appears in a reduced iterate of $\sigma$. By Proposition~\ref{Prop definition CT}~$(11)$, there exist an edge $e$, an integer $s \in \ZZ$ and a closed Nielsen path $w$ such that $\delta=ew^se^{-1}$. Note that $\gamma$ is not contained in $w^s$ since $\gamma \in \mathcal{P}(\mathcal{F}\wedge \mathcal{A}(\phi))$ and $w^s$ is a concatenation of paths in $G_{PG,\mathcal{F}}$ and $\mathcal{N}_{PG,\mathcal{F}}$ by Lemma~\ref{Lem NEG INP in Npg} and the fact that $\phi$ is atoroidal relative to $\mathcal{F}$. Hence if $\gamma$ is contained in $\delta$, it is either an initial or a terminal segment of $\delta$. Let $M_1$ be the maximal integer $|d|$ such that $\gamma$ contains an INP of the form $ew^de^{-1}$. Let $M_2$ be the minimal integer $|d|$ such that $\gamma \cap (ew^de^{-1})$ is either an initial or a terminal segment of $ew^de^{-1}$. Let $M_3$ be the maximal integer $|d|$ such that $ew^de^{-1}$ is contained in $[f(\sigma')]$ with $\sigma'$ a splitting unit which is either an edge in an irreducible stratum or a maximal taken connecting path in a zero stratum. Let $M=\max\{M_1,M_2,M_3\}$. We add one letter for each $ew^de^{-1}$ with $|d|\leq M+1$. We add exactly one letter representing every $ew^de^{-1}$ with $|d| >M+1$.

\medskip

\noindent{$(e)$ }  Let $\delta$ be an exceptional path appearing in a reduced iterate of $\sigma$. There exist edges $e_1,e_2$, a nonzero integer $s$ and a closed Nielsen path $w$ such that $\delta=e_1w^se_2^{-1}$. Note that $\gamma$ is not contained in $w^s$ since $\gamma \in \mathcal{P}(\mathcal{F}\wedge \mathcal{A}(\phi))$ and $w^s$ is a concatenation of paths in $G_{PG,\mathcal{F}}$ and $\mathcal{N}_{PG,\mathcal{F}}$ by Lemma~\ref{Lem NEG INP in Npg} and the fact that $\phi$ is atoroidal relative to $\mathcal{F}$. Let $M_4$ be the maximal integer $|d|$ such that $\gamma$ contains an exceptional path of the form $e_1w^de_2^{-1}$. Let $M_5$ be the minimal integer $|d|$ such that $\gamma \cap e_1w^de_2^{-1}$ is either a proper initial or terminal segment of $e_1w^de_2^{-1}$. Let $M_6$ be the maximal integer $|d|$ such that $e_1w^de_2^{-1}$ is contained in $[f(\sigma')]$ with $\sigma'$ a splitting unit which is either an edge in an irreducible stratum or a maximal taken connecting path in a zero stratum. Let $M'=\max\{M_4,M_5,M_6\}$. We add one letter for each $e_1w^de_2^{-1}$ with $|d|\leq M'+1$. We add one letter representing every $e_1w^de_2^{-1}$ with $|d| >M'+1$.

\medskip

We claim that the alphabet $A_{\gamma}$ is finite. Indeed, since the graph $G$ is finite, so is the number of letters in the first category. By Proposition~\ref{Prop definition CT}~$(3)$, the zero strata of $G_{r-1}$ are exactly the contractible components of $G_{r-1}$. Hence the number of letters in the second category is finite. The number of letters in the third category is finite by Proposition~\ref{Prop definition CT}~$(9)$. The remaining letters of $A_{\gamma}$ are finite by definition. Let $\zeta$ be the following substitution on $A_{\gamma}$. If $a \in A_{\gamma}$ represents a unique path in $G$, we set $\zeta(a)=[f(a)]$. If $a \in A_{\gamma}$ represents several paths in $G$, we set $\zeta(a)=a$. We claim that $\zeta$ is a well-defined substitution. Indeed, by Proposition~\ref{Prop definition CT}~$(6)$, if $a$ is a letter in $A_{\gamma}$ which represents a unique path in $G$, then $[f(a)]$ is completely split and every splitting unit in $[f(a)]$ is represented by a unique letter by the construction of letters in the fourth and fifth category. Moreover, if $a \in A_{\gamma}$ represents several paths, then the definition of $\zeta$ does not depend on the choice of a representative of $a$. Hence $\zeta$ is a well-defined substitution. 

We claim that if $a \in A_{\gamma}$ represents several paths in $G$, then, for every representative $\alpha$ of $a$, the path $[f(\alpha)]$ is represented by $a$. Indeed, the claim is immediate when $a$ represents several INPs, so we focus on the case where $a$ represents several exceptional paths. Let $e_1,e_2$ be edges in $G$, let $w$ be a closed Nielsen path in $G$ and let $d \in \ZZ$ be such that $e_1w^de_2^{-1}$ is represented by the letter $a$. There exist a splitting unit $\sigma'$ of a reduced iterate of $\sigma$ by $[f]$, an integer $N \in \NN^*$ and an integer $d_1 \in \ZZ$ such that $e_1w^{d_1}e_2^{-1}$ is a subpath of $[f^N(\sigma')]$. Thus, using the constants given in $(e)$, we have $|d_1| \leq M_6 \leq M$. By the construction of the alphabet $A_{\gamma}$, there exists a letter $a'$ in $A_{\gamma}$ corresponding to the path $e_1w^{d_1}e_2^{-1}$ and $a'$ represents a unique path. For every $n \in \NN$, let $d_n \in \ZZ$ be such that $[f^n(e_1w^{d_1}e_2^{-1})]=e_1w^{d_n}e_2^{-1}$. Then the sequence $(d_n)_{n \in \NN}$ is monotonic. Let $m_0$ be the minimal integer such that the path $e_1w^{d_{m_0}}e_2^{-1}$ is represented by $a$. Note that $m_0 >1$ as $a'$ represents a unique path. By monotonicity, $d_{m_0} \neq d_1$. Thus, if $d_{m_0}> d_1$, then for every $m \geq m_0$, we have $d_m \geq d_{m_0}$ and if $d_{m_0}< d_1$, then for every $m \geq m_0$, we have $d_m \leq d_{m_0}$. Hence for every $m \geq m_0$, the path $e_1w^{d_{m+1}}e_2^{-1}$ is represented by $a$. This shows that if $\alpha \in a$ then $[f(\alpha)] \in a$. This concludes the proof of the claim. Hence $\zeta$ only depends on the function $[f(.)]$. 

By reordering columns and rows, we may suppose that, if $M$ is the matrix associated with $\zeta$, then columns and rows of $M$ with index greater than $p$ are precisely the letters in $A_{\gamma}$ representing splitting units which are concatenations of paths in $G_{PG,\mathcal{F}}$ and $\mathcal{N}_{PG,\mathcal{F}}$. By Lemma~\ref{Lem iterate of a path in GPG}, iterates by $\zeta$ of letters of $A_{\gamma}$ representing concatenations of paths in $G_{PG,\mathcal{F}}$ and $\mathcal{N}_{PG,\mathcal{F}}$ are words on $A_{\gamma}$ whose letters represent  concatenations of paths in $G_{PG,\mathcal{F}}$ and $\mathcal{N}_{PG,\mathcal{F}}$. Thus, the matrix $M$ is a lower block triangular matrix, where every block of index at most $p$ corresponds to either edges in a common stratum, or the $0$ matrix when the associated letter is a maximal taken connecting path in a zero stratum.

Since $\sigma$ is expanding, it has a reduced iterate which contains splitting units which are edges in EG strata. Hence if $a_{\sigma}$ is the letter in $A_{\gamma}$ corresponding to $\sigma$, the iterates $\zeta^n(a_{\sigma})$ contain letters of $A_{\gamma}$ in a Perron-Frobenius block with eigenvalue greater than $1$. Since the initial direction of $\sigma$ is fixed by Proposition~\ref{Prop number of occurences substitution}, for every word $w$ in the alphabet $A_{\gamma}$, the limit $$\lim_{m \to \infty} \frac{( w,[\zeta^{m}(\sigma)])}{\lvert\lvert\zeta^{m}(\sigma)\rvert\rvert_{(p)}}$$ exists and is finite. Hence the limit $$\lim_{m \to \infty} \frac{\left\langle w,[\zeta^{m}(\sigma)] \right\rangle}{\lvert\lvert\zeta^{m}(\sigma)\rvert\rvert_{(p)}}$$ exists and is finite. 

\medskip 

\noindent{\bf Claim. } There exists a matrix $M'$ obtained from $M$ by multiplying rows and columns by positive scalars and such that, for every $m \in \NN^*$, we have $\ell_{\mathcal{F}}([f^m(\sigma)])=\lVert M'^m(\sigma)\rVert_{(p)}$.

\medskip 

\dem Remark that if $e_1w^se_2^{-1}$ is an exceptional path, and if $e_1w^de_2^{-1}$ is an exceptional path with distinct width, then their $\mathcal{F}$-lengths are equal and at most equal to $2$. Indeed, since $\phi$ is an atoroidal outer automorphism relative to $\mathcal{F}$, every closed Nielsen path of $G$ is contained in $G_p$. Since $w$ is a closed Nielsen path, we see that $w$ is a concatenation of paths in $G_{PG,\mathcal{F}}$ and $\mathcal{N}_{PG,\mathcal{F}}$ by Lemma~\ref{Lem closed nielsen paths in Npg'}. Hence we have $$\ell_{\mathcal{F}}(e_1w^se_2^{-1})=\ell_{\mathcal{F}}(e_1)+\ell_{\mathcal{F}}(e_2) \leq 2.$$ Similarly, if $ew^se^{-1}$ and $ew^de^{-1}$ are INP intersecting the same maximal NEG stratum, then their $\mathcal{F}$-length are equal and at most equal to $2$. Let $M'$ be the matrix obtained from $M$ by multiplying every row correponding to either an exceptional path not contained in $G_p$, an INP not contained in $G_p$, a collection of exceptional paths not contained in $G_p$, a collection of INPs not contained in $G_p$ or a maximal taken connecting path not contained in $G_p$, by the corresponding $\mathcal{F}$-length. Note that, by the above remarks, this does not depend on the choice of a representative when the letter corresponds to a collection of paths. Then for every $m \in \NN^*$, the value $\lVert M'^m(\sigma)\rVert_{(p)}$ corresponds to the sum of the $\mathcal{F}$-length of every splitting unit in $[f^m(\sigma)]$ not contained in $G_p$. By Lemma~\ref{Lem completely split and paths in Npg}, complete splittings are $PG$-relative complete splittings. By Lemma~\ref{Lem compute exp length completely split}~$(2)$, we have $\ell_{\mathcal{F}}([f^m(\sigma)])=\lVert M'^m(\sigma)\rVert_{(p)}$. This proves the claim.
\hfill\qedsymbol

\medskip

By the claim, we see that for every $m \in \NN^*$, there exists a constant $K$ such that we have $$\frac{1}{K}\lvert\lvert\zeta^{m}(\sigma)\rvert\rvert_{(p)}\leq \ell_{\mathcal{F}}([f^{m}(\sigma)])\leq K\lvert\lvert\zeta^{m}(\sigma)\rvert\rvert_{(p)}.$$ Using the claim in the proof of Proposition~\ref{Prop number of occurences substitution} (replacing $\lVert M_{\ell}^n(p)(\alpha) \rVert$ by $\ell_{\mathcal{F}}([f^n(\sigma)])$ which is possible since $\ell_{\mathcal{F}}([f^n(\sigma)])$ is the norm of a matrix by the claim), the limit $$\lim_{m \to \infty} \frac{\left\langle w,[f^{m}(\sigma)] \right\rangle}{\ell_{\mathcal{F}}([f^{m}(\sigma)])}$$ exists and is finite. We now construct a finite set of words $W(\gamma)$ in the alphabet $A_{\gamma}$ such that for every $m \in \NN^*$, there exists a bijection between occurrences of $\gamma$ in $[f^m(\sigma)]$ and occurrences of a word $w \in W(\gamma)$ in $[\zeta^{m}(\sigma)]$. This will conclude the proof of Case~2. Let $W(\gamma)$ be the set of words in $A_{\gamma}$ consisting in every path contained in a reduced iterate of $\sigma$ which contains $\gamma$, which is completely split and which is minimal for these properties. By construction, every occurrence of $\gamma$ in a reduced iterate of $\sigma$ is contained in a word in $W(\gamma)$. We claim that the set $W(\gamma)$ is finite. Indeed, let $w$ be a word in $W(\gamma)$. Then $w$ corresponds to a path in a reduced iterate of $\sigma$ which is a concatenation of splitting units $w=\sigma_1\ldots \sigma_k$. By minimality of $w$, if $w' \in W(\gamma)$ is distinct from $w'$, then the number of splitting units in $w'$ is at most equal to $k$ and $w'$ might differ from $w$ by changing $\sigma_1$ and $\sigma_k$. Thus, $W(\gamma)$ is finite. For every $w \in W(\gamma)$, let $m_w$ be the number of occurrences of $\gamma$ in $w$. Since $\gamma$ is not contained in $G_p$, the value $m_w$ does not depend on the choice of a representative of $w$ if $w$ represents a collection of paths. Therefore, for every $m \in \NN^*$, we have $$\left\langle \gamma, f^{m}(\sigma) \right\rangle= \sum_{w \in W(\gamma)} m_w\left\langle w, f^{m}(\sigma) \right\rangle.$$ This shows that the limit  $$\sigma_{\gamma}=\lim_{m \to \infty} \frac{\left\langle \gamma,f^{m}(\sigma) \right\rangle}{\ell_{\mathcal{F}}(f^{m}(\sigma))}$$ exists and is finite. This proves Assertion~$(1)$ of the proposition when $\sigma$ is an edge in an irreducible stratum.

Suppose now that $\sigma$ is a maximal taken connecting path in a zero stratum. We prove the proposition by induction on the height $r$ of the splitting unit $\sigma$. Suppose first that $\sigma$ is an expanding splitting unit which is a maximal taken connecting path in a zero stratum of minimal height $r$. Then $[f(\sigma)]$ has height $r-1$, hence it does not contain splitting units which are maximal taken connecting path in zero strata. In this case, the proof follows from the above case. Suppose now that $\sigma$ is a maximal taken connecting path in a zero stratum. Then its reduced image is completely split and has height at most $r-1$. In this case the claim follows by induction applied to $[f(\sigma)]$. This concludes the proof of Assertion~$(1)$.

\medskip

\noindent{$(2)$ } Let us prove that for every element $\gamma \in \mathcal{P}(\mathcal{F} \wedge \mathcal{A}(\phi))$, we have:

\medskip

\noindent{$(i)$ } $0 \leq \sigma_{\gamma} <\infty$;

\medskip

\noindent{$(ii)$ } $\sigma_{\gamma}=\sigma_{\gamma^{-1}}$;

\medskip

\noindent{$(iii)$ } $\sigma_{\gamma}=\sum_{e \in E} \sigma_{\gamma e}$, where $E$ is the subset of $\vec{E}G$ consisting in all edges that are incident to the endpoints of $\gamma$ and distinct from the inverse of the last edge of $\gamma$.

\medskip

The point $(i)$ follows from Assertion~$(1)$. The second point follows from the definition of $\left\langle \gamma, f^m(\sigma) \right\rangle$. In order to prove the third point, remark that $\left\langle \gamma,f^m(\sigma) \right\rangle$ and $\sum_{e \in E} \left\langle\gamma e,f^n(\sigma)\right\rangle$ differ only when $[f^m(\sigma)]$ ends with $\gamma$ or $\gamma^{-1}$. Therefore the difference between $\left\langle \gamma,f^m(\sigma) \right\rangle $ and $\sum_{e \in E} \left\langle \gamma e,f^m(\sigma) \right\rangle$ is at most $2$. This implies that 
$$\left\lvert \frac{\left\langle \gamma,f^m(\sigma) \right\rangle}{\ell_{\mathcal{F}}(f^m(\sigma))} - \sum_{e \in E} \frac{\left\langle \gamma e,f^m(\sigma) \right\rangle}{\ell_{\mathcal{F}}(f^m(\sigma))} \right\rvert \to 0 \text{ as } n \to \infty.$$ This proves the third point. By \cite[Lemma~3.2]{Guerch2021currents}, since the map $\gamma \mapsto \sigma_{\gamma}$ satisfies the conditions $(i)-(iii)$, it determines a projective relative current $n_{\sigma} \in \PCurr(F_{\tt n},\mathcal{F})$. This current is unique since a relative current is entirely determined by its set of values on cylinders of finite paths $\gamma \in \mathcal{P}(\mathcal{F} \wedge \mathcal{A}(\phi))$ by Lemma~\ref{Lem path positive exp length cover double boundary}. This concludes the proof.
\hfill\qedsymbol

\bigskip

\begin{defi}\label{Defi attractive convex}
Let ${\tt n} \geq 3$ and let $\mathcal{F}$ be a free factor system of $F_{\tt n}$. Let $\phi \in \Out(F_{\tt n},\mathcal{F})$ be an almost atoroidal outer automorphism relative to $\mathcal{F}$ and let $\mathcal{F}_1$ be a free factor system such that $\mathcal{F} \leq \mathcal{F}_1$ and such that the extension $\mathcal{F}_1 \leq \{F_{\tt n}\}$ is sporadic and such that $\phi|_{\mathcal{F}_1}$ is atoroidal relative to $\mathcal{F}$. In the case that $\phi$ is atoroidal relative to $\mathcal{F}$, we assume that $\mathcal{F}_1=\{[F_{\tt n}]\}$. Let $f \colon G \to G$ be a CT map representing a power of $\phi$ with filtration $$\varnothing=G_0 \subsetneq G_1 \subsetneq \ldots \subsetneq G_k=G,$$ such that there exists $i \in \{1,\ldots,k-1\}$ with $\mathcal{F}(G_i)=\mathcal{F}_1$. We define the \emph{simplex of attraction of $\phi$}, denoted by $\Delta_{+}(\phi)$, as the set of projective classes of nonnegative linear combinations of currents $\mu_{\sigma}$ obtained from Proposition~\ref{Prop existence relative currents atoroidal automorphisms} applied to $\phi|_{\mathcal{F}_1}$ and $f$ and which correspond to splitting units $\sigma$ whose exponential length grows exponentially fast under iteration of $f$. The \emph{simplex of repulsion} of $\phi$, denoted by $\Delta_-(\phi)$, is $\Delta_{+}(\phi^{-1})$.
\end{defi}

\begin{rmq}\label{Rmq attractive repulsive currents}
The definitions of attractive and repulsive currents given in Definition~\ref{Defi attractive convex} rely on the choice of CT maps representing powers of the almost atoroidal outer automorphisms $\phi$ and $\phi^{-1}$. However, it will be a consequence of Proposition~\ref{Prop invariance 3} and Proposition~\ref{Prop North south dynamics outside neighborhood} that the attractive and repulsive currents depend only on $\phi$.
\end{rmq}

We now prove properties of the subsets $\Delta_{\pm}(\phi)$. As explained above Proposition~\ref{Prop existence relative currents atoroidal automorphisms}, there are only finiely many expanding splitting units. Hence the subsets $\Delta_{\pm}(\phi)$ are closed. Since $\PCurr(F_{\tt n},\mathcal{F}\wedge \mathcal{A}(\phi))$ is a Hausdorff, compact space by Lemma~\ref{Lem PCurr compact} and since $\Delta_{\pm}(\phi)$ are closed subsets, we have the following.

\begin{lem}\label{Lem Delta compact}
Let ${\tt n} \geq 3$ and let $\mathcal{F}$ be a free factor system of $F_{\tt n}$. Let $\phi \in \Out(F_{\tt n},\mathcal{F})$ be an atoroidal outer automorphism relative to $\mathcal{F}$. The subsets $\Delta_{\pm}(\phi)$ are compact and contain finitely many extremal points.
\hfill\qedsymbol
\end{lem}

Note that one compute $\left\lVert\mu(\sigma)\right\rVert_{\mathcal{F}}$ by counting the number of occurrences of every $PG$-relative splitting unit of positive $\mathcal{F}$-length in a reduced iterate of $\sigma$ and taking the limit. This is precisely the limit of the $\mathcal{F}$-length of reduced iterates of $\sigma$ by Lemma~\ref{Lem compute exp length completely split}. Hence we have the following result.

\begin{lem}\label{Lem Delta value norm}
Let ${\tt n} \geq 3$ and let $\mathcal{F}$ be a free factor system of $F_{\tt n}$. Let $\phi \in \Out(F_{\tt n},\mathcal{F})$ be an atoroidal outer automorphism relative to $\mathcal{F}$. We have $\left\lVert\mu(\sigma)\right\rVert_{\mathcal{F}}=1$.
\hfill\qedsymbol
\end{lem}

We now prove that the subsets $\Delta_{\pm}(\phi)$ are $\phi$-invariant. We first recall some lemmas.

\begin{lem}\cite[Bounded Cancellation]{Cooper87}\label{Lem Bounded cancellation lemma}
Let ${\tt n} \geq 2$ and let $G$ be a marked graph of $F_{\tt n}$. Let $f \colon G \to G$ be a graph map. There exists a constant $C_f$ such that for any reduced path $\rho=\rho_1\rho_2$ in $G$ we have $$\ell([f(\rho)]) \geq \ell([f(\rho_1)]) + \ell([f(\rho_2)])-2C_f.$$
\end{lem}

\begin{lem}\cite[Lemma~5.7]{LustigUyanik2019} \label{Lem invariance 1}
For any graph $G$ without valence $1$ vertices there exists a constant $K \geq 0$ such that for any finite reduced edge path $\gamma$ in $G$ there exists an edge path $\gamma'$ of length at most $K$ such that the concatenation $\gamma\gamma'$ exists and is a reduced circuit.
\end{lem}

\begin{lem}\label{Lem invariance 2}
Let $f \colon G \to G$ be as in Proposition~\ref{Prop existence relative currents atoroidal automorphisms}. Let $K_1 \geq 0$ be any constant, let $\sigma$ be an expanding splitting unit and let $\eta_{\sigma}$ be the current associated with $\sigma$ given by Proposition~\ref{Prop existence relative currents atoroidal automorphisms}~$(2)$. Let $m \in \NN$ and let $\gamma_m'$ be a reduced edge path of length at most $K_1$. Let $\gamma_m=[f^m(\sigma)]^*\gamma_m'$, where $[f^m(\sigma)]^*$ is obtained from $[f^m(\sigma)]$ by erasing an initial and a terminal subpath of length $K_1$. For every element $\gamma \in \mathcal{P}(\mathcal{F}\wedge \mathcal{A}(\phi))$, we have $$\lim_{m \to \infty} \frac{\left\langle \gamma,\gamma_m \right\rangle}{\ell_{\mathcal{F}}(\gamma_m)}=\left\langle \gamma,\eta_{\sigma} \right\rangle.$$
\end{lem}

\dem The proof follows \cite[Lemma~5.8]{LustigUyanik2019}. Since $\ell(\gamma_m') \leq K_1$, we have $$\ell_{\mathcal{F}}([f^m(\sigma)]^*) \geq \ell_{\mathcal{F}}([f^m(\sigma)])-2K_1.$$ Since $\sigma$ is expanding, we have $\lim_{m \to \infty}\ell_{\mathcal{F}}([f^m(\sigma)])=+\infty$. Hence we have $$\lim_{m \to \infty} \frac{\left\langle \gamma, \gamma_m \right\rangle}{\left\langle \gamma,[f^m(\sigma)] \right\rangle}=1$$ and $$\lim_{m \to \infty} \frac{\ell_{\mathcal{F}}(\gamma_m)}{\ell_{\mathcal{F}}([f^m(\sigma)])}=1.$$ Hence the result follows from Proposition~\ref{Prop existence relative currents atoroidal automorphisms}~$(1)$.
\hfill\qedsymbol

\begin{prop}\label{Prop invariance 3}
Let ${\tt n} \geq 3$ and let $\mathcal{F}$ be a free factor system of $F_{\tt n}$. Let $\phi \in \Out(F_{\tt n},\mathcal{F})$ be an almost atoroidal outer automorphism relative to $\mathcal{F}$. Let $f \colon G \to G$ be as in Proposition~\ref{Prop existence relative currents atoroidal automorphisms}. Let $\sigma$ be an expanding splitting unit and let $\eta_{\sigma}$ be the current associated with $\sigma$ given by Proposition~\ref{Prop existence relative currents atoroidal automorphisms}~$(2)$. There exists $\lambda_{\sigma} >0$ such that $$\phi(\eta_{\sigma})=\lambda_{\sigma}\eta_{\sigma}.$$
\end{prop}

\dem The proof follows \cite[Proposition~5.9]{LustigUyanik2019}. Let $K \geq 0$ be the constant associated with $G$ given by Lemma~\ref{Lem invariance 1}. Let $m \in \NN$, and let $\gamma_m'$ be the path of length at most $K$ given by Lemma~\ref{Lem invariance 1} such that $\gamma_m=[f^m(\sigma)]\gamma_m'$  is a reduced circuit. Since $\lim_{t \to \infty} \ell_{exp}([f^t(\sigma)])=+\infty$, for large values of $m$, we have $\ell_{exp}(\gamma_m) >0$. Let $w_m$ be an element of $F_{\tt n}$ whose conjugacy class is represented by $\gamma_m$. Note that, by Lemma~\ref{Lem psi0 equal exp length}, we have $\ell_{\mathcal{F}}(\gamma_m)=\lVert \eta_{w_m} \rVert_{\mathcal{F}}$. By Proposition~\ref{Prop circuits in Gpg are elements in poly subgroup}, since $\ell_{exp}(\gamma_m)>0$, we see that $w_m$ is $\mathcal{F}\wedge \mathcal{A}(\phi)$-nonperipheral, hence $w_m$ defines a current $\eta_{[w_m]} \in \Curr(F_{\tt n},\mathcal{F}\wedge \mathcal{A}(\phi))$.

Let $\alpha_m=[f^{m+1}(\sigma)][f(\gamma_m')]$. Note that since $\ell(\gamma_m') \leq K$, the value $\ell([f(\gamma_m')])$ is bounded by a constant $K_0$ which only depends on $K$. Let $C'$ be the constant given by Lemma~\ref{Lem Bounded cancellation lemma} and let $K_1=\max\{K_0,C'\}$. Then, with the notations of Lemma~\ref{Lem invariance 2}, the reduced circuit $\gamma_m''=[\alpha_m]$ can be written as a product $\gamma_m''=[f^m(\sigma)]^*\beta_m$ where $\ell(\beta_m) \leq K_1$ and $\ell_{\mathcal{F}}([f^m(\sigma)]^*) \geq \ell_{\mathcal{F}}([f^m(\sigma)])-2K_1$. Applying Lemma~\ref{Lem invariance 2} twice, we see that, for every element $\gamma \in \mathcal{P}(\mathcal{F}\wedge \mathcal{A}(\phi))$, we have
$$\lim_{m \to \infty} \frac{\left\langle \gamma,\gamma_m \right\rangle}{\ell_{\mathcal{F}}(\gamma_m)}=\left\langle \gamma,\eta_{\sigma} \right\rangle$$ and $$\lim_{m \to \infty} \frac{\left\langle \gamma,\gamma_m'' \right\rangle}{\ell_{\mathcal{F}}(\gamma_m'')}=\left\langle \gamma,\eta_{\sigma} \right\rangle.$$ By Lemma~\ref{Lem Caracterisation limites PCurr}, we have $$\lim_{m \to \infty} \frac{\eta_{[w_m]}}{\lVert \eta_{[w_m]} \rVert_{\mathcal{F}}}=\eta_{\sigma}.$$ From the continuity of the $\Out(F_{\tt n})$-action on $\PCurr(F_{\tt n},\mathcal{F}\wedge\mathcal{A}(\phi))$ and from $\phi(\eta_{\eta_{[w_m]}})=\eta_{\phi([w_m])}$, we see that $$\lim_{m \to \infty} \frac{\eta_{\phi([w_m])}}{\lVert \eta_{[w_m]} \rVert_{\mathcal{F}}}=\phi(\eta_{\sigma}).$$ Since the reduced circuit $\gamma_m''$ represents the conjugacy class $\phi([w_m])$, the second of the above equalities implies that $$\lim_{m \to \infty} \frac{\eta_{\phi([w_m])}}{\lVert \eta_{\phi([w_m])} \rVert_{\mathcal{F}}}=\eta_{\sigma}.$$ Recall that $\lim_{m \to \infty} \frac{\ell_{\mathcal{F}}(\gamma_m)}{\ell_{\mathcal{F}}([f^m(\sigma)])}=1$, that $\lim_{m \to \infty} \frac{\ell_{\mathcal{F}}(\gamma_m'')}{\ell_{\mathcal{F}}([f^{m+1}(\sigma)])}=1$, that $\ell_{\mathcal{F}}(\gamma_m)=\lVert \eta_{[w_m]} \rVert_{\mathcal{F}}$ and that $\ell_{\mathcal{F}}(\gamma_m'')=\lVert \eta_{\phi([w_m])} \rVert_{\mathcal{F}}$. Recall from the claim in the proof of Proposition~\ref{Prop existence relative currents atoroidal automorphisms} that $\ell_{\mathcal{F}}([f(\sigma)])$ is the norm of a matrix. The conclusion of Proposition~\ref{Prop invariance 3} then follows from the fact (see~\cite[Remark~3.3]{LustigUyanik2017}) that there exists $\lambda_{\sigma}>0$ such that $$\lim_{m \to \infty}\frac{\ell_{\mathcal{F}}([f^{m+1}(\sigma)])}{\ell_{\mathcal{F}}([f^m(\sigma)])}=\lambda_{\sigma}.\eqno\qedsymbol$$

We now prove a lemma which will be used in \cite{Guerch2021Polygrowth}.

\begin{lem}
Let ${\tt n} \geq 3$ and let $\mathcal{F}$ be a free factor system of $F_{\tt n}$. Let $\phi \in \Out(F_{\tt n},\mathcal{F})$ be an expanding outer automorphism relative to $\mathcal{F}$. Let $f \colon G \to G$ be as in Proposition~\ref{Prop existence relative currents atoroidal automorphisms}. Let $\sigma$ be an expanding splitting unit and let $\eta_{\sigma}$ be the current associated with $\sigma$ given by Proposition~\ref{Prop existence relative currents atoroidal automorphisms}~$(2)$. 

\medskip

\noindent{$(1)$ } There exists a projective current $[\eta] \in \PCurr(F_{\tt n},\mathcal{F} \wedge \mathcal{A}(\phi))$ whose support is contained in the support of $\eta_{\sigma}$ and such that $\Supp(\eta)$ is uniquely ergodic. In particular, the support of every extremal current of $\Delta_{\pm}(\phi)$ contains a closed subset which is uniquely ergodic.

\medskip

\noindent{$(2)$ } There exists only finitely many projective currents $[\eta] \in \PCurr(F_{\tt n},\mathcal{F} \wedge \mathcal{A}(\phi))$ whose support is contained in the support of $\eta_{\sigma}$ and such that $\Supp(\eta)$ is uniquely ergodic.
\end{lem}

\dem {$(1)$ } Note that, since $\phi$ is expanding relative to $\mathcal{F}$, we have $\mathcal{F} \wedge \mathcal{A}(\phi)=\mathcal{A}(\phi)$. Let $r \in \NN$ be the minimal integer such that $H_r$ is an EG stratum and a reduced iterate of $\sigma$ contains an edge of $H_r$. Such a stratum $H_r$ exists since $\sigma$ is expanding. Let $e$ be an edge of $H_r$ with fixed initial direction and let $\eta_e$ be the current in $\PCurr(F_{\tt n},\mathcal{A}(\phi))$ associated with $e$ given by Proposition~\ref{Prop existence relative currents atoroidal automorphisms}~$(2)$.

\medskip

\noindent{\bf Claim. } The support of $\eta_e$ is uniquely ergodic.

\medskip
 
\dem By minimality of $r$, every edge contained in a reduced iterate of $e$ is either in $H_r$ or in $G_{PG}'$. Let $G'$ be the minimal subgraph of $G$ which contains every reduced iterate of $e$ and let $A$ be a subgroup of $F_{\tt n}$ such that $\pi_1(G')$ is a conjugate of $A$ when $\pi_1(G)$ is identified with $F_{\tt n}$. Then $G'$ is $f$-invariant and hence $[A]$ is $\phi$-invariant. Let $G_1',\ldots,G_k'$ be the connected component of $\overline{G'-H_r}$ and let $\mathcal{F}'$ be the free factor system of $F_{\tt n}$ determined by $G_1',\ldots,G_k'$. Let $\Phi \in \phi$ be such that $\Phi(A)=A$. Note that $[\Phi|_{A}] \in \Out(A)$ is fully irreducible relative to $\mathcal{F}'$. For every $i \in \{1,\ldots,k\}$, we have $G_i' \subseteq G_{PG}'$. By Proposition~\ref{Prop definition CT}~$(3)$, for every $i \in \{1,\ldots,k\}$, either $G_i'$ is contractible or $G_i' \subseteq G_{PG}$. By Proposition~\ref{Prop circuits in Gpg are elements in poly subgroup} for every subgroup $H$ of $F_{\tt n}$ such that $[H] \in \mathcal{F}'$, there exists a subgroup $H'$ of $F_{\tt n}$ such that $[H'] \in \mathcal{A}(\phi)$ and $H \subseteq H'$. Hence we have $\mathcal{F}' \leq \mathcal{A}(\phi)$. Moreover by Proposition~\ref{Prop circuits in Gpg are elements in poly subgroup} and Proposition~\ref{Prop definition CT}~$(9)$, if $\gamma$ is a cyclically reduced circuit of $G'$ of height $r$ whose growth under iteration of $f$ is polynomial, then $\gamma$ contains (up to taking inverse) the only height $r$ EG INP $\sigma_r$. As one of the endpoints of $\sigma_r$ is not contained in $G_{r-1}$ by~\cite[Fact~I.1.42]{HandelMosher20}, we see that either $\sigma_{r}$ is not closed and $\gamma$ does not exist or $\sigma_r$ is closed and $\gamma$ is an iterate of $\sigma_r$ or $\sigma_r^{-1}$. Let $b \in F_{\tt n}$ be the (possibly trivial) element associated with $\sigma_r$. Then, we have $$\partial^2 A \cap \partial^2(F_{\tt n},\mathcal{A}(\phi))=\partial^2(A,\mathcal{F}'\cup \{[b]\}).$$ Let $\PCurr(\Supp(\eta_e))$ be the set of projective currents in $\PCurr(F_{\tt n},\mathcal{F} \wedge \mathcal{A}(\phi))$ whose support is contained in $\Supp(\eta_e)$. We now construct an injective application $$\Theta \colon \PCurr(\Supp(\eta_e)) \to \PCurr(A,\mathcal{F}')$$ such that for every projective current $\mu \in \PCurr(\Supp(\eta_e))$ we have $$\Supp(\Theta([\mu]))=\Supp([\mu]) \cap \partial^2 A.$$ Let $\mathcal{P}(\mathcal{F}')$ be the set of paths in a Cayley tree of $F_{\tt n}$ defined above Lemma~\ref{Lem path positive exp length cover double boundary} (replacing $\mathcal{F} \wedge \mathcal{A}(\phi)$ by $\mathcal{F}'$). Let $\mathcal{P}_A(\mathcal{F}')$ be the set of paths in $\mathcal{P}(\mathcal{F}')$ contained in $G'$. By Lemma~\ref{Lem path positive exp length cover double boundary}, the set consisting in elements $C(\gamma)$ with $\gamma \in \mathcal{P}(\mathcal{F}')$ covers $\partial^2(A,\mathcal{F}')$. Thus, by \cite[Lemma~3.2]{Guerch2021currents}, it suffices to prove that for every projective current $\eta \in \PCurr(\Supp(\eta_e))$, we can associate a function $\widetilde{\eta} \colon \mathcal{P}_A(\mathcal{F}') \to \RR $ such that for every $\gamma \in  \mathcal{P}_A(\mathcal{F}')$, we have

\medskip

\noindent{$(i)$ } $0 \leq \widetilde{\eta}(\gamma) <\infty$;

\medskip

\noindent{$(ii)$ } $\widetilde{\eta}(\gamma)=\sigma_{\gamma^{-1}}$;

\medskip

\noindent{$(iii)$ } $\widetilde{\eta}(\gamma)=\sum_{e \in E} \sigma_{\gamma e}$, where $E$ is the subset of $\vec{E}G'$ consisting in all edges that are incident to the endpoints of $\gamma$ and distinct from the inverse of the last edge of $\gamma$.

\medskip

Let $\eta \in \PCurr(\Supp(\eta_e))$. If $\gamma \in \mathcal{P}_A(\mathcal{F}')$ is not contained in the axis of a conjugate of $b$, we may set $\widetilde{\eta}(\gamma)=\eta(C(\gamma))$. Since $\sigma_e$ is $r$-legal, a reduced iterate of $\sigma_e$ cannot contain the only height $r$ EG INP. Thus, we may set, for every path $\gamma \in  \mathcal{P}_A(\mathcal{F}')$ contained in the axis of a conjugate of $b$: $\widetilde{\eta}(\gamma)=0$. The function $\widetilde{\eta}$ satisfies conditions~$(i)-(iii)$ as $\eta$ is a relative currents, hence it defines a unique current in $\PCurr(A,\mathcal{F}')$, which we still denote by $\widetilde{\eta}$. Note that for every element $\gamma \in \mathcal{P}_A(\mathcal{F}')$, we have 
$$\widetilde{\eta}(C(\gamma) \cap \partial^2 A \cap \partial^2(F_{\tt n},\mathcal{A}(\phi)))=\eta(C(\gamma) \cap \partial^2 A \cap \partial^2(F_{\tt n},\mathcal{A}(\phi))),$$ so that the application $\PCurr(\Supp(\eta_e)) \to \PCurr(A,\mathcal{F}')$ is injective. Moreover, we have $\Supp(\widetilde{\eta})=\Supp(\eta)) \cap \partial^2 A$.

Hence $\eta_e$ defines a current $\widetilde{\eta_e} \in \PCurr(A,\mathcal{F}')$. This current coincides with the attractive projective current associated with $[\Phi|_{A}]$ defined by Gupta in~\cite[Proposition~8.12]{gupta2017relative}. By \cite[Lemma~4.17]{gupta18}, the support of $\widetilde{\eta_e}$ is uniquely ergodic. Thus the support of $\eta_e$ is uniquely ergodic.
\hfill\qedsymbol
\medskip

By the claim, it remains to prove that $\Supp(\eta_e) \subseteq \Supp(\eta_{\sigma})$. But an element $\eta \in \partial^2(F_{\tt n },\mathcal{A}(\phi))$ is contained in the support of $\eta_{\sigma}$ if for every element $\gamma \in \mathcal{P}(\mathcal{F}\wedge \mathcal{A}(\phi))$ such that $\beta \in C(\gamma)$, we have $\eta_{\sigma}(C(\gamma))>0$. Thus, the support of $\eta_{\sigma}$ contains all the cylinder sets of the form $C(\gamma)$ where $\gamma \in \mathcal{P}(\mathcal{F}\wedge \mathcal{A}(\phi))$ and $\gamma$ is contained in a reduced iterate of $\sigma$. In particular, since $e$ is contained in a reduced iterate of $\sigma$, we have $\Supp(\eta_e) \subseteq \Supp(\eta_{\sigma})$. This proves Assertion~$(1)$.

\medskip

\noindent{$(2)$ } Suppose towards a contradiction that there exist infinitely many pairwise distinct projective currents $([\eta_m])_{m \in \NN} \in \PCurr(F_{\tt n},\mathcal{F} \wedge \mathcal{A}(\phi))$ whose support is contained in the support of $\eta_{\sigma}$ and such that for every $m \in \NN$, the support $\Supp(\eta_m)$ is uniquely ergodic. By compactness of $\PCurr(F_{\tt n},\mathcal{F} \wedge \mathcal{A}(\phi))$ (see~Lemma~\ref{Lem PCurr compact}) up to passing to a subsequence, there exists a projective current $[\eta] \in \PCurr(F_{\tt n},\mathcal{F} \wedge \mathcal{A}(\phi))$ such that $\lim_{m \to \infty} [\eta_m]=[\eta]$. Let $K \in \NN^*$ be such that $\mathcal{P}(\mathcal{F} \wedge \mathcal{A}(\phi))$ contains reduced edge paths of length equal to $K$. By additivity of $\eta$, there exists $\gamma_,\ldots,\gamma_t \in \mathcal{P}(\mathcal{F} \wedge \mathcal{A}(\phi))$ of length equal to $K$ such that the support $\Supp(\eta)$ is contained in $\bigcup_{j=1}^t C(\gamma_j)$ and for every $j \in \{1,\ldots,m\}$, we have $\eta(C(\gamma_j))>0$. Then, there exists $N \in \NN^*$ such that, for every $m \geq N$ and every $j \in \{1,\ldots,t\}$, we have $\eta_m(C(\gamma_j))>0$. Hence for every $m \geq N$, we have $$\Supp(\eta) \subseteq \bigcup_{j=1}^t C(\gamma_j) \subseteq \Supp(\eta_m).$$ By unique ergodicity, for every $m \geq N$, we have $[\eta]=[\eta_m]$, a contradiction.
\hfill\qedsymbol

\section{North-South dynamics for almost atoroidal outer automorphisms}\label{Section North South dyn 1}

Let $X$ be a compact metric space and let $G$ be a group acting on $X$ by homeomorphisms. We say that an element $g \in G$ acts on $X$ with \emph{generalized north-south dynamics} if the action of $g$ on $X$ has two invariant disjoint closed subsets $\Delta_{-}$ and $\Delta_{+}$ such that, for every open neighborhood $U_{\pm}$ of $\Delta_{\pm}$ and every compact set $K_{\pm} \subseteq X-\Delta_{\mp}$, there exists $M >0$ such that, for every $n \geq M$, we have $$g^{\pm n} K_{\pm} \subseteq U_{\pm}. $$

In this section we prove the following theorem. Recall that a relative expanding outer automorphism is relative atoroidal, hence relative almost atoroidal.

\begin{theo}\label{Theo North south dynamics relative atoroidal}
Let ${\tt n} \geq 3$ and let $\mathcal{F}$ be a free factor system of $F_{\tt n}$. Let $\phi \in \Out(F_{\tt n},\mathcal{F})$ be a relative expanding outer automorphism. Let $\Delta_+(\phi)$ and $\Delta_{-}(\phi)$ be the simplexes of attraction and repulsion of $\phi$. Then $\phi$ acts on $\PCurr(F_{\tt n},\mathcal{F})$ with generalized north-south dynamics with respect to $\Delta_{+}(\phi)$ and $\Delta_{-}(\phi)$. 
\end{theo}

Theorem~\ref{Theo intro 1} in the introduction follows from Theorem~\ref{Theo North south dynamics relative atoroidal} since every exponentially growing element of $\Out(F_{\tt n})$ is expanding relative to its polynomial part. 

\subsection{Relative exponential length and goodness}\label{Section relative exp length}

Let ${\tt n} \geq 3$ and let $\mathcal{F}$ be a free factor system of $F_{\tt n}$. Let $\phi \in \Out(F_{\tt n},\mathcal{F})$ be an almost atoroidal outer automorphism relative to $\mathcal{F}$. In this section we define and prove the properties of the objects needed in order to prove Theorem~\ref{Theo North south dynamics relative atoroidal}. 
Let $f \colon G \to G$ be a CT map representing a power of $\phi$ with filtration $\varnothing=G_0 \subsetneq G_1 \subsetneq \ldots \subsetneq G_k=G$ and let $p \in \{1,\ldots,k\}$ be such that $\mathcal{F}(G_p)=\mathcal{F}$. The proof of Theorem~\ref{Theo North south dynamics relative atoroidal} relies on the study of $PG$-relative completely split edge paths. More precisely, given a reduced circuit $\gamma$ of $G$, we study the proportion of subpaths of $\gamma$ which have $PG$-relative complete splittings. This proportion will be measured using the exponential length. However, the lack of equality in Lemma~\ref{Lem Exponential length less exp length subpaths} shows that the exponential length is not well-adapted to study the exponential length of a path by comparing it with the exponential length of its subpaths. Instead, we define a notion of \emph{exponential length of a subpath relative to $\gamma$}. We first need some preliminary results regarding splittings of edge paths.

\begin{defi}
Let $\gamma$ be a reduced edge path in $G$ and let $\gamma=\gamma_0\gamma_1'\gamma_1\ldots \gamma_k\gamma_k'$ be the exponential decomposition of $\gamma$ (see the beginning of Section~\ref{Section exponential length}). Let $\alpha$ be a subpath of $\gamma$. The \emph{exponential length of $\alpha$ relative to $\gamma$}, denoted by $\ell_{exp}^{\gamma}(\alpha)$ is: $$\ell_{exp}^{\gamma}(\alpha)=\sum_{i=1}^k\ell_{exp}(\alpha \cap \gamma_k').$$ We define the \emph{$\mathcal{F}$-length of $\alpha$ relative to $\gamma$} similarly replacing $\ell_{exp}$ by $\ell_{\mathcal{F}}$ and the exponential decomposition by the $\mathcal{F}$-exponential decomposition.
\end{defi}

Note that, for every reduced edge path $\gamma$ of $G$, we have $\ell_{exp}^{\gamma}(\gamma)=\ell_{exp}(\gamma)$. The exponential length relative to a path $\gamma$ is well-adapted to compute the exponential length of $\gamma$ using its subpaths, as shown by the following lemma.

\begin{lem}\label{Lem compute exponential length optimal factor}
Let $\gamma$ be a reduced edge path and let $\gamma'=\alpha\beta \subseteq \gamma$ be a subpath of $\gamma$. Then $$\ell_{exp}^{\gamma}(\gamma')=\ell_{exp}^{\gamma}(\alpha)+\ell_{exp}^{\gamma}(\beta).$$ In particular, when $\gamma'=\gamma$, we have $$\ell_{exp}(\gamma)=\ell_{exp}^{\gamma}(\alpha)+\ell_{exp}^{\gamma}(\beta). $$ The same statement is true replacing $\ell_{exp}^{\gamma}$ by $\ell_{\mathcal{F}}^{\gamma}$.
\end{lem}

\dem The proof is similar for both $\ell_{exp}^{\gamma}$ and $\ell_{\mathcal{F}}^{\gamma}$, so we only do the proof for $\ell_{exp}^{\gamma}$. Let $\gamma=\gamma_0\gamma_1'\gamma_1\ldots \gamma_k\gamma_k'$ be the exponential decomposition of $\gamma$. Then, for every $i \in \{1,\ldots,k\}$, the paths $\alpha \cap \gamma_i'$ and $\beta \cap \gamma_i'$ do not contain a subpath of a path in $\mathcal{N}_{PG}^{\max}(\gamma)$. In particular, for every $i \in \{1,\ldots,k\}$, one computes $\ell_{exp}(\alpha \cap \gamma_i')$ and $\ell_{exp}(\beta \cap \gamma_i')$ by removing edges from $G_{PG}'$. Since $\ell_{exp}^{\gamma}(\gamma')$ is computed by removing edges in $G_{PG}'$ from every $\gamma_i'$ with $i \in \{1,\ldots,k\}$, the proof follows.
\hfill\qedsymbol

\bigskip

In Lemma~\ref{Lem bound exponential length subpath}, we will show that if $\gamma$ is a reduced edge path in $G$ and that $\alpha$ is a subpath of $\gamma$, then $\ell_{exp}(\alpha)$ and $\ell_{exp}^{\gamma}(\alpha)$ differ by a uniform additive constant. This will allow us to compute directly $\ell_{exp}(\alpha)$ rather than $\ell_{exp}^{\gamma}(\alpha)$.

Let $\gamma$ be a reduced edge path in $G$ and let $\gamma=\gamma_1\ldots\gamma_m$ be a splitting of $\gamma$. Let $J_{CS,PG} \subseteq \{\gamma_1,\ldots,\gamma_m\}$ be the subset consisting in all subpaths which have a $PG$-relative complete splitting. If $\ell_{exp}(\gamma)>0$, let $$\mathfrak{g}_{CT,PG}(\gamma,\gamma_1,\ldots,\gamma_m)=\frac{\sum_{\gamma_i \in J_{CS,PG}} \ell_{exp}^{\gamma}(\gamma_i)}{\ell_{exp}(\gamma)}.$$  The \emph{goodness of $\gamma$}, denoted by $\mathfrak{g}(\gamma)$, is the least upperbound of $\mathfrak{g}_{CT,PG}(\gamma)$ over all splittings of $\gamma$ if $\ell_{exp}(\gamma)>0$, and is equal to $0$ otherwise. When $\gamma$ is a circuit, the value $\mathfrak{g}_{CT,PG}(\gamma)$ is defined using only circuital splittings.

Since there are only finitely many decompositions of a finite edge path into subpaths, the value $\mathfrak{g}(\gamma)$ is realized for some splitting of $\gamma$. A splitting for which $\mathfrak{g}(\gamma)$ is realized is called an \emph{optimal splitting} of $\gamma$, and an {\it optimal circuital splitting} when $\gamma$ is a circuit.

A subpath of $\gamma$ which is the concatenation of consecutive splitting units of an optimal splitting of $\gamma$ is called a \emph{factor} of $\gamma$. When $\ell_{exp}(\gamma)=0$, we use the convention that the only factor of $\gamma$ is $\gamma$ itself. The factors of $\gamma$ that admit a $PG$-relative complete splitting are called \emph{complete factors}. The factors in an optimal splitting which do not admit $PG$-relative complete splittings are said to be \emph{incomplete}. Remark that, by Proposition~\ref{Prop definition CT}~$(6),(8)$ and by Lemma~\ref{Lem iterate of a path in GPG}, the $[f]$-image of a $PG$-relative complete path is $PG$-relative complete, and the reduced iterates of an incomplete factor are eventually $PG$-relative complete.

\bigskip

Using Lemma~\ref{Lem compute exponential length optimal factor}, we have the following result.

\begin{lem}\label{Lem goodness versus decomposition}
Let $\gamma$ be a reduced edge path and let $\gamma=\gamma_0'\gamma_1\gamma_1'\ldots\gamma_m\gamma_m'$ be an optimal splitting of $\gamma$, where, for every $i \in \{0,\ldots,m\}$, the path $\gamma_i'$ is an incomplete factor of $\gamma$ and, for every $i \in \{1,\ldots,m\}$, the path $\gamma_i$ is complete. Then $$\mathfrak{g}(\gamma)= \frac{\sum_{i=1}^m \ell_{exp}^{\gamma}(\gamma_i)}{\sum_{i=1}^m\ell_{exp}^{\gamma}(\gamma_i)+\sum_{j=0}^m\ell_{exp}^{\gamma}(\gamma_i')}.\eqno\qedsymbol
$$ 
\end{lem}

\begin{defi}
Let ${\tt n},\mathcal{F},\phi,f,p$ be as in the beginning of Section~\ref{Section relative exp length}. Let $K \geq 1$. The CT map $f$ is \emph{$3K$-expanding} if for every edge $e$ of $\overline{G-G_{PG}'}$, we have $$\ell_{exp}([f(e)]) \geq 3K.$$
\end{defi}

Note that, by Lemma~\ref{Lem exponential length goes to infinity}, for every $K \geq 1$, the CT map $f$ has a power which is $3K$-expanding. Note that, since $\phi$ is exponentially growing, we have $G \neq G_{PG}'$, so that the definition of $3K$-expanding is not empty.

In the rest of the section, let $K \geq 1$ be a constant such that, for every reduced edge path $\sigma$ which is either in $\mathcal{N}_{PG}$ or a path in a zero stratum, we have $\ell(\sigma) \leq \frac{K}{2}$. Such a $K$ exists since $\mathcal{N}_{PG}$ is finite by Lemma~\ref{Lem Nielsen paths in NPG properties}~$(1)$ and since every zero stratum is contractible by Proposition~\ref{Prop definition CT}~$(3)$. We fix a constant $C_f$ given by Lemma~\ref{Lem Bounded cancellation lemma}. Let \begin{equation}\label{Equation defi C}
C=\max\{K,C_f\}.
\end{equation} 
Recall that, if $\sigma$ is a $PG$-relative splitting unit, $\sigma$ is either an edge in an irreducible stratum, a path in a zero stratum or a path in $\mathcal{N}_{PG}$. Thus, the choice of $K$ implies that for every $PG$-relative splitting unit $\sigma$, we have $\ell_{exp}(\sigma) \leq \frac{K}{2}$. 

\begin{lem}\label{Lem bound exponential length subpath}
Let $\gamma$ be a  reduced edge path in $G$ and let $\gamma'$ be a subpath of $\gamma$. Let $\gamma=\gamma_0\gamma_1'\gamma_1\ldots \gamma_k\gamma_k'$ be the exponential decomposition of $\gamma$. There exist three (possibly empty) subpaths $\delta_1$, $\delta_2$ and $\tau$ of $\gamma$ such that for every $i \in \{1,2\}$, the path $\delta_i$ is a proper subpath of a splitting unit of some $\gamma_j$, we have $\ell_{exp}(\tau)=\ell_{exp}^{\gamma}(\tau)=\ell_{exp}^{\gamma}(\gamma')$ and $\gamma'=\delta_1\tau\delta_2$. In particular, we have $$\ell_{exp}^{\gamma}(\gamma') \leq \ell_{exp}(\gamma') \leq \ell_{exp}^{\gamma}(\gamma')+2C \leq \ell_{exp}(\gamma)+2C.$$
The same statement is true replacing $\ell_{exp}$ by $\ell_{\mathcal{F}}$ and $\ell_{exp}^{\gamma}$ by $\ell_{\mathcal{F}}^{\gamma}$.
\end{lem}

\dem The proof is similar for both $\ell_{exp}$ and $\ell_{\mathcal{F}}$, so we only do the proof for $\ell_{exp}$. Since $\gamma'$ is a subpath of $\gamma$, there exist three (possibly trivial) paths $\delta_1'$, $\tau'$ and $\delta_2'$ such that:

\noindent{$(a)$ } for every $i \in \{1,2\}$, there exists $k_i \in \{0,\ldots,k\}$ such that the path $\delta_i'$ is a subpath of some $\gamma_{k_i}$;

\noindent{$(b)$ } for every $j \in \{0,\ldots,k\}$, either $\gamma_j$ is contained in $\tau'$ or $\gamma_j$ does not contain edges of $\tau'$;

\noindent{$(c)$ } we have $\gamma'=\delta_1'\tau'\delta_2'$. 

The path $\delta_1'$ has a decomposition $\delta_1'=\delta_1f_1$, where $f_1$ is a (possibly trivial) factor of $\gamma_{k_1}$ and $\delta_1$ is properly contained in a splitting unit of $\gamma_{k_1}$ for some fixed choice of optimal splitting of $\gamma_{k_1}$. Similarly, the path $\delta_2'$ has a decomposition $\delta_2'=f_2\delta_2$, where $f_2$ is a (possibly trivial) factor of $\gamma_{k_2}$ and $\delta_2$ is properly contained in a splitting unit of $\gamma_{k_2}$ for some fixed choice of optimal splitting of $\gamma_{k_2}$. Let $\tau=f_1\tau'f_2$. Then $\gamma'=\delta_1\tau\delta_2$. It remains to show that $\ell_{exp}(\tau)=\ell_{exp}^{\gamma}(\tau)=\ell_{exp}^{\gamma}(\gamma')$. Since for every $i \in \{1,2\}$, the path $f_i$ is a path in $\mathcal{N}_{PG}$, we have $\ell_{exp}(\tau)=\ell_{exp}(\tau')$. By $(b)$, one obtains $\ell_{exp}(\gamma')$ by deleting edges in $G_{PG}'$ and every path of $\mathcal{N}_{PG}^{\max}(\gamma)$ contained in $\tau'$. Hence we have $$\ell_{exp}^{\gamma}(\tau')=\sum_{i=1}^k\ell_{exp}(\tau' \cap \gamma_k')=\sum_{i=1}^k\ell_{exp}(\tau \cap \gamma_k')=\ell_{exp}^{\gamma}(\tau).$$ Since $\delta_1$ and $\delta_2$ are contained in paths of $\mathcal{N}_{PG}^{\max}(\gamma)$, we have $\ell_{exp}^{\gamma}(\gamma')=\ell_{exp}^{\gamma}(\tau)$, that is, the second equality holds. 

We now prove the final inequalities in the lemma. The first inequality follows from the fact that every path in $\mathcal{N}_{PG}^{\max}(\gamma')$ is a subpath of some $\gamma_i$ for $i \in \{0,\ldots,k\}$. Thus, we have $\ell_{exp}^{\gamma}(\gamma') \leq \ell_{exp}(\gamma')$. By Lemma~\ref{Lem Exponential length less exp length subpaths}, we have $$\ell_{exp}(\gamma') \leq \ell_{exp}(\delta_1)+\ell_{exp}(\tau)+\ell_{exp}(\delta_2) \leq \ell_{exp}^{\gamma}(\gamma')+\ell(\delta_1)+\ell(\delta_2).$$ By definition of the constant $K$ and the fact that $K \leq C$, we have: $$\ell_{exp}^{\gamma}(\gamma')+\ell(\delta_1)+\ell(\delta_2)\leq \ell_{exp}^{\gamma}(\gamma')+2C \leq \ell_{exp}(\gamma)+2C,$$ where the last inequality follows from Lemma~\ref{Lem compute exponential length optimal factor}.
\hfill\qedsymbol

\begin{lem}\label{Lem p length of a completely split path grows linearly}
Let $f \colon G \to G$ be a $3K$-expanding CT map. Let $\gamma$ be a $PG$-relative completely split edge path of positive exponential length. Then $$\ell_{exp}([f(\gamma)]) \geq 3\;\ell_{exp}(\gamma).$$
\end{lem}

\dem Consider a $PG$-relative complete splitting $\gamma=\gamma_0'\gamma_1\gamma_1'\ldots\gamma_m\gamma_m'$ of $\gamma$, where, for every $i \in \{0,\ldots,m\}$, the path $\gamma_i'$ is either a (possibly trivial) concatenation of paths in $G_{PG}$ and in $\mathcal{N}_{PG}$ or a (possibly trivial) reduced maximal taken connecting path in a zero stratum and, for every $i \in \{1,\ldots,m\}$, the path $\gamma_i$ is an edge in an irreducible stratum of positive exponential length. By Lemma~\ref{Lem splitting units positive exp length}, we have $$\ell_{exp}(\gamma)=\sum_{i=1}^m \ell(\gamma_i).$$ Since $f$ is $3K$-expanding, for every $i \in \{1,\ldots,m\}$, we have $$\ell_{exp}([f(\gamma_i)]) \geq 3K \ell_{exp}(\gamma_i).$$ Since the reduced image of a $PG$-relative complete splitting is a $PG$-relative complete splitting by Lemma~\ref{Lem iterate of a path in GPG}, by Lemma~\ref{Lem compute exp length completely split}~$(2)$, we see that $$\ell_{exp}([f(\gamma)]) \geq \sum_{i=1}^m \ell_{exp}([f(\gamma_i)]) \geq \sum_{i=1}^m 3K \ell_{exp}(\gamma_i) \geq 3\ell_{exp}(\gamma).$$ This concludes the proof.
\hfill\qedsymbol

\begin{lem}\label{Lem concatenation of completely split paths and complete edges}
Let $f \colon G \to G$ be a $3K$-expanding CT map. Let $\gamma=\gamma_1\gamma_2$ be a (not necessarily reduced) edge path of positive exponential length, where $\gamma_1$ and $\gamma_2$ are reduced edge paths. Let $\gamma_1=a_1b_1\ldots a_kb_k$ be an optimal splitting of $\gamma_1$ where for every $i \in \{1,\ldots,k\}$, the path $a_i$ is an incomplete factor and for every $i \in \{1,\ldots,k\}$ the path $b_i$ is complete. For every $i \in \{1,2\}$, let $\gamma_i'$ be the subpath of $\gamma_i$ contained in $[\gamma]$. Let $\gamma_1'=\gamma_1^{-}\gamma_1^{+}$ be a decomposition of $\gamma_1'$ into two subpaths where $\gamma_1^{+}$ is the maximal terminal segment of $\gamma_1'$ such that $\sum_{i=1}^k\ell_{exp}(\gamma_1^{+}\cap b_i)=2C$. Then every $PG$-relative complete factor $b'$ of $\gamma_1$ contained in $\gamma_1^-$ (for the given optimal splitting) is also a $PG$-relative complete factor of $[\gamma]$.
\end{lem}

\begin{figure}
\centering
\captionsetup{justification=centering}
\begin{tikzpicture}[scale=1.5]
\draw (0:0) -- (90:1);
\draw (0:0) -- (330:1);
\draw (0:0) -- (210:1);
\draw (0:0) node {$\bullet$};
\draw (90:1) node {$\bullet$};
\draw (210:1) node {$\bullet$};
\draw (330:1) node {$\bullet$};
\draw[red] (-0.1,0) -- (-0.1,1);
\draw[red] (-0.1,0) -- (206:1);
\draw[blue] (0.1,0) -- (0.1,1);
\draw[blue] (0.1,0) -- (334:1);
\draw[red] (-0.5,0.35) node {$\gamma_1$};
\draw[blue] (0.5,0.35) node {$\gamma_2$};
\draw[green] (0:0) -- (210:0.5);
\draw[green] (240:0.3) node {$\gamma_2^+$};
\draw (0,-0.75) node {$[\gamma]$};
\end{tikzpicture}
\caption{Illustration of Lemma~\ref{Lem concatenation of completely split paths and complete edges}. If a complete factor of $\gamma_1$ contained in $[\gamma]$ is not contained in $\gamma_2^+$, then it is a complete factor of $[\gamma]$.}\label{Figure 2}
\end{figure}
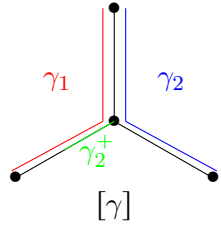

\begin{rmq}\label{Rmq explication lemma}
\noindent{$(1)$ } We emphasize that, in the statement of Lemma~\ref{Lem concatenation of completely split paths and complete edges}, if the path $\gamma_1$ is $PG$-relative completely split, the path $\gamma_1'$ is not necessarily $PG$-relative completely split. Indeed, there might be some identification with the path $\gamma_2$ that might create incomplete factors in $\gamma_1'$.

\medskip

\noindent{$(2)$ } Lemma~\ref{Lem concatenation of completely split paths and complete edges} also implies that if $\gamma_1$ is $PG$-relative completely split, the intersection of an incomplete factor of $[\gamma]$ with $\gamma_1'$ is contained in a terminal segment of $\gamma_1'$ of exponential length at most equal to $2C$ (see Figure~\ref{Figure 2}). Indeed, the claim in the proof of Lemma~\ref{Lem concatenation of completely split paths and complete edges} shows that the path $\gamma_1^-$ is a complete factor of $\gamma_1$, hence a complete factor of $[\gamma]$ by Lemma~\ref{Lem concatenation of completely split paths and complete edges}. Moreover, we have $k=1$, $a_1$ is trivial and $\ell_{exp}(\gamma_1^+)=\ell_{exp}(\gamma_1^+ \cap b_1)$.
\end{rmq}

\dem Let $t \in \{1,\ldots,k\}$ be the minimal integer such that $\gamma_1^-$ is contained in $\delta'=a_1b_1\ldots a_tb_t$. Let $b_t=\delta_1\ldots\delta_{s'}$ be a $PG$-relative complete splitting of $b_t$. Let $s \in \{1,\ldots,s'\}$ be the minimal integer such that $\gamma_1^-$ is contained in $\delta=a_1b_1\ldots a_t\delta_1\ldots\delta_{s}$. The integer $s$ exists since, by maximality of $\gamma_1^+$, for every $i \in \{1,\ldots,k\}$, either $\gamma_1^+ \cap a_i=a_i$ or $\gamma_1^+ \cap a_i=\varnothing$.

\medskip

\noindent{\bf Claim. } We have $\delta=\gamma_1^-$. 

\medskip

\dem By minimality of $t$ and $s$, the path $\gamma_1^-$ contains an edge of $\delta_s$. We claim that $\delta_s$ is contained in $\gamma_1'$. Indeed, it is clear if $\delta_s$ is an edge. Suppose towards a contradiction that $\delta_s$ is not contained in $\gamma_1'$. Then the concatenation point of $\gamma_1'$ and $\gamma_2'$ is contained in $\delta_s$. If $\delta_s$ is a maximal taken connecting path in a zero stratum, then, by the choice of $K$, we have $\ell(\delta_s) \leq \frac{K}{2} \leq \frac{C}{2}$. Since $\ell(\gamma_1^+) \geq 2C$, the path $\delta_s \cap \gamma_1'$ would be contained in $\gamma_1^+$, contradicting the fact that $\gamma_1^-$ contains the first edge of $\delta_s$. Suppose that $\delta_s$ is a concatenation of paths in $G_{PG}$ and $\mathcal{N}_{PG}$. Then $\delta_s$ has a decomposition $\delta_s=\beta_1^{(s)}\alpha_1^{(s)}\beta_1^{(s)}\ldots\alpha_{k_s-1}^{(s)}\beta_{k_s}^{(s)}\alpha_{k_s}^{(s)}$, where for every $i \in \{1,\ldots,k_s\}$, the path $\beta_i^{(s)}$ is contained in $G_{PG}$, for every $i \in \{1,\ldots,k_s-1\}$, the path $\alpha_i^{(s)}$ is contained in $\mathcal{N}_{PG}^{\max}(\sigma)$ and $\alpha_{k_s}^{(s)}$ is a subpath of a path in $\mathcal{N}_{PG}^{\max}(\delta_s)$. By the choice of $K$, we have $\ell_{exp}(\delta_s) \leq \ell(\alpha_{k_s}) \leq \frac{K}{2} \leq \frac{C}{2}$. Since $\ell_{exp}(\gamma_1^{+}) \geq 2C$, the path $\delta_s \cap \gamma_1'$ would be contained in $\gamma_1^+$, contradicting the fact that $\gamma_1^-$ contains the first edge of $\delta_s$. Hence, in every case, the path $\delta_s$ is contained in $\gamma_1'$. Note that, since $\gamma_1^+$ is the maximal subpath of $\gamma_1'$ for the property that $\sum_{i=1}^k\ell_{exp}(\gamma_1^{+} \cap b_i)=2C$, the $PG$-relative splitting unit $\delta_s$ is not a concatenation of paths in $G_{PG}$ and in $\mathcal{N}_{PG}$ or a maximal taken connecting path in a zero stratum. Indeed, otherwise it is properly contained in $\gamma_1^+$, contradicting the fact that $\gamma_1^-$ intersects $\delta_s$. Hence $\delta_s$ is an edge and $\delta=\gamma_1^-$. 
\hfill\qedsymbol

\medskip

By the claim, we see that $\gamma_1^-= a_1b_1\ldots a_t\delta_1\ldots\delta_s$ is an optimal splitting of $\gamma_1^-$. Let $r \in \{1,\ldots,k\}$ be the minimal integer such that $\gamma_1'$ is contained in $a_1b_1\ldots a_rb_r$. The last edge of $\gamma_1'$ is either contained in $a_r$ or in $b_r$. In the first case, for every $i \in \{1,\ldots,k\}$, either $b_i$ is contained in $\gamma_1'$ or $b_i\cap \gamma_1'$ is at most a point. In the second case, it is possible that $b_r \cap \gamma_1' \neq b_r$ and that $b_r \cap \gamma_1'$ contains an edge. Let $\alpha'$ be the (possibly trivial) terminal segment of $\gamma_1^{+}$ which is properly contained in a splitting unit $\sigma$ of $b_r$. If $\sigma$ is a maximal taken connecting path in a zero stratum, then, by the choice of $K$, we have $\ell_{exp}(\alpha') \leq \ell(\alpha') \leq \ell(\sigma) \leq \frac{K}{2} \leq \frac{C}{2}$. Suppose that $\sigma$ is a concatenation of paths in $G_{PG}$ and $\mathcal{N}_{PG}$. Then $\alpha'$ has a decomposition $\alpha'=\beta_1\alpha_1\beta_1\ldots\alpha_{\ell-1}\beta_{\ell}\alpha_{\ell}$, where for every $i \in \{1,\ldots,\ell\}$, the path $\beta_i$ is contained in $G_{PG}$, for every $i \in \{1,\ldots,\ell-1\}$, the path $\alpha_i$ is contained in $\mathcal{N}_{PG}^{\max}(\sigma)$ and $\alpha_{\ell}$ is a subpath of a path in $\mathcal{N}_{PG}^{\max}(\sigma)$. By the choice of $K$, we have $\ell_{exp}(\alpha') \leq \ell(\alpha_{\ell}) \leq \frac{K}{2} \leq \frac{C}{2}$. Since $\ell_{exp}(\gamma_1^{+}) \geq 2C$, there exists a $PG$-relative complete factor $\alpha_0$ of $b_r$ such that $\gamma_1^{+}=\delta_{s+1}\ldots \delta_{s'} a_{t+1}b_{t+1}\ldots a_r\alpha_0\alpha'=\alpha\alpha'$ and $$\sum_{i=1}^k\ell_{exp}\left(\alpha \cap b_i \right) \geq C.$$

We now prove that every $PG$-relative complete factor of $\gamma_1$ contained in $\gamma_1^-$ is a $PG$-relative complete factor of $\gamma$. Note that the decomposition $\gamma_1^-\alpha$ is a splitting. Thus, it suffices to prove that, for every $k \in \NN^*$, the path $[f^k(\gamma_1^-)]$ is contained in $[f^k(\gamma)]$ as any identification in order to obtain $[f^k(\gamma)]$ which involves a path in $f^k(\gamma_1^-)$ will be induced by an identification in order to obtain $[f^k(\gamma_1^-)]$ from $f^k(\gamma_1^-)$. By Lemma~\ref{Lem p length of a completely split path grows linearly} applied to $\delta_{s+1},\ldots, \delta_{s'}$, to the paths $b_i$ with $i \in \{1,\ldots,k\}$ such that $b_i \subseteq \alpha$ and to $\alpha_0$, we have $$\begin{array}{ccl} \sum_{i=1}^k\ell_{exp}([f(\alpha)] \cap [f(b_i)]) &\geq& \sum\limits_{i=s+1}^{s'}\ell_{exp}([f(\delta_i)])+\sum\limits_{i=t+1}^{r-1}\ell_{exp}([f(b_i)])+\ell_{exp}([f(\alpha_0)]) \\
{} &\geq &3\;\sum_{i=1}^k\ell_{exp}\left(\alpha \cap b_i \right) \geq 3C,
\end{array}$$ where the first inequality follows from  the fact that the decomposition $$\alpha=\delta_{s+1}\ldots \delta_{s'} a_{t+1}b_{t+1}\ldots a_r\alpha_0$$ is an optimal splitting of $\alpha$. Note that, since the decomposition $\gamma_1^{-}\alpha$ is a splitting, for every $k \in \NN^*$, the path $[f^k(\alpha)]$ is contained in $[f^k(\gamma_1^{-}\alpha)]$. Remark that Lemma~\ref{Lem Bounded cancellation lemma} implies that the segment of $[f(\gamma_1^{-}\alpha)]$ which is $C$ away from the concatenation point between $[f(\gamma_1^{-}\alpha)]$ and $[f(\alpha'\gamma_2')]$ remains in $[f([\gamma])]$. In particular, the edges of $[f(\gamma_1^{-}\alpha)]$ which are cancelled with edges of $[f(\alpha'\gamma_2')]$ are contained in $[f(\alpha)]$. Recall that $\sum_{i=1}^k\ell_{exp}([f(\alpha)] \cap [f(b_i)]) \geq 3C$ and that the subpath of $[f(\alpha)]$ which is contained in $[f([\gamma])]$ is obtained by the concatenation of at most $C$ edges of $[f(\alpha)]$. Thus, we see that the sum over $i$ of the exponential length of the subpaths of $[f(\alpha)] \cap [f(b_i)]$ which are contained in $[f([\gamma])]$ is at least equal to $2C$. Hence the path $[f(\gamma_1^{-})]$ is a subpath of $[f([\gamma])]$ and $\sum_{i=1}^k\ell_{exp}([f(\gamma_1^{+})]\cap [f(b_i)] \cap [f([\gamma])]) \geq 2C$. Thus, we can apply the same arguments to show that for every $k \geq 1$, the path $[f^k(\gamma_1^{-})]$ is contained in $[f^k([\gamma])]$ and the exponential length of the subpath of $[f^k(\alpha)]$ contained in $[f^k([\gamma])]$ is at least equal to $2C$. Hence every $PG$-relative complete factor of the path $\gamma_1$ contained in $\gamma_1^{-}$ is a complete factor of an optimal splitting of $[\gamma]$.
\hfill\qedsymbol

\begin{lem}\label{Lem identification paths complete exp grow path and Gpg}
\noindent{$(1)$ } Let $\gamma=\alpha\beta$ be a reduced path. Let $N \in \NN^*$ be such that $[f^N(\alpha)]$ has a $PG$-relative complete splitting and that $[f^N(\beta)]$ is a concatenation of paths in $G_{PG}$ and in $\mathcal{N}_{PG}$. For every $m \geq N$, let $\alpha_m$, $\beta_m$ and $\sigma_m$ be paths such that $[f^m(\alpha)]=\alpha_m\sigma_m$ and $[f^m(\beta)]=\sigma_m^{-1}\beta_m$. For every $m \geq N$, we have $\ell_{exp}(\sigma_m) \leq C$, $\ell_{exp}(\alpha_m) \geq \ell_{exp}([f^m(\alpha)])-C$ and $\ell_{exp}(\beta_m) \leq C$.

\medskip

\noindent{$(2)$ } Let $\gamma=\beta^{(1)}\alpha\beta^{(2)}$ be a reduced path. Let $N \in \NN^*$ be such that $[f^N(\alpha)]$ has a $PG$-relative complete splitting and, for every $i \in \{1,2\}$, the path $[f^N(\beta^{(i)})]$ is a concatenation of paths in $G_{PG}$ and in $\mathcal{N}_{PG}$. For every $m \geq N$, let $\alpha_m$, $\beta_m^{(1)}$, $\beta_m^{(2)}$, and $\sigma_m^{(1)}$, $\sigma_m^{(2)}$ be paths such that $[f^m(\alpha)]=\sigma_m^{(1)}\alpha_m\sigma_m^{(2)}$, $[f^m(\beta^{(1)})]=\beta_m^{(1)}\sigma_m^{(1)-1}$ and $[f^m(\beta^{(2)})]=\sigma_m^{(2)-1}\beta_m$. For every $m \geq N$, either $\ell_{exp}(\alpha_m) \leq 2C$ or we have $\ell_{exp}(\sigma_m^{(1)}),\ell_{exp}(\sigma_m^{(2)}) \leq C$, $\ell_{exp}(\alpha_m) \geq \ell_{exp}([f^m(\alpha)])-2C$ and $\ell_{exp}(\beta_m^{(1)}),\ell_{exp}(\beta_m^{(2)}) \leq C$.

\end{lem}

\dem The proof of Assertion~$(2)$ follows from Assertion~$(1)$ by applying Assertion~$(1)$ twice: one with $\gamma=\alpha\beta^{(2)}$ and one with $\gamma=\alpha^{-1}\beta^{(1)}$. If for some $m \in \NN^*$, $\ell_{exp}(\alpha_m) \geq 2C$, there is no identification between $[f^m(\beta^{(1)})]$ and $[f^m(\beta^{(2)})]$, so Assertion~$(2)$ follows from Assertion~$(1)$. Therefore, we focus on the proof of Assertion~$(1)$. Let $m \geq N$. When $\sigma_m$ is reduced to a point, we have $\ell_{exp}(\alpha_m)=\ell_{exp}([f^m(\alpha)])$ and $\ell_{exp}(\beta_m)=\ell_{exp}([f^m(\beta)])=0$ by Lemma~\ref{Lem exponential length paths in Gpg}. This concludes the proof in this case. So we may suppose that $\sigma_m$ is nontrivial. Let $[f^m(\alpha)]=a_1\ldots a_{k}$ be a $PG$-relative complete splitting of $[f^m(\alpha)]$. Suppose that, for every $i \in \{1,\ldots,k\}$ such that $a_i$ is a concatenation of paths in $G_{PG}$ and $\mathcal{N}_{PG}$, the path $a_i$ is a maximal subpath of $[f^m(\alpha)]$ for the property of being a factor which is a concatenation of paths in $G_{PG}$ and $\mathcal{N}_{PG}$. For every $j \in \{1,\ldots,k\}$, let $r_j$ be the height of $a_j$. Let $i \in \{1,\ldots,k\}$ be such that $a_i$ contains the first edge of $\sigma_m$. Let $\sigma' \in \mathcal{N}_{PG}^{\max}(\sigma_m)$. Note that there exists $\sigma'' \in \mathcal{N}_{PG}^{\max}([f^m(\alpha)])$ such that $\sigma' \subseteq \sigma''$. By Lemma~\ref{Lem compute exp length completely split}~$(1)$ applied to $\sigma''$ and $[f^m(\alpha)]$, the path $\sigma''$ is contained in a factor which is a concatenation of paths in $G_{PG}$ and $\mathcal{N}_{PG}$. By the maximality assumption, there exists $j \in \{1,\ldots,k\}$ such that $\sigma' \subseteq \sigma'' \subseteq a_j$. Hence we can compute $\ell_{exp}(\sigma_m)$ by removing, for every $j \in \{1,\ldots,k\}$, paths in the intersection $\sigma_m \cap a_j$. Thus, we have, $$\ell_{exp}(\sigma_m) = \sum_{j>i} \ell_{exp}(a_j)+ \ell_{exp}(a_i \cap \sigma_m).$$

Note that, by Lemma~\ref{Lem iterate of a path in GPG}, the path $[f^m(\beta)]=\sigma_m^{-1}\beta_m$ is a concatenation of paths in $G_{PG}$ and in $\mathcal{N}_{PG}$. Let $j \in \{i,\ldots,k\}$.

\medskip

\noindent{\bf Claim. } If $j>i$, then either $a_j$ is not an edge in an EG stratum and $\ell_{exp}(a_j \cap \sigma_m)=0$, or $\ell_{exp}((a_i \ldots a_j) \cap \sigma_m) \leq C$. If $j=i$, then $\ell_{exp}(a_j \cap \sigma_m) \leq C$. 

\medskip

\dem We distinguish several cases, according to the nature of $a_j$.

\noindent{$(i)$ } Suppose that $a_j$ is maximal taken connecting path in a zero stratum. By definition we have $\ell_{exp}(a_j \cap \sigma_m)=0$.

\noindent{$(ii)$ } Suppose that $a_j$ is a concatenation of paths in $G_{PG}$ and in $\mathcal{N}_{PG}$. If $j>i$, we have $a_j \cap \sigma_m=a_j$. By Lemma~\ref{Lem exponential length paths in Gpg} applied to $a_j$, we have $\ell_{exp}(a_j \cap \sigma_m)=0$. Suppose that $i=j$. Suppose that the first edge of $\sigma_m$ is not contained in a path in $\mathcal{N}_{PG}^{\max}(a_i)$. Then $a_i$ has a decomposition $a_i=a_i^0a_i^1a_i^2$ where $a_i^1$ is a path contained in $G_{PG}$ such that the first edge of $\sigma_m$ is contained in $a_i^1$ and such that, for every path $\delta \in \mathcal{N}_{PG}^{\max}(a_i)$, either $\delta \subseteq a_i^0$ or $\delta \subseteq a_i^2$. Note that a terminal segment of $a_i$ whose first edge is contained in $a_i^1$ is a concatenation of paths in $G_{PG}$ and in $\mathcal{N}_{PG}$. In particular, the path $a_i \cap \sigma_m$ is a concatenation of paths in $G_{PG}$ and in $\mathcal{N}_{PG}$. By Lemma~\ref{Lem exponential length paths in Gpg} applied to $a_i \cap \sigma_m$, we have $\ell_{exp}(a_i \cap \sigma_m)=0$. Suppose now that the first edge of $\sigma_m$ is contained in a path $\delta \in \mathcal{N}_{PG}^{\max}(a_i)$. Then $a_i$ has a decomposition $a_i^1\delta a_i^2$, where the first edge of $\sigma_m$ is contained in $\delta$. Note that $a_i^2$ is a concatenation of paths in $G_{PG}$ and in $\mathcal{N}_{PG}$ which is contained in $\sigma_m$. By Lemma~\ref{Lem Exponential length less exp length subpaths} applied to $a_i \cap \sigma_m=(\delta \cap \sigma_m) a_i^2$, by Lemma~\ref{Lem exponential length paths in Gpg} applied to $a_i^2$ and by definition of the constant $K$, we have $$\ell_{exp}(\sigma_m \cap a_i) \leq \ell_{exp}(\sigma_m \cap \delta)+ \ell_{exp}(a_i^2)=\ell_{exp}(\delta \cap \sigma_m) \leq \ell(\delta) \leq K \leq C.$$

\noindent{$(iii)$ } Suppose that $a_j$ is an edge in an irreducible stratum with positive exponential length. Since $[f^m(\beta)]$ is a concatenation of paths in $G_{PG}$ and in $\mathcal{N}_{PG}$, there exists a path $\gamma' \in \mathcal{N}_{PG}^{\max}([f^m(\beta)])$ such that $a_j$ is contained in $\gamma'$. By Lemma~\ref{Lem compute exp length completely split}~$(1)$, every path in $\mathcal{N}_{PG}^{\max}([f^m(\alpha)])$ is contained in a minimal factor of $[f^m(\alpha)]$ consisting in $PG$-relative splitting units which are concatenation of paths in $G_{PG}$ and $\mathcal{N}_{PG}$. Since $a_j$ is a $PG$-relative splitting unit of $[f^m(\alpha)]$ which is not a concatenation of paths in $G_{PG}$ and in $\mathcal{N}_{PG}$, the path $a_j$ is not contained in a path of $\mathcal{N}_{PG}^{\max}([f^m(\alpha)])$. Hence the path $\gamma'$ is not contained in $\sigma_m$ as otherwise it would be contained in a path of $\mathcal{N}_{PG}^{\max}([f^m(\alpha)])$. Therefore, we see that $(a_i \ldots a_j) \cap \sigma_m \subseteq \gamma'$. Hence, by the choice of $K$, we have $$\ell_{exp}((a_i\ldots a_j) \cap \sigma_m) \leq \ell((a_i\ldots a_j) \cap \sigma_m) \leq \ell(\gamma') \leq C.$$ 

\noindent This proves the claim as we considered all possible $PG$-relative splitting units.
\hfill\qedsymbol

\bigskip

Let $m \in \NN^*$. By the claim, either $\ell_{exp}((a_i \ldots a_j) \cap \sigma_m) \leq C$ or, for every $j >i$, we have $\ell_{exp}(a_j \cap \sigma_m)=0$. In the second case, we have $$\ell_{exp}(\sigma_m) = \sum_{j>i} \ell_{exp}(a_j)+ \ell_{exp}(a_i \cap \sigma_m)=\ell_{exp}(a_i \cap \sigma_m) \leq C,$$ where the las inequality follows from the case $j=i$ of the claim. Hence, for every $m \in \NN^*$, we have $\ell_{exp}(\sigma_m) \leq C$. Note that, by Lemma~\ref{Lem Exponential length less exp length subpaths} applied to $[f^m(\alpha)]=\alpha_m\sigma_m$, we have $$\ell_{exp}(\alpha_m) \geq \ell_{exp}([f^m(\alpha)])-\ell_{exp}(\sigma_m) \geq \ell_{exp}([f^m(\alpha)])-C.$$

It remains to prove that $\ell_{exp}(\beta_m) \leq C$. But $\beta_m$ can be written as $\beta_m=\delta_1\delta_2$ where $\delta_2$ is a concatenation of paths in $G_{PG}$ and in $\mathcal{N}_{PG}$ and $\delta_1$ is a (possibly trivial) path contained in a path of $\mathcal{N}_{PG}^{\max}([f^m(\beta)])$. By Lemma~\ref{Lem exponential length paths in Gpg} applied to $\delta_2$ and by the choice of $K$ (since $\delta_1$ is a subpath of a path in $\mathcal{N}_{PG}$), we have $$\ell_{exp}(\beta_m) \leq \ell_{exp}(\delta_1)+ \ell_{exp}(\delta_2)=\ell_{exp}(\delta_1) \leq \ell(\delta_1) \leq C.$$ This concludes the proof.
\hfill\qedsymbol

\begin{lem}\label{Lem complete splitting rel poly grow}
Let $L \geq 1$. There exists $n_0=n_0(L) \in \NN^*$  which satisfies the following properties. Let $\gamma$ be a reduced edge path of $G$ such that $\ell_{exp}(\gamma) \leq L$. For every $n \geq n_0$ and every optimal splitting of $[f^n(\gamma)]$, either $[f^n(\gamma)]$ is a concatenation of paths in $G_{PG}$ and in $\mathcal{N}_{PG}$ or the following two assertions hold:
\medskip

\noindent{$(a)$ } the path $[f^n(\gamma)]$ contains a complete factor of exponential length at least equal to $10C$;

\medskip

\noindent{$(b)$ } the exponential length of an incomplete factor of $[f^n(\gamma)]$ is at most equal to $8C$.
\end{lem}

\dem By Lemma~\ref{Lem exponential length goes to infinity}, there exists an integer $m' \in \NN^*$ depending only on $f$ such that for every edge $e$ of $\overline{G-G_{PG}'}$ and every $n \geq m'$, we have $\ell_{exp}[f^n(e)] \geq 16C+1$. Let $\gamma=\gamma_0\gamma_1'\gamma_1\ldots \gamma_{\ell}\gamma_{\ell}'$ be the exponential decomposition of $\gamma$. Let $$\gamma=\beta_0\alpha_1\beta_1\ldots \alpha_k\beta_k$$ be a nontrivial decomposition of $\gamma$ such that, for every $i \in \{0,\ldots,k\}$, the path $\beta_i$ is a concatenation of paths in $G_{PG}$ and in $\mathcal{N}_{PG}$ and for every $i \in \{1,\ldots,k\}$, the path $\alpha_i$ is a concatenation of edges in irreducible strata not contained in some $\gamma_j$ with $j \in \{0,\ldots,\ell\}$ and paths in zero strata. The main point of the proof is to show that, up to applying an iterate of $[f]$, there is no cancellation between the subpaths $\alpha_i$. By definition of the exponential length, for every $i \in \{1,\ldots,k\}$, we have $\ell_{exp}(\gamma)=\sum_{i=1}^k \ell_{exp}(\alpha_i)$. Therefore, since $\ell_{exp}(\gamma) \leq L$, for every $i \in \{1,\ldots,k\}$, we have $\ell_{exp}(\alpha_i) \leq L$. Note that, for every $i \in \{1,\ldots,k\}$, we have $\ell_{exp}(\alpha_i)=\ell(\alpha_i)-\ell(\alpha_i \cap \mathcal{Z})$ where $\mathcal{Z}$ is the subgraph of $G$ consisting in all zero strata. By the choice of $C$ the length of every path contained in a zero stratum is at most equal to $C$. Hence for every $i \in \{1,\ldots,k\}$, we have $\ell(\alpha_i) \leq CL$. By Proposition~\ref{Prop definition CT}~$(8)$ there exists $m'' \in \NN^*$ depending only on $L$ such that, for all $i \in \{1,\ldots,k\}$ and $m \geq m''$, the path $[f^m(\alpha_i)]$ is completely split. Let $m=m'+m''$. By Lemma~\ref{Lem compute exp length completely split}~$(2)$, for every $n \geq m$ and every $i \in \{1,\ldots,k\}$, since $[f^{n-m'}(\alpha_i)]$ is completely split, one compute its exponential length by adding the exponential length of all its splitting units. Thus, if $[f^{n-m'}(\alpha_i)]$ contains a splitting unit which is an edge $e$ in $\overline{G-G_{PG}'}$, we have \begin{equation}\label{Equation new lemma}
\ell_{exp}([f^n(\alpha_i)]) \geq \ell_{exp}([f^{m'}(e)]) \geq 16C+1.
\end{equation}

Let $C_m$ be a bounded cancellation constant for $f^m$ given by Lemma~\ref{Lem Bounded cancellation lemma}. Note that, if there exists $i \in \{1,\ldots,k-1\}$ such that $\ell(\beta_i) < C_m$, then there might exist some identifications between $[f^m(\alpha_{i-1})]$ and $[f^m(\alpha_i)]$ when reducing the paths in order to obtain $[f^m(\gamma)]$. This is why we replace the decomposition $\gamma=\beta_0\alpha_1\beta_1\ldots \alpha_k\beta_k$ of $\gamma$ by a new one. This new decompostion is defined as follows. Since every lift of $f^{m}$ to the universal cover of $G$ is a quasi-isometry, there exists $M_{m} >0$ depending only on $m$ such that, for every reduced edge path of length $\ell(\beta) >M_{m}$, we have $\ell([f^{m}(\beta)]) \geq 2C_{m}+1$. Let $\Gamma_m=\{\beta_i \;|\; \ell(\beta_i) \leq M_m\}$. Note that $|\Gamma_m| \leq k+1$. Note that, by Lemma~\ref{Lem No zero path and Nielsen path adjacent} and Proposition~\ref{Prop definition CT}~$(4)$ and Lemma~\ref{Lem No zero path and Nielsen path adjacent}, for every $i \in \{1,\ldots,k\}$, if $\beta_{i-1}$ or $\beta_i$ is not trivial, then $\alpha_i$ is not contained in a zero stratum. In particuliar, we may suppose that, for every $i \in \{1,\ldots,k\}$, we have $\ell_{exp}(\alpha_i)>0$. Thus, since $\ell_{exp}(\gamma)=\sum_{i=1}^k\ell_{exp}(\alpha_i) \leq L$, and, for every $i \in \{1,\ldots,k\}$, we have $\ell_{exp}(\alpha_i)>0$, we see that $k \leq L$. Hence we have $|\Gamma_m| \leq k+1 \leq L+1$.

\medskip

\noindent{\bf Claim. } There exist $m_1 \geq m$ depending only on $|\Gamma_m|$ (and hence on $L$) and a decomposition $\gamma=\beta_0^{(1)}\alpha_1^{(1)}\beta_1^{(1)}\ldots\alpha_{k_1}^{(1)}\beta_{k_1}^{(1)}$ such that:

\noindent{$(a')$ } for every $i \in \{1,\ldots,k_1\}$, the path $[f^{m_1}(\alpha_i^{(1)})]$ is completely split;

\noindent{$(b')$ } for every $i \in \{0,\ldots,k_1\}$, the path $\beta_i^{(1)}$ is a concatenation of paths in $G_{PG}$ and in $\mathcal{N}_{PG}$;

\noindent{$(c')$ } for every $i \in \{0,\ldots,k_1\}$, the subpath of $[f^{m_1}(\beta_i^{(1)})]$ contained in $[f^{m_1}(\gamma)]$ is not reduced to a point;

\noindent{$(d')$ } for every $i \in \{1,\ldots,k_1\}$, for every $n \geq m'$, if $[f^{n-m'}(\alpha_i^{(1)})]$ contains a splitting unit which is an edge in $\overline{G-G_{PG}'}$ then $\ell_{exp}([f^n(\alpha_i^{(1)})]) \geq 16C+1$.

\medskip

\dem The proof is by induction on $|\Gamma_m|$. Suppose first that $\Gamma_m=\varnothing$. By the definition of $|\Gamma_m|$ and $M_m$, for every $i \in \{0,\ldots,k\}$, the path $[f^m(\beta_i)]$ has length at least equal to $2C_m+1$. By Lemma~\ref{Lem Bounded cancellation lemma}, for every $i \in \{0,\ldots,k\}$, the subpath of $[f^m(\beta_i)]$ contained in $[f^m(\gamma)]$ is not reduced to a point. So the integer $m_1=m$ and the decomposition $\gamma=\beta_0\alpha_1\beta_1\ldots \alpha_k\beta_k$ satisfy the assertions of the claim (Assertion~$(d')$ follows from Equation~\eqref{Equation new lemma}).

Suppose now that $\Gamma_m \neq \varnothing$. Then $$\sum_{i=1}^k \ell(\alpha_i) + \sum_{\beta_i \in \Gamma_m} \ell(\beta_i) \leq kCL +M_mL\leq CL^2+M_mL.$$ Let $m_2' \geq m$ be such that for every path $\beta$ of length at most equal to $CL^2+M_mL$ and every $n \geq m_2'$, the path $[f^n(\beta)]$ is completely split. Then $\gamma$ has a decomposition $\gamma=\beta_0^{(2)}\alpha_1^{(2)}\beta_2^{(2)}\ldots\alpha_{k_2}^{(2)}\beta_{k_2}^{(2)}$ such that, for every $i \in \{1,\ldots,k_2\}$, the path $[f^{m_2'}(\alpha_i^{(2)})]$ is completely split and for every $i \in \{0,\ldots,k_2\}$, the path $\beta_i^{(2)}$ is a concatenation of paths in $G_{PG}$ and in $\mathcal{N}_{PG}$ of length greater than $M_m$. Let $m_2=m_2'+m'$. Then for every $i \in \{1,\ldots,k_2\}$, the paths $[f^{m_2}(\alpha_i^{(2)})]$ and $[f^{m_2-m'}(\alpha_i^{(2)})]$ are completely split. Moreover, if $[f^{m_2-m'}(\alpha_i^{(2)})]$ contains a splitting unit which is an edge in $\overline{G-G_{PG}'}$, then $\ell_{exp}([f^m(\alpha_i^{(2)})]) \geq 16C+1$ as in Equation~\eqref{Equation new lemma}. Let $C_{m_2}$ be a bounded cancellation constant associated with $f^{m_2}$ and let $M_{m_2} \geq M_m$ be such that, for every reduced edge path of length $\ell(\beta) >M_{m_2}$, we have $\ell([f^{m_1}(\beta)]) \geq 2C_{m_2}+1$. Let $\Gamma_{m_2}=\{\beta_i^{(2)} \;|\; \ell(\beta_i) \leq M_{m_2}\}$. Note that $|\Gamma_{m_2}|<|\Gamma_m|$. Hence we can apply the induction hypothesis to the decomposition $\gamma=\beta_0^{(2)}\alpha_1^{(2)}\beta_2^{(2)}\ldots\alpha_{k_2}^{(2)}\beta_{k_2}^{(2)}$ to obtain the desired decomposition of $\gamma$. This concludes the proof of the claim.
\hfill\qedsymbol

\bigskip

Let $m_1$ and $\gamma=\beta_0^{(1)}\alpha_1^{(1)}\beta_1^{(1)}\ldots\alpha_{k_1}^{(1)}\beta_{k_1}^{(1)}$ be as in the assertion of the claim. By Assertion $(c')$ of the claim, for every $i \in \{1,\ldots,k_1\}$, there is no identification between edges of $[f^{m_1}(\alpha_i^{(1)})]$, $[f^{m_1}(\alpha_{i-1}^{(1)})]$ and $[f^{m_1}(\alpha_{i+1}^{(1)})]$ when reducing in order to obtain $[f^{m_1}(\gamma)]$. 

For every $i \in \{1,\ldots,k_1\}$, since $[f^{m_1}(\alpha_i^{(1)})]$ is $PG$-relative completely split, we can distinguish three possible cases for $[f^{m_1}(\alpha_i^{(1)})]$:

\medskip

\noindent{$(i)$ } the path $[f^{m_1}(\alpha_i^{(1)})]$ contains a $PG$-relative splitting unit which is an edge in $\overline{G-G_{PG}'}$ (by Lemma~\ref{Lem splitting units positive exp length} this case happens exactly when $\ell_{exp}([f^{m_1}(\alpha_i^{(1)})]) >0$);

\medskip

\noindent{$(ii)$ } $\ell_{exp}([f^{m_1}(\alpha_i^{(1)})])=0$ and the path $[f^{m_1}(\alpha_i^{(1)})]$ is a concatenation of paths in $G_{PG}$ and in $\mathcal{N}_{PG}$;

\medskip

\noindent{$(iii)$ } $\ell_{exp}([f^{m_1}(\alpha_i^{(1)})])=0$ and $[f^{m_1}(\alpha_i^{(1)})]$ contains a maximal taken connecting path in a zero stratum.

\medskip

We claim that if there exists $i \in \{1,\ldots,k_1\}$ such that $[f^{m_1}(\alpha_i^{(1)})]$ satisfies $(iii)$, then $[f^{m_1}(\gamma)]$ is contained in a zero stratum. Indeed, suppose that $[f^{m_1}(\alpha_i^{(1)})]$ satisfies $(iii)$. By Lemma~\ref{Lem splitting units positive exp length} applied to the $PG$-relative completely split edge path $[f^{m_1}(\alpha_i^{(1)})]$, since $\ell_{exp}([f^{m_1}(\alpha_i^{(1)})])=0$ the path $[f^{m_1}(\alpha_i^{(1)})]$ does not contain an edge in $\overline{G-G_{PG}'}$. Therefore, the path $[f^{m_1}(\alpha_i^{(1)})]$ is a concatenation of paths in $G_{PG}'$ and in $\mathcal{N}_{PG}$. By Proposition~\ref{Prop definition CT}~$(4)$ and Lemma~\ref{Lem No zero path and Nielsen path adjacent}, there is no path in a zero stratum which is adjacent to a concatenation of paths in $G_{PG}$ and in $\mathcal{N}_{PG}$. Hence $[f^{m_1}(\alpha_i^{(1)})]=\sigma$, where $\sigma$ is a maximal taken connecting path in a zero stratum not contained in $G_{PG}$. But the endpoints of $\sigma$ are the endpoints of $[f^{m_1}(\beta_{i-1}^{(1)})]$ and $[f^{m_1}(\beta_i^{(1)})]$, which are concatenation of paths in $G_{PG}$ and in $\mathcal{N}_{PG}$. As above, this implies that $[f^{m_1}(\gamma)]=\sigma$. Since zero strata are contractible, there exists $m_3 \in \NN^*$ such that $[f^{m_3}(\gamma)]$ is $PG$-relative completely split. Hence Assertion~$(b)$ of Lemma~\ref{Lem complete splitting rel poly grow} follows. Applying a further power of $[f]$ (which can be chosen uniformly as there are finitely many reduced edge paths contained in a zero stratum), there exists $m_4 \in \NN^*$ such that $[f^{m_4}(\gamma)]$ is a concatenation of paths in $G_{PG}$ and in $\mathcal{N}_{PG}$ or it satisfies Assertion~$(a)$ of Lemma~\ref{Lem complete splitting rel poly grow}. This concludes the proof of Lemma~\ref{Lem complete splitting rel poly grow} in case~$(iii)$.

Hence we may suppose that for every $i \in \{1,\ldots,k_1\}$, the path $[f^{m_1}(\alpha_i^{(1)})]$ satisfies either $(i)$ or $(ii)$. Note that, if $i \in \{1,\ldots,k_1\}$ is such that the path $[f^{m_1}(\alpha_i^{(1)})]$ satisfies $(i)$, then $[f^{m_1}(\alpha_i^{(1)})]$ also satisfies the hypothesis of Assertion~$(d')$ of the claim. Thus $$\ell_{exp}([f^{m_1+m'}(\alpha_i^{(1)})]) \geq 16C+1.$$ Let $m_1'=m_1+m'$ and let $n' \geq m_1'$. Let $\Lambda_{exp}=\{\alpha_i^{(1)} \;|\; \ell_{exp}([f^{n'}(\alpha_i^{(1)})]) \geq 16C+1\}$. For every \mbox{$j \in \{1,\ldots,k_1\}$} and every $n \in \NN^*$, let $\alpha_j^{(n)}$ be the subpath of $[f^n(\alpha_j^{(1)})]$ contained in $[f^n(\gamma)]$. For every $j \in \{0,\ldots,k_1\}$ and every $n \in \NN^*$, let $\beta_j^{(n)}$ be the subpath of $[f^n(\beta_j^{(1)})]$ contained in $[f^n(\gamma)]$. Suppose first that $\Lambda_{exp}$ is not empty and let $\alpha_i^{(1)} \in \Lambda_{exp}$. By Lemma~\ref{Lem identification paths complete exp grow path and Gpg}~$(2)$ applied to $\beta^{(1)}=[f^{n'}(\beta_{i-1}^{(1)})]$, $\alpha=[f^{n'}(\alpha_i^{(1)})]$ and $\beta^{(2)}=[f^{n'}(\beta_i^{(1)})]$, we have $\ell_{exp}(\alpha_i^{(n')}) \geq 14C+1$. By Remark~\ref{Rmq explication lemma}~$(2)$ applied twice (once with $\gamma_1=[f^{n'}(\alpha_i^{(1)})]$ and $\gamma_2=[f^{n'}(\beta_i^{(1)}\ldots\alpha_{k_1}^{(1)}\beta_{k_1}^{(1)})]$, and once with $\gamma_1=[f^{n'}(\alpha_i^{(1)})]^{-1}$ and $\gamma_2=[f^{n'}(\beta_0^{(1)}\ldots\alpha_{i-1}^{(1)}\beta_{i-1}^{(1)})]^{-1}$), the path $\alpha_i^{(n')}$ contains a complete factor of $[f^{n'}(\gamma)]$ of exponential length at least equal to $14C+1-4C=10C+1$. This proves Assertion~$(a)$ of Lemma~\ref{Lem complete splitting rel poly grow}. Moreover, Remark~\ref{Rmq explication lemma}~$(2)$ implies that the intersection of an incomplete factor of $[f^{n'}(\gamma)]$ with $\alpha_i^{(n')}$ is contained in the union of an initial and a terminal segment of $\alpha_i^{(n')}$ of exponential lengths at most $2C$. For every $i \in \{1,\ldots,k_1\}$ such that $\alpha_i^{(1)} \in \Lambda_{exp}$, let $\tau_i^{1}$ be the maximal initial segment of $\alpha_i^{(n')}$ of exponential length equal to $2C$ and let $\tau_i^{2}$ be the maximal terminal segment of $\alpha_i^{(n')}$ of exponential length equal to $2C$. 

We now prove Assertion~$(b)$ of Lemma~\ref{Lem complete splitting rel poly grow}. Suppose that there exists $i \in \{1,\ldots,k_1\}$ such that $\alpha_i^{(1)} \notin \Lambda_{exp}$, so that in particular $[f^{m_1}(\alpha_i^{(1)})]$ does not satisfy $(i)$. Then $[f^{m_1}(\alpha_i^{(1)})]$ satisfies $(ii)$ and is a concatenation of paths in $G_{PG}$ and in $\mathcal{N}_{PG}$. By Lemma~\ref{Lem iterate of a path in GPG}~$(3)$, the path $[f^{n'}(\alpha_i^{(1)})]$ is a concatenation of paths in $G_{PG}$ and in $\mathcal{N}_{PG}$. By Lemma~\ref{Lem concatenation paths Gpg Npg}, the path $[[f^{n'}(\beta_{i-1}^{(1)})][f^{n'}(\alpha_i^{(1)})][f^{n'}(\beta_i^{(1)})]]$ is a concatenation of paths in $G_{PG}$ and in $\mathcal{N}_{PG}$. Thus, the path $\beta_{i-1}^{(n')}\alpha_i^{(n')}\beta_i^{(n')}$ is a subpath of a concatenation of paths in $G_{PG}$ and in $\mathcal{N}_{PG}$. Hence $[f^{n'}(\gamma)]$ has a decomposition $$[f^{n'}(\gamma)]=\epsilon_1\alpha_1^{(n',+)}\epsilon_2\ldots\alpha_{k_2}^{(n',+)}\epsilon_{k_2}$$ where for every $j \in \{1,\ldots,k_2\}$, the path $\alpha_j^{(n',+)}$ is in $\Lambda_{exp}$ and for every $j \in \{0,\ldots,k_2\}$, the path $\epsilon_j$ is contained in a path $\iota_j$ which is a concatenation of paths in $G_{PG}$ and in $\mathcal{N}_{PG}$. Hence, for every $j \in \{0,\ldots,k_2\}$, we have $\ell_{exp}(\iota_j)=0$ by Lemma~\ref{Lem exponential length paths in Gpg} and, by Lemma~\ref{Lem bound exponential length subpath}, we have $\ell_{exp}(\epsilon_j) \leq 2C$.

If $\gamma'$ is an incomplete factor of $[f^{n'}(\gamma)]$, as explained above, there exists $i \in \{1,\ldots,k_2\}$ such that $\gamma'$ is contained in $\tau_{i-1}^2\epsilon_{i-1}\tau_i^1$. By Lemma~\ref{Lem bound exponential length subpath}, we have $$\ell_{exp}(\gamma') \leq \ell_{exp}(\tau_{i-1}^2\epsilon_{i-1}\tau_i^1)+2C.$$ By Lemma~\ref{Lem Exponential length less exp length subpaths}, the exponential length of $\gamma'$ is at most equal to $$\ell_{exp}(\tau_{i-1}^2)+\ell_{exp}(\epsilon_{i-1})+\ell_{exp}(\tau_i^1)+2C \leq 6C+\ell_{exp}(\epsilon_{i-1}) \leq 8C.$$ This proves $(b)$.

Finally, suppose that $\Lambda_{exp}$ is empty. For every $j \in \{1,\ldots,k_1\}$, the path $[f^{m_1}(\alpha_j^{(1)})]$ is a concatenation of paths in $G_{PG}$ and in $\mathcal{N}_{PG}$. By Lemma~\ref{Lem concatenation paths Gpg Npg}, the path $[f^{m_1}(\gamma)]$ is a concatenation of paths in $G_{PG}$ and in $\mathcal{N}_{PG}$. By Lemma~\ref{Lem iterate of a path in GPG}, for every $n' \geq m_1$, the path $[f^{n'}(\gamma)]$ is a concatenation of paths in $G_{PG}$ and in $\mathcal{N}_{PG}$. This concludes the proof.
\hfill\qedsymbol

\begin{lem}\label{Lem uniform bound on bad subpaths}
Let $f \colon G \to G$ be a $3K$-expanding CT map. There exists $N \in \NN^*$ such that for every reduced edge path $\gamma$ and every $m \geq N$, the total exponential length of incomplete factors in any optimal splitting of $[f^m(\gamma)]$ is uniformly bounded by $8C\ell_{exp}(\gamma)$.
\end{lem}

\dem By Proposition~\ref{Prop definition CT}~$(8)$, there exists $N \in \NN^*$ such that, for every reduced edge path $\alpha$ of length at most equal to $C+1$, the path $[f^N(\alpha)]$ is completely split. Suppose first that $\ell_{exp}(\gamma)=0$. Then, by definition of the exponential length, the path $\gamma$ is a concatenation of paths in $G_{PG}'$ and in $\mathcal{N}_{PG}$. By Proposition~\ref{Prop definition CT}~$(4)$, every edge in a zero stratum is adjacent to either an edge in a zero stratum or an edge in an EG stratum. Moreover, by Lemma~\ref{Lem No zero path and Nielsen path adjacent}, there does not exist a subpath of $\gamma$ contained in a zero stratum which is adjacent to a Nielsen path. Hence $\gamma$ is either a concatenation of paths in $G_{PG}$ and in $\mathcal{N}_{PG}$ or a path in a zero stratum. In the first case, the path $\gamma$ is $PG$-relative completely split. In the second case, by the definition of the constant $K$ and Equation~\eqref{Equation defi C}, we have $\ell(\gamma) \leq K \leq C$. By the choice of $N$, for every $m \geq N$, the path $[f^m(\gamma)]$ is completely split. By Lemma~\ref{Lem completely split and paths in Npg}, for every $m \geq N$, the path $[f^m(\gamma)]$ is $PG$-relative completely split. By Lemma~\ref{Lem exponential length paths in Gpg}, for every $m \geq N$, we have $\ell_{exp}([f^m(\gamma)])=0$.

So we may suppose that $\ell_{exp}(\gamma)>0$. Let $\gamma=\gamma_0\gamma_1'\gamma_1\ldots \gamma_{\ell}\gamma_{\ell}'$ be the exponential decomposition of $\gamma$ (see the beginning of Section~\ref{Section exponential length}). By Lemma~\ref{Lem No zero path and Nielsen path adjacent}, there does not exist a subpath of $\gamma$ contained in a zero stratum which is adjacent to a Nielsen path. Therefore, the path $\gamma$ has a decomposition $\alpha_0\beta_1\alpha_1\ldots\beta_k\alpha_k$ where, for every $i \in \{0,\ldots,k\}$, the path $\alpha_i$ is a (possibly trivial) concatenation of paths in $G_{PG}$ and in $\mathcal{N}_{PG}$ and, for every $i \in \{1,\ldots,k\}$, the path $\beta_i$ is a concatenation of a (possibly trivial) maximal reduced path in a zero stratum and an edge in an irreducible stratum not contained in $G_{PG}$ or in some $\gamma_i$. By construction of $K$, for every $i \in \{1,\ldots,k\}$, we have $\ell(\beta_i) \leq C+1$. By the choice of $N$, for every $m \geq N$, the path $[f^m(\beta_i)]$ is completely split. Note that, for every $i \in \{1,\ldots,k\}$, we have $\ell_{exp}(\beta_i)=1$ and that $$\ell_{exp}(\gamma)=\sum_{i=1}^k \ell_{exp}(\beta_i)=k.$$ By Lemma~\ref{Lem iterate of a path in GPG}, for every $i \in \{0,\ldots,k\}$ and every $m \geq M$, the path $[f^m(\alpha_i)]$ is a concatenation of paths in $G_{PG}$ and in $\mathcal{N}_{PG}$. By Lemma~\ref{Lem exponential length paths in Gpg}, for every $m \geq M$, we have $\ell_{exp}([f^m(\alpha_i)])=0$. By Lemma~\ref{Lem bound exponential length subpath}, the exponential length of the subpath of $[f^m(\alpha_i)]$ contained in $[f^m(\gamma)]$ is at most equal to $2C$. For every $i \in \{0,\ldots,k\}$ (resp. $i \in \{1,\ldots,k\})$ and every $m \geq N$, let $\alpha_{i,m}$ (resp. $\beta_{i,m}$) be the subpath of $[f^m(\alpha_i)]$ (resp. $[f^m(\beta_i)]$) contained in $[f^m(\gamma)]$. By Remark~\ref{Rmq explication lemma}~$(2)$, for every $i \in \{1,\ldots,k\}$ and every $m \geq N$, the exponential length of any incomplete factor in $\beta_{i,m}$ is at most equal to $4C$. By Lemma~\ref{Lem Exponential length less exp length subpaths}, for every $m \geq N$, the sum of the exponential lengths of the incomplete factors in $[f^m(\gamma)]$ is at most equal to $$\sum_{i=0}^k \ell_{exp}(\alpha_{i,m})+4Ck \leq 2C(k+1) +4kC \leq 4Ck+4Ck=8Ck=8C\ell_{exp}(\gamma).$$ The conclusion of the lemma follows.
\hfill\qedsymbol

\begin{lem}\label{Lem zero relative exp length}
Let $f \colon G \to G$ be a $3K$-expanding CT map. Let $\gamma$ be a reduced edge path in $G$. Suppose that $\gamma$ has a splitting $\gamma=b_1ab_2$ where, for every $i \in \{1,2\}$, the path $b_i$ is a possibly trivial $PG$-relative completely split. If $\ell_{exp}^{\gamma}(a)=0$ then $\ell_{exp}(a)=0$.
\end{lem}

\dem Let $\gamma=\gamma_0\gamma_1'\gamma_1\ldots \gamma_k\gamma_k'$ be the exponential decomposition of $\gamma$. By Lemma~\ref{Lem bound exponential length subpath}, there exist three (possibly trivial) paths $\delta_1$, $\delta_2$ and $\tau$ such that for every $i \in \{1,2\}$, the path $\delta_i$ is a proper initial or terminal subpath of a splitting unit of some $\gamma_j$, we have $\ell_{exp}(\tau)=\ell_{exp}^{\gamma}(\tau)=\ell_{exp}^{\gamma}(a)$ and $a=\delta_1\tau\delta_2$. Since $\ell_{exp}^{\gamma}(a)=0$, we have $\ell_{exp}(\tau)=0$. Hence $\tau$ is a concatenation of paths in $G_{PG}'$ and in $\mathcal{N}_{PG}$. By Proposition~\ref{Prop definition CT}~$(4)$, every edge in a zero stratum is adjacent to either an edge in a zero stratum or an edge in an EG stratum. Moreover, by Lemma~\ref{Lem No zero path and Nielsen path adjacent}, there does not exist a subpath of $\gamma$ contained in a zero stratum which is adjacent to a Nielsen path. Hence $\tau$ is either a concatenation of paths in $G_{PG}$ and in $\mathcal{N}_{PG}$ or a path in a zero stratum. If $\tau$ is contained in a zero stratum, by Lemma~\ref{Lem No zero path and Nielsen path adjacent}, we see that $\delta_1$ and $\delta_2$ are trivial, that is, $a=\tau$. Thus, we have $\ell_{exp}(a)=\ell_{exp}(\tau)=0$.

So we may suppose that $\tau$ is a concatenation of paths in $G_{PG}$ and in $\mathcal{N}_{PG}$. Suppose towards a contradiction that there exists $i \in \{1,2\}$ such that $\delta_i$ is not trivial. For every $i \in \{1,2\}$ such that $\delta_i \neq \varnothing$, let $\sigma_i$ be the splitting unit of some $\gamma_j$ containing $\delta_i$ and let $r_i$ be the height of $\sigma_i$. By~\cite[Lemma~5.11]{BesHan92}, for every $i \in \{1,2\}$ such that $\delta_i$ is not trivial, there exist two distinct $r_i$-legal paths $\alpha_i$ and $\beta_i$ such that $\sigma_i=\alpha_i\beta_i$ and such that the turn $\{Df(\alpha_i^{-1}),Df(\beta_i)\}$ is the only height $r_i$ illegal turn. Moreover, there exists a path $\tau_i'$ such that $[f(\alpha_i)]=\alpha_i\tau_i'$ and $[f(\beta_i)]=\tau_i'^{-1}\beta_i$. Let $\epsilon_1^{(1)},\epsilon_1^{(2)}$ be two paths such that $\sigma_1=\epsilon_1^{(1)}\epsilon_1^{(2)}$, the path $\epsilon_1^{(1)}$ is contained in $b_1$ and the path $\epsilon_1^{(2)}$ is contained in $a$. Similarly, let $\epsilon_2^{(1)},\epsilon_2^{(2)}$ be two paths such that $\sigma_2=\epsilon_2^{(1)}\epsilon_2^{(2)}$, the path $\epsilon_2^{(2)}$ is contained in $b_2$ and the path $\epsilon_2^{(1)}$ is contained in $a$. 

\medskip

\noindent{\bf Claim. } $(1)$ For every path $b \in \mathcal{N}_{PG}^{\max}(b_1)$ (resp. $b \in \mathcal{N}_{PG}^{\max}(b_2)$), the path $b$ does not contain edges of $\epsilon_1^{(1)}$ (resp. $\epsilon_2^{(2)}$).

\noindent{$(2)$ } The path $\epsilon_1^{(1)}$ is $r_1$-legal and the path $\epsilon_2^{(2)}$ is $r_2$-legal.

\medskip

\dem We prove the claim for $b_1$, the proof for $b_2$ being similar.

\medskip

\noindent{$(1)$ } Let $b \in \mathcal{N}_{PG}^{\max}(b_1)$. There exists $c \in \mathcal{N}_{PG}^{\max}(\gamma)$ such that $b \subseteq c$. Moreover, by Lemma~\ref{Lem Nielsen paths in NPG properties}~$(3)$ applied to $\gamma'=b$ and $\gamma=c$, either $b$ is a concatenation of splitting units of $c$, or $b$ is properly contained in a splitting unit of $c$ and is not an initial or a terminal segment of $c$. Since $b_1$ is an initial segment of $\gamma$, the second case cannot occur. Hence $b$ is a concatenation of splitting units of $c$.  Since $\sigma_1$ is not contained in $b_1$, the path $b$ cannot contain edges of $\sigma_1$. Since $\epsilon_1^{(1)} \subseteq \sigma_1$, the path $b$ cannot contain edges of $\epsilon_1^{(1)}$.

\medskip

\noindent{$(2)$ } Suppose towards a contradiction that $\epsilon_1^{(1)}$ is not $r_1$-legal. Then it contains the illegal turn $\{Df(\alpha_1^{-1}),Df(\beta_2)\}$.  Recall that the path $b_1$ is $PG$-relative completely split. By the description of $PG$-relative splitting units, the illegal turn must be contained in a $PG$-relative splitting unit of $b_1$ which is a concatenation of paths in $G_{PG}$ and in $\mathcal{N}_{PG}$. Since the last edge of $\alpha_1$ is an edge in an EG stratum, the last edge of $\alpha_1$ must be contained in a path contained in $\mathcal{N}_{PG}$. Hence $\epsilon_1$ intersects a path in $\mathcal{N}_{PG}^{\max}(b_1)$. This contradicts Assertion~$(1)$.
\hfill\qedsymbol

\medskip

By Assertion~$(2)$ of the claim, for every $i \in \{1,2\}$ such that $\sigma_i$ is not trivial, the path $\epsilon_i^{(i)}$ is $r_i$-legal. Moreover, by Assertion~$(1)$ of the claim an INP contained in $b_i$ cannot intersect the path $\epsilon_i^{(i)}$. Since the paths $b_1$ and $b_2$ are $PG$-relative completely split, the paths $b_1$ and $b_2$ split respectively at the origin of $\epsilon_1^{(1)}$ and at the end of $\epsilon_2^{(2)}$. So we may suppose that $b_1=\epsilon_1^{(1)}$ and $b_2=\epsilon_2^{(2)}$. Therefore, there exists a (possibly trivial) path $\tau_1$ such that, up to taking a power of $f$ so that the length of $[f(b_1)]$ is greater than $\alpha_1$, we have $[f(b_1)]=\alpha_1\tau_1$ and $[f(\epsilon_1^{(2)})]=\tau_1^{-1}\beta_1$. Similarly, there exists a path $\tau_2$ such that $[f(\epsilon_2^{(1)})]=\alpha_2\tau_2$ and $[f(b_2)]=\tau_2^{-1}\beta_2$. 

Since $\gamma$ splits at the concatenation points of $b_1$, $a$ and $b_2$, the paths $\tau_1^{-1}$ and $\tau_2$ contained in $[f(\epsilon_1^{(2)})][f(\tau)][f(\epsilon_2^{(1)})]$ must be identified when passing to $[f(a)]$. Suppose first that $[f(\tau)]$ is a point. Then since the EG INPs $\sigma_1$ and $\sigma_2$ are uniquely determined by their initial and terminal edges by Proposition~\ref{Prop definition CT}~$(9)$, we see that $\sigma_1=\sigma_2^{-1}$. But then there are some identifications between $b_1$ and $b_2$, which contradicts the fact that $b_1ab_2$ is a splitting. 

Thus, we may suppose that $[f(\tau)]$ is nontrivial. By Lemma~\ref{Lem iterate of a path in GPG}, since $\tau$ is a concatenation of paths in $G_{PG}$ and in $\mathcal{N}_{PG}$ so is $[f(\tau)]$. Note that, since an EG INP is completely determined by its initial and terminal edges by Proposition~\ref{Prop definition CT}~$(9)$, if $[f(\tau)]$ contains the initial or the terminal edge of an EG INP $\sigma$, then $\sigma$ is contained in $[f(\tau)]$. Note that there are identifications between edges of $[f(\epsilon_1^{(2)})]$ and $[f(\tau)]$ or between edges of $[f(\tau)]$ and $[f(\epsilon_2^{(1)})]$. Therefore, $[f(\tau)]$ starts with $\sigma_1^{-1}$ or $[f(\tau)]$ ends with $\sigma_2^{-1}$. Thus, one of the following holds: 

\medskip

\noindent{$(a)$ } $[f(\tau)]=\sigma_1^{-1}\tau'$ with $\tau'$ a (possibly trivial) path which is a concatenation of paths in $G_{PG}$ and in $\mathcal{N}_{PG}$ which does not end by $\sigma_2^{-1}$;

\medskip

\noindent{$(b)$ } $[f(\tau)]=\tau'\sigma_2^{-1}$ with $\tau'$ a (possibly trivial) path which is a concatenation of paths in $G_{PG}$ and in $\mathcal{N}_{PG}$ which does not start  by $\sigma_1^{-1}$; 

\medskip

\noindent{$(c)$ } $[f(\tau)]=\sigma_1^{-1}\tau'\sigma_2^{-1}$ with $\tau'$ a (possibly trivial) path. 

\medskip

Note that $\sigma_1^{-1}\tau'\sigma_2^{-1}$ is reduced, so that there is no identification between $\alpha_1^{-1}$ and $\tau'$ and between $\tau'$ and $\beta_2^{-1}$. Let $e_{\sigma_1}$ be the terminal edge of $\sigma_1$ and let $e_{\sigma_2}$ be the initial edge of $\sigma_2$. By Proposition~\ref{Prop definition CT}~$(9)$, both $e_{\sigma_1}$ and $e_{\sigma_2}$ are edges in EG strata. Since $f$ is $3K$-expanding, for every $i \in \{1,2\}$, the path $[f(e_{\sigma_i})]$ has length at least equal to $3K$. Recall that, for every $i \in \{1,2\}$, by definition of $K$, we have $\ell(\sigma_i) \leq K$, so that $\ell(\alpha_i),\ell(\beta_i) \leq K$. Since $[f(\epsilon_1^{(2)})]=\alpha_1\tau_1$ and $[f(\epsilon_2^{(1)})]=\alpha_2\tau_2$, the path $[f(e_{\sigma_1})]$ contains a nondegenerate terminal segment of $\tau_1^{-1}$ and the path $[f(e_{\sigma_2})]$ contains a nondegenerate initial segment of length $2K$ of $\tau_2$. As $e_{\sigma_1}$ is $r_1$-legal and as $f$ is a relative train track by Proposition~\ref{Prop definition CT}~$(1)$, we see that the last edge of $\tau_1^{-1}$ is not the last edge of $\alpha_1$. Similarly, the first edge of $\tau_2$ is not the first edge of $\beta_2$. Therefore, we have $[\tau_1^{-1}\beta_1\sigma_1^{-1}]=\tau_1^{-1}\alpha_1^{-1}$ and $[\sigma_2^{-1}\alpha_2\tau_2]=\beta_2^{-1}\tau_2$. Thus we have $$[[f(\epsilon_1^{(2)})][f(\tau)][f(\epsilon_2^{(1)})]]=[\tau_1^{-1}\beta_1\sigma_1^{-1}\tau'\sigma_2^{-1}\alpha_2\tau_2]=[\tau_1^{-1}\alpha_1^{-1}\tau'\beta_2^{-1}\tau_2],$$ and there is no identification between $\tau_1^{-1}$ and $\alpha_1^{-1}$, $\alpha_1^{-1}$ and $\tau'$, $\tau'$ and $\beta_2^{-1}$ and $\beta_2^{-1}$ and $\tau_2$. Therefore, if $\tau'$ is not trivial, then we have a contradiction as $\tau_1^{-1}$ and $\tau_2$ are not identified in $[f(a)]$. Suppose that $\tau'$ is trivial. Then the paths $\tau_1^{-1}$ and $\tau_2$ are identified in $[f(a)]$ only if a terminal segment of $\alpha_1^{-1}$ is identified with an initial segment of $\beta_2^{-1}$. Since EG INP are uniquely determined by their initial and terminal edges by Proposition~\ref{Prop definition CT}~$(9)$, we see that $\sigma_1=\sigma_2^{-1}$. Hence $\alpha_1^{-1}=\beta_2$ and either $\tau_1^{-1}$ is an initial segment of $\tau_2^{-1}$ or $\tau_2$ is an initial segment of $\tau_1$. Up to changing the orientation of $\gamma$, we may suppose that $\tau_1^{-1}$ is an initial segment of $\tau_2^{-1}$. If $\tau_1^{-1}=\tau_2^{-1}$, then $[f(a)]$ is a vertex. Moreover, as $\sigma_1=\sigma_2^{-1}$, the segment $b_1=\epsilon_1^{(1)}$ is equal to $b_2^{-1}$. Therefore, a terminal segment of $b_1$ is identified with an initial segment of $b_2$, a contradiction. If $\tau_1^{-1}$ is a proper initial segment of $\tau_2^{-1}$, then $\tau_2$ is identified with edges in $b_1$, a contradiction. As we have considered every case, we see that $\delta_1$ and $\delta_2$ are trivial and $\ell_{exp}(a)=\ell_{exp}(\tau)=0$. 

\hfill\qedsymbol

\begin{lem}\label{Lem Goodness nondecreasing}
Let $f \colon G \to G$ be a $3K$-expanding CT map. There exists $n_0 \in \NN^*$ such that for every $n \geq n_0$, and every closed reduced edge path $\gamma$ of $G$, we have the following relation between the goodness of $\gamma$ and the one of $[f^n(\gamma)]$: $$\mathfrak{g}([f^n(\gamma)]) \geq \mathfrak{g}(\gamma).$$
\end{lem}

\dem By Lemma~\ref{Lem exponential length monotonic for completely split edge path}, there exists $N_0 \in \NN^*$ such that, for every $n \geq N_0$ and every $PG$-relative splitting unit $\sigma$, the exponential length of the path $[f^n(\sigma)]$ is at least equal to the one of $\sigma$. By Lemma~\ref{Lem uniform bound on bad subpaths}, there exists $N_1$ such that for every $n \geq N_1$ and every closed reduced edge path $\gamma$ of $G$, the total exponential length of incomplete segments in any optimal splitting of $[f^n(\gamma)]$ is bounded by $8C\ell_{exp}(\gamma)$. Let $N_2=\lceil \log_3(10C+16C^2) \rceil \in \NN^*$ be such that for every $x,y \geq 0$ such that $(x,y) \neq (0,0)$, we have $$\frac{(3^{N_2}-2C)x}{(3^{N_2}-2C)x+8C(1+2C)y} \geq \frac{x}{x+y}.$$ Let $n_0=\max\{N_0,N_1,N_2\}$.

Let $\gamma$ be a closed reduced edge path in $G$. All splittings of $\gamma$ are circuital splittings in what follows. Let $\gamma=\alpha_0\beta_1\alpha_1\ldots\beta_k\alpha_k$ be an optimal splitting of $\gamma$, where for every $i \in \{0,\ldots,k\}$, the path $\alpha_i$ is an incomplete factor of $\gamma$ and for every $i \in \{1,\ldots,k\}$, the path $\beta_i$ is a $PG$-relative complete factor of $\gamma$. First note that, for every $i \in \{1,\ldots,k\}$, and for every $n \geq 1$, the path $[f^n(\beta_i)]$ is $PG$-relative completely split by Proposition~\ref{Prop definition CT}~$(6)$ and Lemma~\ref{Lem iterate of a path in GPG}. Therefore, if $n \geq n_0 \geq N_0$, the total exponential length of such $PG$-relative complete segments is nondecreasing under $[f^n]$. We now distinguish two cases, according to the growth of the paths $\beta_i$. 

Suppose first that for every $i \in \{1,\ldots,k\}$, the exponential length of $\beta_i$ relative to $\gamma$ is equal to zero. Since the splitting $\gamma=\alpha_0\beta_1\alpha_1\ldots\beta_k\alpha_k$ is optimal and since for every $i \in \{1,\ldots,k\}$, we have $\ell_{exp}^{\gamma}(\beta_i)=0$, we have $\mathfrak{g}(\gamma)=0$. Therefore, for every $n \in \NN^*$, we have $\mathfrak{g}([f^n(\gamma)]) \geq \mathfrak{g}(\gamma)$.

Suppose now that there exists $i \in \{1,\ldots,k\}$ such that the exponential length of $\beta_i$ relative to $\gamma$ is positive. By Lemma~\ref{Lem exponential length goes to infinity}, the sequence $(\ell_{exp}([f^n(\beta_i)]))_{n \in \NN^*}$ grows exponentially with $n$. We can now modify the splitting of $\gamma$ into the following splitting: $\gamma=\alpha_0'\beta_1'\alpha_1'\ldots\beta_m'\alpha_m'$ where:

\noindent{$(a)$ } for every $j \in \{0,\ldots,m\}$, the path $\alpha_i'$ is a concatenation of incomplete factors and complete factors of zero exponential length relative to $\gamma$ of the previous splitting;

\noindent{$(b)$ } for every $j \in \{1,\ldots,m\}$, the path $\beta_i'$ is a complete factor of positive exponential length relative to $\gamma$ of the previous splitting. 
 
Note that, by definition of the exponential length relative to $\gamma$, for every \mbox{$i \in \{1,\ldots,m\}$} and every path $\gamma' \in \mathcal{N}_{PG}^{\max}(\gamma)$, the path $\beta_i'$ is not contained in $\gamma'$. Therefore, if there exists $j \in \{0,\ldots,m\}$ and $\gamma' \in \mathcal{N}_{PG}^{\max}(\gamma)$ such that $\alpha_j'$ intersects $\gamma'$ nontrivially, then $\gamma'$ is contained in $\beta_{j-1}'\alpha_j'\beta_j'$. In particular, Lemma~\ref{Lem zero relative exp length} applies and for every $j \in \{0,\ldots,m\}$, if $\ell_{exp}^{\gamma}(\alpha_j')=0$, then $\ell_{exp}(\alpha_j')=0$. Let $\Lambda$ be the subset of $\{0,\ldots,m\}$ such that for every $j \in \Lambda$, we have $\ell_{exp}^{\gamma}(\alpha_j')>0$. 

By Lemma~\ref{Lem bound exponential length subpath} and Lemma~\ref{Lem p length of a completely split path grows linearly}, for every $j \in \{1,\ldots,m\}$ and every $M \in \NN^*$, we have $$\ell_{exp}^{[f^M(\gamma)]}([f^M(\beta_i')]) \geq \ell_{exp}([f^M(\beta_i')])-2C\geq 3^M\ell_{exp}(\beta_i')-2C' \geq (3^M-2C)\ell_{exp}^{\gamma}(\beta_i').$$ By Lemma~\ref{Lem bound exponential length subpath}, for every $j \in \{0,\ldots,m\}$, we have $\ell_{exp}^{\gamma}(\alpha_j') \leq \ell_{exp}(\alpha_j')$. Note that, for every $i \in \{1,\ldots,m\}$, and every $n \in \NN^*$, the path $[f^n(\beta_i')]$ is $PG$-relative completely split. In particular, for every $n \in \NN^*$, any incomplete factor of $[f^n(\gamma)]$ is contained in a reduced iterate of some $\alpha_i'$. Thus, by Lemma~\ref{Lem uniform bound on bad subpaths}, for every $n \geq n_0 \geq N_1$, the total exponential length of incomplete segments in $[f^n(\gamma)]$ is bounded by $8C\sum_{j=1}^k\ell_{exp}(\alpha_j')=8C \sum_{j \in \Lambda}\ell_{exp}(\alpha_j')$. Note that the function $$x \mapsto \frac{x}{x+8C\sum_{j \in \Lambda} \ell_{exp}(\alpha_j')}$$ is nondecreasing. Recall that, for every $n \in \NN^*$, the goodness function is a supremum over splittings of $[f^n(\gamma)]$. Thus, by Lemma~\ref{Lem goodness versus decomposition}, for every $n \geq n_0$, we have:
$$\mathfrak{g}([f^n(\gamma)]) \geq \frac{(3^{n}-2C)\sum_{i=1}^m \ell_{exp}^{\gamma}(\beta_i')}{(3^{n}-2C)\sum_{i=1}^m \ell_{exp}^{\gamma}(\beta_i')+8C\sum_{j \in \Lambda} \ell_{exp}(\alpha_j')}.$$ By Lemma~\ref{Lem bound exponential length subpath}, we have $$8C\sum_{j \in \Lambda} \ell_{exp}(\alpha_j') \leq 8C\sum_{j \in \Lambda} (\ell_{exp}^{\gamma}(\alpha_j')+2C) \leq 8C(1+2C)\sum_{j \in \Lambda} \ell_{exp}^{\gamma}(\alpha_j'),$$ where the last inequality follows from the fact that, for every $j \in \Lambda$, we have $\ell_{exp}^{\gamma}(\alpha_j')\geq 1$. Therefore, since $n_0 \geq N_2$, for every $n \geq n_0$, we have: 
$$\frac{(3^{n}-2C)\sum_{j=1}^m \ell_{exp}^{\gamma}(\beta_j')}{(3^{n}-2C)\sum_{j=1}^m \ell_{exp}^{\gamma}(\beta_j')+8C(1+2C)\sum_{j \in \Lambda} \ell_{exp}^{\gamma}(\alpha_j')} \geq \frac{\sum_{j=1}^m \ell_{exp}^{\gamma}(\beta_j')}{\sum_{j=1}^m \ell_{exp}^{\gamma}(\beta_j')+\sum_{j \in \Lambda} \ell_{exp}^{\gamma}(\alpha_j')}.$$ 
By Lemma~\ref{Lem compute exponential length optimal factor}, we have $$\ell_{exp}(\gamma)=\sum_{j=1}^m \ell_{exp}^{\gamma}(\beta_j')+\sum_{j=0}^m \ell_{exp}^{\gamma}(\alpha_j')=\sum_{j=1}^m \ell_{exp}^{\gamma}(\beta_j')+\sum_{j \in \Lambda} \ell_{exp}^{\gamma}(\alpha_j').$$ Thus, we see that $$\frac{\sum_{j=1}^m \ell_{exp}^{\gamma}(\beta_j')}{\sum_{j=1}^m \ell_{exp}^{\gamma}(\beta_j')+\sum_{j \in \Lambda} \ell_{exp}^{\gamma}(\alpha_j')}=\mathfrak{g}(\gamma),$$ which gives the result.
\hfill\qedsymbol

\begin{rmq}\label{Rmq Convention for relative atoroidal CT map 2}
In the next lemmas, we will adopt the following conventions.

\noindent{Let $\phi \in \Out(F_{\tt n},\mathcal{F})$ be an almost atoroidal outer automorphism relative to $\mathcal{F}$.
Let \mbox{$f \colon G \to G$} be a CT map representing a power of $\phi$ with filtration $$\varnothing=G_0 \subsetneq \ldots \subsetneq G_k=G.$$ Let $p \in \{1,\ldots,k-1\}$ be such that $\mathcal{F}(G_p)=\mathcal{F}$. By Lemma~\ref{Lem exponential length goes to infinity}, up to taking a power of $f$, we may suppose that $f$ is $3K$-expanding. By Lemma~\ref{Lem Goodness nondecreasing}, up to passing to a power of $f$, we may suppose that for every closed reduced edge path $\gamma$ of $G$, we have $\mathfrak{g}([f(\gamma)]) \geq \mathfrak{g}(\gamma)$. }
\end{rmq}

\begin{lem}\label{Lem control of the goodness}
Let $f \colon G \to G$ be as in Remark~\ref{Rmq Convention for relative atoroidal CT map 2}. 

\medskip

\noindent{$(1)$ } For every $\delta > 0$, there exists $m \in \NN^*$ such that for every reduced edge path $\gamma$ such that $\mathfrak{g}(\gamma)\geq \delta$ and every $n \geq m$, the total exponential length relative to $[f^n(\gamma)]$ of complete factors in $[f^n(\gamma)]$ denoted by $TEL(n,\gamma)$ is at least $$TEL(n,\gamma) \geq \mathfrak{g}(\gamma)\ell_{exp}(\gamma)(3^{n}-2C).$$

\medskip

\noindent{$(2)$ } For every $\delta >0$ and every $\epsilon >0$, there exists $m \in \NN^*$ such that for every cyclically reduced circuit $\gamma$ such that $\ell_{exp}(\gamma)>0$, $\mathfrak{g}(\gamma)\geq \delta$ and every $n \geq m$, we have $\mathfrak{g}([f^n(\gamma)]) \geq 1-\epsilon$.
\end{lem}

\dem Let $\gamma=\alpha_0\beta_1\alpha_1\ldots\alpha_k\beta_k$ be an optimal splitting, where for every $i \in \{0,\ldots,k\}$, the path $\alpha_i$ is an incomplete factor of $\gamma$ and for every $i \in \{1,\ldots,k\}$, the path $\beta_i$ is a $PG$-relative complete factor of $\gamma$. We may assume that $\ell_{exp}(\gamma)>0$, otherwise $\mathfrak{g}(\gamma)=0$ and the result is immediate. Note that, since $\mathfrak{g}(\gamma) \geq \delta>0$, there exists $i \in \{1,\ldots,k\}$ such that $\ell_{exp}^{\gamma}(\beta_i) >0$. Let $\Lambda_{\gamma}$ be the set consisting in all complete factors $\beta_i$ of $\gamma$ whose exponential length relative to $\gamma$ is positive. Let $\ell_{exp}^{\gamma}(\Lambda_{\gamma})$ be the sum of the exponential lengths relative to $\gamma$ of all factors that belongs to $\Lambda_{\gamma}$. Note that $$\ell_{exp}^{\gamma}(\Lambda_{\gamma})=\sum_{\beta_i \in \Lambda_{\gamma}}\ell_{exp}^{\gamma}(\beta_i)=\mathfrak{g}(\gamma)\ell_{exp}(\gamma).$$

Note that, for every $n \in \NN^*$, the value $TEL(n,\gamma)$ is a supremum over all splittings of $[f^n(\gamma)]$. Thus, by Lemma~\ref{Lem bound exponential length subpath} and Lemma~\ref{Lem p length of a completely split path grows linearly}, for every $n \in \NN^*$, we have:

$$
TEL(n, \gamma)  \geq  \sum_{\beta_i \in \Lambda_{\gamma}} \ell_{exp}^{[f^n(\gamma)]}([f^n(\beta_i)]) \geq (3^n-2C)\ell_{exp}^{\gamma}(\Lambda_{\gamma}) \geq (3^n-2C) \mathfrak{g}(\gamma)\ell_{exp}(\gamma).
$$

This proves $(1)$.
We now prove $(2)$. By Lemma~\ref{Lem uniform bound on bad subpaths}, there exists $n_0 \in \NN^*$ such that for every $n \geq n_0$, the total exponential length of incomplete segments in $[f^n(\gamma)]$ is bounded by $8C\ell_{exp}(\gamma)$. By Lemma~\ref{Lem bound exponential length subpath}, the total exponential length relative to $\gamma$ of incomplete segments in $[f^n(\gamma)]$ is hence bounded by $10C\ell_{exp}(\gamma)$. Note that, for every $n \in \NN^*$, the value $\mathfrak{g}([f^n(\gamma)])$ is a supremum over all splittings of $[f^n(\gamma)]$. Thus, by Lemma~\ref{Lem goodness versus decomposition}, for every $n \geq n_0$, we have:

\begin{align*}
\mathfrak{g}([f^n(\gamma)]) & \geq \frac{ \mathfrak{g}(\gamma)\ell_{exp}(\gamma)(3^n-2C)}{ 10C\ell_{exp}(\gamma)+\mathfrak{g}(\gamma)\ell_{exp}(\gamma)(3^n-2C)} & {} \\
{} &  = \frac{ \mathfrak{g}(\gamma)(3^n-2C)}{ 10C+\mathfrak{g}(\gamma)(3^n-2C)} \geq  \frac{ \delta (3^n-2C)}{10C+\delta (3^n-2C)}.
\end{align*}

\noindent The last term is independent of $\gamma$ and converges to $1$ as $n$ goes to infinity. Therefore the conclusion of Lemma~\ref{Lem control of the goodness} holds for some $n$ large enough which does not depend on $\gamma$. This proves~$(2)$ and this concludes the proof.
\hfill\qedsymbol

\medskip

\subsection{North-South dynamics for a relative atoroidal outer automorphisms}

Let ${\tt n} \geq 3$ and let $\mathcal{F}$ be a free factor system of $F_{\tt n}$. Let $\phi \in \Out(F_{\tt n},\mathcal{F})$ be an almost atoroidal automorphism relative to $\mathcal{F}$. In this subsection we prove Theorem~\ref{Theo North south dynamics relative atoroidal}. The proof of Theorem~\ref{Theo North south dynamics relative atoroidal} is inspired by the proof of the same result due to Uyanik (\cite{Uyanik2019}) in the context of an atoroidal outer automorphism for $\Out(F_{\tt n})$, that is, in the special case when $\mathcal{F}=\varnothing$. The proof relies on the study of splittings of reduced edge paths in the graph associated with a CT map representing a power of $\phi$. Indeed, we show that, when a cyclically reduced edge path representing $w \in F_{\tt n}$ has a splitting which is close to a complete splitting, then some iterate of $\phi$ sends $[w]$ into an open neighborhood of $\Delta_+(\phi)$ (see Definition~\ref{Defi attractive convex}), and this iterate can be chosen uniformly (see~Lemma~\ref{Lem conversion goodness closeness PCurr}).

\bigskip

Let $\phi \in \Out(F_{\tt n},\mathcal{F})$ be an almost atoroidal outer automorphism which satisfies Definition~\ref{Defi almost atoroidal outer automorphism}~$(2)$. Let $\mathcal{F} \leq \mathcal{F}_1 \leq \mathcal{F}_2=\{[F_{\tt n}]\}$ be a sequence of free factor system given in this definition. Let $f \colon G \to G$ be a CT map representing a power of $\phi$ with filtration $\varnothing=G_0 \subsetneq G_1 \subsetneq \ldots \subsetneq G_k=G$ and such that there exist $p$ and $i$ in $\{1,\ldots,k\}$ such that $\mathcal{F}(G_p)=\mathcal{F}$ and $\mathcal{F}(G_i)=\mathcal{F}_1$. We denote by $\Curr(\mathcal{F}_1,\mathcal{F}_1 \wedge \mathcal{A}(\phi))$ the set of currents of $\Curr(F_{\tt n},\mathcal{F}_1 \wedge \mathcal{A}(\phi))$ whose support is contained in $\partial^2\mathcal{F}_1$. Note that, since the extension $\mathcal{F}_1 \leq \{[F_{\tt n}]\}$ is sporadic, either $\mathcal{F}_1=\{[H_1],[H_2]\}$ or $\mathcal{F}_1=\{[H]\}$ for some subgroups $H_1,H_2$ of $F_{\tt n}$. Up to assuming that $H_2$ is the trivial group, we may assume that $\mathcal{F}_1=\{[H_1],[H_2]\}$. Moreover, we have $\mathcal{F}_1 \wedge \mathcal{A}(\phi)=\{[A_1],\ldots,[A_s],[B_1],\ldots,[B_{t}]\}$ where, for every $j \in \{1,\ldots,s\}$, the group $A_j$ is contained in $H_1$ and for every $j \in \{1,\ldots,t\}$, the group $B_j$ is contained in $H_2$. Since $\mathcal{F}_1 \wedge \mathcal{A}(\phi)$ is a malnormal subgroup system, the set $\{[A_1],\ldots,[A_s]\}$ is a malnormal subgroup system of $H_1$ and the set $\{[B_1],\ldots,[B_t]\}$ is a malnormal subgroup system of $H_2$.

Let $$X(\mathcal{F}_1)=\Curr(H_1,\{[A_1],\ldots,[A_s]\}) \times \Curr(H_2,\{[B_1],\ldots,[B_t]\}).$$ Let $\mu \in \Curr(\mathcal{F}_1,\mathcal{F}_1 \wedge \mathcal{A}(\phi))$. We set $\psi_1(\mu)=(\mu|_{\partial^2 H_1},\mu|_{\partial^2H_2}) \in X(\mathcal{F}_1)$. Since $\mu$ is $F_{\tt n}$-invariant, $\psi_1(\mu)$ does not depend on the choice of the representatives of the conjugacy classes of $H_1$ and $H_2$. Let $(\mu_1,\mu_2) \in X(\mathcal{F}_1)$. Since the subgroup system $\mathcal{F}_1 \wedge \mathcal{A}(\phi)$ is malnormal, for every $j \in \{1,2\}$, the current $\mu_j$ can be extended in a canonical way to a current $\mu_j^{\ast} \in \Curr(F_{\tt n},\mathcal{F}_1 \wedge \mathcal{A}(\phi))$. The current $\mu_j^{\ast}$ is such that, for every Borel subset $B$ of $\partial^2(F_{\tt n},\mathcal{F}_1 \wedge \mathcal{A}(\phi))$, we have $$\mu_j^{\ast}(B)=\mu_j^{\ast}(B \cap \partial^2H_j)=\mu_j(B \cap \partial^2H_j).$$ We set $\psi_2((\mu_1,\mu_2))=\mu_1^{\ast}+\mu_2^{\ast}$. By the property of $\mu_j^{\ast}$ described above, we see that $\psi_2((\mu_1,\mu_2)) \in \Curr(\mathcal{F}_1,\mathcal{F}_1 \wedge \mathcal{A}(\phi))$. The maps $\psi_1$ and $\psi_2$ are clearly continuous.

\begin{lem}\label{Lem homeo curr F1}
The space $\Curr(\mathcal{F}_1,\mathcal{F}_1 \wedge \mathcal{A}(\phi))$ is homeomorphic to $X(\mathcal{F}_1)$.
\end{lem}

\dem We prove that $\psi_1$ and $\psi_2$ are inverse from each other. Let $\mu \in \Curr(F_{\tt n},\mathcal{F}_1 \wedge \mathcal{A}(\phi))$. Then $\psi_2 \circ \psi_1(\mu)=(\mu|_{\partial^2 H_1})^{\ast}+(\mu|_{\partial^2 H_2})^{\ast}$. Note that $\mu$ and $\psi_2 \circ \psi_1(\mu)$ coincide on Borel subsets contained in $\partial^2\mathcal{F}_1$. Since both have supports contained in $\partial^2\mathcal{F}_1$, they are equal. Conversely, let $(\mu_1,\mu_2) \in X(\mathcal{F}_1)$. Then $$\psi_1\circ\psi_2((\mu_1,\mu_2))=((\mu_1^{\ast}+\mu_2^{\ast})|_{\partial^2 H_1},(\mu_1^{\ast}+\mu_2^{\ast})|_{\partial^2 H_2}).$$ But $\mu_2^{\ast}|_{\partial^2 H_1}=0$ and $\mu_1^{\ast}|_{\partial^2 H_2}=0$. Hence we have $$((\mu_1^{\ast}+\mu_2^{\ast})|_{\partial^2 H_1},(\mu_1^{\ast}+\mu_2^{\ast})|_{\partial^2 H_2})=(\mu_1^{\ast}|_{\partial^2 H_1},\mu_2^{\ast}|_{\partial^2 H_2})=(\mu_1,\mu_2).$$ This concludes the proof.
\hfill\qedsymbol

\bigskip

For every $\phi \in \Out(F_{\tt n},\mathcal{F})$, we refer to the definition of $\mathcal{P}(\mathcal{F} \wedge \mathcal{A}(\phi))$ given above Lemma~\ref{Lem path positive exp length cover double boundary}.

\begin{lem}\label{Rmq almost atoroidal}
Let ${\tt n} \geq 3$ and let $\mathcal{F}$ be a free factor system of $F_{\tt n}$. Let $\phi \in \Out(F_{\tt n},\mathcal{F})$ be an almost atoroidal outer automorphism which satisfies Definition~\ref{Defi almost atoroidal outer automorphism}~$(2)$. Let $\mathcal{F} \leq \mathcal{F}_1 \leq \mathcal{F}_2=\{F_{\tt n}\}$ be a sequence of free factor systems given in this definition. Let $f \colon G \to G$ be a CT map representing a power of $\phi$ with filtration $\varnothing=G_0 \subsetneq G_1 \subsetneq \ldots \subsetneq G_k=G$ and such that there exist $p$ and $i$ in $\{0,\ldots,k-1\}$ such that $\mathcal{F}(G_p)=\mathcal{F}$ and $\mathcal{F}(G_i)=\mathcal{F}_1$.

\medskip

\noindent{$(1)$ } The graph $\overline{G-G_i}$ either is a topological arc whose endpoints are in $G_i$ or it retracts onto a circuit $C$ and there exists exactly one topological arc that connects $C$ and $G_i$. 

\medskip

\noindent{$(2)$ } There do not exist an EG stratum or a zero stratum of height greater than $i$. If $\overline{G-G_i}$ is a topological arc, every edge in $\overline{G-G_i}$ is contained in $G_{PG}$. Otherwise every edge of the circuit $C$ in $\overline{G-G_i}$ is contained in $G_{PG}$.

\medskip

\noindent{$(3)$ } Let $\gamma$ be a path of $G_i$ which is not contained in a concatenation of paths of $G_{PG,\mathcal{F}_1}$ and $\mathcal{N}_{PG,\mathcal{F}_1}$. Then $\gamma$ is not contained in a concatenation of paths in $G_{PG}$ and in $\mathcal{N}_{PG}$.

\medskip

\noindent{$(4)$ } We have $$\partial^2(F_{\tt n},\mathcal{F} \wedge \mathcal{A}(\phi))=\bigcup_{\gamma \in \mathcal{P}(\mathcal{F}_1 \wedge \mathcal{A}(\phi))}C(\gamma).$$ In particular, we have $$\PCurr(F_{\tt n},\mathcal{F}\wedge \mathcal{A}(\phi))=\PCurr(F_{\tt n},\mathcal{F}_1\wedge \mathcal{A}(\phi)).$$

\medskip

\noindent{$(5)$ } For every edge path $\gamma$ in $G$, the value $\ell_{\mathcal{F}_1}(\gamma)-\ell_{exp}(\gamma)$ is the number of edges of $\overline{G-G_i}$ contained in $\gamma$. In particular, for every path $\gamma$ contained in $G_i$, we have $$\ell_{\mathcal{F}_1}(\gamma)=\ell_{exp}(\gamma)$$ and for every current $\mu \in \Curr(F_{\tt n},\mathcal{F} \wedge \mathcal{A}(\phi))$ whose support is contained in $\partial^2\mathcal{F}_1$, we have $$\Psi_0(\mu)=\lVert \mu \rVert_{\mathcal{F}_1}.$$

\medskip

\noindent{$(6)$ } Let $\gamma$ be a circuit in $G$. For every $m \in \NN^*$, we have $$\ell_{\mathcal{F}_1}([f^m(\gamma)])-\ell_{exp}([f^m(\gamma)])=\ell_{\mathcal{F}_1}(\gamma)-\ell_{exp}(\gamma).$$

\medskip

\noindent{$(7)$ } Suppose that $\mathcal{F} \wedge \mathcal{A}(\phi)=\{[A_1],\ldots,[A_r]\}$. One of the following holds. 
\begin{itemize}
\item There exist distinct $i,j \in \{1,\ldots,r\}$ such that $$\mathcal{A}(\phi)=(\mathcal{F} \wedge \mathcal{A}(\phi))-\{[A_i],[A_j]\}) \cup \{[A_i \ast A_j]\}.$$ 

\item There exists $i \in \{1,\ldots,r\}$ and an element $g \in F_{\tt n}$ such that $$\mathcal{A}(\phi)=(\mathcal{F} \wedge \mathcal{A}(\phi))-\{[A_i]\}) \cup \{[A_i \ast \left\langle g \right\rangle] \}.$$ In that case, there exists a subgroup $A$ of $F_{\tt n}$ such that $\mathcal{F}=\{[A]\}$ and $F_{\tt n}=A \ast \left\langle g \right\rangle$.

\item There exists $g \in F_{\tt n}$ such that $\mathcal{A}(\phi)=\mathcal{F} \wedge \mathcal{A}(\phi) \cup \{[\left\langle g \right\rangle]\}$. In that case, there exists a subgroup $A$ of $F_{\tt n}$ such that $\mathcal{F}=\{[A]\}$ and $F_{\tt n}=A \ast \left\langle g \right\rangle$.
\end{itemize}

\end{lem}

\dem $(1)$ It is a consequence of \cite[Lemma~II.2.5]{HandelMosher20}. Note that, in the terminology of \cite[Lemma~2.2.5]{HandelMosher20}, the first case is called a \emph{one-edge extension} and the second case is called a \emph{lollipop extension}.

\medskip

\noindent{$(2)$ } By Proposition~\ref{Prop definition CT}~$(4)$, it suffices to show that there does not exist an EG stratum of height greater than $i$. This follows from \cite[Corollary~3.2.2]{BesFeiHan00} (where the stratum described in it is the whole graph $\overline{G-G_i}$) We now prove the second part of Assertion~$(2)$. Let $w$ be an element of $F_{\tt n}$ represented by $\gamma$. Then there exists a subgroup $A$ of $F_{\tt n}$ such that $[A] \in \mathcal{A}(\phi)$ and $w \in A$. Since $\phi|_{\mathcal{F}_1}$ is expanding relative to $\mathcal{F}$ but $\phi$ is not expanding relative to $\mathcal{F}$ by Definition~\ref{Defi almost atoroidal outer automorphism}~$(2)$, there exists a reduced circuit $\gamma$ in $G$ which is not contained in $G_i$ which has polynomial growth under iterates of $f$. By Proposition~\ref{Prop circuits in Gpg are elements in poly subgroup}, the circuit $\gamma$ is a concatenation of paths in $G_{PG}$ and in $\mathcal{N}_{PG}$. By the first part of Assertion~$(2)$, the intersection $\gamma \cap \overline{G-G_i}$ does not contain EG INPs, hence consists in edges in $G_{PG}$. Hence if $\overline{G-G_i}$ is a lollipop, then the circuit $C$ in $\overline{G-G_i}$ is contained in $\gamma$, hence is contained in $G_{PG}$. If $\overline{G-G_i}$ is a topological arc, the graph $\overline{G-G_i}$ is contained in $\gamma$, hence consists in edges in $G_{PG}$. This proves $(2)$. 

\medskip

\noindent{$(3)$ } Let $\gamma$ be as in Assertion~$(3)$. By Assertion~$(2)$, every edge of $\overline{G-G_i}$ is contained in an NEG stratum. In particular, there does not exist an EG INP of height greater than $i$. Hence $\mathcal{N}_{PG}=\mathcal{N}_{PG,\mathcal{F}_1}$. Since $\gamma$ is contained in $G_i$ and since $G_{PG} \cap G_i=G_{PG,\mathcal{F}_1}$, the path $\gamma$ is not contained in a concatenation of paths in $G_{PG}$ and $\mathcal{N}_{PG}$.

\medskip

\noindent{$(4)$ } Since $\phi|_{\mathcal{F}_1}$ is expanding relative to $\mathcal{F}$, we see that $\mathcal{F}_1 \wedge \mathcal{A}(\phi)=\mathcal{F} \wedge \mathcal{A}(\phi)$. Thus, we have $\partial^2(F_{\tt n},\mathcal{F} \wedge \mathcal{A}(\phi))=\partial^2(F_{\tt n},\mathcal{F}_1 \wedge \mathcal{A}(\phi))$. Assertion~$(4)$ then follows from  Lemma~\ref{Lem path positive exp length cover double boundary} applied to $\mathcal{F}_1 \wedge \mathcal{A}(\phi)$.

\medskip

\noindent{$(5)$ } By Assertion~$(2)$, there does not exist an EG INP of height at least $i+1$. Hence $\ell_{\mathcal{F}_1}(\gamma)$ differs from $\ell_{exp}(\gamma)$ by the number of edges in $G_{PG}$ of height at least $i+1$. Since every edge in $\overline{G-G_i}$ is in $G_{PG}$ by Assertion~$(2)$, the conclusion of the first claim of Assertion~$(5)$ follows. The claim about paths contained in $G_i$ is then a direct consequence. Let $\mu$ be a current in $\Curr(\mathcal{F}_1,\mathcal{F}_1 \wedge \mathcal{A}(\phi))$. By Lemma~\ref{Lem homeo curr F1}, there exists $(\mu_1,\mu_2) \in X(\mathcal{F}_1)$ such that $\mu=\mu_1^{\ast}+\mu_2^{\ast}$. Since rational currents are dense in $\Curr(H_1,\{[A_1],\ldots,[A_s]\})$ and $\Curr(H_2,\{[B_1],\ldots,[B_t]\})$ by Proposition~\ref{Prop density rational currents}, linear combination of rational currents are dense in $\Curr(\mathcal{F}_1,\mathcal{F}_1 \wedge \mathcal{A}(\phi))$. The last claim of Assertion~$(5)$ then follows from the linearity and continuity of $\Psi_0$ and $\lVert. \rVert_{\mathcal{F}_1}$.

\medskip

\noindent{$(6)$ } Let $m \in \NN^*$. By Assertion~$(5)$, it suffices to prove that the number of edges in $\overline{G-G_i}$ contained in $[f^m(\gamma)]$ is equal to the number of edges in $\overline{G-G_i}$ contained in $\gamma$. In the case that $\overline{G-G_i}$ is a lollipop extension and that $\gamma$ is the circuit $C$ in $\overline{G-G_i}$, then $\gamma$ is fixed by $f$ by~\cite[Definition~I.1.29~$(3)$]{HandelMosher20} (that is the filtration associated with $f$ is \emph{reduced}). Hence $[f^m(\gamma)]=\gamma$ and the claim follows. Otherwise, if $\overline{G-G_i}$ is either a one-edge extension or a lollipop extension, the circuit $\gamma$ is not contained in $\overline{G-G_i}$. Moreover, if $\gamma$ or $[f^m(\gamma)]$ contains an edge in $\overline{G-G_i}$, then it contains $\overline{G-G_i}$. Hence it suffices to count the number of occurrences of $\overline{G-G_i}$ in $\gamma$ and $[f^m(\gamma)]$. Since $f$ preserves $G_i$, the result follows from Assertion~$(1)$ and \cite[Corollary~3.2.2]{BesFeiHan00} (where the stratum in it is the graph $\overline{G-G_i}$).

\medskip

\noindent{$(7)$ } Note that since $\phi|_{\mathcal{F}_1}$ is expanding relative to $\mathcal{F}$, we have $\mathcal{F}_1 \wedge \mathcal{A}(\phi)=\mathcal{F} \wedge \mathcal{A}(\phi)$. Recall the definition of the graph $G^{\ast}$ and the map $p_{G^{\ast}} \colon G^{\ast} \to G$ from above Lemma~\ref{Lem Injectivity pi1}. By Proposition~\ref{Prop circuits in Gpg are elements in poly subgroup} and Lemma~\ref{Lem Injectivity pi1}~$(2)$, the malnormal subgroup system $\mathcal{A}(\phi)$ is precisely the subgroup system associated with the fundamental groups of the connected components of $G^{\ast}$. Moreover, the malnormal subgroup system associated with $\mathcal{F}_1 \wedge \mathcal{A}(\phi)=\mathcal{F} \wedge \mathcal{A}(\phi)$ is the subgroup system associated with the connected components of $p_{G^{\ast}}^{-1}(G_i)$. By Assertion~$(1)$, the graph $\overline{G-G_i}$ is either a topological arc or a lollipop. Suppose first that $\overline{G-G_i}$ is a topological arc. By Assertion~$(2)$, the graph $\overline{G-G_i}$ consists in edges in $G_{PG}$. Thus, the graph $G^{\ast}$ is obtained from $p_{G^{\ast}}^{-1}(G_i)$ by adding a topological arc $\tau$. If the endpoints of $\tau$ are in two distinct connected components of $G^{\ast}$, then the first case of Assertion~$(7)$ occurs and otherwise the second case of Assertion~$(7)$ occurs. Moreover, if the second case occurs, the extension $\mathcal{F} \leq \{[F_{\tt n}]\}$ is an HNN extension. Thus there exists a subgroup $A$ of $F_{\tt n}$ such that $\mathcal{F}=\{[A]\}$. By~\cite[Corollary~3.2.2]{BesFeiHan00}, one can obtain an element $g$ of $F_{\tt n}$ such that $F_{\tt n}=A \ast g$ by taking a circuit in the image of $p_{G^{\ast}}$ which contains $\overline{G-G_i}$ exactly once. Suppose now that $\overline{G-G_i}$ is a lollipop extension. By Assertion~$(2)$, the circuit $C$ in $\overline{G-G_i}$ consists in edges in $G_{PG}$. Thus, either $G^{\ast}$ is obtained from $p_{G^{\ast}}^{-1}(G_i)$ by adding a lollipop extension, or $G^{\ast}$ is obtained from $p_{G^{\ast}}^{-1}(G_i)$ by adding a connected component which is homotopy equivalent to a circle. If $G^{\ast}$ is obtained from $p_{G^{\ast}}^{-1}(G_i)$ by adding a lollipop extension, the second case of Assertion~$(7)$ occurs. If $G^{\ast}$ is obtained from $p_{G^{\ast}}^{-1}(G_i)$ by adding a connected component which is homotopy equivalent to a circle, the third case of Assertion~$(7)$ occurs. The proof of the fact about HNN extension is similar to the proof for the one-edge extension case. This concludes the proof.
\hfill\qedsymbol

\begin{rmq}\label{Last rmq}
By Lemma~\ref{Rmq almost atoroidal}~$(1)$, $\overline{G-G_i}$ is either a topological arc or it retracts onto a circuit $C$ and there exists exactly one topological arc that connects $C$ and $G_i$.  In the second case, we will adopt the convention that $\overline{G-G_i}=C$, so that, by Lemma~\ref{Rmq almost atoroidal}~$(2)$, in both cases of Lemma~\ref{Rmq almost atoroidal}~$(1)$, every edge in $\overline{G-G_i}$ is in $G_{PG}$.
\end{rmq}

\begin{lem}\label{Lem conversion goodness closeness PCurr}
Let $\phi \in \Out(F_{\tt n},\mathcal{F})$ and let $f \colon G \to G$ be as in Remark~\ref{Rmq Convention for relative atoroidal CT map 2}. 

\medskip

\noindent{$(1)$ } Let $U$ be an open neighborhood of $\Delta_+(\phi)$, let $V$ be a neighborhood of $K_{PG}(\phi)$ (see~Definition~\ref{Defi Kpg}). There exist $N \in \NN^*$ and $\delta \in (0,1)$ such that for every $m \geq 1$ and every $w \in F_{\tt n}$ with $\mathfrak{g}(\gamma_w)> \delta$ and $\eta_{[w]} \notin V$, we have $$(\phi^N)^m(\eta_{[w]}) \in U.$$

\medskip

\noindent{$(2)$ } Suppose that $\phi$ is an almost atoroidal outer automorphism relative to $\mathcal{F}$ as in Definition~\ref{Defi almost atoroidal outer automorphism}~$(2)$. Let $\mathcal{F} \leq \mathcal{F}_1 \leq \mathcal{F}_2$ be an associated sequence of free factor systems. 

For every $\epsilon >0$ and $L >0$, there exists $\delta \in (0,1)$ and $M >0$ such that, for every $n \geq M$, for every reduced edge path $\gamma \in \mathcal{P}(\mathcal{F} \wedge \mathcal{A}(\phi))$ of length at most $L$ contained in $G_{i}$, for every nonperipheral element $w \in F_{\tt n}$ with $\mathfrak{g}(\gamma_w)> \delta$, there exists $[\mu_w] \in \Delta_+(\phi)$ such that: 
$$\left| \frac{\left\langle \gamma, [f^n(\gamma_w)] \right\rangle}{\ell_{exp}([f^n(\gamma_w)])}-\frac{\left\langle \gamma, [\mu_{w}]) \right\rangle}{\lVert[\mu_{w}]\rVert_{\mathcal{F}_1}} \right| < \epsilon.$$

\end{lem}

\dem The proof is similar to the one of \cite[Lemma~6.1]{LustigUyanik2019}. By Lemma~\ref{Lem compute exponential length optimal factor} and Lemma~\ref{Lem control of the goodness}~$(1)$, up to passing to a power of $f$, we may assume that for every $w \in F_{\tt n}$ such that $\mathfrak{g}(\gamma_w) \geq \frac{1}{2}$, and every $n \in \NN^*$, we have $\mathfrak{g}([f^n(\gamma_w)]) \geq \mathfrak{g}(\gamma_w)$ and \begin{equation}\label{Equation TEL}
\ell_{exp}([f^n(\gamma_w)]) \geq TEL(n,\gamma) \geq (3^{n}-2C)\mathfrak{g}(\gamma_w) \ell_{exp}(\gamma_w).
\end{equation}

Let $N \in \NN^*$ be such that $3^N >2C$. Let $\lambda >0$ be such that, for every edge $e \in \vec{E}G$ and every $n \in \NN^*$, we have \begin{equation}\label{Equation first of lem p65}
\ell([f^n(e)]) \leq \lambda^n.
\end{equation} 

By Lemma~\ref{Lem Caracterisation limites PCurr}, a sequence $([\nu_m])_{m \in \NN}$ of projective relative currents tends to a projective current $[\nu] \in \PCurr(F_{\tt n},\mathcal{F} \wedge \mathcal{A}(\phi))$ if for every $\epsilon >0$ and $R > 0$ there exists $M \in \NN^*$ such that, for every $m \geq M$ and every reduced edge path $\gamma \in \mathcal{P}(\mathcal{F} \wedge \mathcal{A}(\phi))$ with $\ell(\gamma) \leq R$, we have 

\begin{equation}\label{Equation 0}
\left|\frac{\left\langle \gamma, \nu \right\rangle}{\left\lVert \nu \right\rVert_{\mathcal{F}}}-\frac{\left\langle \gamma, \nu_m \right\rangle}{\left\lVert \nu_m \right\rVert_{\mathcal{F}}} \right| < \epsilon.
\end{equation}

For every $\mathcal{F}$-expanding splitting unit $\sigma$, we denote by $\mu(\sigma)$ the corresponding current given by Proposition~\ref{Prop existence relative currents atoroidal automorphisms}. By Lemma~\ref{Lem Delta value norm}, we have $\left\lVert\mu(\sigma)\right\rVert_{\mathcal{F}}=1$. Since $\Delta_+(\phi)$ is compact by Lemma~\ref{Lem Delta compact}, there exist $\epsilon,R >0$ such that for every $m \geq M$, if there exists $\nu \in \Delta_+(\phi)$ such that $\nu_m,\nu,R,\epsilon$ satisfy Equation~\eqref{Equation 0}, then $\nu_m \in U$. Since there are only finitely many expanding splitting units of positive exponential length and finitely many edge paths $\gamma \in \mathcal{P}(\mathcal{F} \wedge \mathcal{A}(\phi))$ such that $\ell(\gamma) \leq R$, there exists $M_0 \in \NN^*$ such that for every $m \geq M_0$, for every expanding splitting unit $\sigma$ and for every reduced edge path $\gamma \in \mathcal{P}(\mathcal{F} \wedge \mathcal{A}(\phi))$ with $\ell(\gamma) \leq R$, we have: 
\begin{equation*}
\left|\frac{\left\langle \gamma, [f^m(\sigma)] \right\rangle}{\ell_{\mathcal{F}}([f^m(\sigma)])}-\left\langle \gamma, \mu(\sigma) \right\rangle \right| < \frac{\epsilon}{6}.
\end{equation*}

Recall that $\left\langle \gamma, \mu(\sigma) \right\rangle $ is equal to $\mu(\sigma)(C(\gamma))$ by definition of the number of occurrences of $\gamma$ in $\mu(\sigma)$. Let $\gamma'$ be a reduced edge path in $G$. By Lemma~\ref{Lem bound exponential length subpath}, for every reduced edge path $\sigma$ of $G$ contained in $\gamma'$, we have $\ell_{\mathcal{F}}(\sigma) \geq \ell_{\mathcal{F}}^{\gamma'}(\sigma) \geq \ell_{\mathcal{F}}(\sigma)-2C$. Hence
there exists $M_1 \in \NN^*$ such that for every $m \geq M_1$, for every expanding splitting unit $\sigma$, for every edge path $\gamma'$ containing $\sigma$ as a splitting unit and for every reduced edge path $\gamma \in \mathcal{P}(\mathcal{F} \wedge \mathcal{A}(\phi))$ with $\ell(\gamma) \leq R$, we have: 

\begin{equation}\label{Equation 1}
\left|\frac{\left\langle \gamma, [f^m(\sigma)] \right\rangle}{\ell_{\mathcal{F}}^{[f^m(\gamma')]}([f^m(\sigma)])}-\left\langle \gamma, \mu(\sigma) \right\rangle \right| < \frac{\epsilon}{6}.
\end{equation}

Recall the definition of the continuous function $\Psi_0 \colon \Curr(F_{\tt n},\mathcal{F} \wedge \mathcal{A}(\phi)) \to \RR$ given above Definition~\ref{Defi Kpg}. Recall that, by Lemma~\ref{Lem polynomially growing currents empty}~$(3)$, for every current $\mu \in \Curr(F_{\tt n},\mathcal{F}\wedge\mathcal{A}(\phi))$, we have $\lVert \mu \rVert_{\mathcal{F}} >0$. Let $$\begin{array}{cccl}
\Psi \colon & \Curr(F_{\tt n},\mathcal{F} \wedge \mathcal{A}(\phi)) & \to & \RR \\
{} & [\nu] & \mapsto & \frac{\Psi_0(\nu)}{\lVert \nu \rVert_{\mathcal{F}}}.
\end{array}$$ Since $\Psi$ is continuous and since $\PCurr(F_{\tt n},\mathcal{F} \wedge \mathcal{A}(\phi))-V$ is compact, there exists $s >0$ such that for every $\nu \in \PCurr(F_{\tt n},\mathcal{F} \wedge \mathcal{A}(\phi))-V$, we have:
$$\Psi([\nu]) \geq s.$$ In particular, by Lemma~\ref{Lem psi0 equal exp length}, for every nonperipheral element $w \in F_{\tt n}$ such that $\eta_{[w]} \notin V$, we have 
\begin{equation}\label{Equation p66}
\frac{\ell_{exp}(\gamma_w)}{\ell_{\mathcal{F}}(\gamma_w)}=\frac{\Psi_0(\eta_{[w]})}{\lVert \eta_{[w]} \rVert_{\mathcal{F}}}=\Psi([\eta_{[w]}]) \geq s.
\end{equation}

Now let $w \in F_{\tt n}$ be a nonperipheral element such that $\mathfrak{g}(\gamma_w) \geq \frac{1}{2}$ and $\eta_{[w]} \notin V$. Let $\gamma_w=\alpha_0\beta_1\alpha_1\ldots\alpha_k\beta_k$ be an optimal splitting of $\gamma_w$, where for every $i \in \{0,\ldots,k\}$, the path $\alpha_i$ is an incomplete factor of $\gamma_w$ and for every $i \in \{1,\ldots,k\}$, the path $\beta_i$ is a complete factor of $\gamma$.  Using this optimal splitting, we construct another decomposition of $\gamma_w$ (which is not necessarily a splitting of $\gamma_w$). Since concatenations of paths in $G_{PG}$ and in $\mathcal{N}_{PG}$ have zero exponential length by Lemma~\ref{Lem paths in Npg' have zero exponential length}, we  change the decomposition in such a way that every subpath of $\gamma_w$ which is a concatenation of paths in $G_{PG}$ and in $\mathcal{N}_{PG}$ is in some $\alpha_i$ for $i \in \{1,\ldots,k\}$. In particular, for every $i \in \{1,\ldots,k\}$, the exponential lengths of $\beta_i$ and $\alpha_i$ are equal to their exponential lengths relative to $\gamma_w$. Let $i \in \{0,\ldots,k\}$. The path $\alpha_i$ has a decomposition $\alpha_i=\alpha_i^{(1)}\alpha_i^{(1')}\ldots \alpha_i^{(k_i)}\alpha_i^{(k_i')}$ where, for every $j \in \{1,\ldots,k_i\}$, the path $\alpha_i^{(j)}$ is a concatenation of paths in $G_{PG}$ and $\mathcal{N}_{PG}$ and, for every $j \in \{1,\ldots,k_i\}$, the path $\alpha_i^{(j')}$ is a path in $\overline{G-G_{PG}}$ such that every edge of $\alpha_i^{(j')}$ either has positive exponential length relative to $\gamma_w$ or is in a zero stratum. Note that, by Proposition~\ref{Prop definition CT}~$(4)$, for every $j \in \{1,\ldots,k_i\}$ and every maximal subpath $\tau$ of $\alpha_i^{(j')}$ contained in some zero stratum, the path $\tau$ is adjacent to a path in $\gamma_w$ of positive exponential length. Suppose that $\tau$ is nontrivial. Since no zero path is adjacent to a path which is a concatenation of paths in $G_{PG}$ and $\mathcal{N}_{PG}$ by Lemma~\ref{Lem No zero path and Nielsen path adjacent} and Proposition~\ref{Prop definition CT}~$(4)$, either $\alpha_i=\tau$ or $\ell_{exp}(\alpha_i^{(j')})>0$. In the first case, we have $\ell(\tau) \leq C$ by definition of $C$. Thus, there exists $n \in \NN^*$ such that $[f^n(\tau)]$ is completely split. Thus, if the first case occurs, we may suppose, up to taking a power of $f$, that $\alpha_i$ is a completely split and is a splitting unit of some $\beta_j$. Let $i \in \{1,\ldots,k\}$. Since $\beta_i$ does not contain splitting units which are concatenation of paths in $G_{PG}$ and $\mathcal{N}_{PG}$, every splitting unit of $\beta_i$ is an edge in $\overline{G-G_{PG}'}$ or a maximal taken connecting path in a zero stratum. By Lemma~\ref{Lem exponential length goes to infinity}, every splitting unit of $\beta_i$ which is an edge in $\overline{G-G_{PG}'}$ is expanding. Let $\sigma'$ be a splitting unit of $\beta_i$ which is a maximal taken connecting path in a zero stratum and which is not expanding. Let $n \in \NN^*$ be such that $[f^n(\sigma')]$ is completely split. By Lemma~\ref{Lem exponential length goes to infinity} and Lemma~\ref{Lem compute exp length completely split}, the path $[f^n(\sigma')]$ does not contain splitting units which are edges in $\overline{G-G_{PG}}$. If $[f^n(\sigma')]$ contains a splitting unit which is contained in a zero stratum, then an inductive argument shows that, up to taking a larger $n$, the path $[f^n(\sigma')]$ is a concatenation of paths in $G_{PG}$ and $\mathcal{N}_{PG}$. Thus, the $\mathcal{F}$-length of $\sigma'$ grows at most polynomially fast under iterates of $f$. Thus, we see that $\gamma_w$ has a decomposition $$\gamma_w=a_0b_0a_1c_1^{(1)}c_2^{(1)}\ldots c_{k_1}^{(1)}a_2b_2\ldots a_{t}c_1^{(t)}c_2^{(t)}\ldots c_{k_t}^{(t)}a_{t+1}b_{t+1}a_{t+2},
$$
where:

\medskip

\noindent{$(a)$ } for every $i \in \{0,\ldots,t+2\}$, the path $a_i$ is either possibly trivial, a concatenation of paths in $G_{PG}$ and in $\mathcal{N}_{PG}$ or a maximal taken connecting path whose $\mathcal{F}$-length grows at most polynomially fast;

\medskip

\noindent{$(b)$ } for every $i \in \{0,\dots,t+1\}$, the path $b_i$ is a subpath of positive exponential length relative to $\gamma_w$ of an incomplete path of $\gamma_w$ such that every edge of $b_i$ either has positive exponential length relative to $\gamma_w$ or is in a zero stratum;

\medskip

\noindent{$(c)$ } for every $i \in \{1,\dots,t\}$ and every $j \in \{1,\ldots,k_i\}$, the path $c_j^{(i)}$ is a (possibly trivial) expanding splitting unit of a complete factor of $\gamma_w$. 

\medskip

Recall that the length of every path in a zero stratum is bounded by $C$. Thus, for every $i \in \{0,\ldots,t+1\}$, we have $$\ell(b_i) \leq C\ell_{exp}(b_i).$$ We claim that the exponential length relative to $\gamma_w$ of one of the edges at the concatenation point of two consecutive nontrivial paths of the form $a_ib_i$, $b_ia_{i+1}$, $a_ic_1^{(i)}$, $c_j^{(i)}c_{j+1}^{(i)}$ or $c_{k_i}^{(i)}a_{i+1}$ is positive. Indeed, for every $i \in \{1,\dots,t\}$(resp. $i \in \{0,\ldots,t+1\}$) and every $j \in \{1,\ldots,k_i\}$, the path $c_j^{(i)}$ (resp. $b_i$) either has positive exponential length relative to $\gamma_w$ or is contained in a zero stratum. Note that by hypothesis, for every $i \in \{0,\dots,t+1\}$, the path $b_i$ is not contained in a zero stratum. Moreover, if $b_i$ is adjacent to a path $a_i$, then the first edge of $b_i$ is not in a zero stratum by Proposition~\ref{Prop definition CT}~$(4)$, Lemma~\ref{Lem No zero path and Nielsen path adjacent} and the fact that the paths in zero strata that we consider in our subdivision are maximal. Hence one of the edges at the concatenation point of every path of the form $a_ib_i$, $b_ia_{i+1}$ has positive exponential length relative to $\gamma_w$. By maximality of the splitting units contained in zero strata, one of the splitting unit in a path $c_j^{(i)}c_{j+1}^{(i)}$ is an edge in $\overline{G-G_{PG}'}$, hence has positive exponential length relative to $\gamma_w$. Since paths in zero strata and concatenations of paths in $G_{PG}$ and $\mathcal{N}_{PG}$ cannot be adjacent by Proposition~\ref{Prop definition CT}~$(4)$ and Lemma~\ref{Lem No zero path and Nielsen path adjacent}, paths of the form $a_ic_1^{(i)}$ and $c_{k_i}^{(i)}a_{i+1}$ have positive exponential length since in this case $c_{1}^{(i)}$ or $c_{k_i}^{(i)}$ is an edge in $\overline{G-G_{PG}'}$. This proves the claim.

Remark that, by construction and the definition of goodness of a reduced path, we have $$\sum_{i=1}^t\sum_{j=1}^{k_i} \ell_{exp}(c_j^{(i)})=\ell_{exp}(\gamma_w)\mathfrak{g}(\gamma_w).$$

Note that the length of reduced iterates of edges in $G_{PG}$ grows at most polynomially fast, hence the $\mathcal{F}$-length of reduced iterates of edges in $G_{PG}$ grows at most polynomially fast. Let $C' >0$ and $k \in \NN^*$ be such that, for every splitting unit $\sigma'$ which is either an edge in $G_{PG}$ or a maximal taken connecting path in a zero stratum whose $\mathcal{F}$-length grows at most polynomially fast, and every $m \in \NN^*$, we have: $$\ell_{\mathcal{F}}([f^m(\sigma')]) \leq C'm^k\ell_{\mathcal{F}}(\sigma').$$ The constants $C'$ and $k$ exist by the claim in Proposition~\ref{Prop circuits in Gpg are elements in poly subgroup}. Let $i \in \{0,\ldots,t+2\}$ and let $a_i=\alpha_0\ldots\alpha_{\ell_i}$ be a decomposition of $a_i$ such that, for every $j \in \{0,\ldots,\ell_i\}$, $\alpha_{\ell_i}$ is either an edge in $G_{PG}$, a path in $\mathcal{N}_{PG}^{\max}(a_i)$ or a maximal taken connecting path in a zero stratum whose $\mathcal{F}$-length grows at most polynomially fast. By Lemma~\ref{Lem Exponential length less exp length subpaths}, for every $m \in \NN^*$, we have $$\ell_{\mathcal{F}}([f^m(a_i)]) \leq \sum_{j=0}^{\ell_i} \ell_{\mathcal{F}}([f^m(\alpha_j)]) \leq C'm^k\sum_{j=1}^{\ell_i} \ell_{\mathcal{F}}(\alpha_j) = C'm^k\ell_{\mathcal{F}}(a_i),$$ where the last equality follows from the fact that a path in $\mathcal{N}_{PG}$ is contained in some subpath $\alpha_j$ by hypothesis. In particular, 
\begin{equation}\label{Equation p67}
\sum\limits_{i=0}^{t+2} \ell_{\mathcal{F}}([f^n(a_i)]) \leq C'm^k\sum_{i=0}^{t+2} \ell_{\mathcal{F}}(a_i) \leq C'\ell_{\mathcal{F}}(\gamma_w)n^k,
\end{equation}
where the last inequality follows from the fact that, by hypothesis, every path in $\mathcal{N}_{PG}^{\max}(\gamma)$ is contained in some $a_i$. Thus, if $\mathfrak{g}(\gamma_w) \geq \frac{1}{2}$, there exists $C''>0$ such that, for every $n \geq N$, by Equations~\eqref{Equation TEL},~\eqref{Equation p67} and~\eqref{Equation p66}, we have:

$$
\frac{\sum_{i=0}^{t+2}\ell_{\mathcal{F}}([f^n(a_i)])}{\ell_{exp}([f^n(\gamma_w)])} \leq \frac{C'\ell_{\mathcal{F}}(\gamma_w)n^k}{(3^n-2C)\mathfrak{g}(\gamma_w)\;\ell_{exp}(\gamma_w)} 
\leq  \frac{C'\frac{1}{s}\ell_{exp}(\gamma_w)n^k}{(3^n-2C)\mathfrak{g}(\gamma_w)\;\ell_{exp}(\gamma_w)} 
\leq C''\frac{n^k}{(3^n-2C)\mathfrak{g}(\gamma_w)}.
$$

Recall that, for every reduced edge path $\gamma$ of $G$, we have $$\ell_{exp}(\gamma) \leq \ell_{\mathcal{F}}(\gamma).$$ Up to taking a larger $N \in \NN^*$, we may suppose that, for every $n \geq N$, we have 

\begin{equation}\label{New equation p78}
C''\frac{n^k}{(3^n-2C)\mathfrak{g}(\gamma_w)} \leq \frac{\epsilon}{48\mathfrak{g}(\gamma_w)R}.
\end{equation}

For every $n \geq N$ and every nonperipheral element $w \in F_{\tt n}$ such that $\mathfrak{g}(\gamma_w) \geq \frac{1}{2}$, by Equation~\eqref{Equation TEL}, we have

$$
\frac{2R\;\ell_{exp}(\gamma_w)}{\ell_{\mathcal{F}}([f^n(\gamma_w)])}
\leq \frac{2R\;\ell_{exp}(\gamma_w)}{(3^n-2C)\mathfrak{g}(\gamma_w)\ell_{exp}(\gamma_w)}
=\frac{2R}{(3^n-2C)\mathfrak{g}(\gamma_w)}\;.
$$

Up to taking a larger $N$, we may assume that for every $n \geq N$ and every $w \in F_{\tt n}$ such that $\mathfrak{g}(\gamma_w) \geq \frac{1}{2}$, we have:

\begin{equation}\label{Equation 3bis}
\frac{2R\;\ell_{exp}(\gamma_w)}{\ell_{\mathcal{F}}([f^n(\gamma_w)])} \leq \frac{2R}{(3^n-2C)\mathfrak{g}(\gamma_w)} \leq \frac{\epsilon}{12\mathfrak{g}(\gamma_w)}.
\end{equation}

Let $$\delta =\max\left\{\frac{1}{1+\frac{\epsilon}{6}},\frac{1}{1+\frac{2RC\epsilon\lambda^N}{(3^N-2C)6}},\frac{1}{2}\right\}.$$

Thus, in order to prove the first assertion of Lemma~\ref{Lem conversion goodness closeness PCurr}, it suffices to show that for every $m \geq N$ and every $w \in F_{\tt n}$ such that $\mathfrak{g}(\gamma_w) > \delta$ and $\eta_{[w]} \notin V$, the projective current $[\nu_m]=\phi^m([\eta_w])$ is close to an element $[\nu]$ in $\Delta_+(\phi)$ in the sense of Equation~\eqref{Equation 0}. Since the goodness function is monotone by Remark~\ref{Rmq Convention for relative atoroidal CT map 2}, it suffices to prove it for $m=N$. Let $w \in F_{\tt n}$ such that $\mathfrak{g}(\gamma_w) > \delta$ and $\eta_{[w]} \notin V$.

By Equation~\eqref{New equation p78} and the fact that $\mathfrak{g}(\gamma_w) \geq \delta \geq \frac{1}{2}$, we have

\begin{align}\label{Equation 2}
\frac{\sum_{i=0}^{t+2}\ell_{\mathcal{F}}([f^N(a_i)])}{\ell_{\mathcal{F}}([f^N(\gamma_w)])} &\leq  \frac{\sum_{i=0}^{t+2}\ell_{\mathcal{F}}([f^N(a_i)])}{\ell_{exp}([f^n(\gamma_w)])}\nonumber \\ 
{} & \leq  C''\frac{n^k}{(3^N-2C)\mathfrak{g}(\gamma_w)} \leq C''\frac{n^k}{(3^N-2C)\delta} \leq \frac{\epsilon}{24R}.
\end{align}

Moreover, by Equation~\eqref{Equation 3bis} and the fact that $\mathfrak{g}(\gamma_w) \geq \delta \geq \frac{1}{2}$, we have
\begin{equation}\label{Equation 3}
\frac{2R\;\ell_{exp}(\gamma_w)}{\ell_{\mathcal{F}}([f^N(\gamma_w)])} \leq \frac{\epsilon}{6}.
\end{equation}

Note that, for every $w \in F_{\tt n}$ such that $\mathfrak{g}(\gamma_w) > \delta$ and $\eta_{[w]} \notin V$, we have:

\begin{align}\label{Equation 4}
\frac{2RC\lambda^N(1-\mathfrak{g}(\gamma_w))\ell_{exp}(\gamma_w)}{(3^N-2C)\mathfrak{g}(\gamma_w)\ell_{exp}(\gamma_w)}&= 2RC\frac{\lambda^N}{3^N-2C}\left(\frac{1}{\mathfrak{g}(\gamma_w)}-1 \right)\nonumber \\
{} & \leq  2RC\frac{\lambda^N}{3^N-2C}\left(\frac{1}{\delta}-1 \right) \leq \frac{\epsilon}{6},
\end{align}
where the last inequality follows from the definition of $\delta$. 

Let $\gamma \in \mathcal{P}(\mathcal{F} \wedge \mathcal{A}(\phi))$ be of length at most $R$. By the triangle inequality, we have 

\begin{equation}\label{Equation triangle inequality}
{\Large \begin{array}{l}
\left| \frac{\left\langle \gamma, [f^N(\gamma_w)] \right\rangle}{\ell_{\mathcal{F}}([f^N(\gamma_w)])}-\frac{\left\langle \gamma, \sum_{i=1}^t\sum_{j=1}^{k_i}\ell_{\mathcal{F}}^{[f^N(\gamma_w)]}([f^N(c_j^{(i)})])\mu(c_i^{(j)}) \right\rangle}{\sum_{i=1}^t\sum_{j=1}^{k_i}\ell_{\mathcal{F}}^{[f^N(\gamma_w)]}([f^N(c_j^{(i)})])} \right| \\

\leq \left| \frac{\left\langle \gamma, [f^N(\gamma_w)] \right\rangle}{\ell_{\mathcal{F}}([f^N(\gamma_w)])} - \sum_{i=1}^t\sum_{j=1}^{k_i} \frac{\left\langle \gamma, [f^N(c_j^{(i)})] \right\rangle}{\ell_{\mathcal{F}}([f^N(\gamma_w)])}\right| \\

+ \left| \sum_{i=1}^t\sum_{j=1}^{k_i} \frac{\left\langle \gamma, [f^N(c_j^{(i)})] \right\rangle}{\ell_{\mathcal{F}}([f^N(\gamma_w)])} - \frac{\sum_{i=1}^t\sum_{j=1}^{k_i} \left\langle \gamma, [f^N(c_j^{(i)})] \right\rangle}{\sum_{i=1}^t\sum_{j=1}^{k_i} \ell_{\mathcal{F}}^{[f^N(\gamma_w)]}([f^N(c_j^{(i)})])}\right| \smallskip \\

+\left|\frac{\sum_{i=1}^t\sum_{j=1}^{k_i} \left\langle \gamma, [f^N(c_j^{(i)})] \right\rangle}{\sum_{i=1}^t\sum_{j=1}^{k_i} \ell_{\mathcal{F}}^{[f^N(\gamma_w)]}([f^N(c_j^{(i)})])} - \frac{\left\langle \gamma, \sum_{i=1}^t\sum_{j=1}^{k_i}\ell_{\mathcal{F}}^{[f^N(\gamma_w)]}([f^N(c_j^{(i)})])\mu(c_i^{(j)}) \right\rangle}{\sum_{i=1}^t\sum_{j=1}^{k_i}\ell_{\mathcal{F}}^{[f^N(\gamma_w)]}([f^N(c_j^{(i)})])} \right|.
\end{array} }
\end{equation}

Note that an occurrence of $\gamma$ or $\gamma^{-1}$ in $[f^N(\gamma_w)]$ might happen either in some $[f^N(c_j^{(i)})]$ or in some $[f^N(a_i)]$ or in some $[f^N(b_i)]$ or it might cross over the concatenation points. Recall that one of the edges at the concatenation point of paths of the form $a_ib_i$, $b_ia_{i+1}$, $a_ic_1^{(i)}$, $c_j^{(i)}c_{j+1}^{(i)}$ or $c_{k_i}^{(i)}a_{i+1}$ has positive exponential length relative to $\gamma_w$. Recall also that the length of $\gamma$ is at most equal to $R$. Thus the number of such crossings is at most $2R\ell_{exp}(\gamma_w)$. Thus:

{\small $$
\left| \frac{\left\langle \gamma, [f^N(\gamma_w)] \right\rangle}{\ell_{\mathcal{F}}([f^N(\gamma_w)])} - \sum_{i=1}^t\sum_{j=1}^{k_i} \frac{\left\langle \gamma, [f^N(c_j^{(i)})] \right\rangle}{\ell_{\mathcal{F}}([f^N(\gamma_w)])}\right| 
\leq \frac{2R \ell_{exp}(\gamma_w)}{\ell_{\mathcal{F}}([f^N(\gamma_w)])}+ \sum\limits_{i=0}^{t+2} \frac{\left\langle \gamma, [f^N(a_i)] \right\rangle}{\ell_{\mathcal{F}}([f^N(\gamma_w)])}+ \sum\limits_{i=0}^{t+1} \frac{\left\langle \gamma, [f^N(b_i)] \right\rangle}{\ell_{\mathcal{F}}([f^N(\gamma_w)])}.
$$}

Since $\gamma$ is not contained in a concatenation of paths in $G_{PG,\mathcal{F}}$ and $\mathcal{N}_{PG,\mathcal{F}}$, if $\gamma$ is contained in $[f^N(a_i)]$ for $i \in \{1,\ldots,t+1\}$, then $\gamma$ contains an edge of $[f^N(a_i)]$ of positive $\mathcal{F}$-length relative to $[f^N(a_i)]$. Hence we have $\left\langle \gamma, [f^N(a_i)] \right\rangle \leq \ell_{\mathcal{F}}([f^N(a_i)])$. By Equations~\eqref{Equation 3} and~\eqref{Equation 2} with $n=N$, we have 
$$\frac{2R \ell_{exp}(\gamma_w)}{\ell_{\mathcal{F}}([f^N(\gamma_w)])}+ \sum\limits_{i=0}^{t+2} \frac{\left\langle \gamma, [f^N(a_i)] \right\rangle}{\ell_{\mathcal{F}}([f^N(\gamma_w)])} \leq \frac{2R \ell_{exp}(\gamma_w)}{\ell_{\mathcal{F}}([f^N(\gamma_w)])} + \frac{\sum_{i=0}^{t+1}\ell_{\mathcal{F}}([f^N(a_i)])}{\ell_{\mathcal{F}}([f^N(\gamma_w)])}\leq \frac{\epsilon}{4}.$$

Moreover, since for every $i \in \{0,\ldots,t+1\}$, we have $\ell(b_i) \leq C\ell_{exp}(b_i)$ and by Equations~\eqref{Equation TEL},~\eqref{Equation p66} and~\eqref{Equation 4}, we see that:

$$
{\Large \begin{array}{ccl}
\sum\limits_{i=0}^{t+1} \frac{\left\langle \gamma, [f^N(b_i)] \right\rangle}{\ell_{\mathcal{F}}([f^N(\gamma_w)])} & \leq & \sum\limits_{i=0}^{t+1}\frac{\ell([f^N(b_i)])}{\ell_{\mathcal{F}}([f^N(\gamma_w)])} \leq \sum\limits_{i=0}^{t+1} \frac{C\lambda^N\ell_{exp}(b_i)}{(3^N-2C)\mathfrak{g}(\gamma_w)\ell_{exp}(\gamma_w)} \\ 
{} & \leq & \frac{C\lambda^N(1-\mathfrak{g}(\gamma_w))\ell_{exp}(\gamma_w)}{(3^N-2C)\mathfrak{g}(\gamma_w)\ell_{exp}(\gamma_w)} \leq \frac{\epsilon}{6}.
\end{array}}
$$

For the third term of Inequality~\eqref{Equation triangle inequality}, note that, since $\gamma \in \mathcal{P}(\mathcal{F} \wedge \mathcal{A}(\phi))$, it is not contained in a concatenation of paths in $G_{PG,\mathcal{F}}$ and in $\mathcal{N}_{PG,\mathcal{F}}$. Hence an occurrence of $\gamma$ always appear with an edge $e$ of $c$ such that $\ell_{\mathcal{F}}^c(e)=1$. Since $\ell(\gamma) \leq R$, such an edge $e$ can be crossed by at most $R$ occurrences of $\gamma$ in $c$. Thus, for every reduced edge path $c$ in $G$, we have $\left\langle \gamma,c \right\rangle \leq 2R\ell_{\mathcal{F}}(c)$. Hence we have $$\left|\frac{\sum_{i=1}^t\sum_{j=1}^{k_i} \left\langle \gamma, [f^N(c_j^{(i)})] \right\rangle}{\sum_{i=1}^t\sum_{j=1}^{k_i} \ell_{\mathcal{F}}^{[f^N(\gamma_w)]}([f(c_j^{(i)})])} \right| \leq 2R.$$ Since $$\ell_{\mathcal{F}}([f^N(\gamma_w)]) =\sum_{i=1}^t\sum_{j=1}^{k_i} \ell_{\mathcal{F}}^{[f^N(\gamma_w)]}([f^N(c_j^{(i)})]) + \sum_{i=0}^{t+1}\ell_{\mathcal{F}}^{[f^N(\gamma_w)]}([f^N(a_ib_ia_{i+1})]),$$ using Lemma~\ref{Lem compute exponential length optimal factor} and Lemma~\ref{Lem bound exponential length subpath} for the last inequality we have:

$$
\begin{array}{l}
\left| \sum_{i=1}^t\sum_{j=1}^{k_i} \frac{\left\langle \gamma, [f^N(c_j^{(i)})] \right\rangle}{\ell_{\mathcal{F}}([f^N(\gamma_w)])} - \frac{\sum_{i=1}^t\sum_{j=1}^{k_i} \left\langle \gamma, [f^N(c_j^{(i)})] \right\rangle}{\sum_{i=1}^t\sum_{j=1}^{k_i} \ell_{\mathcal{F}}^{[f^N(\gamma_w)]}([f^N(c_j^{(i)})])}\right| \\

= \left|\frac{\left(\sum_{i=1}^t\sum_{j=1}^{k_i} \left\langle \gamma, [f^N(c_j^{(i)})] \right\rangle \right)\left(\sum_{i=0}^{t+1}\ell_{\mathcal{F}}^{[f^N(\gamma_w)]}([f(a_ib_ia_{i+1})]) \right)}{\left( \sum_{i=1}^t\sum_{j=1}^{k_i} \ell_{\mathcal{F}}^{[f^N(\gamma_w)]}([f(c_j^{(i)})]) \right)\left(\sum_{i=1}^t\sum_{j=1}^{k_i} \ell_{\mathcal{F}}^{[f^N(\gamma_w)]}([f^N(c_j^{(i)})]) + \sum_{i=0}^{t+1}\ell_{\mathcal{F}}^{[f^N(\gamma_w)]}([f^N(a_ib_ia_{i+1})])\right)} \right| \\

\leq \left|\frac{\left(\sum_{i=1}^t\sum_{j=1}^{k_i} \left\langle \gamma, [f^N(c_j^{(i)})] \right\rangle \right)\left(\sum_{i=0}^{t+1}\ell_{\mathcal{F}}^{[f^N(\gamma_w)]}([f(a_ib_ia_{i+1})]) \right)}{\left( \sum_{i=1}^t\sum_{j=1}^{k_i} \ell_{\mathcal{F}}^{[f^N(\gamma_w)]}([f(c_j^{(i)})]) \right)\left(\sum_{i=1}^t\sum_{j=1}^{k_i} \ell_{\mathcal{F}}^{[f^N(\gamma_w)]}([f^N(c_j^{(i)})])\right)} \right| \\

\leq \left|\frac{\left(\sum_{i=1}^t\sum_{j=1}^{k_i} \left\langle \gamma, [f^N(c_j^{(i)})] \right\rangle \right)\left(\sum_{i=0}^{t+1}\ell_{\mathcal{F}}([f^N(b_i)]) + 2\sum\limits_{i=0}^{t+2}\ell_{\mathcal{F}}([f^N(a_i)]) \right)}{\left( \sum_{i=1}^t\sum_{j=1}^{k_i} \ell_{\mathcal{F}}^{[f^N(\gamma_w)]}([f(c_j^{(i)})]) \right)\left(\sum_{i=1}^t\sum_{j=1}^{k_i} \ell_{\mathcal{F}}^{[f^N(\gamma_w)]}([f^N(c_j^{(i)})])\right)} \right| \\

\leq 2R\left|\frac{\sum_{i=0}^{t+1}\ell_{\mathcal{F}}([f^N(b_i)]) + 2\sum\limits_{i=0}^{t+2}\ell_{\mathcal{F}}([f^N(a_i)])}{\sum_{i=1}^t\sum_{j=1}^{k_i} \ell_{\mathcal{F}}^{[f^N(\gamma_w)]}([f^N(c_j^{(i)})])} \right|. 
\end{array}
$$

Recall that we have $$\sum_{i=1}^t\sum_{j=1}^{k_i} \ell_{exp}(c_j^{(i)})=\ell_{exp}(\gamma_w)\mathfrak{g}(\gamma_w)$$ and, for every $i \in \{1,\ldots,t\}$ and every $j \in \{1,\ldots,k_i\}$, we have either $\ell_{exp}(c_j^{(i)})=1$ or $\ell_{exp}(c_j^{(i)})=0$. Hence, we have: $${\begin{array}{ccl} \sum_{i=1}^t\sum_{j=1}^{k_i} \ell_{\mathcal{F}}^{[f^N(\gamma_w)]}([f^N(c_j^{(i)})]) & \geq & \sum_{i=1}^t\sum_{j=1}^{k_i} (\ell_{\mathcal{F}}([f^N(c_j^{(i)})])-2C) \geq \sum_{i=1}^t\sum_{j=1}^{k_i}(3^N-2C) \\
{} &\geq & (3^N-2C)\mathfrak{g}(\gamma_w)\ell_{exp}(\gamma_w), 
\end{array}}$$ where the first inequality follows from Lemma~\ref{Lem bound exponential length subpath} and the second inequality follows from the fact that $f$ is $3K$-expanding and $K \geq 1$. Thus, we have

$${\large
\begin{array}{l}
2R\left|\frac{\sum_{i=0}^{t+1}\ell_{\mathcal{F}}([f^N(b_i)]) + 2\sum\limits_{i=0}^{t+2}\ell_{\mathcal{F}}([f^N(a_i)])}{\sum_{i=1}^t\sum_{j=1}^{k_i} \ell_{\mathcal{F}}^{[f^N(\gamma_w)]}([f^N(c_j^{(i)})])} \right| \\

\leq 2R\left|\frac{\sum_{i=0}^{t+1}\ell_{\mathcal{F}}([f^N(b_i)]) +2 \sum_{i=0}^{t+2}\ell_{\mathcal{F}}([f^N(a_i)])}{\sum_{i=1}^t\sum_{j=1}^{k_i} \ell_{\mathcal{F}}^{[f^N(\gamma_w)]]}([f^N(c_j^{(i)})])} \right| \\

\leq 2R\left| \frac{\sum_{i=0}^{t+1}\ell_{\mathcal{F}}([f^N(b_i)])}{(3^N-2C)\mathfrak{g}(\gamma_w)\ell_{exp}(\gamma_w)} \right| + 2R\left| \frac{2\sum_{i=0}^{t+2}\ell_{\mathcal{F}}([f^N(a_i)])}{(3^N-2C)\delta\ell_{exp}(\gamma_w)}\right|.
\end{array}}
$$

By Equation~\eqref{Equation first of lem p65}, we have $$\sum_{i=0}^{t+1} \ell_{\mathcal{F}}([f^N(b_i)]) \leq \sum_{i=0}^{t+1} \ell([f^N(b_i)]) \leq \lambda^N \sum_{i=0}^{t+1} \ell(b_i)\leq C\lambda^N \sum_{i=0}^{t+1} \ell_{exp}(b_i) \leq C\lambda^N\ell_{exp}(\gamma_w)(1-\mathfrak{g}(\gamma_w)).$$ Hence we have:
 
$$
\begin{array}{l}
2R\left| \frac{\sum_{i=0}^{t+1}\ell_{\mathcal{F}}([f^N(b_i)])}{(3^N-2C)\mathfrak{g}(\gamma_w)\ell_{exp}(\gamma_w)} \right| + 2R\left|  \frac{2\sum_{i=0}^{t+2}\ell_{\mathcal{F}}([f^N(a_i)])}{(3^n-2C)\delta\ell_{exp}(\gamma_w)}\right| \\

\leq 2R\left| \frac{C \lambda^N(1-\mathfrak{g}(\gamma_w))\ell_{exp}(\gamma_w)}{(3^N-2C)\mathfrak{g}(\gamma_w)\ell_{exp}(\gamma_w)} \right| + 2R\left|  \frac{2C'\ell_{\mathcal{F}}(\gamma_w)n^k}{(3^N-2C)\delta\ell_{exp}(\gamma_w)}\right| \text{ by Equation~\eqref{Equation p67} } \\

\leq 2R\left| \frac{C \lambda^N(1-\mathfrak{g}(\gamma_w))\ell_{exp}(\gamma_w)}{(3^N-2C)\mathfrak{g}(\gamma_w)\ell_{exp}(\gamma_w)} \right| + 2R\left|  \frac{2C''n^k}{(3^N-2C)\delta}\right| \\

\leq \frac{2\epsilon}{6} \text{ by Equation~\eqref{Equation 2} and \eqref{Equation 4}.}
\end{array}
$$

Finally, using Equation~\eqref{Equation 1} and the fact that for every $i \in \{1,\ldots,t\}$ and every $j \in \{1,\ldots,k_i\}$, the splitting unit $c_j^{(i)}$ is expanding, we have:

$$
\begin{array}{l}
\left|\frac{\sum_{i=1}^t\sum_{j=1}^{k_i} \left\langle \gamma, [f^N(c_j^{(i)})] \right\rangle}{\sum_{i=1}^t\sum_{j=1}^{k_i} \ell_{\mathcal{F}}^{[f^N(\gamma_w)]}([f^N(c_j^{(i)})])} - \frac{\left\langle \gamma, \sum_{i=1}^t\sum_{j=1}^{k_i}\ell_{\mathcal{F}}^{[f^N(\gamma_w)]}([f^N(c_j^{(i)})])\mu(c_i^{(j)}) \right\rangle}{\sum_{i=1}^t\sum_{j=1}^{k_i}\ell_{\mathcal{F}}^{[f^N(\gamma_w)]}([f^N(c_j^{(i)})])} \right| \\

= \left|\frac{\sum_{i=1}^t\sum_{j=1}^{k_i} \ell_{\mathcal{F}}^{[f^N(\gamma_w)]}([f^N(c_j^{(i)})])\left(\frac{\left\langle \gamma, [f^N(c_j^{(i)})] \right\rangle}{\ell_{\mathcal{F}}^{[f^N(\gamma_w)]}([f^N(c_j^{(i)})])}-\left\langle \gamma, \mu(c_i^{(j)}) \right\rangle\right)}{\sum_{i=1}^t\sum_{j=1}^{k_i} \ell_{\mathcal{F}}^{[f^N(\gamma_w)]}([f^N(c_j^{(i)})])} \right| \\

\leq \frac{\frac{\epsilon}{6}\sum_{i=1}^t\sum_{j=1}^{k_i} \ell_{\mathcal{F}}^{[f^N(\gamma_w)]}([f^N(c_j^{(i)})])}{\sum_{i=1}^t\sum_{j=1}^{k_i} \ell_{\mathcal{F}}^{[f^N(\gamma_w)]}([f^N(c_j^{(i)})])}=\frac{\epsilon}{6}.
\end{array}
$$

Combining all inequalities, we have
$$
\left| \frac{\left\langle \gamma, [f^N(\gamma_w)] \right\rangle}{\ell_{\mathcal{F}}([f^N(\gamma_w)])}-\frac{\left\langle \gamma, \sum_{i=1}^t\sum_{j=1}^{k_i}\ell_{\mathcal{F}}^{[f^N(\gamma_w)]}([f^N(c_j^{(i)})])\mu(c_i^{(j)}) \right\rangle}{\sum_{i=1}^t\sum_{j=1}^{k_i}\ell_{\mathcal{F}}^{[f^N(\gamma_w)]}([f^N(c_j^{(i)})])} \right| \leq \frac{\epsilon}{4}+\frac{\epsilon}{6}+\frac{2\epsilon}{6}+\frac{\epsilon}{6} \leq \epsilon.
$$
This concludes the proof of Assertion~$(1)$ of Lemma~\ref{Lem conversion goodness closeness PCurr} since for every $i \in \{1,\ldots,t\}$ and every $j \in \{1,\ldots,k_i\}$, we have $\mu(c_j^{(i)}) \in \Delta_+(\phi)$.

The proof of Assertion~$(2)$ is the same one as the proof of Assertion~$(1)$, replacing $\ell_{\mathcal{F}}$ and $\ell_{\mathcal{F}}^{\gamma}$ by $\ell_{exp}$ and $\ell_{exp}^{\gamma}$, adding the following arguments. Let $\gamma$ and $w \in F_{\tt n}$ be as in Assertion~$(2)$. Then $\gamma$ is not contained in a contenation of paths in $G_{PG}$ and in $\mathcal{N}_{PG}$ by Lemma~\ref{Rmq almost atoroidal}~$(3)$. If $$\gamma_w=a_0b_0a_1c_1^{(1)}c_2^{(1)}\ldots c_{k_1}^{(1)}a_2b_2\ldots a_{t}c_1^{(t)}c_2^{(t)}\ldots c_{k_t}^{(t)}a_{t+1}b_{t+1}a_{t+2},
$$ is the same decomposition of $\gamma_w$ as in the proof of Assertion~$(1)$, then for every $m \in \NN$ and every $i \in \{1,\ldots,t+2\}$, the path $\gamma$ is not contained in $[f^m(a_i)]$ by Lemma~\ref{Lem iterate of a path in GPG}. Similarly, for every $m \in \NN^*$ and every $i \in \{1,\ldots,t+2\}$, we have $\ell_{exp}([f^m(a_i)])=0$. Hence we do not need Equation~\eqref{Equation 2}. By Lemma~\ref{Rmq almost atoroidal}~$(5)$, we have $$\ell_{exp}(\gamma) = \ell_{\mathcal{F}_1}(\gamma).$$ Moreover, by Lemma~\ref{Rmq almost atoroidal}~$(5)$, for every current $[\mu] \in \Delta_+(\phi)$, we have $\Psi_0(\mu)=\lVert \mu \rVert_{\mathcal{F}_1}$. Replacing $\ell_{\mathcal{F}}$ and $\ell_{\mathcal{F}}^{\gamma}$ by $\ell_{exp}$ and $\ell_{exp}^{\gamma}$ in the equations in the proof of Assertion~$(1)$ concludes the proof.
\hfill\qedsymbol

\bigskip

For the next lemma, we need to compute the exponential length of incomplete segments in a circuit $\gamma$ in $G$. Let $\ell_{exp}(\mathrm{Inc}(\gamma))$ be the sum of the exponential lengths of the incomplete segments of an optimal splitting of $\gamma$. Let $\ell_{exp}^{\gamma}(\mathrm{Inc}(\gamma))$ be the sum of the exponential lengths relative to $\gamma$ of the incomplete segments of an optimal splitting of $\gamma$. Note that $\ell_{exp}^{\gamma}(\mathrm{Inc}(\gamma))$ do not depend on the choice of an optimal splitting. Note that $$\ell_{exp}^{\gamma}(\mathrm{Inc}(\gamma))=(1-\mathfrak{g}(\gamma))\ell_{exp}(\gamma) \leq \ell_{exp}(\gamma).$$

\begin{lem}\label{Lem control goodness versus decrease bad length}
Let $\phi \in \Out(F_{\tt n},\mathcal{F})$ and let $f \colon G \to G$ be as in Remark~\ref{Rmq Convention for relative atoroidal CT map 2}. Let $\delta \in (0,1)$, and let $R > 1$. There exists $n_0 \in \NN^*$ such that for every $n \geq n_0$ and every nonperipheral element $w \in F_{\tt n}$ such that $\eta_{[w]} \notin K_{PG}(\phi)$, we either have $$\mathfrak{g}([f^n(\gamma_w)]) \geq \delta$$ or 
$$\ell_{exp}^{[f^n(\gamma_w)]}(\mathrm{Inc}([f^n(\gamma_w)])) \leq \frac{10C}{R}\ell_{exp}^{\gamma_w}(\mathrm{Inc}(\gamma_w)) \textit{ and } \ell_{exp}([f^n(\gamma_w)]) \leq \frac{10C}{(1-\delta)R}\ell_{exp}(\gamma_w). $$
\end{lem}

\dem Let $w \in F_{\tt n}$ be a nonperipheral element such that $\eta_{[w]} \notin K_{PG}(\phi)$. Suppose that $n \in \NN^*$ is such that $\mathfrak{g}([f^n(\gamma_w)])<\delta$. Assuming for now that we have proved that $$\ell_{exp}^{[f^n(\gamma_w)]}(\mathrm{Inc}([f^n(\gamma_w)])) \leq \frac{10C}{R}\ell_{exp}^{\gamma_w}(\mathrm{Inc}(\gamma_w)),$$ we deduce that $\ell_{exp}([f^n(\gamma_w)]) \leq \frac{10C}{(1-\delta)R}\ell_{exp}(\gamma_w)$. Indeed, we have $$\ell_{exp}^{[f^n(\gamma)]}(\mathrm{Inc}([f^n(\gamma)]))=(1-\mathfrak{g}([f^n(\gamma)]))\ell_{exp}([f^n(\gamma)]) \geq (1-\delta)\ell_{exp}([f^n(\gamma)]).$$ Thus we have 
$${\large \begin{array}{ccl}
\ell_{exp}([f^n(\gamma_w)]) & \leq & \frac{1}{1-\delta}\ell_{exp}^{[f^n(\gamma_w)]}(\mathrm{Inc}([f^n(\gamma_w)])) \leq \frac{10C}{(1-\delta)R}\ell_{exp}^{\gamma_w}(\mathrm{Inc}(\gamma_w)) \\
{} & \leq & \frac{10C}{(1-\delta)R}\ell_{exp}(\gamma_w).
\end{array}}$$ 

Therefore, it suffices to prove that there exists $n_0 \in \NN^*$ such that for every $n \geq n_0$, if $\mathfrak{g}([f^n(\gamma_w)])<\delta$, then $$\ell_{exp}^{[f^n(\gamma_w)]}(\mathrm{Inc}([f^n(\gamma_w)])) \leq \frac{10C}{R}\ell_{exp}^{\gamma_w}(\mathrm{Inc}(\gamma_w)).$$ Consider an optimal splitting $\gamma_w=\alpha_0'\beta_1'\alpha_1'\ldots\alpha_m'\beta_m'$, where for every $i \in \{0,\ldots,m\}$, the path $\alpha_i'$ is an incomplete factor of $\gamma_w$ and for every $i \in \{0,\ldots,m\}$, the path $\beta_i'$ is a $PG$-relative complete factor of $\gamma_w$. We can modify the splitting of $\gamma_w$ in a new splitting $\gamma_w=\alpha_0\beta_1\alpha_1\ldots\beta_k\alpha_k$ where:

\medskip

\noindent{$(i)$ } for every $i \in \{0,\ldots,k\}$, the path $\alpha_i$ is a concatenation of incomplete factors and complete factors of zero exponential length relative to $\gamma_w$ of the old splitting;

\medskip

\noindent{$(ii)$ } for every $i \in \{1,\ldots,k\}$, the path $\beta_i$ is a complete factor of positive exponential length relative to $\gamma_w$ of the old splitting. 

\medskip

In the remainder of the proof, we still refer to the paths $\alpha_i$ as incomplete factors. By the last claim of Remark~\ref{Rmq Convention for relative atoroidal CT map 2}, we may suppose that $\mathfrak{g}(\gamma_w)<\delta$, that is: \begin{equation}\label{Equation bad}
\ell_{exp}^{\gamma_w}(\mathrm{Inc}(\gamma_w))=\sum_{i=0}^k\ell_{exp}^{\gamma_w}(\alpha_i) \geq (1-\delta)\ell_{exp}(\gamma_w).
\end{equation} 

\medskip

\noindent{\bf Claim. } For every $i \in \{0,\ldots,k\}$ and every $m \in \NN^*$, we have  $$\ell_{exp}^{[f^m(\gamma_w)]}(\mathrm{Inc}([f^m(\alpha_i)])) \leq 24C^2\; \ell_{exp}^{\gamma_w}(\alpha_i).$$ Similarly, for every $m \in \NN^*$, we have $$\ell_{exp}^{[f^m(\gamma_w)]}(\mathrm{Inc}([f^m(\gamma_w)])) \leq 24C^2\; \ell_{exp}(\gamma_w).$$

\medskip

\dem Since a reduced iterate of a complete factor is complete, every incomplete factor of $[f^m(\gamma_w)]$ is contained in a reduced iterate of some $\alpha_i$. Thus, we have $$\ell_{exp}^{[f^m(\gamma_w)]}(\mathrm{Inc}([f^m(\gamma_w)])) \leq \sum_{i=0}^k\ell_{exp}^{[f^m(\gamma_w)]}(\mathrm{Inc}([f^m(\alpha_i)])). $$ Hence it suffices to prove the result for the paths $\alpha_i$ with $i \in \{0,\ldots,k\}$. By Property~$(ii)$ for every $i \in \{1,\ldots,k\}$, the path $\beta_i$ has positive exponential length relative to $\gamma_w$. Therefore, if there exists $\gamma' \in \mathcal{N}_{PG}^{\max}(\gamma_w)$ such that $\alpha_i$ intersects $\gamma'$ nontrivially, then $\gamma'$ is contained in $\beta_{i}\alpha_i\beta_{i+1}$. In particular, Lemma~\ref{Lem zero relative exp length} applies and for every $i \in \{0,\ldots,k\}$, if $\ell_{exp}^{\gamma_w}(\alpha_i)=0$, then $\ell_{exp}(\alpha_i)=0$. 

Let $i \in \{0,\ldots,k\}$. Suppose first that $\ell_{exp}^{\gamma_w}(\alpha_i)=0$. By the above, we have $\ell_{exp}(\alpha_i)=0$. By Lemma~\ref{Lem uniform bound on bad subpaths}, there exists $N \in \NN^*$ such that for every $m \geq N$, such that the total exponential length of incomplete factors in any optimal splitting of $[f^m(\alpha_i)]$ is equal to $0$. Hence for every $m \geq N$, the path $[f^m(\alpha_i)]$ is $PG$-relative completely split. Up to taking a power of $f$, we may assume that $N=1$. So this concludes the proof of the claim in the case when $\ell_{exp}^{\gamma_w}(\alpha_i)=0$. 

So we may assume that $\ell_{exp}^{\gamma_w}(\alpha_i)>0$.  By Lemma~\ref{Lem uniform bound on bad subpaths}, for every $m \in \NN^*$, the total exponential length of incomplete factors in $[f^m(\alpha_i)]$ is at most equal to $8C \ell_{exp}(\alpha_i)$. By Lemma~\ref{Lem bound exponential length subpath}, for every $i \in \{1,\ldots,k\}$, we have $$\ell_{exp}(\alpha_i) \leq \ell_{exp}^{\gamma_w}(\alpha_i)+2C \leq 3C\ell_{exp}^{\gamma_w}(\alpha_i).$$ Hence again by Lemma~\ref{Lem bound exponential length subpath}, we have $$\ell_{exp}^{[f^m(\gamma_w)]}(\mathrm{Inc}([f^m(\alpha_i)])) \leq \ell_{exp}(\mathrm{Inc}([f^m(\alpha_i)])) \leq 24C^2\ell_{exp}^{\gamma_w}(\alpha_i).$$ This proves the claim.
\hfill\qedsymbol

\bigskip

Let $\Lambda_{\gamma_w}$ be the set consisting in all incomplete factors $\alpha_i$ of $\gamma_w$ whose exponential length relative to $\gamma_w$ is at least equal to $(3.10^8)R^6C^{12}+1$. Let $\Lambda_{\gamma_w}'$ be the set consisting in all incomplete factors $\alpha_i$ of $\gamma_w$ which are not in $\Lambda_{\gamma_w}$. Let $\ell_{exp}^{\gamma_w}(\Lambda_{\gamma_w})$ (resp. $\ell_{exp}^{\gamma_w}(\Lambda_{\gamma_w}')$) be the sum of the exponential lengths relative to $\gamma_w$ of all incomplete factors of $\gamma$ that belongs to $\Lambda_{\gamma_w}$ (resp. $\Lambda_{\gamma_w}'$). We distinguish between two cases, according to the proportion of $\ell_{exp}^{\gamma_w}(\Lambda_{\gamma_w})$ in the exponential length relative to $\gamma_w$ of incomplete factors in $\gamma_w$.

\medskip

\noindent{\bf Case 1 } Suppose that 
$$
\frac{\ell_{exp}^{\gamma_w}(\Lambda_{\gamma_w})}{\ell_{exp}^{\gamma_w}(\mathrm{Inc}(\gamma_w))} < \frac{1}{(24C^2R)^2}.
$$

This implies that \begin{equation}\label{Equation p72}
\frac{\ell_{exp}^{\gamma_w}(\Lambda_{\gamma_w}')}{\ell_{exp}^{\gamma_w}(\mathrm{Inc}(\gamma_w))} \geq \frac{(24C^2R)^2-1}{(24C^2R)^2}.
\end{equation}

Note that, by Lemma~\ref{Lem bound exponential length subpath}, every path in $\Lambda_{\gamma_w}'$ has exponential length at most equal to $(3.10^8)C^{12}R^6+1+2C$. By Lemma~\ref{Lem complete splitting rel poly grow}, there exists $n_0 \in \NN^*$ such that, for every edge path $\beta$ of exponential length at most equal to $(3.10^8)R^6C^{12}+1+2C$ and every $n \geq n_0$ either $[f^n(\beta)]$ is a concatenation of paths in $G_{PG}$ and in $\mathcal{N}_{PG}$ or $[f^{n_0}(\beta)]$ contains a complete factor of exponential length at least equal to $10C$. By Lemma~\ref{Lem bound exponential length subpath}, in the second case, the path $[f^{n_0}(\beta)]$ has a complete factor of positive exponential length relative to $[f^{n_0}(\beta)]$. Let $\Gamma_{\gamma_w}$ be the set consisting in all incomplete paths $\alpha_i$ of $\gamma_w$ such that $\alpha_i \in \Lambda_{\gamma_w}'$ and $[f^{n_0}(\alpha_i)]$ is a concatenation of paths in $G_{PG}$ and in $\mathcal{N}_{PG}$. Let $\Gamma_{\gamma_w}'$ be the set consisting in all incomplete paths $\alpha_i$ of $\gamma_w$ such that $\alpha_i \in \Lambda_{\gamma_w}'$ and $[f^{n_0}(\alpha_i)]$ has at least one complete factor of positive exponential length relative to $[f^{n_0}(\alpha_i)]$. Note that $\Lambda_{\gamma_w}'=\Gamma_{\gamma_w} \cup \Gamma_{\gamma_w}'$. Let $\ell_{exp}^{\gamma_w}(\Gamma_{\gamma_w})$ (resp.  $\ell_{exp}^{\gamma_w}(\Gamma_{\gamma_w}')$) be the sum of the exponential lengths relative to $\gamma_w$ of paths in $\Gamma_{\gamma_w}$ (resp. $\Gamma_{\gamma_w}'$). 

\medskip

\noindent{\bf Subcase~1 } Suppose that $$\frac{\ell_{exp}^{\gamma_w}(\Gamma_{\gamma_w})}{\ell_{exp}^{\gamma_w}(\Lambda_{\gamma_w}')} \geq \frac{24C^2R}{24C^2R+1}.$$ Then $$\ell_{exp}^{\gamma_w}(\Gamma_{\gamma_w}) \geq \frac{24C^2R}{24C^2R+1}\ell_{exp}^{\gamma_w}(\Lambda_{\gamma_w}') \geq \frac{24C^2R-1}{24C^2R}\ell_{exp}^{\gamma_w}(\mathrm{Inc}(\gamma_w)).$$

Note that, for every $n \geq n_0$ and every path $\alpha_i \in \Gamma_{\gamma_w}$, we have $\ell_{exp}([f^n(\alpha_i)])=0$ by Lemma~\ref{Lem exponential length paths in Gpg}. By the claim, for every path $\alpha_i$ such that $\alpha_i \in \Lambda_{\gamma_w}'$ and $\alpha_i \notin \Gamma_{\gamma_w}$, and for every $n \in \NN^*$, the total exponential length of incomplete factors in $[f^n(\alpha_i)]$ relative to $[f^n(\alpha_i)]$ is at most equal to $24C^2 \ell_{exp}^{\gamma_w}(\alpha_i)$. Thus, for every $n \geq n_0$, we have:

$$\begin{array}{ccl}
\ell_{exp}^{[f^n(\gamma_w)]}(\mathrm{Inc}([f^n(\gamma_w)])) & \leq & \sum\limits_{\alpha_i \in \Lambda_{\gamma_w} \cup \Lambda_{\gamma_w}'} \ell_{exp}^{[f^n(\gamma_w)]}(\mathrm{Inc}([f^n(\alpha_i)])) \\
{} &\leq & \sum\limits_{\alpha_i \in \Lambda_{\gamma_w} \cup (\Lambda_{\gamma_w}'-\Gamma_{\gamma_w})} 24C^2\ell_{exp}^{\gamma_w}(\alpha_i) \\
{} & \leq & 24C^2\ell_{exp}^{\gamma_w}(\mathrm{Inc}(\gamma_w))-24C^2\frac{24C^2R-1}{24C^2R}\ell_{exp}^{\gamma_w}(\mathrm{Inc}(\gamma_w)) \\
{} & \leq & \frac{1}{R}\ell_{exp}^{\gamma_w}(\mathrm{Inc}(\gamma_w)).
\end{array}$$ This concludes the proof of Lemma~\ref{Lem control goodness versus decrease bad length} when Subcase~1 occurs.

\medskip

\noindent{\bf Subcase~2 } Suppose that $$\frac{\ell_{exp}(\Gamma_{\gamma_w})}{\ell_{exp}(\Lambda_{\gamma_w}')} < \frac{24C^2R}{24C^2R+1}.$$

Note that the assumption of Subcase~$2$ and Equation~\eqref{Equation p72} imply that $$\ell_{exp}^{\gamma_w}(\Gamma_{\gamma_w}') \geq \frac{1}{24C^2R+1}\ell_{exp}^{\gamma_w}(\Lambda_{\gamma_w}') \geq \frac{(24C^2R)^2-1}{(24C^2R)^2}\frac{1}{24C^2R+1}\ell_{exp}^{\gamma_w}(\mathrm{Inc}(\gamma_w)).$$

Since every path in $\Gamma_{\gamma_w}'$ has exponential length at most equal to $(3.10^8)R^6C^{12}+1+2C$, by Lemma~\ref{Lem p length of a completely split path grows linearly}, up to taking a larger $n_0$, for every path $\alpha_i \in \Gamma_{\gamma_w}'$ such that $\ell_{exp}(\alpha_i)>0$ and every $n \geq n_0$, the exponential length of a complete factor in $[f^n(\alpha_i)]$ is at least equal to $3^{n-n_0}\ell_{exp}(\alpha_i)$. Moreover, for every path $\alpha_i \in \Gamma_{\gamma_w}'$ such that $\ell_{exp}(\alpha_i)=0$ and every $n \geq n_0$, the exponential length of a complete factor in $[f^n(\alpha_i)]$ is at least equal to $3^{n-n_0}$. By Lemma~\ref{Lem bound exponential length subpath}, for every $n \geq n_0$ and every path $\alpha_i \in \Gamma_{\gamma_w}'$ such that $\ell_{exp}(\alpha_i)>0$, the exponential length relative to $[f^n(\alpha_i)]$ of a complete factor in $[f^n(\alpha_i)]$ is at least equal to $$3^{n-n_0}\ell_{exp}(\alpha_i)-2C \geq (3^{n-n_0}-2C)\ell_{exp}(\alpha_i).$$ Thus, for every $n \geq n_0$ and every path $\alpha_i \in \Gamma_{\gamma_w}'$, the exponential length relative to $[f^n(\alpha_i)]$ of a complete factor in $[f^n(\alpha_i)]$ is at least equal to $$(3^{n-n_0}-2C)\ell_{exp}(\alpha_i).$$ Therefore, for every $n \geq n_0$, the sum of the exponential lengths of complete factors in $[f^n(\gamma_w)]$ is at least equal to \begin{equation}\label{Equation p86}
(3^{n-n_0}-2C)\ell_{exp}^{\gamma_w}(\Gamma_{\gamma_w}') \geq (3^{n-n_0}-2C)\frac{(24C^2R)^2-1}{(24C^2R)^2}\frac{1}{24C^2R+1}\ell_{exp}^{\gamma_w}(\mathrm{Inc}(\gamma_w)).
\end{equation}

By the claim, for every $n \in \NN^*$, we have $\ell_{exp}^{[f^n(\gamma_w)]}(\mathrm{Inc}([f^n(\gamma_w)])) \leq 24C^2 \ell_{exp}^{\gamma_w}(\mathrm{Inc}(\gamma_w))$. Recall that the goodness function is a supremum over splittings of the considered path. Thus, by Equation~\eqref{Equation p86} for every $n \geq n_0$, since the maps $t \mapsto \frac{t}{t+a}$ are nonincreasing for every $a>0$, we have

$$\begin{array}{ccl} 
\mathfrak{g}([f^n(\gamma_w)]) & \geq & \frac{(3^{n-n_0}-2C)\frac{(24C^2R)^2-1}{(24C^2R)^2}\frac{1}{24C^2R+1}\ell_{exp}^{\gamma_w}(\mathrm{Inc}(\gamma_w))}{(3^{n-n_0}-2C)\frac{(24C^2R)^2-1}{(24C^2R)^2}\frac{1}{24C^2R+1}\ell_{exp}^{\gamma_w}(\mathrm{Inc}(\gamma_w))+ \ell_{exp}^{[f^n(\gamma_w)]}(\mathrm{Inc}([f^n(\gamma_w)])} \\ 
{} & \geq & \frac{(3^{n-n_0}-2C)\frac{(24C^2R)^2-1}{(24C^2R)^2}\frac{1}{24C^2R+1}\ell_{exp}^{\gamma_w}(\mathrm{Inc}(\gamma_w))}{(3^{n-n_0}-2C)\frac{(24C^2R)^2-1}{(24C^2R)^2}\frac{1}{24C^2R+1}\ell_{exp}^{\gamma_w}(\mathrm{Inc}(\gamma_w))+ 24C^2\ell_{exp}^{\gamma_w}(\mathrm{Inc}(\gamma_w))} \\
{} & \geq & \frac{(3^{n-n_0}-2C)\frac{(24C^2R)^2-1}{(24C^2R)^2}\frac{1}{24C^2R+1}}{(3^{n-n_0}-2C)\frac{(24C^2R)^2-1}{(24C^2R)^2}\frac{1}{24C^2R+1}+ 24C^2},
\end{array}$$ which goes to $1$ as $n$ goes to infinity. Hence there exists $n_1 \in \NN$ which is independent of $\gamma_w$, such that, for every path $\gamma_w$ as in Subcase~2 and every $n \geq n_1$, we have: $\mathfrak{g}([f^n(\gamma_w)]) \geq \delta$. This concludes the proof of Lemma~\ref{Lem control goodness versus decrease bad length} when Case~1 occurs.

\medskip

\noindent{\bf Case 2 } Suppose that, contrarily to Case~1, we have  
$$
\frac{\ell_{exp}^{\gamma_w}(\Lambda_{\gamma_w})}{\ell_{exp}^{\gamma_w}(\mathrm{Inc}(\gamma_w))} \geq \frac{1}{(24C^2R)^2}.
$$

Let $\alpha \in \Lambda_{\gamma_w}$ and consider the decomposition of the reduced path $\alpha$ into maximal subsegments $\alpha^{(1)}\ldots\alpha^{(k_{\alpha})}$ of exponential length relative to $\gamma_w$ equal to $2000R^3C^6$, except possibly the last one of exponential length relative to $\gamma_w$ less than or equal to $2000R^3C^6$. Let

$$\Lambda_{\gamma_w}^{(1)}=\left\{\alpha^{(j)} \; | \; \alpha \in \Lambda_{\gamma_w}, j \in \{1,\ldots,k_{\alpha}\}, \ell_{exp}^{\gamma_w}(\alpha^{(j)})=2000R^3C^6 \right\},$$

$$\Lambda_{\gamma_w}^{(2)}=\left\{\alpha^{(j)} \; | \; \alpha \in \Lambda_{\gamma_w}, j \in \{1,\ldots,k_{\alpha}\}, \ell_{exp}^{\gamma_w}(\alpha)<2000R^3C^6 \right\}.$$

Note that, since for every $\alpha \in \Lambda_{\gamma_w}$, we have $\ell_{exp}^{\gamma_w}(\alpha) \geq (3.10^8)R^6C^{12}+1$, we see that \begin{equation}\label{Equation p74 first one}
|\Lambda_{\gamma_w}^{(1)}| \geq 120000R^3C^6|\Lambda_{\gamma_w}^{(2)}|.
\end{equation}

Note that every element in $\Lambda_{\gamma_w}^{(1)} \cup \Lambda_{\gamma_w}^{(2)}$ has exponential length at most equal to $2000R^3C^6+1+2C$ by Lemma~\ref{Lem bound exponential length subpath}. By Lemma~\ref{Lem complete splitting rel poly grow}, there exists $M \in \NN^*$ depending only on $f$ such that for every $n \geq M$ and every reduced edge path $\alpha$ of exponential length at most equal to $(3.10^8)R^6C^{12}+1+2C$, either $[f^n(\alpha)]$ is a concatenation of paths in $G_{PG}$ and in $\mathcal{N}_{PG}$ or the following holds:

\medskip

\noindent{$(a)$ } there exists a complete factor of $[f^n(\alpha)]$ whose exponential length is at least equal to $10C$;

\medskip

\noindent{$(b)$ } the exponential length of an incomplete factor of $[f^n(\alpha)]$ is at most equal to $8C$.

\medskip

This applies in particular to every element $\alpha \in \Lambda_{\gamma_w}^{(1)} \cup \Lambda_{\gamma_w}^{(2)}$ and to every element $\alpha \in \Lambda_{\gamma_w}'$. For every $\alpha^{(j)} \in \Lambda_{\gamma_w}^{(1)}$ and every $n \geq M$, let $\alpha^{(j,n)}$ be the (possibly degenerate) subpath of $[f^n(\alpha^{(j)})]$ contained in $[f^n(\alpha)]$. Let $\Lambda_{\gamma_w}^{(3)}$ be the subset of $\Lambda_{\gamma_w}^{(1)}$ consisting in all $\alpha^{(j)} \in \Lambda_{\gamma_w}^{(1)}$ such that $\ell_{exp}(\alpha^{(j,M)}) \leq 80C^2$, and let $\Lambda_{\gamma_w}^{(4)}=\Lambda_{\gamma_w}^{(1)}-\Lambda_{\gamma_w}^{(3)}$. 

Suppose first that 
\begin{equation}\label{Equation p74 second equation}
|\Lambda_{\gamma_w}^{(4)}|>\frac{1}{30000R^3C^6}|\Lambda_{\gamma_w}^{(3)}|.
\end{equation}

Therefore, as $|\Lambda_{\gamma_w}^{(1)}|=|\Lambda_{\gamma_w}^{(3)}|+|\Lambda_{\gamma_w}^{(4)}|$, by Equation~\eqref{Equation p74 first one}, we have $$|\Lambda_{\gamma_w}^{(2)}| \leq \frac{30001R^3C^6}{120000R^3C^6}|\Lambda_{\gamma_w}^{(4)}|=K_0|\Lambda_{\gamma_w}^{(4)}|,$$ where $K_0$ is a constant depending only on $C$ and $R$. Note that $\Lambda_{\gamma_w}=\Lambda_{\gamma_w}^{(2)} \cup \Lambda_{\gamma_w}^{(3)} \cup \Lambda_{\gamma_w}^{(4)}$ and for every $j \in \{2,3,4\}$, every path in $\Lambda_{\gamma_w}^{(j)}$ has exponential length at most equal to $2000R^3C^6$. Thus, we see that $$\ell_{exp}^{\gamma_w}(\Lambda_{\gamma_w}) \leq 2000R^3C^6(|\Lambda_{\gamma_w}^{(2)}|+|\Lambda_{\gamma_w}^{(3)}|+|\Lambda_{\gamma_w}^{(4)}|) \leq K_0'|\Lambda_{\gamma_w}^{(4)}|$$ for some constant $K_0'$ depending only on $C$ and $R$.

Recall that if $\alpha^{(j)} \in \Lambda_{\gamma_w}^{(4)}$, then $\ell_{exp}(\alpha^{(j,M)}) > 80C^2$. Suppose towards a contradiction that $[f^M(\alpha^{(j)})]$ is a concatenation of paths in $G_{PG}$ and in $\mathcal{N}_{PG}$. Since $\alpha^{(j,M)}$ is a subpath of $[f^M(\alpha^{(j)})]$, we have $\ell_{exp}^{[f^M(\alpha^{(j)})]}(\alpha^{(j,M)})=0$. By Lemma~\ref{Lem bound exponential length subpath}, we see that $\ell_{exp}(\alpha^{(j,M)}) \leq \ell_{exp}^{[f^M(\alpha^{(j)})]}(\alpha^{(j,M)})+2C=2C$, which leads to a contradiction. Hence $[f^M(\alpha^{(j)})]$ satisfies $(a)$ and $(b)$. Note that $\alpha^{(j,M)}$ is a subpath of $[f^M(\alpha^{(j)})]$. Since $\ell_{exp}(\alpha^{(j,M)}) > 80C^2$, since every incomplete factor of $[f^M(\alpha^{(j)})]$ has exponential length at most equal to $8C$ by $(b)$ and since an incomplete factor of $[f^M(\alpha^{(j)})]$ is followed by a complete factor of $[f^M(\alpha^{(j)})]$, we see that $\alpha^{(j,M)}$ contains a subpath of a complete factor of $[f^M(\alpha^{(j)})]$. Since $\ell_{exp}(\alpha^{(j,M)})> 80C^2$ and since every incomplete subpath of $[f^M(\alpha^{(j)}]$ has exponential length at most equal to $8C$, the path $\alpha^{(j,M)}$ must contain a subpath $\alpha^{(j,M)'}$ such that the total exponential length of complete factors of $\alpha^{(j,M)'}$ is at least equal to $10C$. Let $\alpha_0^{(j,M)}$ be the minimal concatenation of splittings of a fixed optimal splittings of $[f^m(\alpha^{(j)})]$ which contains $\alpha^{(j,M)'}$. Let $\tau_1^{(j,M)}$ and $\tau_2^{(j,M)}$ be paths such that $[f^M(\alpha^{(j)})]=\tau_1^{(j,M)}\alpha_0^{(j,M)}\tau_2^{(j,M)}$. By Lemma~\ref{Lem concatenation of completely split paths and complete edges} applied twice (once with $\gamma=\alpha_0^{(j,M)}\tau_2^{(j,M)}[f^M(\alpha^{(j+1)}\ldots\alpha_k^{(k_{\alpha_k})})]$ and $\gamma_1=\alpha_0^{(j,M)}$ and once with $\gamma^{-1}=[f^M(\alpha_1^{(1)}\ldots \alpha^{(j-1)})]\tau_1^{(j,M)}\alpha_0^{(j,M)}$ and $\gamma_1^{-1}=\alpha_0^{(j,M)}$), we see that $\alpha^{(j,M)}$ contains a complete factor of $[f^M(\gamma_w)]$ of exponential length at least equal to $10C-4C=6C$. By Lemma~\ref{Lem bound exponential length subpath}, the path $\alpha^{(j,M)}$ contains a complete factor of $[f^M(\gamma_w)]$ of exponential length relative to $[f^M(\gamma_w)]$ at least equal to $C$. By Lemma~\ref{Lem p length of a completely split path grows linearly} (with $\gamma$ a complete factor contained in $\alpha^{(j,M)}$), for every $n \geq M$ and every $\alpha^{(j)} \in \Lambda_{\gamma_w}^{(4)}$, the path $\alpha^{(j,n)}$ contains a complete subpath of $[f^n(\gamma_w)]$ of exponential length at least equal to $3^{n-M}C$. By Lemma~\ref{Lem bound exponential length subpath}, for every $n \geq M$ and every $\alpha^{(j)} \in \Lambda_{\gamma_w}^{(4)}$, the path $\alpha^{(j,n)}$ contains a complete subpath of $[f^n(\gamma_w)]$ of exponential length relative to $[f^n(\gamma_w)]$ at least equal to $3^{n-M}C-2C$. Hence for every $n \geq M$, the sum of the exponential length relative to $[f^n(\gamma_w)]$ of complete factors contained in $[f^n(\gamma_w)]$ is at least equal to $(3^{n-M}C-2C)|\Lambda_{\gamma_w}^{(4)}|$.

By the claim, for every $n \geq M$, we have $$\ell_{exp}^{[f^n(\gamma_w)]}(\mathrm{Inc}([f^n(\gamma_w)])) \leq 24C^2\ell_{exp}^{\gamma_w}(\gamma_w) \leq 24C^2\frac{1}{1-\delta}\ell_{exp}^{\gamma_w}(\mathrm{Inc}(\gamma_w)),$$ where the last inequality holds by Equation~\eqref{Equation bad}. Using the above equations and the assumptions of Case~2, we see that
$$\begin{array}{ccl}
\ell_{exp}^{[f^n(\gamma_w)]}(\mathrm{Inc}([f^n(\gamma_w)])) & \leq & 24C^2\frac{1}{1-\delta}\ell_{exp}^{\gamma_w}(\mathrm{Inc}(\gamma_w)) \\
{} & \leq & 24C^2\frac{1}{1-\delta}(24C^2R)^2\ell_{exp}^{\gamma_w}(\Lambda_{\gamma_w}) \\
{} & \leq & 24C^2\frac{1}{1-\delta}(24C^2R)^2K_0'|\Lambda_{\gamma_w}^{(4)}|=K_1|\Lambda_{\gamma_w}^{(4)}|,
\end{array}$$ where $K_1$ is a constant depending only on $C$, $R$ and $\delta$. Thus, since the goodness function is a supremum over all splittings of the considered path, for every $n \geq M$, we have:

$$
\begin{array}{ccl}
\mathfrak{g}([f^n(\gamma_w)]) & \geq & \frac{(3^{n-M}C-2C)|\Lambda_{\gamma_w}^{(4)}|}{(3^{n-M}C-2C)|\Lambda_{\gamma_w}^{(4)}|+\ell_{exp}^{[f^n(\gamma_w)]}(\mathrm{Inc}([f^n(\alpha)]))} \\
{} & \geq & \frac{(3^{n-M}C-2C)|\Lambda_{\gamma_w}^{(4)}|}{(3^{n-M}C-2C)|\Lambda_{\gamma_w}^{(4)}|+K_1|\Lambda_{\gamma_w}^{(4)}|} \\
{} & = & \frac{3^{n-M}C-2C}{3^{n-M}C-2C+K_1}, 
\end{array}
$$
which converges to $1$ as $n$ goes to infinity. Hence there exists $M' \in \NN^*$ depending only on $f$ such that for every $n\geq M$, we have $\mathfrak{g}([f^n(\gamma_w)]) \geq \delta$. This proves Lemma~\ref{Lem control goodness versus decrease bad length} in this case.

\medskip

Suppose now that contrarily to Equation~\eqref{Equation p74 second equation}, we have 
\begin{equation}\label{Equation p75}
|\Lambda_{\gamma_w}^{(4)}| \leq\frac{1}{30000R^3C^6}|\Lambda_{\gamma_w}^{(3)}|.
\end{equation}
Then $$|\Lambda_{\gamma_w}^{(1)}|=|\Lambda_{\gamma_w}^{(3)}|+|\Lambda_{\gamma_w}^{(4)}| \leq \left(1+\frac{1}{30000R^3C^6}\right)|\Lambda_{\gamma_w}^{(3)}|.$$  

\medskip

\noindent{\bf Claim 2 } Let $n \geq M$, let $\alpha^{(j)} \in \Lambda_{\gamma_w}^{(2)} \cup \Lambda_{\gamma_w}^{(4)}$. The total exponential length of incomplete factors of $[f^n(\gamma_w)]$ contained in $\alpha^{(j,n)}$ is at most equal to $12C\ell_{exp}(\alpha^{(j)})$.

\medskip

\dem Let $\sigma$ be an incomplete factor of $[f^n(\gamma_w)]$ which is contained in $\alpha^{(j,M)}$. Then one of the following holds:

\medskip

\noindent{$(i)$ } the path $\sigma$ is an incomplete factor of $[f^n(\alpha^{(j)})]$;

\medskip

\noindent{$(ii)$ } the path $\sigma$ contains a subpath which is complete in $[f^n(\alpha^{(j)})]$. 

\medskip

Note that the total exponential length of incomplete factors of $[f^n(\gamma_w)]$ which satisfy $(i)$ is bounded by the total exponential length of incomplete factors of $[f^n(\alpha^{(j)})]$. Thus, by Lemma~\ref{Lem uniform bound on bad subpaths}, the total exponential length of incomplete factors of $[f^n(\gamma_w)]$ which satisfy $(i)$ is bounded by $8C\ell_{exp}(\alpha^{(j)})$. Suppose that $\sigma$ satisfies~$(ii)$. Let $\alpha^{(j,n)}=a_1ca_2$ be a decomposition of $\alpha^{(j,n)}$ where for every $i \in \{1,2\}$, the total exponential length of complete factors of $[f^n(\alpha^{(j)})]$ contained in $a_i$ is equal to $2C$. By Lemma~\ref{Lem concatenation of completely split paths and complete edges} applied to $\gamma=[f^n(\alpha^{(j)})][f^n(\alpha^{(j+1)}\ldots\alpha_k^{(k_{\alpha_k})})]$ and $\gamma_1=[f^n(\alpha^{(j)})]$ and to $\gamma^{-1}=[f^n(\alpha_1^{(1)}\ldots\alpha^{(j-1)})][f^n(\alpha^{(j)})]$ and $\gamma_1^{-1}=[f^n(\alpha^{(j)})]$, the path $\sigma$ is contained in either $a_1$ or $a_2$. For every $t \in \{1,2\}$, let $a_t=b_1^{(t)}b_1^{(t)'}\ldots b_s^{(t)}b_{s_t}^{(t)'}$ be a decomposition of $a_t$ where, for every $i \in \{1,\ldots,s_t\}$, the path $b_i^{(t)}$ is an incomplete factor of $[f^n(\alpha^{(j)})]$ and for every $i \in \{1,\ldots,s_t\}$, the path $b_i^{(t)'}$ is a complete factor of $[f^n(\alpha^{(j)})]$ contained in $a_t$. 

Suppose that there exists $i \in \{1,\ldots,s_1\}$ such that $b_i^{(1)'}$ is a complete factor of $[f^n(\gamma_w)]$. We claim that for every $j \geq i+1$, the path $b_j^{(1)'}$ is a complete factor of $[f^n(\gamma_w)]$. Indeed, let $n' \geq n$  and let $j \geq i+1$. Then there is no identification between an initial segment of $[f^{n'}(b_i^{(1)'})]$ and an initial segment of $[f^{n}(\gamma_w)]$ not intersecting $\alpha^{(j,n')}$ as otherwise there would exist identifications with $[f^{n'}(b_i^{(1)'})]$, contradicting the fact that $b_i^{(1)'}$ is complete. Similarly, there is no identification between a terminal segment of $[f^{n'}(b_i^{(1)'})]$ and a terminal segment of $[f^{n}(\gamma_w)]$ not intersecting $\alpha^{(j,n')}$ as otherwise there would exist identifications with $[f^{n'}(c)]$. The claim follows. Similarly, if there exists $i \in \{1,\ldots,s_2\}$ such that $b_i^{(2)'}$ is a complete factor of $[f^n(\gamma_w)]$, then for every $j<i$, the path $b_j^{(2)'}$ is a complete factor of $[f^n(\gamma_w)]$. Hence we may assume that for every $t \in \{1,2\}$ and every $s \in \{1,\ldots,s_t\}$, the path $b_s^{(t)'}$ is incomplete in $[f^n(\gamma_w)]$. Therefore, for every $t \in \{1,2\}$, the whole path $a_t$ is incomplete in $[f^n(\gamma_w)]$. Therefore, in order to prove the claim, it suffices to bound the exponential lengths of $a_1$ and $a_2$. Let $t \in \{1,2\}$. By Lemma~\ref{Lem Exponential length less exp length subpaths}, we have $$\ell_{exp}(a_t) \leq \sum\limits_{i=1}^{s_t} \ell_{exp}(b_i^{(t)})+\ell_{exp}(b_i^{(t)'}).$$ For every $i \in \{1\ldots,s_t\}$, the path $b_i^{(t)}$ satisfies $(i)$ and we already have a bound on the total exponential length of such paths. Moreover, since the total exponential length of complete factors of $\alpha^{(j,n)}$ contained in $a_t$ is at most equal to $2C$, we have $$\sum\limits_{i=1}^{s_t} \ell_{exp}(b_i^{(t)'}) \leq 2C.$$ 
Thus, the total exponential length of incomplete factors of $[f^n(\gamma_w)]$ contained in $\alpha^{(j,M)}$ is at most equal to $$8C\ell_{exp}(\alpha^{(j)})+ \sum\limits_{t=1}^2\sum\limits_{i=1}^{s_t} \ell_{exp}(b_i^{(t)'}) \leq 8C\ell_{exp}(\alpha^{(j)})+4C\leq 12C\ell_{exp}(\alpha^{(j)}),$$ where the last inequality follows from the fact that every element of $\Lambda_{\gamma_w}^{(2)} \cup \Lambda_{\gamma_w}^{(4)}$ has positive exponential length.
\hfill\qedsymbol

\medskip

By Claim~2 and Lemma~\ref{Lem bound exponential length subpath}, for every $n \geq M$ and every $\alpha^{(j)} \in \Lambda_{\gamma_w}^{(2)} \cup \Lambda_{\gamma_w}^{(4)}$, the total exponential length relative to $[f^n(\gamma_w)]$ of incomplete factors in the subpath of $[f^n(\gamma_w)]$ contained in $[f^n(\alpha^{(j)})]$ is at most equal to $12C\ell_{exp}^{\gamma_w}(\alpha^{(j)})+2C \leq 14C\ell_{exp}^{\gamma_w}(\alpha^{(j)})$. Hence by definition, for every $n \geq M$ and every path $\alpha^{(j)} \in \Lambda_{\gamma_w}^{(2)} \cup \Lambda_{\gamma_w}^{(4)}$, we have $$\ell_{exp}^{[f^n(\gamma_w)]}(\mathrm{Inc}([f^n(\gamma_w)]) \cap \alpha^{(j,n)}) \leq 14C\ell_{exp}(\alpha^{(j)}).$$ We claim that, for every $n \geq M$, every element in $\Lambda_{[f^n(\gamma_w)]}$ is contained in an iterate of an element in $\Lambda_{\gamma_w}$. Indeed, note that, by the choice of $M$ (in the above application of Lemma~\ref{Lem complete splitting rel poly grow}), for every element $\alpha \in \Lambda_{\gamma_w}'$, the exponential length of an incomplete factor in $[f^n(\alpha)]$ is at most equal to $8C$. Hence an incomplete factor of $[f^n(\alpha)]$ whose exponential length is at least equal to $(3.10^8)R^6C^{12}+1$ cannot be contained in an iterate of an element of $\Lambda_{\gamma_w}$. The claim follows. Therefore, using Equation~\eqref{Equation p75} for the third inequality, the value of $\ell_{exp}(\Lambda_{[f^M(\gamma_w)]})$ is at most equal to 

$$
\begin{array}{l}
\sum\limits_{\alpha^{(j)} \in \Lambda_{\gamma_w}^{(3)}} \ell_{exp}(\alpha^{(j,M)}) + \sum\limits_{\alpha^{(j)} \in \Lambda_{\gamma_w}^{(4)}}\ell_{exp}^{[f^M(\gamma_w)]}(\mathrm{Inc}([f^M(\gamma_w)]) \cap \alpha^{(j,M)}) \\
+\sum\limits_{\alpha^{(j)} \in \Lambda_{\gamma_w}^{(2)}}\ell_{exp}^{[f^M(\gamma_w)]}(\mathrm{Inc}([f^M(\gamma_w)]) \cap \alpha^{(j,M)}) \\
\leq 80C^2|\Lambda_{\gamma_w}^{(3)}|+ 14C\sum\limits_{\beta \in \Lambda_{\gamma_w}^{(4)}}\ell_{exp}(\beta)+14C\sum\limits_{\alpha \in \Lambda_{\gamma_w}^{(2)}} \ell_{exp}(\alpha) \\
\leq 80C^2|\Lambda_{\gamma_w}^{(3)}| + 14C(2000R^3C^6)|\Lambda_{\gamma_w}^{(4)}|+14C\sum\limits_{\alpha \in \Lambda_{\gamma_w}^{(2)}}\ell_{exp}(\alpha) \\
\leq 80C^2|\Lambda_{\gamma_w}^{(3)}| + C|\Lambda_{\gamma_w}^{(3)}|+14C\sum\limits_{\alpha \in \Lambda_{\gamma_w}^{(2)}}\ell_{exp}(\alpha) \\
\leq 81C^2|\Lambda_{\gamma_w}^{(3)}|+14C\sum\limits_{\alpha \in \Lambda_{\gamma_w}^{(2)}}\ell_{exp}(\alpha).
\end{array}
$$

Since by Equation~\eqref{Equation p74 first one} $$\left(1+\frac{1}{30000R^3C^6}\right)|\Lambda_{\gamma_w}^{(3)}| \geq |\Lambda_{\gamma_w}^{(1)}| \geq 120000R^3C^6|\Lambda_{\gamma_w}^{(2)}|,$$ we have $|\Lambda_{\gamma_w}^{(3)}| \geq 60000R^3C^6|\Lambda_{\gamma_w}^{(2)}|$. Hence we have

$$
\begin{array}{ccl}
\ell_{exp}^{[f^M(\gamma_w)]}(\Lambda_{[f^M(\gamma_w)]}) & \leq & 81C^2|\Lambda_{\gamma_w}^{(3)}|+14C\sum_{\alpha \in \Lambda_{\gamma_w}^{(2)}}\ell_{exp}(\alpha) \\
{} & \leq & 81C^2|\Lambda_{\gamma_w}^{(3)}|+(14C)(2000R^3C^6)|\Lambda_{\gamma_w}^{(2)}| \\
{} & \leq & 81C^2|\Lambda_{\gamma_w}^{(3)}|+2C|\Lambda_{\gamma_w}^{(3)}|=83C^2|\Lambda_{\gamma_w}^{(3)}|.
\end{array}
$$

Let $n \geq M$. Suppose first that $$
\frac{\ell_{exp}^{[f^n(\gamma_w)]}(\Lambda_{[f^n(\gamma_w)]})}{\ell_{exp}^{[f^n(\gamma_w)]}(\mathrm{Inc}([f^n(\gamma_w)]))} < \frac{1}{(24C^2R)^2}.
$$
Then we can apply Case~1 to conclude the proof of Lemma~\ref{Lem control goodness versus decrease bad length}. Otherwise, we have $$(24C^2R)^2\ell_{exp}^{[f^n(\gamma_w)]}(\Lambda_{[f^n(\gamma_w)]}) \geq \ell_{exp}^{[f^n(\gamma_w)]}(\mathrm{Inc}([f^n(\gamma_w)])).
$$
By Lemma~\ref{Lem uniform bound on bad subpaths} and Lemma~\ref{Lem bound exponential length subpath}, we have $$\begin{array}{ccl} 
\ell_{exp}^{[f^n(\gamma_w)]}(\mathrm{Inc}([f^n(\gamma_w)]) & \leq & \ell_{exp}(\mathrm{Inc}([f^n(\gamma_w)]) \leq 8C\ell_{exp}(\mathrm{Inc}([f^M(\gamma_w)]) \\
{} & \leq & 10C\ell_{exp}^{[f^M(\gamma_w)]}(\mathrm{Inc}([f^M(\gamma_w)]).
\end{array}$$ Hence we have
$$
\begin{array}{ccl}
\frac{\ell_{exp}^{[f^n(\gamma_w)]}(\mathrm{Inc}([f^n(\gamma_w)]))}{\ell_{exp}^{\gamma_w}(\mathrm{Inc(\gamma_w)})} &\leq & \frac{(24C^2R)^2\ell_{exp}^{[f^n(\gamma_w)]}(\mathrm{Inc}([f^n(\gamma_w)]))}{\ell_{exp}^{[f^M(\gamma_w)]}(\mathrm{Inc}([f^M(\gamma_w)]))}\frac{\ell_{exp}^{[f^M(\gamma_w)]}(\mathrm{Inc}([f^M(\gamma_w)]))}{\ell_{exp}^{\gamma_w}(\mathrm{Inc(\gamma_w)})} \\ 
{} & \leq & \frac{10C(24C^2R)^2\ell_{exp}^{[f^n(\gamma_w)]}(\Lambda_{[f^n(\gamma_w)]})}{\ell_{exp}^{[f^n(\gamma_w)]}(\Lambda_{\gamma_w})} \\
{} & \leq & \frac{10C(24C^2R)^2(83C^2|\Lambda_{\gamma_w}^{(3)}|)}{2000R^3C^6|\Lambda_{\gamma_w}^{(3)}|} \\
{} & \leq & \frac{10C}{R}.
\end{array}
$$

This concludes the proof of Lemma~\ref{Lem control goodness versus decrease bad length}.
\hfill\qedsymbol

\bigskip

In the next proposition, we need to work with CT maps that represent both an almost atoroidal outer automorphism and its inverse. We therefore introduce the following conventions: 

\medskip

\noindent{\it Let $f' \colon G' \to G'$ be a CT map representing $\phi^{-M}$, which exists by Theorem~\ref{Theo existence CT}. We denote by $K'$ the constant similar to the constant $K$ given above~Lemma~\ref{Lem bound exponential length subpath} and by $C_{f'}$ the bounded cancellation constant given by Lemma~\ref{Lem Bounded cancellation lemma}. We set $C'=\max\{K',C_{f'}\}$ as in Equation~\eqref{Equation defi C}. We denote by $G_{p'}$ the invariant subgraph of $G'$ such that $\mathcal{F}(G_{p'})=\mathcal{F}$, by $\ell_{\mathcal{F}'}$ the corresponding $\mathcal{F}$-length  and by $\ell_{exp'}$ the corresponding exponential length. Let $\mathfrak{g}'$ be the corresponding goodness function. If $w \in F_{\tt n}$, we denote by $\gamma_w'$ the corresponding circuit in $G'$. } 

\medskip

We also need a result which shows that the exponential length is invariant by $F_{\tt n}$-equivariant quasi-isometry. In order to prove this, we need some additional definitions. Let $G$ be a connected (pointed) graph whose fundamental group is isomorphic to $F_{\tt n}$ and let $\widetilde{G}$ be the universal cover of $G$. Let $\phi \in \Out(F_{\tt n})$ be an exponentially growing outer automorphism. Let $\widehat{G}$ be the graph obtained from $\widetilde{G}$ as follows. We add one vertex $v_{gA}$ for every left class $gA$, with $g \in F_{\tt n}$ and $A$ is a subgroup of $F_{\tt n}$ such that $[A] \in \mathcal{A}(\phi)$ and we add one edge between $v_{gA}$ and a vertex $v$ of $\widetilde{G}$ if and only if the vertex $v$ is contained in the tree $T_{gAg^{-1}}$. The graph $\widehat{G}$ is known as the \emph{electrification of $\widetilde{G}$} (see for instance~\cite{Bowditch2012}). For a path $\gamma$ in $G$, we denote by $\widetilde{\gamma}$ a lift of $\gamma$ in $\widetilde{G}$. Let $\widehat{\gamma}$ be the path in $\widehat{G}$ constructed as follows. Let $\widetilde{\gamma}=a_1b_1\ldots a_kb_k$ be the decomposition of $\widetilde{\gamma}$ such that, for every $i \in \{1,\ldots,k\}$, the path $b_i$ is contained in some tree $T_{g_iA_ig_i^{-1}}$ with $g_i \in F_{\tt n}$, $A_i$ a subgroup of $F_{\tt n}$ such that $[A_i] \in \mathcal{A}(\phi)$ and $b_i$ is maximal for the property of being contained in such a tree $T_{g_iA_ig_i^{-1}}$. Then $\widehat{\gamma}$ is a path $\widehat{\gamma}=a_1c_1\ldots a_kc_k$ where, for every $i \in \{1,\ldots,k\}$, the path $c_i$ is the two-edge path whose endpoints are the endpoints of $b_i$ and the middle vertex of $c_i$ is $v_{g_iA_i}$. Note that the path $\widehat{\gamma}$ is not uniquely determined. Indeed, it is possible that there exists $i \in \{1,\ldots,k\}$ such that $b_i$ is contained in two distinct trees $T_A$ and $T_B$ with $[A],[B] \in \mathcal{A}(\phi)$. However, if $\widehat{\gamma}$ and $\widehat{\gamma}'$ are two such paths associated with $\widetilde{\gamma}$, then $\ell(\widehat{\gamma})=\ell(\widehat{\gamma}')$.

\begin{prop}\label{Prop explength quasi isometry invariant}
Let ${\tt n} \geq 3$, let $\phi \in \Out(F_{\tt n})$ and let $f \colon G \to G$ be a CT map representing a power of $\phi$.

\noindent{$(1)$ } There exists a constant $B_0\geq 1$ such that, for every element $w \in F_{\tt n}$ with $\ell_{exp}(\gamma_w)>0$, we have: $$ \frac{1}{B_0}\ell_{exp}(\gamma_w) \leq \ell(\widehat{\gamma_w}) \leq B_0\;\ell_{exp}(\gamma_w).$$

\noindent{$(2)$ } Let $f' \colon G' \to G'$ be a CT map representing a power of $\phi^{-1}$. There exists a constant $B >0$ such that, for every element $w \in F_{\tt n}$, we have:
$$ \frac{1}{B}\ell_{exp'}(\gamma_w') \leq \ell_{exp}(\gamma_w) \leq B\ell_{exp'}(\gamma_w').$$
\end{prop}

\dem $(1)$ Recall the definition of the graph $G^{\ast}$ from just above Lemma~\ref{Lem Injectivity pi1}. We can turn the graph $G^{\ast}$ into a metric graph by assigning, to every edge $e \in \vec{E}G^{\ast}$, the length equal to the length of the path $p_{G^{\ast}}(e)$ in $G$. Since the graph $G^{\ast}$ is finite, there exists a constant $B'$ such that the diameter of every maximal subtree of $G^{\ast}$ is at most $B'$. Let $B_0=2B'+2$. 

Let $w \in F_{\tt n}$. Let $\gamma_w=a_1b_1\ldots a_kb_k$ be the decomposition of $\gamma_w$ with $a_1$ and $b_k$ possibly empty such that, for every $i \in \{1,\ldots,k\}$, the path $b_i$ is a maximal concatenation of paths in $G_{PG}'$ and in $\mathcal{N}_{PG}$ and, for every $i \in \{1,\ldots,k\}$ and every edge $e$ of $a_i$, we have $\ell_{exp}^{\gamma_w}(e)=1$. Note that by the definition of the exponential length we have $$\ell_{exp}(\gamma_w)=\sum_{i=1}^k \ell(a_i). $$ Let $A$ be a subgroup of $F_{\tt n}$ such that $[A] \in \mathcal{A}(\phi)$. Let $i \in \{1,\ldots,k\}$ and let $\alpha$ be a subpath of $a_i$ whose lift is contained in $T_A$. By Proposition~\ref{Prop circuits in Gpg are elements in poly subgroup}, the subpath $\alpha$ is contained in a concatenation of paths in $G_{PG}$ and in $\mathcal{N}_{PG}$. Since $a_i$ does not contain any concatenation of paths in $G_{PG}$ and $\mathcal{N}_{PG}$, the path $\alpha$ is a proper subpath of an EG INP. By the definition of $C$ (see~Equation~\eqref{Equation defi C}), we see that $\ell(\alpha) \leq C$. Thus, we have: $\ell(a_i) \leq C\ell(\widehat{a}_i)$ and $$\ell_{exp}(\gamma_w) \leq C\sum_{i=1}^k\ell(\widehat{a}_i).$$ 

\medskip

\noindent{\bf Claim. } Let $A$ be a subgroup of $F_{\tt n}$ such that $[A] \in \mathcal{A}(\phi)$. Let $\beta$ be a subpath of $\gamma_w$ such that a lift of $\beta$ is contained in $T_A$. There does not exist $i \in \{1,\ldots,k\}$ such that both $\beta \cap b_i$ and $\beta \cap b_{i+1}$ are not reduced to a point.

\medskip

\dem Suppose towards a contradiction that such an element $i \in \{1,\ldots,k\}$ exists. Then $a_{i+1}$ is contained in $\beta$. By the above, the path $a_{i+1}$ is contained in an EG INP $\sigma$. Since both $b_i$ and $b_{i+1}$ are concatenations of paths in $G_{PG}'$ and $\mathcal{N}_{PG}$, the path $a_{i+1}$ must contain the initial or the terminal segment of $\sigma$. Since $\beta$ is contained in a concatenation of paths in $G_{PG}$ and in $\mathcal{N}_{PG}$ by Proposition~\ref{Prop circuits in Gpg are elements in poly subgroup}, the EG INP $\sigma$ must be contained in $\beta$ and $\beta \cap a_{i+1} \subseteq \sigma$. This contradicts the maximality of the paths $b_i$ and $b_{i+1}$. 
\hfill\qedsymbol

\medskip

Hence $\beta$ is either contained in $b_ia_{i+1}$ or in $a_{i+1}b_{i+1}$. Let $i \in \{1,\ldots,k\}$ and let $\beta$ be a maximal subpath of $\gamma_w$ containing edges of $a_i$ and such that a lift of $\beta$ is contained in some $T_A$ with $A$ a subgroup of $F_{\tt n}$ such that $[A] \in \mathcal{A}(\phi)$. By the claim, the path $a_i$ has a decomposition $a_i=c_i^+d_ic_i^-$ such that $c_i^+$ and $c_i^-$ are possibly trivial, lifts of $c_i^+$ and $c_i^-$ are contained in trees $T_{A_+}$ and $T_{A_-}$ with $A_+$ and $A_-$ subgroups of $F_{\tt n}$ such that $[A_+],[A_-] \in \mathcal{A}(\phi)$ and one of the following holds:

\noindent{$(a)$ } $\beta \subseteq d_i$;

\noindent{$(b)$ } $\beta \cap a_i \neq \beta$ and $\beta \cap a_i \in \{c_i^+,c_i^-\}$. 

Note that for every $i \in \{1,\ldots,k\}$, we have $\ell(\widehat{a_i}) \leq \ell(\widehat{d_i})+4$. Then $$\ell(\widehat{\gamma_w}) \geq \sum_{i=1}^k\ell(\widehat{d_i}) \geq \sum_{i=1}^k(\ell(\widehat{a_i})-4)=\sum_{i=1}^k\ell(\widehat{a}_i)-4k.$$
Moreover, if $\beta$ is a path which satisfies the hypothesis of the claim, then there exists at most one $i \in \{1,\ldots,k\}$ such that $\beta \cap b_i$ is not reduced to a point. Therefore, we see that $\ell(\widehat{\gamma_w}) \geq k$. Thus, we have $$\ell_{exp}(\gamma_w) \leq C\sum_{i=1}^k\ell(\widehat{a}_i) \leq C(\ell(\widehat{\gamma}_w)+4k) \leq 5C\ell(\widehat{\gamma}_w). $$ 

This proves the first inequality of Assertion~$(1)$. We now prove the second inequality. For every $i \in \{1,\ldots,k\}$, there exists a unique edge path $b_i^{\ast} \subseteq G^{\ast}$ such that $p^{\ast}(b_i^{\ast})=b_i$. Let $i \in \{1,\ldots,k\}$. Since $G^{\ast}$ is a finite graph, there exist (possibly trivial) reduced paths $\beta_i^{\ast},\delta_i^{\ast}$ and $\delta_i^{\ast '}$ such that:

\noindent{$(i)$ } the path $\beta_i^{\ast}$ is a circuit;

\noindent{$(ii)$ } the paths $\delta_i^{\ast}$ and $\delta_i^{\ast '}$ are contained in maximal trees of $G^{\ast}$;

\noindent{$(iii)$ } we have $b_i^{\ast}=\delta_i^{\ast}\beta_i^{\ast}\delta_i^{\ast '}$.

By Lemma~\ref{Lem Injectivity pi1}~$(1)$, the paths $p^{\ast}(\delta_i^{\ast})$, $p^{\ast}(\beta_i^{\ast})$ and $p^{\ast}(\delta_i^{\ast '})$ are reduced edge paths of $G$. By definition of $B'$, we have $\ell(\delta_i^{\ast}),\ell(\delta_i^{\ast '}) \leq B'$. Since $p^{\ast}(\beta_i^{\ast})$ is a circuit which is a concatenation of paths in $G_{PG}$ and in $\mathcal{N}_{PG}$, by Proposition~\ref{Prop circuits in Gpg are elements in poly subgroup}, there exists a subgroup $H_i$ of $F_{\tt n}$ such that $[H_i] \in \mathcal{A}(\phi)$ and the conjugacy classes of elements of $F_{\tt n}$ represented by $p^{\ast}(\beta_i^{\ast})$ are contained in $[H_i]$. Hence the length of $\widehat{p^{\ast}(\beta_i^{\ast})}$ is bounded by $2$. Hence the length of the path $\widehat{b}_i$ is bounded by $2+2B'=B_0$. Therefore, since $\ell_{exp}(\gamma_w)>0$, we have $$\ell(\widehat{\gamma}_w)=\sum_{i=1}^k \ell(a_i)+\ell(\widehat{b}_i) \leq \sum_{i=1}^k (\ell(a_i) +B_0) \leq (B_0+1)\sum_{i=1}^k\ell(a_i) = (B_0+1)\ell_{exp}(\gamma_w).$$ This proves Assertion~$(1)$.

\medskip

\noindent{$(2)$ } Let $f'$ be as in Assertion~$(2)$ and let $w \in F_{\tt n}$. Suppose first that $\ell_{exp}(\gamma_w)=0$. Then $\gamma_w$ is a concatenation of paths in $G_{PG}'$ and in $\mathcal{N}_{PG}$. By Proposition~\ref{Prop definition CT}~$(4)$ and Lemma~\ref{Lem No zero path and Nielsen path adjacent}, there does not exist an edge in a zero stratum which is adjacent to a concatenation of paths in $G_{PG}$ and in $\mathcal{N}_{PG}$. Since zero strata are contractible by Proposition~\ref{Prop definition CT}~$(3)$, it follows that $\gamma_w$ is a concatenation of paths in $G_{PG}$ and in $\mathcal{N}_{PG}$. By Proposition~\ref{Prop circuits in Gpg are elements in poly subgroup}, there exists a subgroup $A$ of $F_{\tt n}$ such that $[A] \in \mathcal{A}(\phi)$ and $w \in A$. Since $\mathcal{A}(\phi)=\mathcal{A}(\phi^{-1})$ by Equation~\eqref{Equation p6}, by Proposition~\ref{Prop circuits in Gpg are elements in poly subgroup}, we have $\ell_{exp'}(\gamma_w')=0$. So we may suppose that $\ell_{exp}(\gamma_w)>0$ and that $\ell_{exp'}(\gamma_w')>0$. By Assertion~$(1)$, in order to prove Assertion~$(2)$, it suffices to prove that $\widehat{G}$ and $\widehat{G}'$ are $F_{\tt n}$-equivariantly quasi-isometric. Since $\mathcal{A}(\phi)$ is a malnormal subgroup system, this follows from \cite[Theorem~7.11]{Bowditch2012} and~\cite[proof of Theorem~5.1]{Hruska10}.
\hfill\qedsymbol

\begin{prop}\label{Prop control goodness in f and inverse}
Let $\phi \in \Out(F_{\tt n},\mathcal{F})$ and let $f \colon G \to G$ be as in Remark~\ref{Rmq Convention for relative atoroidal CT map 2}. Let $f' \colon G' \to G'$ be as in the above convention. Let $\delta \in (0,1)$ and let $W$ be a neighborhood of $K_{PG}(\phi)$ in $\PCurr(F_{\tt n},\mathcal{F} \wedge \mathcal{A}(\phi))$. There exists $n_0 \in \NN^*$ such that for every $n \geq n_0$ and every nonperipheral element $w \in F_{\tt n}$ such that $\eta_{[w]} \notin W$, one of the following holds:
$$\mathfrak{g}([f^n(\gamma_w)]) \geq \delta$$
or
$$\mathfrak{g}'([f'^n(\gamma_w')]) \geq \delta.$$
\end{prop}

\dem Let $w \in F_{\tt n}$ be a nonperipheral element such that $\eta_{[w]} \notin W$. Let $R=\frac{10C}{(1-\delta)^2}8C'B^2$. We use the alternative given by Lemma~\ref{Lem control goodness versus decrease bad length} with the constants $\delta$ and $R$. If the first alternative of Lemma~\ref{Lem control goodness versus decrease bad length} occurs, then we are done. Suppose that $\mathfrak{g}([f^n(\gamma_w)]) < \delta$. There exists $n_0 \in  \NN^*$ depending only on $f$ such that for every $n \geq n_0$, we have $$\ell_{exp}^{[f^n(\gamma_w)]}(\mathrm{Inc}([f^n(\gamma_w)])) \leq \frac{10C}{R}\ell_{exp}^{\gamma_w}(\mathrm{Inc}(\gamma_w)). $$ By Lemma~\ref{Lem Goodness nondecreasing}, since $\mathfrak{g}([f^n(\gamma_w)]) < \delta$, we have $\mathfrak{g}(\gamma_w) < \delta$. Thus, we see that $$\ell_{exp}^{\gamma_w}(\mathrm{Inc}(\gamma_w)) \geq (1-\delta)\ell_{exp}(\gamma_w).$$ Let $\gamma''$ be the reduced circuit in $G$ such that $[f^{n_0}(\gamma'')]=\gamma_w$. Since $\mathfrak{g}(\gamma_w) <\delta$ and $[\eta_{[w]}] \notin K_{PG}(\phi)$, by Lemma~\ref{Lem control goodness versus decrease bad length}, we see that  
$$\ell_{exp}^{\gamma_w}(\mathrm{Inc}(\gamma_w)) \leq \frac{10C}{R}\ell_{exp}^{\gamma''}(\mathrm{Inc}(\gamma'')). $$
We have
$$
\begin{array}{ccl}
\ell_{exp'}([f'^{n_0}(\gamma_w')]) & \geq & \frac{1}{B}\ell_{exp}(\gamma'') \geq \frac{1}{B}\ell_{exp}^{\gamma''}(\mathrm{Inc}(\gamma'')) \\
{} & \geq & \frac{1}{B}\frac{R}{10C}\ell_{exp}^{\gamma_w}(\mathrm{Inc}(\gamma_w))  \geq \frac{1}{B}\frac{(1-\delta)R}{10C}\ell_{exp}(\gamma_w) \\
{} & \geq & \frac{1}{B^2} \frac{(1-\delta)R}{10C}\ell_{exp'}(\gamma_w') = 8C'\frac{1}{1-\delta}\ell_{exp'}(\gamma_w').
\end{array}
$$

But by Lemma~\ref{Lem uniform bound on bad subpaths}, we have:
$$
\ell_{exp'}^{[f'^{n_0}(\gamma_w')]}(\mathrm{Inc}(f'^{n_0}(\gamma_w')) \leq \ell_{exp'}(\mathrm{Inc}(f'^{n_0}(\gamma_w'))\leq  8C'\ell_{exp'}(\gamma_w').
$$

Therefore, we see that 
$$\mathfrak{g}'([f'^{n_0}(\gamma_w')])=1-\frac{\ell_{exp'}^{[f'^{n_0}(\gamma_w')]}(\mathrm{Inc}([f'^{n_0}(\gamma_w')])}{\ell_{exp'}([f'^{n_0}(\gamma_w')])} \geq 1-(1-\delta)=\delta>0.$$

By Lemma~\ref{Lem control of the goodness}, we see that there exists $n_1 \geq n_0$ depending only on $f'$ such that for every $n \ge n_1$, $$\mathfrak{g}'([f'^n(\gamma_w')]) \geq \delta.$$ This concludes the proof.
\hfill\qedsymbol

\begin{prop}\label{Prop North south dynamics outside neighborhood}
Let $\phi \in \Out(F_{\tt n},\mathcal{F})$ and let $f \colon G \to G$ be as in Remark~\ref{Rmq Convention for relative atoroidal CT map 2}. Let $U_+$ be a neighborhood of $\Delta_+(\phi)$, let $U_-$ be a neighborhood of $\Delta_-(\phi)$, let $V$ be a neighborhood of $K_{PG}(\phi)$. There exists $N \in \NN^*$ such that for every $n \geq 1$ and every $\mathcal{F} \wedge \mathcal{A}(\phi)$-nonperipheral $w \in F_{\tt n}$ such that $\eta_{[w]} \notin V$, one of the following holds $$\phi^{Nn}(\eta_{[w]}) \in U_+ \text{~~~~ or ~~~~} \phi^{-Nn}(\eta_{[w]}) \in U_-.$$
\end{prop}

\dem Let $\delta \in (0,1)$ and let $w \in F_{\tt n}$ be a nonperipheral element with $\eta_{[w]} \notin V$. By Proposition~\ref{Prop control goodness in f and inverse}, there exists $n_0 \in \NN^*$ such that for every $n \geq n_0$, we have $\mathfrak{g}([f^n(\gamma_w)]) \geq \delta$ or $\mathfrak{g}'([f'^n(\gamma_w')]) \geq \delta$. By Lemma~\ref{Lem conversion goodness closeness PCurr}~$(1)$, there exists $n_1 \geq n_0$ such that for every $n \geq n_1$, we have $$\phi^{Nn}(\eta_{[w]}) \in U_+ \text{~~~~ or ~~~~} \phi^{-Nn}(\eta_{[w]}) \in U_-.$$ This concludes the proof.
\hfill\qedsymbol

\bigskip

The above proposition gives a result of North-South dynamics outside of a neighborhood of $K_{PG}(\phi)$. As $K_{PG}(\phi)$ is empty for a relative expanding outer automorphism by Lemma~\ref{Lem polynomially growing currents empty}~$(1)$, we can now prove Theorem~\ref{Theo North south dynamics relative atoroidal}.

\bigskip

\noindent{\it Proof of Theorem~\ref{Theo North south dynamics relative atoroidal}.} Let $\phi \in \Out(F_{\tt n},\mathcal{F})$ be an expanding outer automorphism relative to $\mathcal{F}$. By Lemma~\ref{Lem polynomially growing currents empty}, we have $K_{PG}(\phi)=\varnothing$. Let $U_+$ be a neighborhood of $\Delta_+(\phi)$ and let $U_-$ be a neighborhood of $\Delta_-(\phi)$. By Proposition~\ref{Prop North south dynamics outside neighborhood}, there exists $N \in \NN^*$ such that for every $n \geq 1$ and every nonperipheral element $w \in F_{\tt n}$, we have $$\phi^{Nn}(\eta_{[w]}) \in U_+\text{~~~~ or ~~~~}\phi^{-Nn}(\eta_{[w]}) \in U_-.$$ Recall that, by Proposition~\ref{Prop density rational currents}, the rational currents are dense in $\PCurr(F_{\tt n},\mathcal{F} \wedge \mathcal{A}(\phi))$. Hence we can apply \cite[Proposition~3.3]{LustigUyanik2019} to see that $\phi^{2N}$ has generalized North-South dynamics. Then, using \cite[Proposition~3.4]{LustigUyanik2019}, we conclude that $\phi$ has generalized North-South dynamics.
\hfill\qedsymbol

\section{North-South dynamics for almost atoroidal relative outer automorphism}\label{Section North SOuth dyn almost}

Let ${\tt n} \geq 3$ and let $\mathcal{F}$ be a free factor system of $F_{\tt n}$. Let $\phi \in \Out(F_{\tt n},\mathcal{F})$ be an almost atoroidal outer automorphism which satisfies Definition~\ref{Defi almost atoroidal outer automorphism}~$(2)$. Let $\mathcal{F} \leq \mathcal{F}_1 \leq \mathcal{F}_2=\{[F_{\tt n}]\}$ be a sequence of free factor system given in this definition. We use the convention of Remark~\ref{Last rmq}. We will show a result of North-South type dynamics for $\phi$ (see~Theorem~\ref{Theo North-South dynamics almost atoroidal}). Note that, if $\mathcal{A}(\phi)\neq \{[F_{\tt n}]\}$ the simplices $\Delta_{\pm}(\phi)$ are still defined. Note that, by Lemma~\ref{Lem polynomially growing currents empty}~$(3)$ and Lemma~\ref{Rmq almost atoroidal}~$(4)$, for every current $\mu \in \Curr(F_{\tt n},\mathcal{F} \wedge \mathcal{A}(\phi))$, we have $\lVert \mu \rVert_{\mathcal{F}_1}>0$. Let $K_{PG}(\phi)$ be the set of polynomially growing currents of $\phi$. Note that, combining Lemma~\ref{Lem Delta value norm} and Lemma~\ref{Rmq almost atoroidal}~$(5)$, we have $K_{PG}(\phi) \cap \Delta_{\pm}(\phi)=\varnothing$. Let $$\widehat{\Delta}_{\pm}(\phi)=\left\{[t\mu +(1-t)\nu]\;\left|\; t \in [0,1], [\mu] \in \Delta_{\pm}(\phi), [\nu] \in K_{PG}(\phi), \lVert \mu \rVert_{\mathcal{F}_1}=\lVert \nu \rVert_{\mathcal{F}_1}=1 \right.\right\}$$ be the \emph{convexes of attraction and repulsion of $\phi$}. 

In order to promote a global North-South type dynamics, we need to construct contracting neighborhoods of $\widehat{\Delta}_+(\phi)$. To this end, following \cite{clay2019atoroidal}, we introduce a notion of \emph{goodness for currents of $\PCurr(F_{\tt n},\mathcal{F} \wedge \mathcal{A}(\phi))$}.

Let $f \colon G \to G$ be as in Remark~\ref{Rmq Convention for relative atoroidal CT map 2}. By Lemma~\ref{Lem exponential length goes to infinity}, let $N \in \NN^*$ be such that, for every edge $e$ of $\overline{G-G_{PG}'}$, we have $\ell_{exp}([f^N(e)]) \geq 4C+1$. Let $C_N=C_{f^N}$ be a constant associated with $f^N$ given by Lemma~\ref{Lem Bounded cancellation lemma}. Let $L>0$ be such that for every path $\gamma$ of $G$ of length at least $L$, we have $\ell([f^N(\gamma)]) \geq C_N+1$. The constant $L$ exists since $f^N$ lifts to a quasi-isometry on the universal cover of $G$. Let $\mathcal{P}_{cs}$ be the finite set of paths of the form $\gamma=\gamma_1 e \gamma_2$, where, for every $i \in \{1,2\}$, the path $\gamma_i$ has length equal to $L$, the path $e$ is an edge in $\overline{G-G_{PG}'}$ and $\gamma_1 e \gamma_2$ is a splitting of $\gamma$. In Lemma~\ref{Lem properties of goodness of currents}~$(2)$, we prove in particular that $\mathcal{P}_{cs}$ is not empty. We will denote by $\widehat{\gamma}$ the edge $e$. 

Let $[\mu] \in \PCurr(F_{\tt n},\mathcal{F} \wedge \mathcal{A}(\phi))$. Recall the definition of $\Psi_0$ just above Definition~\ref{Defi Kpg}. By Lemma~\ref{Lem polynomially growing currents empty}~$(1)$,~$(2)$, we have $\phi(K_{PG}(\phi))=K_{PG}(\phi)$. Hence, for every current $[\mu] \notin K_{PG}(\phi)$, we have $\Psi_0(\phi(\mu))>0$. Thus, for every current $[\mu] \in \PCurr(F_{\tt n},\mathcal{F} \wedge \mathcal{A}(\phi))-K_{PG}(\phi)$, we can define the \emph{completely split goodness $\overline{\mathfrak{g}}(\mu)$ of $\mu$} by
$$\overline{\mathfrak{g}}(\mu)=\frac{1}{\Psi_0(\phi^N(\mu))}\sum\limits_{\gamma \in  \mathcal{P}_{cs}}\left\langle \gamma,\mu \right\rangle.$$
Observe that the function $\overline{\mathfrak{g}}$ is continuous and that it defines a well-defined continuous function $\PCurr(F_{\tt n},\mathcal{F} \wedge \mathcal{A}(\phi))-K_{PG}(\phi)  \to \RR$. 

\begin{lem}\label{Lem properties of goodness of currents}
Let $f \colon G \to G$ be as in Remark~\ref{Rmq Convention for relative atoroidal CT map 2}. 

\noindent{$(1)$ } Let $w \in F_{\tt n}$ be such that $\ell_{exp}(\gamma_w) >0$. We have $\mathfrak{g}([f^N(\gamma_w)]) \geq \overline{\mathfrak{g}}(\eta_{[w]})$.

\medskip

\noindent{$(2)$ } For every $[\mu] \in \Delta_+(\phi)$, we have $\overline{\mathfrak{g}}([\mu])>0$

\end{lem}

\dem $(1)$ The proof of this assertion is similar to the one of \cite[Lemma~4.9~(2)]{clay2019atoroidal}. Let $\gamma \in \mathcal{P}_{cs}$ be such that $\left\langle \gamma,\eta_{[w]} \right\rangle>0$. Then $\gamma \subseteq \gamma_w$. For every occurrence of $\gamma$ in $\gamma_w$, by the choice of $L$, $C_N$ and by Lemma~\ref{Lem Bounded cancellation lemma}, the path $[f^N(\gamma_w)]$ contains $[f^N(\widehat{\gamma})]$, which has exponential length at least equal to $4C_N+1$. Therefore, Lemma~\ref{Lem concatenation of completely split paths and complete edges} implies that the path $[f^N(\gamma_w)]$ contains a subpath of $[f^N(\widehat{\gamma})]$ of exponential length at least $1$ which is a complete factor of $[f^N(\gamma_w)]$ relative to $G_{PG}$. Hence we have: $$\ell_{exp}([f^N(\gamma_w)])\mathfrak{g}([f^N(\gamma_w)]) \geq \sum\limits_{\gamma \in  \mathcal{P}_{cs}}\left\langle \gamma,\eta_{[w]} \right\rangle. $$ By Lemma~\ref{Lem psi0 equal exp length}, we have $$\Psi_0(\phi^N(\eta_{[w]}))=\ell_{exp}([f^N(\gamma_w)])=\Psi_0(\eta_{[\phi^N(w)]})=\ell_{exp}(\gamma_{\phi^N([w])}).$$ Therefore, we have $$\mathfrak{g}([f^N(\gamma_w)]) \geq \overline{\mathfrak{g}}(\eta_{[w]}).$$

\medskip

\noindent{$(2)$ } Let $[\mu] \in \Delta_+(\phi)$. Since $[\mu]$ is a convex combination of extremal points of $\Delta_+(\phi)$ and since for every element $\gamma \in \mathcal{P}_{cs}$, the application $\left\langle \gamma, .\right\rangle$ is linear, it suffices to prove the result for every extremal point of $\Delta_+(\phi)$. So we may suppose that $[\mu]$ is an extremal point of $\Delta_+(\phi)$. Let $G_i$ be the minimal subgraph of $G$ such that $\mathcal{F}(G_i)=\mathcal{F}_1$. Since $[\mu]$ is extremal and since $\phi|_{\mathcal{F}_1}$ is expanding relative to $\mathcal{F}$, by Proposition~\ref{Prop existence relative currents atoroidal automorphisms}, there exists an expanding splitting unit $\sigma$ in $G_i$ whose initial direction is fixed by $f$ and such that, for every path $\gamma \in \mathcal{P}(\mathcal{F}_1 \wedge \mathcal{A}(\phi))$, we have $$\left\langle \gamma, \mu \right\rangle =\mu(C(\gamma))=\lim_{n \to \infty} \frac{\left\langle \gamma,[f^n(\sigma)] \right\rangle }{\ell_{\mathcal{F}_1}([f^n(\sigma)])}.$$ By Lemma~\ref{Rmq almost atoroidal}~$(5)$, since the path $[f^n(\sigma)]$ is contained in $G_i$ and, for every path $\gamma \in \mathcal{P}(\mathcal{F} \wedge \mathcal{A}(\phi))$, the above limit is finite, we have $$\lim_{n \to \infty} \frac{\left\langle \gamma,[f^n(\sigma)] \right\rangle }{\ell_{\mathcal{F}_1}([f^n(\sigma)])}=\lim_{n \to \infty} \frac{\left\langle \gamma,[f^n(\sigma)] \right\rangle }{\ell_{exp}([f^n(\sigma)])}.$$ Hence it suffices to prove that there exists $\gamma \in \mathcal{P}_{cs}$ such that $$\lim_{n \to \infty} \frac{\left\langle \gamma,[f^n(\sigma)] \right\rangle }{\ell_{exp}([f^n(\sigma)])} >0.$$ Let $e$ be an edge of $\overline{G-G_{PG}'}$. Note that, since $\sigma$ is a splitting unit, for every $m \in \NN^*$, the path $[f^m(\sigma)]$ is completely split. Hence an occurrence of $e$ in $\lim_{m \to \infty} [f^m(\sigma)]$ is contained in a splitting unit of $\lim_{m \to \infty} [f^m(\sigma)]$ which is either an INP or is equal to $e$. By Lemma~\ref{Lem NEG INP in Npg} if an INP $\gamma'$ contains $e$, there exists $\gamma_0' \in \mathcal{N}_{PG}$ such that $e \subseteq \gamma_0' \subseteq \gamma'$. For every $m \in \NN^*$, we denote by $N(m,e)$ the number of occurrences of $e$ or $e^{-1}$ in $[f^m(\sigma)]$ which are splitting units of $[f^m(\sigma)]$ and by $EGINP(e)$ the set of all EG INPs containing $e$. Note that, since the set $\mathcal{N}_{PG}$ is finite by Lemma~\ref{Lem Nielsen paths in NPG properties}, so is the limit $$ \lim_{n \to \infty} \sum\limits_{\gamma \in EGINP(e)} \frac{\left\langle \gamma,[f^n(\sigma)] \right\rangle }{\ell_{exp}([f^n(\sigma)])}.$$ Since for every $m \in \NN^*$, we have $$\left\langle e,[f^m(\sigma)] \right\rangle = N(m,e)+\sum\limits_{\gamma \in EGINP(e)}\left\langle \gamma,[f^n(\sigma)] \right\rangle,$$ we see that the limit $$\lim_{m \to \infty} \frac{N(m,e)}{\ell_{exp}([f^m(\sigma)])}$$ exists. We claim that there exists an edge $e$ of $\overline{G-G_{PG}'}$ such that $$\lim_{m \to \infty} \frac{N(m,e)}{\ell_{exp}([f^m(\sigma)])}>0.$$ Indeed, note that, by Lemma~\ref{Lem splitting units positive exp length} since $[f^m(\sigma)]$ is $PG$-relative completely split, we have $$\ell_{exp}([f^m(\sigma)])=\sum_{e \in \vec{E}\left(\overline{G-G_{PG}'}\right)} N(m,e).$$ Hence $$\sum_{e \in \vec{E}\left(\overline{G-G_{PG}'}\right)} \lim_{m \to \infty} \frac{N(m,e)}{\ell_{exp}([f^m(\sigma)])}=1,$$ which implies the claim. Let $e_0$ be an edge of $\overline{G-G_{PG}'}$ which satisfies the claim. Since, for every $m \in \NN^*$, the path $[f^m(\sigma)]$ is completely split, if an occurrence of $e_0$ or $e_0^{-1}$ in $[f^m(\sigma)]$ is a splitting unit and if $\gamma$ is a path in $[f^m(\sigma)]$ of the form $\gamma=\gamma_1 e_0\gamma_2$ or $\gamma=\gamma_1 e_0^{-1}\gamma_2$, then such a decomposition of $\gamma$ is a splitting of $\gamma$. Thus, if $\ell(\gamma_1)=\ell(\gamma_2)=L$, then the path $\gamma$ is in $\mathcal{P}_{cs}$ and it contains the occurrence of $e_0$. Hence since $\mu=\mu(\sigma)$, we have $$\lim_{m \to \infty} \frac{N(m,e)}{\ell_{exp}([f^m(\sigma)])}=\sum_{\gamma \in\mathcal{P}_{cs}, e_0 \subseteq \gamma} \left\langle \gamma, \mu \right\rangle>0.$$ Therefore, there exists $\gamma \in \mathcal{P}_{cs}$ such that $\left\langle \gamma, \mu \right\rangle>0$ and $\overline{\mathfrak{g}}([\mu])>0$.
\hfill\qedsymbol

\begin{lem}\label{Lem Contracting neigh for exponential simplex}
Let $f \colon G \to G$ be as in Remark~\ref{Rmq Convention for relative atoroidal CT map 2}. Let $U_{\pm}$ be open neighborhoods of $\Delta_{\pm}(\phi)$. There exist open neighborhoods $U_{\pm}' \subseteq U_{\pm}$ of $\Delta_{\pm}(\phi)$ such that $\phi^{\pm 1}(U_{\pm}') \subseteq U_{\pm}'$.
\end{lem}

\dem The proof is similar to the one of \cite[Lemma~4.13]{clay2019atoroidal}. We prove the result for $\Delta_+(\phi)$, the proof for $\Delta_-(\phi)$ being symmetric. By Lemma~\ref{Lem properties of goodness of currents}~$(2)$, for every $[\mu] \in \Delta_+(\phi)$, we have $\overline{\mathfrak{g}}([\mu])>0$. By compactness of $\Delta_+(\phi)$ and continuity of $\overline{\mathfrak{g}}$, there exists $\delta_0 >0$ such that, for every $\mu \in \Delta_+(\phi)$, we have $\overline{\mathfrak{g}}(\mu) \geq \delta_0$.
Let $\delta \in (0,\delta_0)$. Let $U_+$ be a neighborhood of $\Delta_+(\phi)$. Since the function $\overline{\mathfrak{g}}$ is continuous, there exists an open neighborhood $U_+^0 \subseteq U_+$ of $\Delta_+(\phi)$ such that, for every $[\mu] \in U_+^0$, we have $\overline{\mathfrak{g}}([\mu]) > \delta$. Up to taking a smaller $U_+^0$, we may suppose that $K_{PG}(\phi) \cap U_+^0=\varnothing$ (this is possible since $K_{PG}(\phi)$ is compact and $\Delta_+(\phi) \cap K_{PG}(\phi)=\varnothing$). In particular, by Lemma~\ref{Lem psi0 equal exp length}, for every nonperipheral element $w \in F_{\tt n}$ such that $\eta_{[w]} \in U_0^+$, we have $\ell_{exp}(\gamma_w)>0$.

Let $w \in F_{\tt n}$ be a nonperipheral element such that $\eta_{[w]} \in U_0^+$. By Lemma~\ref{Lem properties of goodness of currents}~$(1)$, we have
$$
\mathfrak{g}([f^N(\gamma_w)]) \geq \overline{\mathfrak{g}}(\eta_{[w]}) >\delta.$$

By Lemma~\ref{Lem conversion goodness closeness PCurr}~$(1)$, there exists $M \geq N$ such that, for every $w \in F_{\tt n}$ such that $\eta_{[w]} \in U_+^0$, we have $\phi^M([\eta_{[w]}]) \in U_+^0$. Let $$U_+'=\bigcap_{i=0}^{M-1} \phi^i(U_+^0).$$ Since $\phi(\Delta_+(\phi))=\Delta_+(\phi)$ by Proposition~\ref{Prop invariance 3}, the set $U_+'$ is an open neighborhood of $\Delta_+(\phi)$ which is stable by $\phi$ by density of rational currents (see~Proposition~\ref{Prop density rational currents}) and continuity of $\phi$. This concludes the proof.
\hfill\qedsymbol

\begin{lem}\label{Lem contractible neighborhood for convex}
Let $f \colon G \to G$ be as in Remark~\ref{Rmq Convention for relative atoroidal CT map 2}. Suppose that the outer automorphism $\phi$ is of type~$(2)$ in Definition~\ref{Defi almost atoroidal outer automorphism}. Let $\mathcal{F} \leq \mathcal{F}_1 \leq \mathcal{F}_2=\{F_{\tt n}\}$ be as in the beginning of Section~\ref{Section North SOuth dyn almost}. Let $i \in \{1,\ldots,k-1\}$ be such that $\mathcal{F}(G_i)=\mathcal{F}_1$. Let $\widehat{V}_{\pm}$ be open neighborhoods of $\widehat{\Delta}_{\pm}(\phi)$. There exist open neighborhoods $\widehat{V}_{\pm}'$ of $\widehat{\Delta}_{\pm}(\phi)$ contained in $\widehat{V}_{\pm}$ such that $\phi^{\pm}(\widehat{V}_{\pm}') \subseteq \widehat{V}_{\pm}'$.
\end{lem}

\dem The proof follows \cite[Lemma~4.14]{clay2019atoroidal}. We prove the result for $\widehat{\Delta}_+(\phi)$, the proof for $\widehat{\Delta}_-(\phi)$ being symmetric. Given $[\mu] \in \PCurr(F_{\tt n},\mathcal{F} \wedge \mathcal{A}(\phi))-K_{PG}(\phi)$, a finite set of reduced edge paths $\mathcal{P}$ in $G$ and some $\epsilon >0$ determine an open neighborhood $N([\mu],\mathcal{P},\epsilon)$ of $[\mu]$ in $\PCurr(F_{\tt n},\mathcal{F}\wedge \mathcal{A}(\phi))-K_{PG}(\phi)$ as follows:

$$
N([\mu],\mathcal{P},\epsilon)=\left\{[\nu] \in \PCurr(F_{\tt n},\mathcal{F} \wedge \mathcal{A}(\phi))-K_{PG}(\phi) \;\left|\; \forall \gamma \in \mathcal{P}, \; \left| \frac{\left\langle \gamma,\nu \right\rangle}{\Psi_0(\nu)}-\frac{\left\langle \gamma,\mu \right\rangle}{\Psi_0(\mu)}\right|< \epsilon \right\}\right..
$$

Since $K_{PG}(\phi)$ is compact, if $\epsilon$ is small enough, this defines an open neighborhood of $[\mu]$ in $\PCurr(F_{\tt n},\mathcal{F} \wedge \mathcal{A}(\phi))$. For a subset $X \subseteq \PCurr(F_{\tt n},\mathcal{F}\wedge \mathcal{A}(\phi))-K_{PG}(\phi)$, let $$N(X,\mathcal{P},\epsilon)=\bigcup_{[\mu] \in X} N([\mu],\mathcal{P},\epsilon).$$

For $L >0$, let $\mathcal{P}_+(L)$ be the set of reduced edge paths in $G_i$ of length at most equal to $L$ which are not contained in any concatenation of paths in $G_{PG,\mathcal{F}_1}$ and $\mathcal{N}_{PG,\mathcal{F}_1}$. By Lemma~\ref{Rmq almost atoroidal}~$(3)$, the set $\mathcal{P}_+(L)$ is also the set of reduced edge paths in $G_i$ of length at most equal to $L$ which are not contained in any concatenation of paths in $G_{PG}$ and $\mathcal{N}_{PG}$.
Let $[\mu] \in \Delta_+(\phi)$ and let $t \in [0,1]$. Let
$$
K_{PG}([\mu],t)=\{[(1-t)\nu+t\mu]\;|\; [\nu] \in K_{PG}(\phi), \; \lVert \nu \rVert_{\mathcal{F}_1}=\lVert \mu \rVert_{\mathcal{F}_1}=1\}.
$$
Remark that $$\widehat{\Delta}_+(\phi)=\bigcup_{[\mu] \in \Delta_+(\phi),\; t \in [0,1]} K_{PG}([\mu],t).$$ 
Let $\epsilon >0$. Let $V_{poly}(\epsilon)=[\Psi_0^{-1}((-\epsilon,\epsilon))]$. It is clear, by the continuity of $\Psi_0$ and the definition~\ref{Defi Kpg} of $K_{PG}(\phi)$, that $\bigcap_{\epsilon > 0} V_{poly}(\epsilon)=K_{PG}(\phi)$. Let $t \in (0,1]$ and $[\mu] \in \Delta_+(\phi)$ and let $\mu$ be such that $\lVert \mu \rVert_{\mathcal{F}_1}=1$. By Lemma~\ref{Rmq almost atoroidal}~$(5)$, we have $\Psi_0(\mu)=1$. Let 
$$V_{poly}([\mu],t,\epsilon)=\left\{[\nu] \in \PCurr(F_{\tt n},\mathcal{F} \wedge \mathcal{A}(\phi)) \;\left|\; \begin{array}{c} \lVert \nu \rVert_{\mathcal{F}_1}=\lVert \mu \rVert_{\mathcal{F}_1}=1, \\
 t(1+\epsilon) > \Psi_{0}(\nu) >t(1-\epsilon)\end{array} \right.\right\}.$$ Note that, since $\Psi_0(\mu)=1$, we have $[\nu] \in V_{poly}([\mu],t,\epsilon)$ if for $[\nu]$ such that $\lVert \nu \rVert_{\mathcal{F}_1}=1$, we have $$t\Psi_{0}(\mu)(1+\epsilon) > \Psi_{0}(\nu) >t\Psi_{0}(\mu)(1-\epsilon).$$
Let 
$$V_{\infty}([\mu],t)=\bigcap_{L \to \infty, \epsilon \to 0} N(K_{PG}([\mu],t),\mathcal{P}_+(L),\epsilon) \cap V_{poly}([\mu],t,\epsilon).
$$

\noindent{\bf Claim. } For every $[\mu ] \in \Delta_+(\phi)$ and every $t \in (0,1]$, we have $V_{\infty}([\mu],t])=K_{PG}([\mu],t)$.

\medskip

\dem The inclusion $K_{PG}([\mu],t) \subseteq V_{\infty}([\mu],t])$ being immediate since $\Psi_0$ is linear and vanishes on $K_{PG}(\phi)$, we prove the converse inclusion. Let $\nu \in V_{\infty}([\mu],t)$. By definition~\ref{Defi attractive convex} of $\Delta_+(\phi)$, for every $[\mu'] \in \Delta_+(\phi)$ and for every reduced edge path $\gamma$ not contained in $G_{i}$, we have $\left\langle \gamma, \mu' \right\rangle=0$. Hence, by Lemma~\ref{Rmq almost atoroidal}~$(4)$, the current $[\mu]$ is entirely determined by the cylinder sets determined by reduced edge paths contained in $G_i$ which are not contained in concatenation of paths in $G_{PG,\mathcal{F}_1}$ and $\mathcal{N}_{PG,\mathcal{F}_1}$. By Lemma~\ref{Rmq almost atoroidal}~$(3)$, the current $[\mu]$ is entirely determined by the cylinder sets determined by reduced edge paths contained in $G_i$ which are not contained in concatenation of paths in $G_{PG}$ and $\mathcal{N}_{PG}$. Let $\gamma$ be a reduced edge path which is contained in $G_i$ and which is not contained in a concatenation of paths in $G_{PG}$ and $\mathcal{N}_{PG}$. By Lemma~\ref{Lem polynomially growing currents empty}, for every projective current $[\nu'] \in K_{PG}(\phi)$, the support of $\nu'$ is contained in $\partial^2\mathcal{A}(\phi)$. By Proposition~\ref{Prop circuits in Gpg are elements in poly subgroup}, if $g \in F_{\tt n}$ is such that there exists a subgroup $A$ of $F_{\tt n}$ such that $[A] \in \mathcal{A}(\phi)$ and $g \in A$, then $\gamma_g$ is a concatenation of paths in $G_{PG}$ and $\mathcal{N}_{PG}$. In particular, if $\gamma'$ is a path of $G$ such that $\{g^{+\infty},g^{-\infty}\} \in C(\gamma')$, then $\gamma'$ is contained in a concatenation of paths in $G_{PG}$ and in $\mathcal{N}_{PG}$. In particular, since $\gamma$ is not contained in a concatenation of paths in $G_{PG}$ and in $\mathcal{N}_{PG}$, for every projective current $[\nu'] \in K_{PG}(\phi)$, we have $\left\langle \gamma,\nu' \right\rangle=0$. 

Suppose that $\lVert \nu \rVert_{\mathcal{F}_1}=\lVert \mu \rVert_{\mathcal{F}_1}=1$. By Lemma~\ref{Rmq almost atoroidal}~$(5)$, we also have $\Psi_0(\mu) =1$. There exists $\lambda>0$ such that for every path $\gamma$ which is contained in $G_i$ and which is not contained in a concatenation of paths in $G_{PG}$ and $\mathcal{N}_{PG}$, we have $\left\langle \gamma, \nu \right\rangle=\left\langle \gamma, \lambda t\mu \right\rangle$. We claim that $\nu-\lambda t\mu \in \Curr(F_{\tt n},\mathcal{F} \wedge \mathcal{A}(\phi))$ and that $[\nu-\lambda t\mu] \in K_{PG}(\phi)$. Indeed, for the first part, it suffices to show that for every path $\gamma \in \mathcal{P}(\mathcal{F}_1 \wedge \mathcal{A}(\phi))$, we have $(\nu-\lambda t \mu)(C(\gamma)) \geq 0$. This follows from the fact that, for every path $\gamma \in \mathcal{P}(\mathcal{F}_1 \wedge \mathcal{A}(\phi))$ such that $\gamma \subseteq G_i$, the path $\gamma$ is not contained in a concatenation of paths in $G_{PG}$ and in $\mathcal{N}_{PG}$. Hence we have $\left\langle \gamma, \nu \right\rangle=\left\langle \gamma, \lambda t\mu \right\rangle$. Moreover, if $\gamma \in \mathcal{P}(\mathcal{F}_1 \wedge \mathcal{A}(\phi))$, then we have $\mu(C(\gamma))=0$. This shows that $\nu-\lambda t\mu \in \Curr(F_{\tt n},\mathcal{F} \wedge \mathcal{A}(\phi))$. 

We now prove that $[\nu-\lambda t\mu] \in K_{PG}(\phi)$. Otherwise, by Lemma~\ref{Lem polynomially growing currents empty}, the support of $\nu-\lambda t\mu$ is not contained in $\partial^2\mathcal{A}(\phi)$. By Proposition~\ref{Prop circuits in Gpg are elements in poly subgroup}, there exists a path $\gamma$ which is not contained in a concatenation of paths in $G_{PG}$ and in $\mathcal{N}_{PG}$ such that $$\left\langle \gamma, \nu-\lambda t\mu \right\rangle >0.$$ Consider a decomposition of $\gamma=a_1b_1\ldots a_k b_k$ where, for every $j \in \{1,\ldots,k\}$, the path $a_j$ is contained in $\overline{G-G_i}$ and, for every $j \in \{1,\ldots,k\}$, the path $b_j$ is contained in $G_i$ with $a_1$ and $b_k$ possibly empty. By Lemma~\ref{Rmq almost atoroidal}~$(1)$,~$(2)$ and Remark~\ref{Last rmq}, up to taking a larger path $\gamma$, we may suppose that $b_1$ is nontrivial. By Lemma~\ref{Rmq almost atoroidal}~$(2)$ and Remark~\ref{Last rmq}, for every $j \in \{1,\ldots,k\}$, the path $a_j$ is contained in $G_{PG}$. Since $\gamma$ is not contained in a concatenation of paths in $G_{PG}$ and $\mathcal{N}_{PG}$, there exists $j \in \{1,\ldots,k\}$ such that $b_j$ is not contained in a concatenation of paths in $G_{PG}$ and $\mathcal{N}_{PG}$. But then $\left\langle b_j, \nu \right\rangle=\left\langle b_j, \lambda t\mu \right\rangle$, that is $\left\langle b_j, \nu-\lambda t\mu \right\rangle=0$. By additivity of $\nu-\lambda t\mu$, we have $$\left\langle \gamma, \nu-\lambda t\mu \right\rangle \leq \left\langle b_j, \nu-\lambda t\mu \right\rangle=0.$$ This contradicts the choice of $\gamma$. Hence $[\nu-\lambda t\mu] \in K_{PG}(\phi)$. Therefore, we have $\Psi_0(\nu-\lambda t\mu)=0$. Since $[\nu]\in V_{\infty}([\mu],t)$ and since $\lVert \nu \rVert_{\mathcal{F}_1}=\lVert \mu \rVert_{\mathcal{F}_1}=1$, we see that $$\Psi_{0}(\nu)=t\Psi_{0}(\mu).$$ By linearity of $\Psi_{0}$ and the fact that $\Psi_0(\mu)=1$, we have $$t=t\Psi_{0}(\mu)=\Psi_{0}(\nu)=\lambda t\Psi_{0}(\mu)=\lambda t.$$ Since $t >0$ and $\Psi_{0}(\mu)=1$, we have $\lambda=1$. Suppose first that $t \neq 1$. Let $\nu'=\frac{1}{1-t}(\nu-t\mu)$, so that $[\nu'] \in K_{PG}(\phi)$ and $\lVert \nu' \rVert_{\mathcal{F}}=1$. Then $[\nu]=[(1-t)\nu'+t\mu] \in K_{PG}([\mu],t)$. Thus, we have $V_{\infty}([\mu],t)=K_{PG}([\mu],t)$. 

Suppose now that $t=1$. Then $\Psi_0(\nu)=1=\lVert \nu \rVert_{\mathcal{F}}$. We claim that if $\gamma \in \mathcal{P}(\mathcal{F}_1 \wedge \mathcal{A}(\phi))$ is such that $\nu(C(\gamma))>0$, then $\gamma \subseteq G_i$. Indeed, otherwise there would exist an edge $e$ contained in $\overline{G-G_i}$ such that $\nu(C(e))>0$. By the description of $\overline{G-G_i}$ given in Lemma~\ref{Rmq almost atoroidal}~$(1)$,~$(2)$ and additivity of the current $\nu$, we can choose the edge $e \in \overline{G-G_i}$ in such a way that $e \in G_{PG}$. This would imply that $\lVert \nu \rVert_{\mathcal{F}_1} >\Psi_0(\nu)$, a contradiction. The claim follows. But, since for every path $\gamma \in \mathcal{P}(\mathcal{F}_1 \wedge \mathcal{A}(\phi))$ such that $\gamma \subseteq G_i$, we have $\nu(C(\gamma))=\mu(C(\gamma))$, we see that $\nu=\mu$ and that $\nu \in K_{PG}([\mu],1)$. This concludes the proof of the claim.
\hfill\qedsymbol

\medskip

Since $\widehat{\Delta}_+(\phi)$ is compact, there exist $L>0$ and $\epsilon >0$ such that, for every $[\mu] \in \Delta_+(\phi)$ and every $t \in (0,1]$, we have $$V([\mu],t,L,\epsilon)=N(K_{PG}([\mu],t),\mathcal{P}_+(L),\epsilon) \cap V_{poly}(t,\epsilon) \subseteq \widehat{V}_+.$$
When $t=0$, there exists $\epsilon>0$ such that $V_{poly}(\epsilon) \subseteq \widehat{V}_+$.
Let $s \in (0,1)$, and let $V$ be an open neighborhood of $K_{PG}(\phi)$ such that, for every $[\nu] \in V$ with $\lVert \nu \rVert_{\mathcal{F}_1}=1$, we have: 
\begin{equation}\label{Equation p85}
\Psi_{0}(\nu)<s.
\end{equation} 
For every $[\mu] \in \left(N(\widehat{\Delta}_+(\phi),\widehat{\mathcal{P}}_+(L),\epsilon)-V\right) \cap \widehat{\Delta}_+(\phi)$, there exist $[\mu_{poly}] \in K_{PG}(\phi)$, $[\mu_{exp}] \in \Delta_+(\phi)$ and $t \in (0,1]$ such that $$[\mu]=[t\mu_{exp}+(1-t)\mu_{poly}].$$

By Lemma~\ref{Lem properties of goodness of currents}~$(2)$, for every $[\mu] \in \Delta_+(\phi)$, we have $\overline{\mathfrak{g}}([\mu])>0$. By compactness of $\Delta_+(\phi)$ and continuity of $\overline{\mathfrak{g}}$, there exists $\delta_1 >0$ such that, for every $\mu \in \Delta_+(\phi)$, we have $\overline{\mathfrak{g}}(\mu) \geq \delta_1$. Since $\overline{N(\widehat{\Delta}_+(\phi),\widehat{\mathcal{P}}_+(L),\epsilon)-V} \cap \widehat{\Delta}_+(\phi)$ is compact, and since the function $\overline{\mathfrak{g}}$ is continuous, there exists $\delta_0' >0$ such that the set $U=\overline{\mathfrak{g}}^{-1}((\delta_0',+\infty))$ is an open neighborhood of $(N(\widehat{\Delta}_+(\phi),\widehat{\mathcal{P}}_+(L),\epsilon)-V) \cap \widehat{\Delta}_+(\phi)$ intersecting $V$. Note that $U \cap K_{PG}(\phi)=\varnothing$.
We set $$\widehat{V}_+'=\left( \bigcup_{[\mu] \in \Delta_+(\phi),\; t \in (0,1]} V([\mu],t,L,\epsilon) \cup V_{poly}(\epsilon) \right) \cap \left(U \cup V \right).$$ 
Let $\delta_0$ and $M_0$ be the constants given by Lemma~\ref{Lem conversion goodness closeness PCurr}~$(2)$ for the above choice of $\epsilon>0$ and $L>0$. By replacing $\delta_0$ with a smaller constant and $M_0$ with a larger one, we may suppose that $\delta_0$ and $M_0$ also satisfy the conclusion of Lemma~\ref{Lem conversion goodness closeness PCurr}~$(1)$ for $U$ as well (where the open neighborhood $W$ of $K_{PG}(\phi)$ needed in Lemma~\ref{Lem conversion goodness closeness PCurr}~$(1)$ is such that $W \subseteq V-U$). 

\medskip

\noindent{\bf Claim~2 } There exists $N \in \NN^*$ such that $\phi^N(\widehat{V}_+') \subseteq \widehat{V}_+'$.

\medskip

\dem Let $w \in F_{\tt n}$ be a nonperipheral element such that $\eta_{[w]} \in \widehat{V}_+'$. Suppose first that $\eta_{[w]} \in U \cap \widehat{V}_+'$. Since $\eta_{[w]} \notin K_{PG}(\phi)$, by Lemma~\ref{Lem psi0 equal exp length}, we have $\ell_{exp}(\gamma_w)>0$. By Lemma~\ref{Lem properties of goodness of currents}~$(1)$, we have:
$$
\mathfrak{g}([f^N(\gamma_w)]) \geq \overline{\mathfrak{g}}(\eta_{[w]}) >\delta_0'.$$

By Lemma~\ref{Lem conversion goodness closeness PCurr}~$(1)$, there exists $M \geq M_0+N$ such that, for every $w \in F_{\tt n}$ such that $\eta_{[w]} \in U \cap \widehat{V}_+'$ and every $n \geq 1$, we have $\phi^{Mn}([\eta_{[w]}]) \in U \cap \widehat{V}_+' \subseteq \widehat{V}_+'$.

Suppose now that $\eta_{[w]} \in V \cap \widehat{V}_+'$. By Lemma~\ref{Lem polynomially growing currents empty}~$(3)$ and Lemma~\ref{Rmq almost atoroidal}~$(4)$ for every projective current $[\mu] \in \PCurr(F_{\tt n},\mathcal{F} \wedge \mathcal{A}(\phi))$, we have $\lVert \mu \rVert_{\mathcal{F}_1}>0$. For a projective current $[\mu] \in \PCurr(F_{\tt n},\mathcal{F} \wedge \mathcal{A}(\phi))$, let $$\Psi_{\mathcal{F}_1}([\mu])=\frac{\Psi_0(\mu)}{\lVert \mu \rVert_{\mathcal{F}_1}}.$$ Then, by definition of $V$ and by Lemma~\ref{Lem psi0 equal exp length}, we have 
$$\Psi_{\mathcal{F}_1}([\eta_{[w]}])=\frac{\ell_{exp}(\gamma_w)}{\ell_{\mathcal{F}_1}(\gamma_w)} <s.$$ If $[\eta_{[w]}] \in K_{PG}(\phi)$, then since $\phi(K_{PG}(\phi))=K_{PG}(\phi)$, we are done. Therefore, we may suppose that $[\eta_{[w]}] \notin K_{PG}(\phi)$ and, by Lemma~\ref{Lem psi0 equal exp length}, for every $n \in \NN^*$, we have $\ell_{exp}([f^n(\gamma_w)]) \geq 1$. Let $R >1$ be such that $\frac{1}{1+\frac{R(1-\delta_0)}{10C}(1-s)} \leq \epsilon$. By Lemma~\ref{Lem control goodness versus decrease bad length}, one of the following assertion holds:

\medskip

\noindent{$(1)$ } $\mathfrak{g}([f^M(\gamma_w)]) > \delta_0$,

\medskip

\noindent{$(2)$ } $\ell_{exp}([f^M(\gamma_w)]) < \frac{10C}{(1-\delta_0)R}\ell_{exp}(\gamma_w)$.

\medskip

First assume that Assertion~$(1)$ holds. Let $[\mu_{[\phi^M([w])]}] \in \Delta_+(\phi)$ be the projective current associated with $[\phi^M([w])]$ given by Lemma~\ref{Lem conversion goodness closeness PCurr}~$(2)$. Let $$t=\Psi_{\mathcal{F}_1}([\eta_{[\phi^M([w])]}]).$$ We claim that $[\eta_{[\phi^M([w])]}]\in V([\mu_{[\phi^M([w])]}],t,L,\epsilon)$. Indeed, we clearly have $$[\eta_{\phi^M([w])}] \in V_{poly}([\mu_{[\phi^M([w])]}],t,\epsilon).$$ By Lemma~\ref{Lem conversion goodness closeness PCurr}~$(2)$, for every reduced edge path $\gamma \in \mathcal{P}_+(L)$, we have
$$\left| \frac{\left\langle \gamma, \eta_{[\phi^M([w])]} \right\rangle}{\Psi_0(\eta_{[\phi^M([w])]})}-\frac{\left\langle \gamma, \mu_{[\phi^M([w])]} \right\rangle}{\Psi_0( \mu_{[\phi^M([w])]})} \right| < \epsilon.$$ Therefore we have $[\eta_{[\phi^M([w])]}] \in N(K_{PG}([\mu_{[\phi^M([w])]}],t),\mathcal{P}_+(L),\epsilon)$. The claim follows by Equation~\eqref{Equation p85}. By definition of $\widehat{V}_+'$, we see that $\phi^M([\eta_{[w]}])=[\eta_{[\phi^M([w])]}] \in \widehat{V}_+'$.

Suppose now that Assertion~$(2)$ holds. We claim that $[\eta_{[\phi^M([w])]}] \in V_{poly}(\epsilon)$. By Lemma~\ref{Rmq almost atoroidal}~$(1)$,$(2)$ and Remark~\ref{Last rmq}, the graph $\overline{G-G_i}$ consists in edges in $G_{PG}$. By Lemma~\ref{Rmq almost atoroidal}~$(6)$, we have $$\ell_{\mathcal{F}_1}([f^M(\gamma_w)])-\ell_{exp}([f^M(\gamma_w)])=\ell_{\mathcal{F}_1}(\gamma_w)-\ell_{exp}(\gamma_w).$$ Hence we have

$$
\begin{array}{ccl}
\Psi_{\mathcal{F}_1}([\eta_{[\phi^M(\gamma_w)]}])=\frac{\ell_{exp}([f^M(\gamma_w)])}{\ell_{\mathcal{F}_1}([f^M(\gamma_w)])} & = &  \frac{\ell_{exp}([f^M(\gamma_w)])}{\ell_{exp}([f^M(\gamma_w)])+\ell_{\mathcal{F}_1}([f^M(\gamma_w)])-\ell_{exp}([f^M(\gamma_w)])} \\

{} & = &  \frac{\ell_{exp}([f^M(\gamma_w)])}{\ell_{exp}([f^M(\gamma_w)])+\ell_{\mathcal{F}_1}(\gamma_w)-\ell_{exp}(\gamma_w)} \\

{} & = & \frac{1}{1+\frac{\ell_{\mathcal{F}_1}(\gamma_w)-\ell_{exp}(\gamma_w)}{\ell_{exp}([f^M(\gamma_w)])}} \leq \frac{1}{1+\frac{R(1-\delta_0)}{10C}\frac{\ell_{\mathcal{F}_1}(\gamma_w)-\ell_{exp}(\gamma_w)}{\ell_{exp}(\gamma_w)}} \\

{} & \leq &  \frac{1}{1+\frac{R(1-\delta_0)}{10C}\frac{\ell_{\mathcal{F}_1}(\gamma_w)-\ell_{exp}(\gamma_w)}{\ell_{\mathcal{F}_1}(\gamma_w)}} \leq \frac{1}{1+\frac{R(1-\delta_0)}{10C}(1-s)} \leq \epsilon.
\end{array}
$$

Note that $\psi_{\mathcal{F}_1}^{-1}((0,\epsilon)) \subseteq V_{poly}(\epsilon)$. Thus, we have $\Phi^M([\eta_{[w]}])=[\eta_{[\phi^M([w])]} ]\in V_{poly}(\epsilon) \subseteq \widehat{V}_+'$. Therefore, by density of the rational currents (see~Proposition~\ref{Prop density rational currents}) and continuity of $\phi$, we have $\phi^M(\widehat{V}_+') \subseteq \widehat{V}_+'$. This proves Claim~2.
\hfill\qedsymbol

\bigskip

Let $$\widehat{V}_+''=\bigcap_{i=0}^{M-1} \phi^i(\widehat{V}_+').$$ Since $\phi(\widehat{\Delta}_+(\phi))=\widehat{\Delta}_+(\phi)$, the set $\widehat{V}_+''$ is an open neighborhood of $\widehat{\Delta}_+(\phi)$ which is stable by $\phi$ by construction. This concludes the proof.
\hfill\qedsymbol

\begin{theo}\label{Theo North-South dynamics almost atoroidal}
Let ${\tt n} \geq 3$ and let $\mathcal{F}$ be a free factor system of $F_{\tt n}$. Let $\mathcal{F} \leq \mathcal{F}_1 \leq \mathcal{F}_2$ be a sequence of free factor systems such that the extension $\mathcal{F}_1 \leq \mathcal{F}_2$ is sporadic. Let $\phi \in \Out(F_{\tt n},\mathcal{F})$ be such that $\phi$ preserves $\mathcal{F} \leq \mathcal{F}_1 \leq \mathcal{F}_2$ and $\phi|_{\mathcal{F}_1}$ is an expanding automorphism relative to $\mathcal{F}$. Let $\widehat{\Delta}_{\pm}(\phi)$ be the convexes of attraction and repulsion of $\phi$ and $\Delta_{\pm}(\phi)$ be the simplices of attraction and repulsion of $\phi$. Let $U_{\pm}$ be open neighborhoods of $\Delta_{\pm}(\phi)$ in $\PCurr(F_{\tt n},\mathcal{F} \wedge \mathcal{A}(\phi))$ and $\widehat{V}_{\pm}$ be open neighborhoods of $\widehat{\Delta}_{\pm}(\phi)$ in $\PCurr(F_{\tt n},\mathcal{F} \wedge \mathcal{A}(\phi))$. There exists $M \in \NN^*$ such that for every $n \geq M$, we have 
$$
\phi^{\pm n}(\PCurr(F_{\tt n},\mathcal{F} \wedge \mathcal{A}(\phi))-\widehat{V}_{\mp}) \subseteq U_{\pm}.
$$
\end{theo}

\dem The proof is similar to \cite[Theorem~4.15]{clay2019atoroidal}. We replace $\phi$ by a power so that $\phi$ satisfies Remark~\ref{Rmq Convention for relative atoroidal CT map 2}. By Lemmas~\ref{Lem Contracting neigh for exponential simplex} and~\ref{Lem contractible neighborhood for convex}, we may suppose that $\phi(U_+) \subseteq U_+$ and that $\phi(\widehat{V}_+) \subseteq \widehat{V}_+$. Let $M$ be the exponent given by Proposition~\ref{Prop North south dynamics outside neighborhood} by using $U_+=U_+$ and $U_-=V=\widehat{V}_-$. For every current $[\mu] \in \PCurr(F_{\tt n},\mathcal{F} \wedge \mathcal{A}(\phi))-\phi^M(\widehat{V}_{\mp})$, we have $\phi^M([\mu]) \in U_+$ since $\phi^{-M}([\mu]) \notin \widehat{V}_-$. Therefore, for every $[\mu] \in \PCurr(F_{\tt n},\mathcal{F} \wedge \mathcal{A}(\phi))-\widehat{V}_-$, we have $\phi^{2M}([\mu]) \in U_+$ and for every $n \geq M$, we have $\phi^{2n}([\mu]) \in U_+$ since $\phi(U_+) \subseteq U_+$. Therefore for every $n \geq M$, we see that $$\phi^{2n}(\PCurr(F_{\tt n},\mathcal{F} \wedge \mathcal{A}(\phi))-\widehat{V}_-) \subseteq U_+.$$ A symmetric argument for $\phi^{-1}$ shows that $\phi^2$ acts with generalized North-South dynamics. By \cite[Proposition~3.4]{LustigUyanik2019}, we see that $\phi$ acts with generalized North-South dynamics. This concludes the proof.
\hfill\qedsymbol

\begin{coro}\label{Coro p length iterates currents}
For every open neighborhood $\widehat{V}_{-} \subseteq \PCurr(F_{\tt n},\mathcal{F} \wedge \mathcal{A}(\phi))$ of $\widehat{\Delta}_-$, there exist $M \in \NN^*$ and a constant $L_0$ such that, for every current $[\mu] \in \PCurr(F_{\tt n},\mathcal{F} \wedge \mathcal{A}(\phi))-\widehat{V}_-$, and every $m \geq M$, we have $$\lVert\phi^m(\mu) \rVert_{\mathcal{F}} \geq 3^{m-M}L_0 \lVert\mu \rVert_{\mathcal{F}}.$$
\end{coro}

\dem Let $f \colon G \to G$ be as in Remark~\ref{Rmq Convention for relative atoroidal CT map 2}. By Lemma~\ref{Lem properties of goodness of currents}~$(2)$, there exist a constant $\delta>0$ and an open neighborhood $U$ of $\Delta_+(\phi)$ such that, for every projective current $[\mu] \in U$, we have $\overline{\mathfrak{g}}([\mu]) \geq \delta$. We first prove Corollary~\ref{Coro p length iterates currents} for currents $[\mu] \in U$. By Proposition~\ref{Prop density rational currents}, it suffices to prove the result for rational currents. By Lemma~\ref{Lem properties of goodness of currents}~$(1)$, since $U \cap K_{PG}(\phi)=\varnothing$, for every element $w \in F_{\tt n}$ such that $[\eta_{[w]}] \in U$, we have $\mathfrak{g}([f^N(\gamma_w)]) \geq \delta$. By Lemma~\ref{Lem control of the goodness}~$(1)$ and Lemma~\ref{Lem compute exponential length optimal factor}, there exists a constant $K_1>0$ depending on $\delta$ such that for every $m \geq N$ and for every element $w \in F_{\tt n}$ such that $[\eta_w] \in U$, we have $\ell_{exp}([f^m(\gamma_w)]) \geq TEL(m-N,[f^N(\gamma_w)]) \geq 3^{m-N}K_1\ell_{exp}([f^N(\gamma_w)])$. Since $\PCurr(F_{\tt n},\mathcal{F} \wedge \mathcal{A}(\phi))-\widehat{V}_-$ is compact and since $K_{PG}(\phi) \subseteq \widehat{V}_-$, by Lemma~\ref{Lem psi0 equal exp length} and Lemma~\ref{Lem polynomially growing currents empty}~$(3)$, there exists a constant $K_2>0$ such that such that for every $m \geq N$ and for every element $w \in F_{\tt n}$ such that $[\eta_{[w]}] \in U$, we have $\frac{\ell_{exp}([f^N(\gamma_w)])}{\ell_{\mathcal{F}}([f^N(\gamma_w)])} \geq K_2$. Thus, we have $$ \ell_{\mathcal{F}}([f^m(\gamma_w)]) \geq \ell_{exp}([f^m(\gamma_w)]) \geq 3^{m-N}K_1\ell_{exp}([f^N(\gamma_w)]) \geq 3^{n-M}K_1K_2 \ell_{\mathcal{F}}([f^N(\gamma_w)]).$$ We set $K_3=K_1K_2$. By compactness of $\PCurr(F_{\tt n},\mathcal{F}\wedge \mathcal{A}(\phi))$ and Lemma~\ref{Lem polynomially growing currents empty}~$(3)$, there exists $L>0$ such that for every current $[\mu] \in \PCurr(F_{\tt n},\mathcal{F}\wedge \mathcal{A}(\phi))$, we have $\frac{\lVert\phi^N(\mu)\rVert_{\mathcal{F}}}{\lVert \mu \rVert_{\mathcal{F}}} \geq L$. Hence for every $m \geq N$ and for every element $w \in F_{\tt n}$ such that $[\eta_{[w]}] \in U$, we have $$\ell_{\mathcal{F}}([f^m(\gamma_w)]) \geq 3^{m-N}K_3L\ell_{\mathcal{F}}(\gamma_w).$$ Hence the proof follows when $[\mu] \in U$.

By Theorem~\ref{Theo North-South dynamics almost atoroidal}, there exists $M_1 \in \NN^*$ such that, for every $m \geq M_1$ and every $[\mu] \in \PCurr(F_{\tt n},\mathcal{F} \wedge \mathcal{A}(\phi))-\widehat{V}_-$, we have $\phi^m([\mu]) \in U$. Let $M=M_1+N$. By the above, Lemma~\ref{Lem psi0 equal exp length}, the density of rational currents (see~Proposition~\ref{Prop density rational currents}) and continuity of $\phi$, for every current $[\mu] \notin \widehat{V}_-$, for every $n \geq M$, we have $$\lVert \phi^n(\mu) \rVert_{\mathcal{F}} \geq 3^{n-M}K_3L\lVert \phi^{M_1}(\mu) \rVert_{\mathcal{F}}. $$ By compactness of $\PCurr(F_{\tt n},\mathcal{F}\wedge \mathcal{A}(\phi))$ and Lemma~\ref{Lem polynomially growing currents empty}~$(3)$, there exists $L'>0$ such that for every current $[\mu] \in \PCurr(F_{\tt n},\mathcal{F}\wedge \mathcal{A}(\phi))$, we have $\frac{\lVert\phi^{M_1}(\mu)\rVert_{\mathcal{F}}}{\lVert \mu \rVert_{\mathcal{F}}} \geq L'$. Hence for every $n \geq M$, we have $$\lVert \phi^n(\mu) \rVert_{\mathcal{F}} \geq 3^{n-M}K_3LL'\lVert \mu \rVert_{\mathcal{F}}. $$ This concludes the proof.
\hfill\qedsymbol

\bibliographystyle{alphanum}
\bibliography{bibliographie}

\begin{thebibliography}{BFH2}

\bibitem[BFH1]{BesFeiHan00}
M.~Bestvina, M.~Feighn, and M.~Handel.
\newblock {\it The Tits alternative for Out$(\mathrm{F}_n)$ I: Dynamics of
  exponentially growing automorphisms}.
\newblock {Ann. of Math. {\bf 151} (2000) 517--623}.

\bibitem[BFH2]{BesFeiHan04}
M.~Bestvina, M.~Feighn, and M.~Handel.
\newblock {\it Solvable subgroups of $\operatorname{Out}(F_n)$ are virtually
  abelian}.
\newblock {Geom. Dedicata {\bf 104} (2004) 71--96}.

\bibitem[BH]{BesHan92}
M.~Bestvina and M.~Handel.
\newblock {\it Train tracks and automorphisms of free groups}.
\newblock {Ann. of Math. {\bf 135} (1992) 1--51}.

\bibitem[Bow]{Bowditch2012}
B.~Bowditch.
\newblock {\it Relatively hyperbolic groups}.
\newblock {Internat. J. Algebra Comput. (3) {\bf 22} (2012)}.

\bibitem[Coo]{Cooper87}
D.~Cooper.
\newblock {\it Automorphisms of free groups have finitely generated fixed point
  sets}.
\newblock {J. Algebra {\bf 111} (1987) 453--456}.

\bibitem[Cou]{Coulon2019}
R.~Coulon.
\newblock {\it Examples of groups whose automorphisms have exotic growth}.
\newblock {Preprint {\tt [arXiv:1908.11668]}}.

\bibitem[CU]{clay2019atoroidal}
M.~Clay and C.~Uyanik.
\newblock {\it Atoroidal dynamics of subgroups of $\mathrm{Out} ({F}_N)$}.
\newblock {J. London Math. Soc. (2) {\bf 102} (2020) 818--845}.

\bibitem[CV]{Vogtmann1986}
M.~Culler and K.~Vogtmann.
\newblock {\it Moduli of graphs and automorphisms of free groups.}
\newblock {Invent. Math. {\bf 84} (1986) 91--120}.

\bibitem[DS]{DahmaniKrish2020}
F.~Dahmani and S.~Krishna~M S.
\newblock {\it Relative hyperbolicity of hyperbolic-by-cyclic groups}.
\newblock {Preprint {\tt [arXiv:2006.07288]}, to appear in Groups, Geom. ,Dyn.}

\bibitem[FH]{FeiHan06}
M.~Feighn and M.~Handel.
\newblock {\it The recognition theorem for $\operatorname{Out}(F_n)$}.
\newblock {Groups Geom. Dyn. (1) {\bf 5} (2011) 39--106}.

\bibitem[GH]{Guirardelhorbez19laminations}
V.~Guirardel and C.~Horbez.
\newblock {\it Algebraic laminations for free products and arational trees}.
\newblock {Alg. Geom. Topol. (5) {\bf 19} (2019) 2283--2400}.

\bibitem[GL]{GabLev95}
D.~Gaboriau and G.~Levitt.
\newblock {\it The rank of actions on $\RR$-trees}.
\newblock {Ann. Scien. Ec. Norm. Sup. (4) {\bf 28} (1995) 549--570}.

\bibitem[Gue1]{Guerch2021currents}
Y.~Guerch.
\newblock {\it Currents relative to a malnormal subgroup system}.
\newblock {Preprint {\tt [arXiv:2112.01112]}}.

\bibitem[Gue2]{Guerch2021Polygrowth}
Y.~Guerch.
\newblock {\it Polynomial growth and subgroups of $\mathrm{Out}(F_n)$}.
\newblock {In preparation}.

\bibitem[Gup1]{gupta2017relative}
R.~Gupta.
\newblock {\it Relative currents}.
\newblock {Conform. Geom. Dyn. {\bf 21} (2017) 319--352}.

\bibitem[Gup2]{gupta18}
R.~Gupta.
\newblock {\it Loxodromic elements for the relative free factor complex}.
\newblock {Geom. Dedicata {\bf 196} (2018) 91--121}.

\bibitem[Hel]{Helfgott2015}
H.A Helfgott.
\newblock {\it Growth in groups: ideas and perspectives}.
\newblock {Bull. Amer. Math. Soc. (N.S.) {\bf 52} (2015) 357–-413}.

\bibitem[HM]{HandelMosher20}
M.~Handel and L.~Mosher.
\newblock {\it Subgroup Decomposition in $\mathrm{Out}(F_n)$}.
\newblock {Mem. Amer. Math. Soc. {\bf 264} (2020)}.

\bibitem[Hru]{Hruska10}
G.C. Hruska.
\newblock {\it Relative hyperbolicity and relative quasiconvexity for countable
  groups}.
\newblock {Algebr. Geom. Topol. {\bf 10} (2010) 1807--1856}.

\bibitem[Iva]{Ivanov92}
N.~V. Ivanov.
\newblock {\it Subgroups of Teichm\"uller modular groups}.
\newblock {Trans. Math. Mono. {\bf 115} Amer. Math. Soc. (1992)}.

\bibitem[Kap]{Kapovich2006}
I.~Kapovich.
\newblock {\it Currents on free groups}.
\newblock {In Topological and asymptotic aspects of group theory, Contemp.
  Math. {\bf 394} (2006) Amer. Math. Soc., 149--176}.

\bibitem[Lev]{Levitt09}
G.~Levitt.
\newblock {\it Counting growth types of automorphisms of free groups}.
\newblock {Geom. Funct. Anal. {\bf 19} (2009) 1119--1146}.

\bibitem[LL]{LevLus03}
G.~Levitt and M.~Lustig.
\newblock {\it Irreducible automorphisms of $F_n$ have north-south dynamics on
  compactified outer space}.
\newblock {J. Inst. Math. Jussieu {\bf 2} (2003) 59--72}.

\bibitem[LS]{LubotzkySegal2003}
A.~Lubotzky and D.~Segal.
\newblock {\it Subgroup growth}.
\newblock {Progress in Mathematics, {\bf 212} Birkh{\"a}user Verlag, Basel,
  2003}.

\bibitem[LU1]{LustigUyanik2017}
M.~Lustig and C.~Uyanik.
\newblock {\it Perron-Frobenius theory and frequency convergence for reducible
  substitutions}.
\newblock {Discrete Contin. Dyn. Syst. (1) {\bf 37} (2017) 355--385}.

\bibitem[LU2]{LustigUyanik2019}
M.~Lustig and C.~Uyanik.
\newblock {\it North-South dynamics of hyperbolic free group automorphisms on
  the space of currents}.
\newblock {J. Topol. Anal. (2) {\bf 11} (2019) 427--466}.

\bibitem[Man]{Mann12}
A.~Mann.
\newblock {\it How groups grow}.
\newblock {London Math.~Soc.~Lect.~Note Ser. {\bf 395}, Cambridge Univ.~Press,
  2012}.

\bibitem[Mar]{Martin95}
R.~Martin.
\newblock {\it Non-uniquely ergodic foliations of thin-type, measured currents
  and automorphisms of free groups}.
\newblock {PhD thesis, University of California, Los Angeles, 1995}.

\bibitem[McC]{McCarthy85}
J.~McCarthy.
\newblock {\it A ``Tits-alternative" for subgroups of surface mapping class
  groups}.
\newblock {Trans. Amer. Math. Soc. {\bf 291} (1985) 583--612}.

\bibitem[MM]{MasMin99}
H.~Masur and Y.~Minsky.
\newblock {\it Geometry of the complex of curves I: Hyperbolicity}.
\newblock {Invent. Math. {\bf 138} (1999) 103--149}.

\bibitem[Que]{Queffelec1987}
M.~Queff{\'e}lec.
\newblock {\it Substitution dynamical systems-spectral analysis}.
\newblock {Lect. Notes in Math. {\bf 1294}, Springer Verlag, 1987}.

\bibitem[Ser]{Serre83}
J.-P. Serre.
\newblock {\it Arbres, amalgames, SL$_2$}.
\newblock {3\`eme \'ed. corr., Ast\'erisque {\bf 46}, Soc. Math. France, 1983}.

\bibitem[SW]{ScoWal79}
G.~P. Scott and C.~T.~C. Wall.
\newblock {\it Topological methods in group theory}.
\newblock {dans ``Homological group theory'', C.~T.~C.~Wall ed., Lond. Math.
  Soc. Lect. Notes {\bf 36}, Cambridge Univ. Press (1979) 137--203}.

\bibitem[Thu]{Thurston88}
W.~Thurston.
\newblock {\it On the geometry and dynamics of diffeomorphisms of surfaces}.
\newblock {Bull. Amer. Math. Soc. {\bf 19} (1988) 417-432}.

\bibitem[Tit]{Tits72}
J.~Tits.
\newblock {\it Free subgroups in linear goups}.
\newblock {J. Algebra {\bf 20} (1972) 250--270}.

\bibitem[Uya1]{Uyanik2015}
C.~Uyanik.
\newblock {\it Generalized north-south dynamics on the space of geodesic
  currents}.
\newblock {Geom. Dedicata {\bf 177} (2015) 129--148}.

\bibitem[Uya2]{Uyanik2019}
C.~Uyanik.
\newblock {\it Hyperbolic extensions of free groups from atoroidal ping-pong}.
\newblock {Algebr. Geom. Topol. (3) {\bf 19} (2019) 1385--1411}.

\end{thebibliography}

\noindent \begin{tabular}{l}
Laboratoire de mathématique d'Orsay\\
UMR 8628 CNRS \\
Université Paris-Saclay\\
91405 ORSAY Cedex, FRANCE\\
{\it e-mail: yassine.guerch@universite-paris-saclay.fr}
\end{tabular}

\end{document}